\documentclass[preprint,11pt]{elsarticle}
\usepackage[bookmarks=true,colorlinks=true,linkcolor=blue]{hyperref}

\pdfoutput=1

\biboptions{sort&compress}

\usepackage{color}
\definecolor{jwbGreen}{rgb}{0, .6, 0}

\usepackage{soul}

\usepackage[usenames,dvipsnames,svgnames,table]{xcolor}

\usepackage{amssymb,amsmath}
\usepackage{amsthm}

\usepackage{calc}

\usepackage{calc}
\usepackage[margin=1.in]{geometry}

\renewcommand{\url}[1]{}
\newcommand{\citeCount}[1]{}
\newtheorem{theorem}{Theorem}

\newtheorem{definition}{Definition}

\newcommand{\esup}{{(e)}}
\newcommand{\psup}{{(p)}}
\newcommand{\ksup}{{(k)}}

\newcommand{\dx}{\Delta x}
\newcommand{\dy}{\Delta y}
\newcommand{\dt}{\Delta t}

\newcommand{\dn}{\Delta n}

\newcommand{\Dp}{D_{+}}
\newcommand{\Dm}{D_{-}}

\newcommand{\bogus}[1]{{}}

\newenvironment{myIndent}%
 {\list{}{\leftmargin=0.1in\rightmargin=0.1in}\item[]}%
  {\endlist}

\newenvironment{myProcedure}[1]%
{\noindent\textbf{Procedure~}{#1}\begin{myIndent}\em}
{\end{myIndent}}

\newenvironment{centre}{\begin{center}}{\end{center}}

\newcommand{\eqdef}{\overset{{\rm def}}{=}}

\newcommand{\av}{\mathbf{ a}}

\newcommand{\ev}{\mathbf{ e}}

\newcommand{\iv}{\mathbf{ i}}

\newcommand{\nv}{\mathbf{ n}}

\newcommand{\tv}{\mathbf{ t}}

\newcommand{\vv}{\mathbf{ v}}

\newcommand{\xv}{\mathbf{ x}}

\newcommand{\Fv}{\mathbf{ F}}
\newcommand{\Gv}{\mathbf{ G}}

\newcommand{\Real}{{\mathbb R}}
\newcommand{\Imag}{\text{Im}}

\newcommand{\Bc}{{\mathcal B}}

\newcommand{\Dc}{{\mathcal D}}

\newcommand{\Gc}{{\mathcal G}}
\newcommand{\Gs}{{\mathcal G}}
\newcommand{\Ic}{{\mathcal I}}

\newcommand{\Nc}{{\mathcal N}}

\newcommand{\Sc}{{\mathcal S}}

\newcommand{\sigmav}{\boldsymbol{\sigma}}

\newcommand{\grad}{\nabla}

\newcommand{\tableFont}{\scriptsize}

\newcommand{\num}[2]{#1e#2} 

\clearpage

\newcommand{\rhos}{\rho_b}
\newcommand{\rhob}{\rho_b}

\newcommand{\normalss}{\sffamily}

\newcommand{\Hs}{\bar{H}}

\newcommand{\uHat}{\hat{u}}
\newcommand{\vHat}{\hat{v}}
\newcommand{\pHat}{\hat{p}}

\newcommand{\OmegaF}{\Omega}

\newcommand{\OmegaB}{{\Omega_b}}
\newcommand{\GammaB}{\Gamma_b}

\newcommand{\blue}{\color{blue}}
\newcommand{\red}{\color{red}}

\newcommand{\strutt}{\rule{0pt}{9pt}}

\newcommand\Omegaf{{\Omega}}

\newcommand\tnv{{\mathbf{t}}}

\renewcommand\tv{\tnv}

\newcommand{\mrb}{m_{b}}
\newcommand{\vb}{v_b}
\newcommand{\xvb}{\xv_b}
\newcommand{\vvb}{\vv_b}
\newcommand{\avb}{\av_b}

\newcommand{\mb}{{m_b}}

\newcommand{\adc}{\beta_d}

\newcommand{\Da}{\Dc}

\newcommand{\ampRB}{AMP-RB}

\newcommand{\alphas}{\bar{\alpha}}

\newcommand{\Duu}{\Da^{u}}
\newcommand{\Dun}{\Da^{u}_n}
\newcommand{\Duua}{{{\cal D}^\omega}}

\newcommand{\dr}{\Delta r}

\newcommand{\ub}{u_b}

\newcommand{\ab}{\abu}

\newcommand{\gvnp}{g_v(t^{n+1})}
\newcommand{\gunp}{g_u(t^{n+1})}

\newcommand{\Dyh}{D_{yh}}

\newcommand{\mbBar}{\bar{m}_b}

\newcommand{\bodyStepIComment}{Preliminary body evolution step}

\newcommand{\wHat}{{\hat w}}
\def\aa{r_1}
\def\bb{r_2}

\newcommand{\cxi}{{C_\xi}}
\newcommand{\ceta}{{C_\eta}}

\newcommand{\NAb}{\Nc_b}
\newcommand{\NAv}{\Nc_v}
\newcommand{\SRb}{\Sc_{b}}
\newcommand{\SRv}{\Sc_{v}}
\newcommand{\NAbt}{\widetilde\Nc_b}
\newcommand{\NAvt}{\widetilde\Nc_v}

\newcommand{\cone}{{\tilde C_1}}
\newcommand{\ctwo}{{\tilde C_2}}
\newcommand{\cmmm}{{\tilde C_m}}
\newcommand{\IbBar}{{\bar I_b}}

\newcommand{\abu}{a_u}
\newcommand{\abv}{a_v}
\newcommand{\tpRB}{{TP-RB}}

\usepackage{tikz}

\newlength{\tfwidth}
\newlength{\tfheight}
\newlength{\tfxa}
\newlength{\tfxb}
\newlength{\tfya}
\newlength{\tfyb}
\newcommand{\trimFigWithBox}[6]{%
\setlength\fboxsep{0pt}%
\setlength\fboxrule{1.0pt}
\fbox{\includegraphics[width=#2, clip, trim=#3 #4 #5 #6]{#1}}%
}
\newcommand{\trimFigNoBox}[6]{%
\setlength\fboxsep{1pt}
\setlength\fboxrule{0.0pt}
\fbox{\includegraphics[width=#2, clip, trim=#3 #4 #5 #6]{#1}}%
}
\newcommand{\trimFigHeightWithBox}[6]{%
\setlength\fboxsep{0pt}%
\setlength\fboxrule{1.0pt}
\fbox{\includegraphics[height=#2, clip, trim=#3 #4 #5 #6]{#1}}%
}
\newcommand{\trimFigHeightNoBox}[6]{%
\setlength\fboxsep{1pt}
\setlength\fboxrule{0.0pt}
\fbox{\includegraphics[height=#2, clip, trim=#3 #4 #5 #6]{#1}}%
}
\newcommand{\trimFigb}[6]{%
\setlength{\tfwidth}{(#2+#2*\real{#3})+#2*\real{#4}}
\setlength{\tfheight}{(#2+#2*\real{#5})+#2*\real{#6}}%
\setlength{\tfxa}{\tfwidth*\real{#3}}%
\setlength{\tfxb}{\tfwidth*\real{#4}}%
\setlength{\tfya}{\tfheight*\real{#5}}%
\setlength{\tfyb}{\tfheight*\real{#6}}%
\trimFigWithBox{#1}{#2}{\tfxa}{\tfya}{\tfxb}{\tfyb}%
}
\newcommand{\trimFig}[6]{%
\setlength{\tfwidth}{(#2+#2*\real{#3})+#2*\real{#4}}
\setlength{\tfheight}{(#2+#2*\real{#5})+#2*\real{#6}}%
\setlength{\tfxa}{\tfwidth*\real{#3}}%
\setlength{\tfxb}{\tfwidth*\real{#4}}%
\setlength{\tfya}{\tfheight*\real{#5}}%
\setlength{\tfyb}{\tfheight*\real{#6}}%
\trimFigNoBox{#1}{#2}{\tfxa}{\tfya}{\tfxb}{\tfyb}%
}

\newsavebox\figBox

\newcommand{\trimw}[6]{%
\sbox\figBox{\includegraphics{#1}}
\setlength{\tfwidth}{\the\wd\figBox}
\setlength{\tfheight}{\the\ht\figBox}
\setlength{\tfxa}{\tfwidth*\real{#3}}%
\setlength{\tfxb}{\tfwidth*\real{#4}}%
\setlength{\tfya}{\tfheight*\real{#5}}%
\setlength{\tfyb}{\tfheight*\real{#6}}%
\trimFigNoBox{#1}{#2}{\tfxa}{\tfya}{\tfxb}{\tfyb}%
}

\newcommand{\trimwb}[6]{%

\sbox\figBox{\includegraphics{#1}}
\setlength{\tfwidth}{\the\wd\figBox}
\setlength{\tfheight}{\the\ht\figBox}
\setlength{\tfxa}{\tfwidth*\real{#3}}%
\setlength{\tfxb}{\tfwidth*\real{#4}}%
\setlength{\tfya}{\tfheight*\real{#5}}%
\setlength{\tfyb}{\tfheight*\real{#6}}%
\trimFigWithBox{#1}{#2}{\tfxa}{\tfya}{\tfxb}{\tfyb}%
}

\newcommand{\trimh}[6]{%
\sbox\figBox{\includegraphics{#1}}
\setlength{\tfwidth}{\the\wd\figBox}
\setlength{\tfheight}{\the\ht\figBox}
\setlength{\tfxa}{\tfwidth*\real{#3}}%
\setlength{\tfxb}{\tfwidth*\real{#4}}%
\setlength{\tfya}{\tfheight*\real{#5}}%
\setlength{\tfyb}{\tfheight*\real{#6}}%
\trimFigHeightNoBox{#1}{#2}{\tfxa}{\tfya}{\tfxb}{\tfyb}%
}

\newcommand{\trimhb}[6]{%

\sbox\figBox{\includegraphics{#1}}
\setlength{\tfwidth}{\the\wd\figBox}
\setlength{\tfheight}{\the\ht\figBox}
\setlength{\tfxa}{\tfwidth*\real{#3}}%
\setlength{\tfxb}{\tfwidth*\real{#4}}%
\setlength{\tfya}{\tfheight*\real{#5}}%
\setlength{\tfyb}{\tfheight*\real{#6}}%
\trimFigHeightWithBox{#1}{#2}{\tfxa}{\tfya}{\tfxb}{\tfyb}%
}
\usepackage{algorithm} 
\usepackage{algpseudocode} 

\usepackage{setspace}

\usepackage[font=footnotesize]{caption}

\begin{document}

\small

\begin{frontmatter}

\title{
A stable partitioned FSI algorithm for rigid bodies and incompressible flow.
       Part I: Model problem analysis
}

\author[rpi]{J.~W.~Banks\fnref{DOEThanks,PECASEThanks}}
\ead{banksj3@rpi.edu}

\author[rpi]{W.~D.~Henshaw\corref{cor1}\fnref{DOEThanks,NSFgrantNew}}
\ead{henshw@rpi.edu}

\author[rpi]{D.~W.~Schwendeman\fnref{DOEThanks,NSFgrantNew}}
\ead{schwed@rpi.edu}

\author[rpi]{Qi Tang\fnref{QiThanks}}
\ead{tangq3@rpi.edu}

\address[rpi]{Department of Mathematical Sciences, Rensselaer Polytechnic Institute, Troy, NY 12180, USA.}

\cortext[cor1]{Department of Mathematical Sciences, Rensselaer Polytechnic Institute, 110 8th Street, Troy, NY 12180, USA.}

\fntext[QiThanks]{Research supported by the Eliza Ricketts Postdoctoral Fellowship.}

\fntext[DOEThanks]{This work was supported by contracts from the U.S. Department of Energy ASCR Applied Math Program.}

\fntext[NSFgrantNew]{Research supported by the National Science Foundation under grant DMS-1519934.}

\fntext[PECASEThanks]{Research supported by a U.S. Presidential Early Career Award for Scientists and Engineers.}

\begin{abstract}

A stable partitioned algorithm is developed for fluid-structure interaction (FSI) problems involving viscous incompressible flow and rigid bodies. This {\em added-mass partitioned} (AMP) algorithm remains stable, without sub-iterations, for light and even zero mass rigid bodies when added-mass and viscous added-damping effects are large.  The scheme is based on a generalized Robin interface condition for the fluid pressure that includes terms involving the linear acceleration and angular acceleration of the rigid body.  Added-mass effects are handled in the Robin condition by inclusion of a boundary integral term that depends on the pressure.  Added-damping effects due to the viscous shear forces on the body are treated by inclusion of added-damping tensors that are derived through a linearization of the integrals defining the force and torque.  Added-damping effects may be important at low Reynolds number, or, for example, in the case of a rotating cylinder or rotating sphere when the
rotational moments of inertia are small.  In this first part of a two-part series, the properties of the AMP scheme are motivated and evaluated through the development and analysis of some model problems.  The analysis shows when and why the traditional partitioned scheme becomes unstable due to either added-mass or added-damping effects.  The analysis also identifies the proper form of the added-damping which depends on the discrete time-step and the grid-spacing normal to the rigid body.  The results of the analysis are confirmed with numerical simulations that also demonstrate a second-order accurate implementation of the AMP scheme.

\end{abstract}


\end{frontmatter}

\clearpage
\tableofcontents

\clearpage

\section{Introduction} \label{sec:intro}

We consider fluid-structure interaction (FSI) problems involving the motion of
rigid bodies in an incompressible fluid.  
FSI problems of this type occur in a wide variety
of applications, such as ones involving particulate flows (suspensions, sedimentation, fluidized beds),
valves and moving appendages, buoy structures, ship maneuvering, and underwater vehicles, to name a few. 
A wide range of numerical techniques have been developed to simulate such FSI problems, including arbitrary
Lagrangian--Eulerian (ALE) methods~\cite{takashi1992, HuPatankarZhu2001, vierendeels2005analysis}, methods
based on level-sets~\cite{coquerelle2008vortex, gibou2012efficient}, fictitious domain methods
\cite{glowinski1999distributed, glowinski2000distributed, GlowinskiPanHelsaJosephPeriaux2001},
embedded boundary methods \cite{costarelli2016embedded} and immersed boundary methods
\cite{kajishima2002interaction, uhlmann2005immersed, kim2006immersed, lee2008immersed,
  borazjani2008curvilinear, breugem2012second, kempe2012improved, yang2012simple, bhalla2013unified,
  yang2015non, wang2015strongly, kim2016penalty, lacis2016stable}.

Numerical schemes for simulating FSI problems can be classified as monolithic or partitioned.
Monolithic schemes couple the governing equations for all domains into a single large system.
Partitioned schemes make use of separate solvers for the fluid and solid domains with the separate
solvers being coupled at fluid-solid interfaces. The choice of discrete interface conditions is one of
the primary challenges in developing accurate and stable partitioned schemes. 
Many partitioned schemes for FSI problems coupling rigid bodies and incompressible
flow suffer from numerical instabilities for light bodies.
These instabilities usually originate from added-mass 
and added-damping effects as noted in~\cite{Conca1997,HuPatankarZhu2001,uhlmann2005immersed,kempe2012improved,breugem2012second,yang2015non,lacis2016stable} for example.\footnote{At the analytic level the notions of added-mass and added-damping have been well known for decades; 
in a more recent discussion, for example, the form of the added-mass and added-damping is nicely described
in Conca~et~al.~\cite{Conca1997}, where they consider some model problems for which the contributions of added-mass or added-damping in the solutions can be determined explicitly.}
The dominant instability is often due to added-mass effects that arise
from the pressure-induced force and torque on the body.
As illustrated in Figure~\ref{fig:addedMassCartoon}, when a force is applied instantaneously to 
a rigid body in an incompressible fluid, the resultant 
acceleration of the body depends both on the mass of the body and on the {\em added-mass} of a portion
of the surrounding fluid that is accelerated due to the instantaneous pressure field created by the motion of the body.
If the body is heavy, then the effect of the added-mass is small and its associated instability is weak.  For light
bodies, the added-mass can be significant and lead to added-mass instabilities for partitioned solvers.

{
\newcommand{\lbfont}{\small}
\def\ysb{-.15} 
\def\ysa{-2}   
\def\ya{0} 
\def\yb{2} 
\def\xL{5}
\def\rad{2.}
\newcommand{\plotDisk}[1]{
\fill[fill=#1!20,draw=red,line width=2pt] 
      plot[samples=100, domain=0.:360] ( {\rad*cos(\x)} , {\rad*sin(\x)} ) -- cycle ;
}
%
%
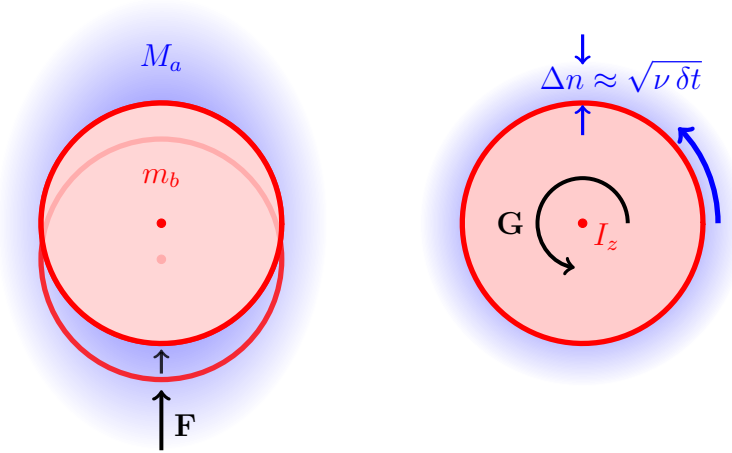
\begin{figure}[hbt]
\newcommand{\textFont}{\normalss}
\begin{center}
  \begin{tikzpicture}[scale=.8]
  \useasboundingbox (0,0) rectangle (13,6.75);  
  %
  \begin{scope}[xshift=3cm,yshift=3cm]
    \begin{scope}[yshift=0cm]
      \shadedraw[shading=radial,inner color=blue!70,draw=white] (0,0) ellipse (2.75cm and 3.75cm);
    \end{scope}
    \plotDisk{white}
    \begin{scope}[yshift=-.6cm,opacity=.8]
      \draw[draw=red,line width=2pt] (0,0) ellipse (2cm and 2cm);
      \draw[->,very thick,black] (0,-1.9) -- (0,-1.5); 
      \draw[fill=red,draw=red]  (0,0) circle (2pt);
    \end{scope}
    \begin{scope}[opacity=.8]
      \plotDisk{red}
    \end{scope}
    \draw[fill=red,draw=red]  (0,0) circle (2pt);
    \begin{scope}[xshift=0pt,yshift=-2.6cm,black]
      \draw[->,line width=1.5pt,yshift=-5pt] (0,-1) -- (0,.0); %
      \draw (0,-.75) node[anchor=west,xshift=1pt] {\large$\Fv$}; 
    \end{scope}
    \draw (0,.75) node[] {\large\red$\mb$}; 
    \draw (0,2.75) node[] {\large\blue $M_a$}; 
  \end{scope}
  \begin{scope}[yshift=3cm,xshift=10cm]
    \begin{scope}[yshift=0cm]
      \shadedraw[shading=radial,inner color=blue!90,draw=white] (0,0) circle (2.7cm);
    \end{scope}
    \plotDisk{red}
    \draw[,line width=1.5pt,->] (.75,0) arc (0:260:.75cm); %
    \draw (-.75,0) node[anchor=east,xshift=-1pt] {\large$\Gv$}; 
    \draw[fill=red,draw=red]  (0,0) circle (2pt);
    \draw (0,0) node[xshift=9pt,yshift=-5pt] {\large\red $I_z$}; 
    \begin{scope}[yshift=2cm]
      \draw (0,0) node[anchor=south west,xshift=-20pt] {\large\blue$\Delta n \approx \sqrt{\nu\, \delta t}$}; 
      \draw[->,very thick,draw=blue,yshift=-1pt] (0,-.5) -- (0,0);
      \draw[->,very thick,draw=blue,yshift=4pt] (0,1.) -- (0,.5);
    \end{scope}
    \draw[,line width=2pt,->,blue] (2.25cm,0) arc (0:45:2.25cm); %
  \end{scope}
  %
  \end{tikzpicture}
\end{center}
\caption{Cartoon illustrating added-mass and added-damping. 
Left: a force $F$ applied to a rigid-body of mass $m_b$ immersed in an incompressible fluid 
accelerates a nearby volume of fluid whose mass, the {\em added-mass} denoted by $M_a$, contributes to an effective mass of the body. 
Right: a torque $G$ applied to a rigid disk with moment of inertia $I_z$ over a time interval $\delta t$ results in a viscous boundary layer of fluid of approximate width $\Delta n\approx \sqrt{\nu\,\delta t}$, where $\nu$ is the
kinematic viscosity of the fluid.  The drag on the surface of the disk associated with this
layer leads to an {\em added-damping} on the rotation of the disk.
} \label{fig:addedMassCartoon}
\end{figure}
}

A secondary {\em added-damping} instability may also occur due to the choice of the numerical
treatment of the force and torque on the body associated with the viscous shear stress on the surface of the body.
Again referring to
Figure~\ref{fig:addedMassCartoon}, the resultant angular motion of a rigid disk due to a torque
applied about its centre of mass depends both on the moment of inertia of the disk and
on the viscous drag created by the fluid in a boundary layer near the surface of the disk.
The width of this fluid layer grows with time so that it has a damping effect on the motion
over an interval of time rather than an instantaneous effect on the acceleration as in the case of added-mass.
An implication of this is
that instabilities due to added-damping can usually be overcome by taking a sufficiently small time-step (except for
the case of a body of zero mass), while added-mass instabilities generally 
cannot be overcome by taking a smaller time-step if the mass of the body is too small relative to the added-mass.

Implementations of partitioned schemes in the literature based
on the traditional coupling between the fluid and solid\footnote{The traditional coupling approach, which
uses the velocity of the body as a boundary condition on the fluid, and the fluid stress as a boundary
condition on the rigid body, is described
in Section~\ref{sec:TBRB}.}
 usually become unstable
when the ratio of the density of the solid to that of the fluid is approximately
one or less~\cite{uhlmann2005immersed,borazjani2008curvilinear}.  The common approach to overcome these
instabilities
is to apply multiple under-relaxed sub-iterations
per time-step, often coupled with an acceleration technique such as Aitken's method 
(e.g., see~\cite{borazjani2008curvilinear, vierendeels2005analysis, wang2015strongly}).
However, these sub-iterations can be expensive since a recomputation of the fluid velocity and pressure is
required at each iteration, and since a large number of sub-iterations is required when the body is light.
In recent years, improving the numerical stability for simulations involving light rigid bodies has attracted considerable attention.
For instance, the stability of the immersed boundary method in Uhlmann~\cite{uhlmann2005immersed} 
was improved in the work of Kempe and Fr{\"o}hlich~\cite{kempe2012improved}
by explicitly evaluating a volume integral related to the artificial flow field inside the body, which
resulted in an
improved stability bound to 
a  density ratio of $0.3$ or greater.
Yang and Stern~\cite{yang2015non} improved this bound further to $0.1$ by relying on a  non-iterative immersed boundary approach.
They predicated the velocity field in a temporary non-inertial frame that enabled an implicit direct solution of the rigid-body motion
and an exact match of the velocity of the fluid and rigid body at the interface.
Lacis~et~al.~\cite{lacis2016stable} recently proposed a different non-iterative immersed boundary method that 
is stable for a density ratio as low as $10^{-4}$. 
Their approach is similar to the block-factorization approach used for similar FSI problems
by Fedkiw~et~al.~\cite{Robinson-MosherSchroederFedkiw2011} and Badia~et~al.~\cite{BadiaQuainiQuarteroni2008}. 
In this approach, an approximate block LU factorization
of the discrete monolithic equations is performed that exposes
one set of decoupled equations for the fluid velocity and another set of equations that couples the pressure 
with the motion for the body.
This treatment 
is similar in spirit to the manner that added-mass
effects are handled by our new scheme (as described herein), 
except that our coupling conditions are derived at the continuous level; these continuous conditions are conducive to the
development of high-order accurate approximations.
Our approach for treating added-damping appears to be new.  

The present paper is the first of a two-part series in which we describe a new
added-mass/added-damping partitioned approach, referred to as the~\ampRB~scheme, for FSI problems
involving the motion of rigid bodies in an incompressible flow.  The principal feature of
the~\ampRB~scheme is that it remains stable, without the need for sub-time-step iterations, for FSI
problems with light bodies (and even zero-mass bodies) when added-mass and added-damping effects are
large.  The time-stepping scheme is based on a fractional-step approach for the fluid in which the velocity
is advanced in one stage, and the pressure is determined in a second stage~\cite{ICNS,splitStep2003}.  
The viscous terms in the stress
tensor are handled implicitly so that the velocity can be advanced with a larger stable time step.
 
The key ingredients of the~\ampRB~scheme are contained in the AMP interface conditions which couple the equations of
motion of the rigid body to a compatibility boundary condition for the pressure on the surface of
the body.  These conditions are derived at a continuous level by matching the {\em acceleration} of
the body to that of the fluid.  As a result, the integration of the equations of motion of the
rigid-body are coupled strongly to the update of the fluid pressure;  this 
ensures the proper balances of forces at the interface thereby suppressing instabilities
due to added-mass effects.  Suppressing instabilities due to added-damping is more subtle.
The fluid forces on the body depend on the viscous shear stresses which, in turn, implicitly depend on the 
velocity of the body. This implicit dependence of the fluid forces on the body velocity
is explicitly exposed and, after some simplifying approximations,
is expressed in terms of {\em added-damping tensors} which
are incorporated into the AMP interface conditions as a means to overcome added-damping
instabilities.

The work here in Part~I focuses on the development of the scheme for two relatively simple model
problems, one involving a rectangular geometry and the other posed in an annular domain.
The~\ampRB~scheme is described fully for the two model problems and the stability of the scheme is
analyzed.  For the two model problems, it is shown that the scheme is stable (without sub-time-step
iterations) for any ratio of the density of the body to that of the fluid (including the case of a
rigid body with zero mass).  The present implementation of the numerical scheme is second-order
accurate, and this is verified (along with its stability) through a series of test calculations.  In
Part~II~\cite{rbins2016r}, the~\ampRB~scheme is extended to FSI problems involving more general
geometric configurations.  This is done using moving overlapping grids~\cite{mog2006} to handle more
complex fluid domains involving the motion of one or more rigid bodies.  Several FSI simulations are
performed which also confirm the stability and accuracy of the approach for the general case.

Finally, it is noted that the~\ampRB~scheme described here, and its extension in Part~II, 
follows the development of previous AMP schemes for other FSI regimes.
These schemes include ones for inviscid compressible
fluids coupled to rigid solids~\cite{lrb2013}, linearly elastic solids~\cite{fsi2012} and nonlinear
elastic solids~\cite{flunsi2016}, and ones for incompressible fluids coupled to elastic bulk
solids~\cite{fib2014}, elastic structural shells~\cite{fis2014} and deforming beams
\cite{beamins2016}. For all cases, the various AMP schemes are stable for light solids without requiring
sub-time-step iterations.

The reminder of the paper is organized as follows.
Section~\ref{sec:governing} defines the governing equations for a rectangular-domain FSI model problem that
contains the essential elements needed to address added-mass and added-damping effects. 
Section~\ref{sec:TBRB} defines, in detail, a traditional-partitioned (TP-RB) algorithm in terms
of a number of procedures to lay the ground work for describing the changes needed for the new {\ampRB}~scheme
described in the subsequent Section~\ref{sec:ampRB}. 
A normal-mode stability analysis of the TP-RB and {\ampRB}~schemes is performed in Section~\ref{sec:stability},
while Section~\ref{sec:analysisDisk} extends the analysis to a second FSI model problem described in an annular geometry.
Numerical verification of the stability and accuracy of the {\ampRB}~scheme is provided in Section~\ref{sec:numericalVerification}.
Conclusions of the work in Part~I are given in the final Section~\ref{sec:conclusions}.

\section{Governing equations for a FSI model problem} \label{sec:governing}

As noted in the introduction, added-mass and added-damping effects arise from the response of the fluid to forces and torques applied to the rigid body.  To understand the essential features of these effects with the aim of designing stable numerical methods, it is useful to identify key FSI model problems and then undertake a detailed study of these problems.  Consider, for example, the local region near the curved boundary of the rigid body immersed in an incompressible fluid as shown on the left of Figure~\ref{fig:addedMassModelProblemCartoon}.  To derive a suitable FSI model problem for this configuration, we employ a standard mathematical technique and consider a PDE boundary-value problem in a local region near the curved boundary and transform this region to a half plane problem. The transformed half-plane problem is then linearized, the coefficients are frozen, and the lower-order terms are dropped. The result is a model that consists of a rectangular body in the lower half-plane adjacent to an incompressible fluid in the upper-half plane, as shown in the right side of Figure~\ref{fig:addedMassModelProblemCartoon}. For this FSI problem, added-mass effects arise from a vertical acceleration of the body while added-damping effects occur from a horizontal acceleration.  The problem can be simplified further by considering a finite domain and assuming, for example, periodic boundary conditions for the fluid in the horizontal direction and a boundary condition for velocity or pressure of the fluid at a height $H$ in the vertical direction. 

{
\newcommand{\lbfont}{\small}
\def\ysb{.0} 
\def\ysa{-2}   
\def\ya{0} 
\def\yb{2} 
\def\xL{5}
\def\rad{2.}
\newcommand{\plotDisk}{
\fill[fill=red!20,draw=red,line width=2pt] 
      plot[samples=100, domain=0.:360] ( {\rad*cos(\x)} , {\rad*sin(\x)} ) -- cycle ;
}
%
%
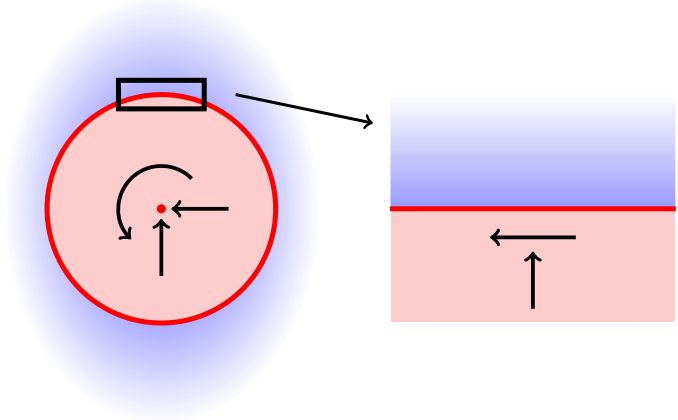
\begin{figure}[hbt]
\newcommand{\textFont}{\normalss}
\begin{center}
 \resizebox{10cm}{!}{
 \begin{tikzpicture}[scale=.8]
  \useasboundingbox (0,.1) rectangle (13,6.5);  
  %
  \begin{scope}[xshift=3cm,yshift=3cm]
    \begin{scope}[yshift=0cm]
      \shadedraw[shading=radial,inner color=blue!70,draw=white] (0,0) ellipse (2.75cm and 3.75cm);
    \end{scope}
    \plotDisk
    \draw[->,line width=1.5pt,yshift=-5pt] (0,-1.) -- (0,0); %
    \draw[line width=1.5pt,->] (.530,.530) arc (45:225:.75cm); %
    \draw[<-,line width=1.5pt,xshift=5pt] (0,0) -- (1,0); %
    \draw[fill=red,draw=red]  (0,0) circle (2pt);
    \begin{scope}[xshift=0cm,yshift=2cm]
      \draw[-,line width=2pt]  (-.75,-.25) rectangle (.75,.25);
    \end{scope}
    \draw[->,very thick] (1.3,2) -- (3.7,1.5); 
    \begin{scope}[xshift=4cm,yshift=0cm]
      \shade[bottom color=blue!40,top color=white,draw=white] (0,\ya) rectangle (\xL,\yb);
      \draw[thick,fill=red!20,draw=white] (0,\ysa) rectangle (\xL,\ysb);
      \draw[draw=red,line width=2pt] (0,\ysb) -- (\xL,\ysb);
      \draw[line width=1.5pt,<-,yshift=-1.5cm] (1.75,1) -- (3.25,1); 
      \draw[line width=1.5pt,->,yshift=-2.25cm] (2.5,.5) -- (2.5,1.5); 
    \end{scope}
  \end{scope}
  %
  \end{tikzpicture}
  }
\end{center}
\caption{Prototypical added-mass and added-damping problems are derived from localizing the problem on the left near the boundary of the rigid body. 
This leads to the simplified rectangular-geometry problem on the right for which added-mass effects arise due to vertical motions
while added-damping effects arise from horizontal translations.
} \label{fig:addedMassModelProblemCartoon}
\end{figure}
}

Following the prescription outlined above, we consider a first model FSI problem in which the fluid occupies the rectangular domain, $\OmegaF = [0, L] \times [0,H]$ with coordinates $\xv=[x, y]^T$, and lies above a rigid body occupying the region $\OmegaB = [0, L] \times [-\bar{H}, 0]$, see Figure~\ref{fig:rectangularModelProblem}.  The fluid in~$\OmegaF$ has velocity, $\vv(\xv,t)=[u(\xv,t), v(\xv,t)]^T$, and pressure, $p(\xv,t)$, and is taken to be viscous and incompressible.  It is assumed that the rigid body undergoes small displacements in the $x$ or $y$-direction, but does not rotate, and that the interface separating the fluid and the rigid body is linearized about the fixed portion of the $x$-axis given by $\GammaB=[0, L] \times \{ 0 \}$.  The velocity and acceleration of the body are denoted by $\vvb(t)=[u_b(t), v_b(t)]^T$ and $\dot{\bf v}_b(t)=\avb(t)=[\abu(t), \abv(t)]^T$, respectively, and the motion of the body is connected to that of the fluid by conditions on velocity and stress along $\GammaB$.  The full set of equations governing the model problem are
\begin{alignat}{2}
 \text{Fluid:}\quad & \rho \frac{\partial \vv}{\partial t} + \grad p = \mu \Delta\vv , \quad&& \xv\in\OmegaF, \label{eq:velocityEquation}\\
  &  \grad\cdot\vv = 0   , \quad&&\xv\in\OmegaF, \label{eq:continuityEquation}\\
\text{Rigid body:}\quad  &   m_b\, \abu = \int_0^L \mu \frac{\partial u}{\partial y}(x,0,t)\,dx + g_u(t) , \quad && \label{eq:bodyEquation1}\\
  &   m_b\, \abv = -\int_0^L p(x,0,t)\,dx + g_v(t) ,   \label{eq:bodyEquation2} \\
 \text{Interface:}\quad 
&  \vv(x,0, t) = \vvb(t), \quad&& x\in[0,L] , \label{eq:matchV} \\
 \text{Top fluid BCs:} \quad 
&  u(x,H,t) =u_H(x,t), \quad&& x\in[0,L] ,\\
&  p(x,H,t)=p_H(x,t), \quad&& x\in[0,L] ,\\
\text{Side fluid BCs:}\quad
& \vv(0,y,t)=\vv(L,y,t),\quad&& y\in[0,H] . \label{eq:modelProblemLast}
\end{alignat}
Here, $\rho$ is the constant fluid density, $\mu$ is the constant fluid viscosity, $\mb$ is the mass of the
rigid body, and  $g_u(t)$ and $g_v(t)$ are external forces on the rigid body in the horizontal and vertical directions, respectively.  These equations are augmented with appropriate initial conditions. 

{
\newcommand{\lbfont}{\small}
\def\ysb{-.15} 
\def\ysa{-2}   
\def\ya{0} 
\def\yb{3} 
\def\xL{8}
\begin{figure}[hbt]
	\newcommand{\textFont}{\normalss}
	\begin{center}
           \resizebox{6cm}{!}{
		\begin{tikzpicture}[scale=.9]
		\useasboundingbox (0,\ysa) rectangle (\xL,3.2);  
		\draw[thick,fill=red!20,draw=red,line width=2pt] (0,\ysa) rectangle (\xL,\ysb);
		\draw[thick,black] (4,-1.2) node {rigid body: $\OmegaB$};
		\draw[thick,fill=blue!20,draw=blue,line width=2pt] (0,\ya) rectangle (\xL,\yb);
		\draw[thick,blue] (0,\ya) rectangle (\xL,\yb);
		\draw[thick,black] (4,1.75) node {fluid: $\OmegaF$};
		\draw[thick,black] (4,\ya) node[anchor=south] {interface: $\GammaB$};

		\draw[-,thick,red] (0  ,\ysa) -- (\xL,\ysa);
		\draw[-,thick,red] (0  ,\ysb) -- (\xL,\ysb);
		\draw[-,thick,red] (0  ,\ya) -- (0  ,\ya) node[anchor=east,black,yshift=-4pt] {\lbfont$y=0$};
		\draw[-,thick,red] (\xL,\ysa) -- (\xL,\ysb);
		\draw[-,thick,blue] (0,\yb) -- (\xL,\yb);
		\draw[-,thick,blue] (\xL,0) -- (\xL,\yb);
		\draw (0  ,\ysa) node[anchor=north,black,yshift=-4pt] {\lbfont$x=0$};
		\draw (\xL,\ysa) node[anchor=north,black,yshift=-4pt] {\lbfont$x=L$};
		\draw (0,  \ysa) node[anchor=east,black] {\lbfont$y=-\Hs$};
		\draw (0,   \yb) node[anchor=east,black] {\lbfont$y=H$};
		%
		%
		\end{tikzpicture}
            }
	\end{center}
	\caption{The geometry for the rectangular geometry FSI model problems.} \label{fig:rectangularModelProblem}
\end{figure}
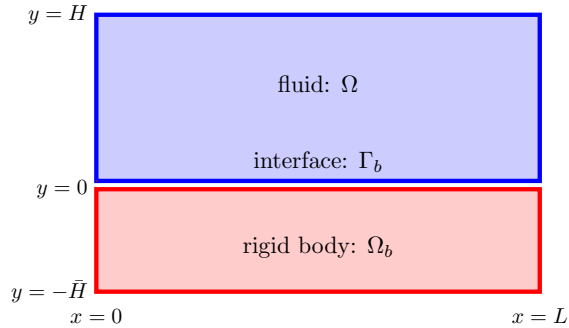
}

An elliptic equation for the fluid pressure can be derived by taking the divergence of the momentum equation in~\eqref{eq:velocityEquation} and using the continuity equation in~\eqref{eq:continuityEquation} to give
\begin{alignat}{3}
  &  \Delta p &= 0 ,  \qquad&& \xv\in\OmegaF . \label{eq:fluidPressure}
\end{alignat}
In the velocity-pressure form of the Stokes equations, the momentum equation for the velocity in~\eqref{eq:velocityEquation} is used with the Laplace equation for the pressure in~\eqref{eq:fluidPressure} instead of the divergence equation in~\eqref{eq:continuityEquation}.  To retain consistency with the original velocity-divergence form of the equations, zero divergence
is specified as an additional boundary condition, i.e.~$\grad\cdot\vv=0$ for
$\xv\in\partial\OmegaF$, see~\cite{ICNS} for additional details.  

In the discussion to follow we use a fractional-step method based on the velocity-pressure form of the fluid equations to
define partitioned FSI algorithms. In Section~\ref{sec:TBRB} a traditional partitioned scheme,
referred to as the \tpRB~algorithm, is described. Then in Section~\ref{sec:ampRB} we present our new
{\it added-mass partitioned} scheme which will be referred to as the \ampRB~algorithm. The important
difference in the two approaches lies in the coupling of the fractional-step method for the fluid with the
equations governing the motion of the body.

\section{Traditional Partitioned/Rigid-Body  algorithm \tpRB} \label{sec:TBRB}

We begin our discussion of partitioned algorithms for the rectangular model FSI problem by considering a {\em
  traditional partitioned} scheme.  This is done to establish the basic elements of a partitioned
FSI algorithm based on a fractional-step time-stepping scheme for the velocity-pressure
formulation of the fluid equations and to lay the groundwork for a later discussion of the new {\em
  added-mass partitioned} scheme.

To set some notation, let $\xv_\iv$ denote the grid-point coordinates of a mesh covering the fluid
domain, $\xv\in\OmegaF$, where $\iv=(i_1,i_2)$ is a multi-index, and let $\iv\in\OmegaF_h$ denote
the set of indices of points in the interior of the fluid grid while $\iv\in\Gamma_h$ are the
indices of points on the interface, $\xv\in\GammaB$.  Grid functions for the fluid velocity and
pressure are given by $\vv_\iv^n \approx \vv(\xv_\iv,t^n)$ and $p_{\iv}^n \approx p(\xv_{\iv},t^n)$,
respectively, where $t^n=n\dt$, $n=0,1,2,\ldots$, denote discrete times with time step $\dt$.  The
rigid-body displacement, velocity, and acceleration at discrete times are given by
  $\xvb^n\approx\xvb(t^n)$, $\vvb^n\approx \vvb(t^n)$, and $\avb^n\approx \avb(t^n)$ respectively.
Second-order accurate approximations are used for the spatial operators
in~\eqref{eq:velocityEquation}, \eqref{eq:continuityEquation} and~\eqref{eq:fluidPressure}, and we
denote these discrete approximations using a subscript~$h$.  Precise definitions for these
approximations are not important in the present discussion. In describing the algorithm,
  $\vv_\iv^n$, $p_\iv^n$, $\vvb^n$ and~$\avb^n$ are assumed to be known at~$t\sp n$ and~$t\sp{n-1}$,
  and the output is the discrete solution at~$t\sp{n+1}$.

The \tpRB~time-stepping scheme is described in Algorithm~\ref{alg:tpRB}, and uses four procedures which are defined below.
The resulting scheme is formally second-order accurate in both space and time, using either the predicted values as the solution at $t\sp{n+1}$, or the discrete solution after the correction steps.  Generally, the predictor-corrector scheme has favorable stability properties in comparison to the scheme using the prediction steps alone, but in either case the algorithm suffers from instabilities when added-mass or added-damping effects are significant.  This occurs when the ratio of the mass of the rigid body to that of the fluid contained in $\OmegaF$ is too small as shown in Sections~\ref{sec:amSubmodel} and~\ref{sec:adSubmodel}.  The added-mass/added-damping instabilities can be suppressed by coupling the calculation of the fluid pressure and the rigid-body acceleration as discussed in the next section.

\begin{algorithm}\caption{Traditional partitioned (\tpRB) scheme}\small
\medskip
\[
\begin{array}{l}
\hbox{// \textsl{{\bodyStepIComment}}}\smallskip\\
1.\quad (\,\avb\sp\esup\,,\,\vvb\sp\esup\,)=\hbox{\bf predictAndExtrapRigidBody}(\,\avb^n\,,\,\avb^{n-1}\,,\,\vvb^{n-1}\,)\bigskip\\
\hbox{// \textsl{Prediction steps}}\smallskip\\
2.\quad (\,\vv_\iv\sp\psup\,)=\hbox{\bf advanceFluidVelocity}(\,\vvb\sp\esup\,,\,2p_{\iv}\sp{n}-p_{\iv}\sp{n-1}\,,\,p_{\iv}^n\,,\,\vv_\iv\sp{n}\,)\smallskip\\
3.\quad (\,p_\iv\sp\psup\,)=\hbox{\bf updateFluidPressure}(\,\vv_\iv\sp\psup\,,\,\avb\sp\esup\,)\smallskip\\
4.\quad (\,\vvb\sp\psup\,,\,\avb\sp\psup\,)=\hbox{\bf advanceBodyVelocityAndAcceleration}(\,\vvb^{n}\,,\,\avb^{n}\,,\,\vv_\iv\sp\psup\,,\,p_\iv\sp\psup\,)\bigskip\\
\hbox{// \textsl{Correction steps}}\smallskip\\
5.\quad (\,\vv_\iv\sp{n+1}\,)=\hbox{\bf advanceFluidVelocity}(\,\vvb\sp\psup\,,\,p_{\iv}\sp\psup\,,\,p_{\iv}^n\,,\,\vv_\iv\sp{n}\,)\smallskip\\
6.\quad (\,p_\iv\sp{n+1}\,)=\hbox{\bf updateFluidPressure}(\,\vv_\iv\sp{n+1}\,,\,\avb\sp\psup\,)\smallskip\\
7.\quad (\,\vvb\sp{n+1}\,,\,\avb\sp{n+1}\,)=\hbox{\bf advanceBodyVelocityAndAcceleration}(\,\vvb^{n}\,,\,\avb^{n}\,,\,\vv_\iv\sp{n+1}\,,\,p_\iv\sp{n+1}\,)
\end{array}
\]
\label{alg:tpRB}
\end{algorithm}

\newcommand{\stageI}{BVA(p)}
\newcommand{\stageII}{FV(p)}
\newcommand{\stageIII}{FP+BA}
\newcommand{\stageIV}{BV}
\newcommand{\stageV}{FV}

\newcommand{\TPstageIII}{FP}
\newcommand{\TPstageIV}{BVA}

\vskip\baselineskip
\begin{myProcedure}{$(\,\avb\,,\,\vvb\,)=\hbox{\bf predictAndExtrapRigidBody}(\,\avb^n\,,\,\avb^{n-1}\,,\,\vvb^{n-1}\,)$}
  Compute provisional values $(\avb,\vvb)$ for the acceleration and velocity of the rigid body at $t^{n+1}$ using a linear extrapolation in time for the acceleration followed by a leap-frog integration for the velocity: 
\begin{align*}
&    \avb = 2\avb^n - \avb^{n-1} , \\
&    \vvb = \vvb^{n-1} + 2\dt\, \avb^{n}. 
\end{align*}
\end{myProcedure}

\vskip\baselineskip
\begin{myProcedure}{$(\,\vv_\iv\,)=\hbox{\bf advanceFluidVelocity}(\,\vvb\sp{*}\,,\,p_{\iv}\sp{*}\,,\,p_{\iv}^n\,,\,\vv_\iv\sp{n}\,)$}
Advance the fluid velocity $\vv_\iv$ to $t^{n+1}$ using a second-order accurate semi-implicit scheme for the momentum equation in~\eqref{eq:velocityEquation}: 
\begin{alignat*}{3}
&   \rho\frac{\vv_{\iv} -\vv^n_{\iv}}{\dt} + \frac{1}{2}\grad_h p_{\iv}^* +\frac{1}{2} \grad_h p_{\iv}^{n} 
         = \frac{\mu}{2}\Big( \Delta_h \vv_{\iv} + \Delta_h \vv_\iv^n \Big)  , \quad&&\iv\in\OmegaF_h, \\
&   \vv_\iv = \vvb\sp{*} ,\quad \grad_h\cdot\vv_\iv=0, \quad&&\iv\in\Gamma_h. 
\end{alignat*}
In this procedure, the viscous term for the velocity is computed using an implicit trapezoidal method, while the pressure gradient uses a predicted value for the pressure at $t^{n+1}$ given by $p_{\iv}^*$ so that this term is explicit.
A discrete approximation of the continuity equation is used to determine the vertical (normal) component of the velocity in the ghost points associated with the interface $\iv\in\Gamma_h$, while the horizontal (tangential) component of the velocity is determined in the ghost points using extrapolation.  Discrete boundary conditions are also applied at the boundary of the mesh corresponding to $y=H$, as well as periodic conditions along the sides corresponding to $x=0$ and $x=L$.
\end{myProcedure}

\vskip\baselineskip
\begin{myProcedure}{$(\,p_\iv\,)=\hbox{\bf updateFluidPressure}(\,\vv_\iv\sp{*}\,,\,\avb\sp{*}\,)$}
 Update the fluid pressure $p_\iv$ at $t\sp{n+1}$ by solving a discrete Laplace equation for the pressure
\begin{equation}
 \Delta_h p_\iv = 0 , \qquad \iv \in \Omega_h\cup\Gamma_h,
\label{eq:discretePressure}
\end{equation}
together with the interface condition
\begin{equation}
\Dyh \, p_\iv = -\rho \abv\sp{*}
               -\mu \grad_h\times\grad_h\times v_\iv^{*} , \qquad\qquad\iv\in\Gamma_h, \label{eq:pressureInterfaceEquationTP}
\end{equation}
which is derived from the $y$-component of the momentum equation for the fluid
velocity.\footnote{The term $\Delta v$ from the $y$-component of~\eqref{eq:velocityEquation} has
  been replaced by $-\grad\times\grad\times v$ to avoid a viscous time-step restriction in the
  application of the discrete interface condition in~\eqref{eq:pressureInterfaceEquationTP},
  see~\cite{Petersson00} for additional details.}  Here, $D_{yh}$ is a second-order accurate
approximation to the $y$-derivative at the interface, and $\vv_\iv\sp{*}$ and $\avb\sp{*}$ are
predicted values for the fluid velocity and the acceleration of the body, respectively, at
$t\sp{n+1}$.  Periodic boundary conditions are applied on $x_{\iv}=0$ and $L$, and appropriate
boundary conditions are applied on $y_{\iv}=H$.
\end{myProcedure}

\vskip\baselineskip
\begin{myProcedure}{$(\,\vvb\,,\,\avb\,)=\hbox{\bf advanceBodyVelocityAndAcceleration}(\,\vvb^{n}\,,\,\avb^{n}\,,\,\vv_\iv\sp{*}\,,\,p_\iv\sp{*}\,)$}
Compute the rigid-body acceleration $\avb$ at $t\sp{n+1}$ by solving discrete approximations of the rigid-body equations in~\eqref{eq:bodyEquation1} and~\eqref{eq:bodyEquation2} given by
\begin{align}
&   a_u = {1\over\mrb}\left\{ \sum_{\iv\in\Gamma_h}  \mu \bigl(\Dyh\,u^{*}_\iv\bigr)\dx + g_u(t^{n+1})\right\}, \label{eq:TBRBau}\\
&   \abv = {1\over\mrb}\left\{-\sum_{\iv\in\Gamma_h}  p^{*}_\iv \dx + g_v(t^{n+1})\right\}, \label{eq:TBRBav}
\end{align}
and then advance the rigid-body velocity $\vvb$ using the trapezoidal method
\begin{alignat*}{3}
&   
   \vvb = \vvb^{n} + \frac{\dt}{2}\big( \avb + \avb^{n}\big)   .
\end{alignat*}
Here, $\vv_\iv\sp{*}$ and $p_\iv\sp{*}$ are computed grid functions for the fluid velocity and pressure, respectively, at  $t\sp{n+1}$.
\end{myProcedure}

\section{Added-Mass Partitioned/Rigid-Body algorithm \ampRB}\label{sec:ampRB}

In this section we develop the \ampRB~algorithm for the FSI model problem
given in~\eqref{eq:velocityEquation}--\eqref{eq:modelProblemLast}. Before defining the algorithm, we
first describe the key interface conditions that are used to account for the added-mass and
added-damping effects. The simplified interface conditions described here
 are a special case of the more general
conditions discussed in Part~II~\cite{rbins2016r}. The simplified conditions 
motivate the more general conditions and also 
illustrate the key ideas behind the AMP approach for incompressible flows
and rigid bodies. The resulting AMP algorithm for the model problem
retains the essential characteristics of the general AMP scheme described in Part~II, 
but is simpler to describe and analyze.

\subsection{Added-mass interface conditions}

In the model problem described in Section~\ref{sec:governing}, 
added-mass effects arise in~\eqref{eq:bodyEquation2} from the force on the rigid body due to the fluid pressure 
through the integral term
\begin{equation}
  \Ic_p(p, \abv) \eqdef -\int_0^L p(x,0,t)\,dx.
\label{eq:icp}
\end{equation}
Notice that the integral depends on the fluid pressure, $p$, but the pressure, in turn, is an implicit function of the vertical acceleration of the body, $\abv$, through the coupled  motion of the fluid and the body.  
Therefore one may regard the integral in~\eqref{eq:icp} to be a function of $p$ and $\abv$.  
The contribution to the added-mass from this integral term is revealed by matching the vertical {\em acceleration} of the fluid and body on the interface,
\begin{equation}
    \frac{\partial v}{\partial t} = \frac{d v_b}{dt} = \abv,  \quad \xv\in\Gamma ,\label{eq:matchA}
\end{equation}
which follows from the time derivative of the vertical component of the interface condition in~\eqref{eq:matchV}. 
Combining~\eqref{eq:matchA} with the $y$-component of the fluid momentum equation in~\eqref{eq:velocityEquation} implies the two coupled interface conditions given by
\begin{align}
 &     \rho \abv + \frac{\partial p}{\partial y} =  \mu \Delta v, \qquad \xv\in\Gamma, \label{eq:ADBC1} \\
 &  \mb\, \abv = \Ic_p(p, \abv) + g_v(t) . \label{eq:ADBC2}
\end{align}
In the \ampRB~algorithm described below, equations~\eqref{eq:ADBC1} and~\eqref{eq:ADBC2} are used as interface conditions when solving the Laplace equation for the fluid pressure (with $v$ on the left-hand side of~\eqref{eq:ADBC1} assumed to be given by a previous step in the algorithm).  We note that the vertical acceleration of body, $\abv$, in~\eqref{eq:ADBC1} can be viewed as a Lagrange multiplier that is adjusted to enforce the scalar constraint given by~\eqref{eq:ADBC2}.  In this formulation, the vertical component of the acceleration of the body is determined with the fluid pressure in the solution of the elliptic problem, and the key idea is that this can be done when $\mb$ is small, or even zero.

Note that when $\mb\ne0$, the vertical component of the acceleration, $\abv$, can be eliminated from~\eqref{eq:ADBC1} and~\eqref{eq:ADBC2} to give a generalized Robin condition on the pressure
\begin{equation}
  \frac{\mb}{\rho}  \frac{\partial p}{\partial y} + \int_0^L p\,dx
        =  \frac{\mb}{\rho} \, \mu \Delta v +  g_v(t),\qquad  \xv\in\Gamma. \label{eq:ADBC3}
\end{equation}
Use of~\eqref{eq:ADBC3}, however, is not valid as the only boundary condition for the pressure
when $\mb=0$ since only one independent condition
remains; furthermore~\eqref{eq:ADBC3} may be poorly conditioned in the numerical scheme 
when $\mb$ is very small.

\subsection{Added-damping interface conditions} \label{sec:adinterface}

Added-damping effects arise primarily due to viscous shear forces on the body, and 
in the model problem of Section~\ref{sec:governing} these are
associated with the integral of the shear stress given by
\begin{align}
   \Ic_u(u,\ub) \eqdef \int_0^L \mu \frac{\partial u}{\partial y}(x,0,t)\,dx, \label{eq:shearIntegral}
\end{align}
which determines the horizontal force on the body in~\eqref{eq:bodyEquation1}.  Here $\Ic_u$ can be
considered to be a function of the horizontal components of the velocity of the fluid and the body
since $u$ is an implicit function of $u_b$.  In order to suppress added-damping effects in
the \ampRB~algorithm, it is important to reveal the dependence of $\Ic_u$ on the velocity $u_b$ (or
more specifically the acceleration $\abu$), at least in an approximate sense.  Let us assume
that predicted values for the velocities, $u^\psup$ and $\ub^{\psup}$, are available and consider
the linear approximation
\begin{align}
\Ic_u(u,\ub) &\approx \Ic_u(u^\psup,\ub^\psup) + \left(\frac{\partial\Ic_u}{\partial\ub}\right)\sp\psup\big( \ub - \ub^\psup \big), \nonumber\\
             & = \Ic_u(u^\psup,\ub^\psup)  -   \Duu \big( \ub - \ub^\psup \big),
\label{eq:addedDampingConditionU}
\end{align}
where $\Duu$ is the {\em added-damping coefficient}\footnote{For more general FSI problems, $\Duu$ takes the form of an {\em added-damping tensor}, see Part II~\cite{rbins2016r}.} defined by
\begin{equation}
\Duu \eqdef -\frac{\partial\Ic_u}{\partial\ub}.
\label{eq:addedDampingDu}
\end{equation}
A formula for $\Duu$ can be found by considering a suitable variational problem associated with the FSI model problem as is discussed later in Section~\ref{sec:addedDampingCoefficient}.  

In practice, it is convenient to work with the acceleration $\abu$, instead of the velocity $\ub$, which suggests an alternate linearization of the form
\begin{align}
    \Ic_u(u,\ub) &\approx \Ic_u(u^\psup,\ub^\psup) -  \adc \, \dt\, \Duu\, \big( \abu - \abu^\psup ),
     \label{eq:addedDampingCondition}
\end{align}
where $\abu^\psup$ is a predicted acceleration of the body (note the factor of $\dt$) and 
where we have introduced the {\em added-damping parameter} $\adc$, whose value plays an important role in the stability of the~\ampRB~scheme. 
\begin{definition} {\bf Added-damping parameter $\adc$.} The parameter
$\adc$ is called the added-damping parameter and is a free parameter 
that multiplies the added-damping coefficient $\Duu$ in the linearization of the integral of the shear stress.
\end{definition}
A stability analysis of the \ampRB~algorithm applied to an {\em added-damping} model problem reveals a
range of stable values for $\adc$ (see Section~\ref{sec:stabilityMP-AD}).  For now, we note that the
linearization in~\eqref{eq:addedDampingCondition} is an important element in the \ampRB~algorithm
since it enables a stable calculation of the horizontal acceleration of the body, $\abu$,
in~\eqref{eq:bodyEquation1} even when $\mb$ is small, or even zero, i.e.~when added-damping effects
are significant.

\subsection{\ampRB~algorithm}

We are now in a position to describe the \ampRB~algorithm for the FSI model problem.
As in the \tpRB~algorithm, the fluid is advanced using a fractional-step method based on the velocity-pressure form of the equations.  While this approach is common to both algorithms, there are important differences.  The principal difference is that the \ampRB~algorithm uses the interface conditions in~\eqref{eq:ADBC1} and~\eqref{eq:ADBC2} for added-mass, and the condition in~\eqref{eq:addedDampingCondition} for added-damping in the stage involving the elliptic problem for the fluid pressure.  Thus, the pressure and the acceleration of the rigid body are determined {\em together} in one stage of the \ampRB~algorithm, while the velocity of the body is advanced in time in a subsequent stage.

\begin{algorithm}\caption{Added-mass partitioned (\ampRB) scheme}\small
\medskip
\[
\begin{array}{l}
\hbox{// \textsl{{\bodyStepIComment}}}\smallskip\\
1.\quad (\,\avb\sp\esup\,,\,\vvb\sp\esup\,)=\hbox{\bf predictAndExtrapRigidBody}(\,\avb^n\,,\,\avb^{n-1}\,,\,\vvb^{n-1}\,)\bigskip\\
\hbox{// \textsl{Prediction steps}}\smallskip\\
2.\quad (\,\vv_\iv\sp\psup\,)=\hbox{\bf advanceFluidVelocity}(\,\vvb\sp\esup\,,\,2p_{\iv}\sp{n}-p_{\iv}\sp{n-1}\,,\,p_{\iv}^n\,,\,\vv_\iv\sp{n}\,)\smallskip\\
3.\quad (\,p_\iv\sp\psup\,,\,\avb\sp\psup\,)=\hbox{\bf updateFluidPressureAndBodyAcceleration}(\,\vv_\iv\sp\psup\,,\,\avb\sp\esup\,)\smallskip\\
4.\quad (\,\vvb\sp\psup\,)=\hbox{\bf advanceBodyVelocityOnly}(\,\avb^\psup\,,\,\avb^{n}\,,\,\vvb^{n}\,)\bigskip\\
\hbox{// \textsl{Correction steps}}\smallskip\\
5.\quad (\,\vv_\iv\sp{n+1}\,)=\hbox{\bf advanceFluidVelocity}(\,\vvb\sp\psup\,,\,p_{\iv}\sp\psup\,,\,p_{\iv}^n\,,\,\vv_\iv\sp{n}\,)\smallskip\\
6.\quad (\,p_\iv\sp{n+1}\,,\,\avb\sp{n+1}\,)=\hbox{\bf updateFluidPressureAndBodyAcceleration}(\,\vv_\iv\sp{n+1}\,,\,\avb\sp\psup\,)\smallskip\\
7.\quad (\,\vvb\sp{n+1}\,)=\hbox{\bf advanceBodyVelocityOnly}(\,\avb^{n+1}\,,\,\avb^{n}\,,\,\vvb^{n}\,)\bigskip\\
\hbox{// \textsl{Fluid-velocity correction step (optional)}}\smallskip\\
8.\quad (\,\vv_\iv\sp{n+1}\,)=\hbox{\bf advanceFluidVelocity}(\,\vvb\sp{n+1}\,,\,p_{\iv}\sp{n+1}\,,\,p_{\iv}^n\,,\,\vv_\iv\sp{n}\,)
\end{array}
\]
\label{alg:ampRB}
\end{algorithm}

The \ampRB~algorithm is given in Algorithm~\ref{alg:ampRB}, and uses two procedures introduced in Section~\ref{sec:TBRB} as well as two additional procedures defined below. The \ampRB~algorithm uses a similar predictor-corrector time-stepping approach as in the \tpRB~algorithm, and it is also second-order accurate in space and time.  For the \ampRB~algorithm, we have included an optional fluid {\em velocity correction} (VC) step which serves to adjust the fluid velocity so that it exactly matches the velocity of the body along the interface at $t\sp{n+1}$.  In practice, we have found that this step, which requires an additional implicit solve for the velocity, is only needed for difficult problems with strong added-mass and added-damping effects as we show in Section~\ref{sec:stabilityMP-AD}.  Also, it is possible to repeat the correction group of procedures multiple times, but in practice we only use one correction group.

\vskip\baselineskip
\begin{myProcedure}{$(\,p_\iv\,,\,\avb\,)=\hbox{\bf updateFluidPressureAndBodyAcceleration}(\,\vv_\iv\sp{*}\,,\,\avb\sp{*}\,)$}
Update the fluid pressure $p_\iv$ and the rigid-body acceleration $\abv$ at $t\sp{n+1}$ by solving a discrete Laplace equation for the pressure in~\eqref{eq:discretePressure}
together with the discrete AMP interface conditions
\begin{align}
&   \Dyh \, p_\iv + \rho \abv  
               =  -\mu \grad_h\times\grad_h\times v_\iv^{*} , \qquad\qquad\iv\in\Gamma_h, \label{eq:pressureInterfaceEquation} \\
&   \mrb \abv + \sum_{\iv\in\Gamma_h}  p_\iv \dx = g_v(t^{n+1}) , \label{eq:pressureInterfaceEquationII}\\
&   \left[ \mrb + \adc \, \dt\, \Duu\right] a_u = 
       \sum_{\iv\in\Gamma_h}  \mu \bigl(\Dyh \, u^{*}_\iv\bigr)\dx + \adc\, \dt\, \Duu\, \abu^{*}  + g_u(t^{n+1}) ,
        \label{eq:pressureInterfaceEquationIII}
\end{align}
which are approximations of the interface conditions in~\eqref{eq:ADBC1}, \eqref{eq:ADBC2} and~\eqref{eq:addedDampingCondition}.  
Here $\Duu$ is the added-damping coefficient and $\adc$ is the added-damping parameter introduced in Section~\ref{sec:adinterface}.
As before, $D_{yh}$ is a second-order accurate approximation to the $y$-derivative, and $\vv_\iv\sp{*}$ and $\avb\sp{*}$ are predicted values for the fluid velocity and the acceleration of the body, respectively, at $t\sp{n+1}$.  Periodic boundary conditions are applied on $x_{\iv}=0$ and $L$, and appropriate boundary conditions are applied on $y_{\iv}=H$.
\end{myProcedure}

For the model problem considered here, we note that the discrete pressure, $p_{\iv}$, and the
vertical component  of the rigid body acceleration, $\abv$, are coupled and
determined by the solution to the elliptic problem consisting of
the discrete Laplace equation in~\eqref{eq:discretePressure} with interface conditions
in~\eqref{eq:pressureInterfaceEquation} and~\eqref{eq:pressureInterfaceEquationII}, while the
horizontal component of the rigid body, $\abu$, is determined independently by the interface condition
in~\eqref{eq:pressureInterfaceEquationIII}.  However, for FSI problems with general curved
interfaces, these equations for the pressure and body accelerations are fully coupled, see Part II for more details.

\vskip\baselineskip
\begin{myProcedure}{$(\,\vvb\,)=\hbox{\bf advanceBodyVelocityOnly}(\,\avb^{*}\,,\,\avb^{n}\,,\,\vvb^{n}\,)$}
Advance the rigid-body velocity $\vvb$ to $t\sp{n+1}$ using the trapezoidal method
\begin{alignat*}{3}
&   
   \vvb = \vvb^{n} + \frac{\dt}{2}\big( \avb^{*} + \avb^{n}\big) ,
\end{alignat*}
using a computed value for the rigid-body $\avb^{*}$ at $t\sp{n+1}$.
\end{myProcedure}

\newcommand{\amplanar}{{MP-AM}}
\newcommand{\adplanar}{{MP-AD}}
\section{Stability analyses of the \ampRB~and \tpRB~schemes for FSI model problems} \label{sec:stability}

The stability properties of the~{\ampRB} and \tpRB~algorithms can be understood by considering two
sub-problems derived from the FSI model problem given
in~\eqref{eq:velocityEquation}--\eqref{eq:modelProblemLast}.  The first problem focuses on the
effects of added-mass, while the second problem reveals issues associated with added-damping.  A
derivation of the two problems begins with a Fourier analysis of the periodic solution for the fluid
variables.  Let $\uHat_k(y,t)$, $\vHat_k(y,t)$ and $\pHat_k(y,t)$ denote coefficient functions of
the Fourier series in the $x$-direction for the components of the fluid velocity and the pressure.
For example, $\uHat_k(y,t)$ is defined by
\begin{equation}
u(x,y,t) = \sum_{k=-\infty}^\infty \uHat_k(y,t) e^{2\pi i k x/L} ~, 
\label{eq:fourier}
\end{equation}
with similar formulas defining $\vHat_k(y,t)$ and $\pHat_k(y,t)$.  The equations governing the coefficient functions for the $k=0$ mode of the fluid, together with the equations and interface conditions for the rigid-body velocity, decouple into the two sub-problems, one for the pressure and the components of fluid and rigid-body velocity in the vertical direction, and the other the components of the velocity in the horizontal direction.  It is the problem for the pressure and vertical components of velocity that forms the model for added-mass effects, with the other problem being the model for added damping.  The problem for added mass, called {\amplanar}, is discussed in detail in Section~\ref{sec:amSubmodel}, while a discussion of the problem for added damping, called {\adplanar}, is given in Section~\ref{sec:adSubmodel}.

We note that the problems for the other Fourier-modes with $k\ne 0$ satisfy homogeneous Dirichlet
conditions at the interface and do not couple to the motion of the rigid body.  Since these modes
are not coupled to the motion of the rigid body, they are therefore not
  germane to our discussion of instabilities related to added-mass or added-damping effects and will
  not be considered further here.

\subsection{Analysis of a model problem for added mass} \label{sec:amSubmodel}

A model problem describing the added-mass effect is first derived. 
Subsequently, stability analyses of the TP-RB and {\ampRB}~schemes to solve this model problem are given.

\subsubsection{Added-mass model problem {\amplanar}}
The added-mass model problem 
is obtained from the $k=0$ mode in the Fourier expansion of the fluid variables
together with the equation for the vertical component of the rigid body.
This model problem for added mass ({\amplanar}) is given by
\begin{equation}
 \text{{\amplanar}:}\;  \left\{ 
   \begin{alignedat}{3}
  & \rho \frac{\partial \vHat}{\partial t} +\frac{\partial \pHat}{\partial y} = 0 , \quad&&  y\in(0,H),  \\
  &  \frac{\partial \vHat}{\partial y} = 0   , \quad&& y\in(0,H), \\
  &   m_b \frac{d v_b}{dt} = - L\, \pHat(0,t)+ g_v(t), \\
  &  \vHat(0,t)=v_b, \quad \pHat(H,t)=\pHat_H(t),
   \end{alignedat}  \right. 
  \label{eq:MP-AM}
\end{equation}
where the subscript on the coefficient functions corresponding to $k=0$ has been dropped for notational convenience.  The solution of the model problem is
\begin{align}
	v_b(t)&={1\over m_b+M_a}\int_0\sp{t} \Big[ g_v(\tau)-L\pHat_H(\tau) \Big] \,d\tau \, +v_b(0), \label{eq:vbsoln} \\
\pHat(y,t)&=\pHat_H(t)+\left(1-{y\over H}\right){g_v(t)/L-\pHat_H(t)\over m_b/M_a+1}, \label{eq:pHatsoln} \\
\vHat(t)&=v_b(t), \label{eq:vHatsoln}
\end{align}
where
\begin{equation}
M_a=\rho LH,
\label{eq:addedmassplanar}
\end{equation}
is the added mass for this problem and corresponds to the mass of the fluid in the whole rectangular domain.  
Note that the fluid velocity, $\vHat$, which only depends on time and not space,
is equal to the velocity of the rigid body according to the condition at the interface, $y=0$.
Meanwhile, the fluid pressure, $\pHat$, is a linear function of $y$ with time-dependent forcing.  
The acceleration of the body is given by
\begin{equation}
\abv(t)={dv_b\over dt}= {g_v(t)-L\pHat_H(t)\over m_b+M_a},
\label{eq:absoln}
\end{equation}
which is the total (time-dependent) force applied to the body divided by the sum of its mass, $m_b$, and the added mass, $M_a$. 
Importantly, the acceleration is seen to be bounded, even when $m_b=0$, assuming that the external forcing is bounded and that $M_a$ is bounded away from zero.
 
\subsubsection{Analysis of the \tpRB~algorithm for {\amplanar}}

A stability analysis of the \tpRB~time-stepping scheme in Algorithm~\ref{alg:tpRB} for the {\amplanar} model problem 
is now presented.
For this added-mass model problem, the fully discretized scheme simplifies to the one
described in detail in Algorithm~\ref{alg:tpRBam}.  Here, $\vHat_j\sp{n}$ and $\pHat_j\sp{n}$
approximate $\vHat(y,t)$ and $\pHat(y,t)$, respectively, at $y_j=j\dy$ and $t\sp{n}=n\dt$, and
$v_b\sp{n}$ and $a_v\sp{n}$ approximate $v_b(t)$ and $a_v(t)$, respectively, at $t\sp{n}$.  The
symbols $D_0$, $D_+$ and $D_-$ denote centred, forward and backward divided difference operators,
respectively, defined by
\begin{equation}
D_0\pHat_j\sp{n}\eqdef{\pHat_{j+1}\sp{n}-\pHat_{j-1}\sp{n}\over2\dy},\qquad 
D_+\pHat_j\sp{n}\eqdef{\pHat_{j+1}\sp{n}-\pHat_{j}\sp{n}\over\dy},\qquad 
D_-\pHat_j\sp{n}\eqdef{\pHat_{j}\sp{n}-\pHat_{j-1}\sp{n}\over\dy},
\label{eq:ddopts}
\end{equation}
while $D_{yh}$ is a second-order accurate one-sided approximation of the first derivative given by
\begin{equation}
D_{yh}\pHat_j\sp{n}={-3\pHat_{j}\sp{n}+4\pHat_{j+1}\sp{n}-\pHat_{j+2}\sp{n}\over2\dy}.
\label{eq:funkyD}
\end{equation}
For the stability analysis of this scheme, 
it is convenient to consider the homogeneous equations with zero forcing so that the discrete variables appearing in the various steps of the algorithm are regarded as {\em perturbations} to the exact discrete solution.
The objective of the analysis is then to determine under what conditions the perturbations can be bounded.

\begin{algorithm}\caption{\tpRB~scheme for {\amplanar} 
                         \\ \protect\hphantom{Algorithm 5~~} {\em Added-mass model problem in a rectangular geometry.}}\small
\medskip
\[
\begin{array}{l}
\hbox{// \textsl{{\bodyStepIComment}}}\smallskip\\
1.\quad a_v^\esup = 2 a_v^n - a_v^{n-1} , \qquad \vb^\esup = \vb^{n-1} + 2\dt\, a_v^{n},\bigskip\\
\hbox{// \textsl{Prediction steps}}\smallskip\\
2.\quad \vHat_j\sp\psup = \vHat_j^n - \frac{3\dt}{2\rho} D_0 \pHat_{j}^n +\frac{\dt}{2\rho} D_0 \pHat_j^{n-1}, \quad j=1,2,\ldots,N, \qquad \vHat_0^\psup = \vb^\esup, \quad D_0\vHat_{N}^\psup = 0, \smallskip\\
3.\quad D_+D_- \pHat_j^\psup = 0, \quad j=1,2,\ldots,N-1, \qquad  D_{yh} \pHat_0^\psup= -\rho a_v^\esup, \quad \pHat_N^\psup=\pHat_H(t\sp{n+1}), \smallskip\\
4.\quad a_v^\psup =  (-L\pHat_0^\psup +\gvnp)/\mb , \qquad  \vb^\psup = \vb^{n} + \frac{\dt}{2}( a_v^\psup + a_v^n ),\bigskip\\
\hbox{// \textsl{Correction steps}}\smallskip\\
5.\quad \vHat_j\sp{n+1} = \vHat_j^n - \frac{\dt}{2\rho}( D_0 \pHat_{j}^\psup + D_0 \pHat_j^{n}), \quad j=1,2,\ldots,N, \qquad \vHat_0^{n+1} = \vb^\psup, \quad D_0\vHat_{N}^{n+1} = 0,\smallskip\\
6.\quad D_+D_- \pHat_j^{n+1} = 0, \quad j=1,2,\ldots,N-1, \qquad  D_{yh} \pHat_0^{n+1}= -\rho a_v^\psup, \quad \pHat_N^{n+1}=\pHat_H(t\sp{n+1}), \smallskip\\
7.\quad a_v^{n+1} =  (-L\pHat_0^{n+1} +\gvnp)/\mb , \qquad  \vb^{n+1} = \vb^{n} + \frac{\dt}{2}( a_v^{n+1} + a_v^n ),
\end{array}
\]
\label{alg:tpRBam}
\end{algorithm}

The first two steps in Algorithm~\ref{alg:tpRBam} result in preliminary
values for the velocity and acceleration of the rigid body as well as a predicted value for the fluid velocity. Of primary interest is the extrapolated value for the vertical component of the acceleration given by
\begin{equation}
\abv^\esup=2\abv^{n}-\abv^{n-1}.
\label{eq:abvp}
\end{equation}
This extrapolated value for the acceleration appears in the boundary condition of the elliptic problem for $\pHat_j\sp\psup$ in Step~3, whose solution is linear in $y$ and takes the form
\begin{equation}
\pHat_j\sp\psup=-\rho \abv^\esup\bigl(y_j-H\bigr),\qquad j=0,1,\ldots,N,
\label{eq:phattp}
\end{equation}
where $\abv^\esup$ is given in~\eqref{eq:abvp}.  The solution for the pressure perturbation in~\eqref{eq:phattp} evaluated at $y_0=0$ is used in Step~4 to compute the vertical component of the acceleration of the body.  With the forcing set to zero, we find
\begin{equation}
\abv^\psup=-{L\pHat_0\sp\psup\over \mb}=-\left({\rho LH\over\mb}\right)\abv\sp\esup,
\label{eq:tpaccelupdate}
\end{equation}
where \eqref{eq:phattp} has been used to eliminate $\pHat_0\sp\psup$.  The velocity of the body, $v_b\sp\psup$, is also computed at this step of the algorithm, but its value does not play an important role in the stability analysis.

The final group of correction steps in Algorithm~\ref{alg:tpRBam} consists of a fluid velocity correction, Step~5, a second fluid pressure update, Step~6, and a final correction of the vertical components of the rigid-body velocity and acceleration, Step~7.  Note that since $\pHat_j\sp\psup$ is linear in $y_j$, the pressure gradient terms in the discrete momentum equation in Step~5 is constant in~$j$ so that this step reduces to setting the (spatially uniform) fluid velocity, $\vHat_j\sp{n+1}$, equal to the body velocity, $v_b\sp\psup$, computed in Step~4.   The solution for the discrete fluid pressure, $\pHat_j\sp{n+1}$, in Step~6 gives the solution in~\eqref{eq:phattp} but with $\abv^\esup$ replaced by the predicted value, $\abv^\psup$, given in~\eqref{eq:tpaccelupdate}.  Using this solution for the pressure and the value for $\abv^\esup$ in~\eqref{eq:abvp}, the vertical component of the acceleration of the body at $t\sp{n+1}$ becomes
\begin{equation}
\abv^{n+1}={2\abv^{n}-\abv^{n-1}\over M_r\sp2},
\label{eq:tpaccel}
\end{equation}
where
\begin{equation}
M_r={\mb\over\rho LH}={\mb\over M_a},
\label{eq:mr}
\end{equation}
is the ratio of the mass of the rigid body to the added mass, $M_a$.  The vertical component of the velocity of the body is computed at this final step to be
\begin{equation}
\vb^{n+1}=\vb\sp{n}+{\dt\over2}\bigl(\abv\sp{n+1}+\abv\sp{n}\bigr).
\label{eq:tpvel}
\end{equation}

The behaviour of the perturbations in the rigid-body acceleration and velocity for the \tpRB~algorithm applied to the model problem, MP-MA, are given by~\eqref{eq:tpaccel} and~\eqref{eq:tpvel}, respectively, with the perturbations in the fluid velocity and pressure determined by the rigid-body perturbations through the velocity equilibration in Step~5 and the pressure solve in Step~6.  Thus, the stability of the algorithm is determined by~\eqref{eq:tpaccel} and~\eqref{eq:tpvel} alone.  Solutions of these two equations have the form $\abv^{n}=A\sp{n}\abv\sp0$ and $\vb\sp{n}=A\sp{n}\vb\sp0$, where the amplification factor, $A$, solves the quadratic
\[
M_r\sp2A\sp2-2A+1=0,
\]
whose roots are given by
\[
A={1\over M_r\sp2}\left[\,1\pm\sqrt{1-M_r\sp2}\;\right].
\]
Perturbations to the components of the solution given by the time-stepping described by Algorithm~\ref{alg:tpRBam} are bounded if $\vert A\vert<1$, and this occurs when $M_r>1$, i.e.~if the mass of the rigid body is larger than the added mass. 
This result is encapsulated in the following theorem:
\begin{theorem} The traditional-partitioned (\tpRB) algorithm given in Algorithm~\ref{alg:tpRBam}  for the added-mass
model problem~{\amplanar} is stable if and only if the mass of the body is greater than the total mass of the fluid
\begin{align*}
    M_r={\mb\over\rho LH}={\mb\over M_a} > 1. 
\end{align*}
\end{theorem}

\subsubsection{Analysis of the \ampRB~algorithm for {\amplanar}}

In order to show that the added-mass instability found in the traditional approach is suppressed in
the AMP scheme, consider the \ampRB~algorithm outlined in Algorithm~\ref{alg:ampRB}.  For the
{\amplanar} model problem, the fully discretized scheme reduces to the one described in detail in
Algorithm~\ref{alg:ampRBam}.  In Steps~3 and~6 of Algorithm~\ref{alg:ampRBam}
the calculation of the body acceleration is coupled
to the elliptic equation for the fluid pressure through the added-mass interface conditions.  These
are the key steps which suppress the added-mass instability.

\begin{algorithm}\caption{\ampRB~scheme for {\amplanar} 
                         \\ \protect\hphantom{Algorithm 5~~} {\em Added-mass model problem in a rectangular geometry.}}\small
\[
\begin{array}{l}
\hbox{//\textsl{{\bodyStepIComment}}}\smallskip\\
1.\quad a_v^\esup = 2 a_v^n - a_v^{n-1} , \qquad \vb^\esup = \vb^{n-1} + 2\dt\, a_v^n,\bigskip\\
\hbox{// \textsl{Prediction steps}}\smallskip\\
2.\quad \vHat_j\sp\psup = \vHat_j^n - \frac{3\dt}{2\rho} D_0 \pHat_{j}^n +\frac{\dt}{2\rho} D_0 \pHat_j^{n-1}, \quad j=1,2,\ldots,N-1, \qquad \vHat_0^\psup = \vb^\esup, \quad D_0\vHat_{N}^\psup = 0, \smallskip\\
3.\quad D_+D_- \pHat_j^\psup = 0, \quad j=1,2,\ldots,N-1, \qquad
\left\{\begin{array}{l} D_{yh} \pHat_0^\psup + \rho a_v^\psup=0 \smallskip\\ \mb a_v^\psup + L\pHat_0^\psup = \gvnp \end{array}\right., \quad \pHat_N^\psup=\pHat_H(t\sp{n+1}), \smallskip\\
4.\quad \vb^\psup = \vb^{n} + \frac{\dt}{2}( a_v^\psup + a_v^n ),\bigskip\\
\hbox{// \textsl{Correction steps}}\smallskip\\
5.\quad \vHat_j\sp{n+1} = \vHat_j^n - \frac{\dt}{2\rho}( D_0 \pHat_{j}^\psup + D_0 \pHat_j^{n}), \quad j=1,2,\ldots,N-1, \qquad \vHat_0^{n+1} = \vb^\psup, \quad D_0\vHat_{N}^{n+1} = 0,\smallskip\\
6.\quad D_+D_- \pHat_j^{n+1} = 0, \quad j=1,2,\ldots,N-1, \qquad  \left\{\begin{array}{l} D_{yh} \pHat_0^{n+1} + \rho a_v^{n+1}=0 \smallskip\\ \mb a_v^{n+1} + L\pHat_0^{n+1} = \gvnp \end{array}\right., \quad \pHat_N^{n+1}=\pHat_H(t\sp{n+1}), \smallskip\\
7.\quad \vb^{n+1} = \vb^{n} + \frac{\dt}{2}( a_v^{n+1} + a_v^n ).
\end{array}
\]
\label{alg:ampRBam}
\end{algorithm}

The first two steps in Algorithm~\ref{alg:ampRBam} result in preliminary values for the velocity and
acceleration of the rigid body as well as a predicted value for the discrete grid function holding the
fluid velocity.  However, these values are not used in subsequent steps of the algorithm for
the {\amplanar} model problem, and so the first important calculation occurs in Step~3.  The
solution of the discrete elliptic equation defined at this step is given by
\[
\pHat_j\sp{\psup}=\pHat(y_j,t\sp{n+1}),\qquad \abv^{\psup}=\abv(t\sp{n+1}),
\]
where $\pHat(y,t)$ and $\abv(t)$ are the exact fluid pressure and vertical body acceleration given in~\eqref{eq:pHatsoln} and~\eqref{eq:absoln}, respectively.  The predicted value for the vertical velocity of the rigid body in Step~4 is a trapezoidal rule quadrature of the body accelerations at $t\sp{n}$ and $t\sp{n+1}$, and this value sets the (spatially uniform) fluid velocity in Step~5 to be
\begin{equation}
\vHat_j\sp{n+1}=v_b\sp{\psup}=v_b\sp{n}+{\dt\over2}\bigl(\abv(t\sp{n+1})+\abv(t\sp{n})\bigr).
\label{eq:vhatamp}
\end{equation}
Step~6 repeats the calculation in Step~3 to give
\begin{equation}
\pHat_j\sp{n+1}=\pHat(y_j,t\sp{n+1}),\qquad \abv^{n+1}=\abv(t\sp{n+1}),
\label{eq:exactdiscretesolutions}
\end{equation}
and thus the result of Step~7 is
\begin{equation}
v_b\sp{n+1}=v_b\sp{n}+{\dt\over2}\bigl(\abv(t\sp{n+1})+\abv(t\sp{n})\bigr).
\label{eq:vbamp}
\end{equation}
Assuming that $ m_b+M_a$ is bounded away from zero, 
the discrete solutions for pressure and acceleration in~\eqref{eq:exactdiscretesolutions} are
readily bounded by the forcing function, $g_v(t)-Lp_H(t)$.  Also, since the discrete solutions for the
velocities in~\eqref{eq:vhatamp} and~\eqref{eq:vbamp} are simply trapezoidal rule quadratures of the exact
solution, we conclude that the discrete solution of the {\amplanar} model problem given by the
\ampRB~algorithm is unconditionally stable, spatially exact, and second-order accurate in time.
This result is encapsulated in the following theorem:
\begin{theorem} The \ampRB~algorithm given in Algorithm~\ref{alg:ampRBam} for the model problem~{\amplanar}
is unconditionally stable provided there exists a constant $K>0$ such that 
\begin{align*}
    m_b+M_a\ge K.
\end{align*}
\end{theorem}
Note that this result is valid independent of the optional velocity correction (VC) step of the \ampRB~algorithm 
which plays no important role for the added-mass case.

\subsection{Analysis of a model problem for added damping} \label{sec:adSubmodel}

We derive a model problem for added-damping and then 
derive an approximation to the added-damping coefficient~\eqref{eq:addedDampingDu} by considering 
the solution to a discrete variational problem. 
With an approximation of this coefficient in hand, we then perform a stability analysis of the \ampRB~algorithm and compare the results to that for the \tpRB~algorithm.

\subsubsection{Added-damping model problem {\adplanar}} \label{sec:mpad}

The added-damping model is obtained from the $k=0$ mode in the Fourier expansion of the fluid variables together
with the equation for the horizontal component of the rigid body.  This results in the {\adplanar} model given by
\begin{equation}
 \text{{\adplanar}:}\;  \left\{ 
   \begin{alignedat}{3}
  & \rho \frac{\partial \uHat}{\partial t} = \mu \frac{\partial^2\uHat}{\partial y^2} , \quad&& y\in(0,H), \\
  &   m_b \frac{d u_b}{dt} =  \mu L \, \frac{\partial\uHat}{\partial y}(0,t) + g_u(t), \\
  &  \uHat(0, t) = u_b(t), \quad \uHat(H, t) = \uHat_H(t),
   \end{alignedat}  \right. 
  \label{eq:MP-AD}
\end{equation}
where the subscript on the coefficient function for the horizontal component of the fluid velocity corresponding to $k=0$ has been dropped for notational convenience.  

\subsubsection{Derivation of the added-damping coefficient} \label{sec:addedDampingCoefficient}

In order to derive an approximation of the added-damping coefficient~\eqref{eq:addedDampingDu}
which appears in the interface condition~\eqref{eq:addedDampingCondition}, consider an implicit (monolithic) discretization of the {\adplanar} model problem of the form
\begin{align}
 & \rho \frac{\uHat_j - \uHat_j^n}{\dt} = 
      \mu \left[ \alpha \Dp\Dm \uHat_j + (1-\alpha) \Dp\Dm \uHat_j^{n}\right], \quad j=1,2,\ldots N-1, \label{eq:AD1}  \\
 &  m_b \ab = \mu L \, \Dyh \uHat_0 + g_u(t^{n+1}) , \label{eq:ADab} \\
 &  \uHat_0 = \ub, \qquad \uHat_N = \uHat_H(t\sp{n+1}), \label{eq:AD2}
\end{align}
where $\uHat_j$, $\ub$ and $\ab$ denote discrete approximations at $t\sp{n+1}$.  It is convenient to introduce the acceleration of the body, $\ab$, in~\eqref{eq:ADab} and relate it to the velocity of the body by
\begin{equation}
\frac{\ub - \ub^n}{\dt} = \alphas \ab + (1-\alphas) \ab^n.
\label{eq:ADvb}
\end{equation}
The divided difference operators, $\Dp$, $\Dm$ and $\Dyh$, appearing in~\eqref{eq:AD1} and~\eqref{eq:ADab} are the same as those defined in~\eqref{eq:ddopts} and~\eqref{eq:funkyD}, although it is informative to also consider a first-order approximation for which $\Dyh=\Dp$.  The constants $\alpha$ and $\alphas$ are included in~\eqref{eq:AD1} and~\eqref{eq:ADvb} to allow some flexibility in the choice of the time-stepping scheme.  Typical choices for the implicit scheme are $\alpha=\alphas=1/2$ corresponding to a second-order accurate trapezoidal method and $\alpha=\alphas=1$ for a first-order accurate backward Euler method.

We recall that the added-damping coefficient, $\Duu$, is related to the derivative of the integral of the shear stress on the body, $\Ic_u$, in~\eqref{eq:shearIntegral} with respect to the horizontal component of the velocity of the body, $\ub$.  In the context of the discretization of the model problem in~\eqref{eq:AD1}--\eqref{eq:ADvb}, the added-damping coefficient is given by
\begin{equation}
\Duu=-\mu L\Dyh \wHat_0,
\label{eq:duuapprox}
\end{equation}
where
\[
\wHat_j={\partial\uHat_j\over\partial\ub},\qquad j=0,1,\ldots,N.
\]
The grid function $\wHat_j$ is obtained from the solution of the discrete variational problem given by
\begin{align}
 & {\rho\over\dt}\wHat_j = \mu \alpha \Dp\Dm \wHat_j, \quad j=1,2,\ldots N-1, \label{eq:ADw}  \\
 &  \wHat_0 = 1, \qquad \wHat_N = 0, \label{eq:ADwbcs}
\end{align}
which is derived from~\eqref{eq:AD1} and~\eqref{eq:AD2} upon differentiation with respect to $\ub$.  The discrete Helmholtz equation in~\eqref{eq:ADw} can be written in the form
\begin{equation}
\wHat_{j+1} - 2 \wHat_j + \wHat_{j-1} = \delta^2 \wHat_j,
\label{eq:ADw2}
\end{equation}
where
\begin{equation}
\delta\eqdef{\dy\over\sqrt{\alpha\nu\dt}},
\label{eq:delta}
\end{equation}
is a dimensionless parameter defined by the ratio of the grid spacing in the vertical direction (normal to the interface) and a viscous length $\sqrt{\alpha\nu\dt}$ determined by the time step.  The solution of the difference equation in~\eqref{eq:ADw2} with boundary conditions in~\eqref{eq:ADwbcs} can be written in the form
\begin{equation}
\wHat_j = \frac{\eta^j - \eta^{2 N-j}}{1-\eta^{2N}},
\label{eq:wHatSoln}
\end{equation}
where
\begin{equation}
\eta = 1 + \delta^2/2 - \delta \sqrt{1 + \delta^2/4} = \frac{1}{1 + \delta^2/2 + \delta \sqrt{1 + \delta^2/4}}<1.
\label{eq:eta}
\end{equation}

The quantity $\Duu$ given by~\eqref{eq:duuapprox} can be obtained using various choices of the divided difference operator $\Dyh$ and the solution for $\wHat_j$ given in~\eqref{eq:wHatSoln}.  For example, if $\Dyh = \Dp$, then
\begin{equation}
\Duu_1=-\mu L\Dp\wHat_0 = \mu L\frac{(1-\eta)}{\dy} ~ \frac{(1+\eta^{2N-1})}{(1-\eta^{2N})} . \label{eq:addDampingCoeff1}
\end{equation}
The terms involving $\eta^{2N}$ arise from the boundary condition at $y=H$, and assuming $\eta^{2N}\ll1$ 
is usually a reasonable approximation\footnote{The assumption $\eta^{2N}\ll1$ is usually reasonable 
unless the fluid layer is very thin, i.e.~of the size of the viscous length $\sqrt{\alpha\nu\dt}$, or $N$ is small.
When $\eta^{2N}\ll1$ is not a good approximation, equation~\eqref{eq:addDampingCoeff1}
indicates that a larger value for the added-damping coefficient would result.}.
Neglecting the terms involving $\eta^{2N}$ gives the approximation
\begin{equation}
\Duu_1\approx \mu L \, \frac{1-\eta}{\dy}.
\label{eq:addDampingCoeff1Approx}
\end{equation}
Similarly if $\Dyh$ is taken to be the second-order accurate one-sided operator defined in~\eqref{eq:funkyD} and if terms involving $\eta^{2N}$ are neglected as before, then
\begin{equation}
\Duu_2\approx \mu L\frac{ 3 - 4 \eta + \eta^2}{2 \dy} .
\label{eq:addDampingCoeff2Approx}
\end{equation}
It is therefore clear that different approximations for $\Dyh$ lead to different approximations to $\Duu$.
In practice it is apparently sufficient to use a simple and low-order approximation to $\Duu$ since the~\ampRB~algorithm
is designed so that the scheme is second-order accurate for any choice of $\Duu$. To this end we 
base our added-damping coefficient $\Duu$ on the first-order accurate approximation~\eqref{eq:addDampingCoeff1Approx}.
Furthermore, notice that when the viscous terms are treated implicitly, it is typical for $\delta \le 1$ since 
the time-step $\dt$ will normally be chosen to 
exceed the usual explicit time-step restriction by some factor $C_\nu$ larger than one (e.g. 10 or 100). Thus in the typical case 
\[ 
  \frac{\alpha \nu\dt}{\dy^2} = \frac{1}{\delta^2} > C_\nu >1  ~\implies  \delta < \frac{1}{\sqrt{C_\nu}} < 1 .
\]
Note that for $\delta$ small, $\eta$ is a second-order accurate approximation to $e^{-\delta}$, while for large $\delta$
both $\eta$ and $e^{-\delta}$ decay to zero. As a result, the sensible
approximation $\eta\approx e^{-\delta}$ leads to the following definition:
%
\begin{definition} The approximate added-damping coefficient $\Dun$ and {\em added-damping length-scale} $\dn$ are
defined by 
\begin{align}
\Dun &\eqdef \mu L \, \frac{1- e^{-\delta}}{\dy} =    \mu L \frac{1}{\dn} , \label{eq:addDampingCoeffn} \\
   \dn &\eqdef \frac{\dy}{1- e^{-\delta}}.\label{eq:dnDef}
\end{align}
\end{definition}
The {\em added-damping length-scale} $\dn$ is a distance in the normal direction 
to the boundary with limiting behaviours
\[
  \dn \sim \begin{cases}
                 \sqrt{\alpha \nu \dt} & \text{for $\delta\rightarrow 0$} , \\
                 \dy           & \text{for $\delta\rightarrow \infty$} .
           \end{cases}
\]
Thus the added-damping length-scale $\dn$ varies in size between 
the viscous length $\sqrt{\alpha\nu\dt}$ and 
the grid-spacing $\dy$ normal to the boundary 
as $\delta$ varies between zero and infinity. 
It is important that $\dn$ has approximately 
the correct behaviour as $\delta$ varies since this makes it possible to 
choose the coefficient of the added-damping once, without requiring it to vary 
as a function of other problem parameters, see Section~\ref{sec:stabilityMP-AD}.

\subsubsection{Analysis of the {\tpRB} and {\ampRB}~algorithms for {\adplanar}} \label{sec:stabilityMP-AD}

Consider now the stability of the {\tpRB} and {\ampRB}~algorithms given in
Algorithms~\ref{alg:tpRB} and~\ref{alg:ampRB}, respectively, for the {\adplanar} model problem.
Algorithm~\ref{alg:ampRBad} gives the specific steps for both algorithms, where the {\tpRB}~scheme
is obtained for the choice $\beta_d=0$ (so that the added-damping term vanishes) and the
{\ampRB}~scheme corresponds to the choice $\beta_d>0$.  The added-damping coefficient $\Dun$ 
appearing in Step 3 of the algorithm is defined in~\eqref{eq:addDampingCoeffn}.  
The horizontal component of the fluid velocity satisfies a discrete diffusion equation
in Steps~2 and~5 of the algorithm (and optionally in Step~8). For convenience the boundary at
$y=H$ is taken to be at infinity where the solution is assumed to be bounded. This choice is made to simplify the algebraic
details of the analysis somewhat, but it has little effect on the main results assuming $H$ is not
small in relation to the grid spacing $\dy$.
 As before, we consider the homogeneous case with $g_u(t)$ set to zero, and regard the
discrete variables as perturbations to the exact solution of the discrete equations.  The stability
analysis considers the growth or decay of these perturbations.  We proceed with the analysis by
first assuming that $\beta_d>0$, and then consider the limiting case $\beta_d=0$ in the end.

\begin{algorithm}\caption{\tpRB~scheme ($\beta_d=0$) and \ampRB~scheme ($\beta_d>0$) for {\adplanar} 
                         \\ \protect\hphantom{Algorithm 5~~} {\em Added-damping model problem in a rectangular geometry.}}\small
\medskip
\[
\begin{array}{l}
\hbox{// \textsl{{\bodyStepIComment}}}\smallskip\\
1.\quad a_u^\esup = 2 a_u^n - a_u^{n-1} , \qquad \ub^\esup = \ub^{n-1} + 2\dt\, a_u^n,\bigskip\\
\hbox{// \textsl{Prediction steps}}\smallskip\\
2.\quad \uHat_j\sp\psup = \uHat_j^n + \frac{\nu\dt}{2}\bigl(D_+D_-\uHat_j\sp\psup+D_+D_-\uHat_j\sp{n}\bigr), \quad j=1,2,\ldots, \qquad \uHat_0^\psup = \ub^\esup, \quad \displaystyle{\lim_{j\rightarrow\infty}\vert\uHat_j\sp\psup\vert<\infty}, \smallskip\\
3.\quad (m_b+\beta_d\dt\Dun)  a_u^\psup=\mu L\bigl(D_{yh}\uHat_0\sp\psup\bigr)+\beta_d\dt\Dun a_u\sp\esup+\gunp, \medskip\\
4.\quad \ub^\psup = \ub^{n} + \frac{\dt}{2}( a_u^\psup + a_u^n ),\bigskip\\
\hbox{// \textsl{Correction steps}}\smallskip\\
5.\quad \uHat_j\sp{n+1} = \uHat_j^n + \frac{\nu\dt}{2}\bigl(D_+D_-\uHat_j\sp{n+1}+D_+D_-\uHat_j\sp{n}\bigr), \quad j=1,2,\ldots, \qquad \uHat_0^{n+1} = \ub^\psup, \quad \displaystyle{\lim_{j\rightarrow\infty}\vert\uHat_j\sp{n+1}\vert<\infty}, \smallskip\\
6.\quad (m_b+\beta_d\dt\Dun)  a_u^{n+1}=\mu L\bigl(D_{yh}\uHat_0\sp{n+1}\bigr)+\beta_d\dt\Dun\, a_u\sp\psup+\gunp, \medskip\\
7.\quad \ub^{n+1} = \ub^{n} + \frac{\dt}{2}( a_u^{n+1} + a_u^n ), \bigskip\\
\hbox{// \textsl{Fluid-velocity correction step (optional)}}\smallskip\\
8.\quad \uHat_j\sp{n+1} = \uHat_j^n + \frac{\nu\dt}{2}\bigl(D_+D_-\uHat_j\sp{n+1}+D_+D_-\uHat_j\sp{n}\bigr), \quad j=1,2,\ldots, \qquad \uHat_0^{n+1} = \ub^{n+1}, \quad \displaystyle{\lim_{j\rightarrow\infty}\vert\uHat_j\sp{n+1}\vert<\infty},
\end{array}
\]
\label{alg:ampRBad}
\end{algorithm}

To address the stability of the schemes we proceed along the classical lines of GKS stability theory~\cite{GKSII}. 
As such, consider solutions to the homogeneous equations of the form
\begin{equation}
  \uHat_j\sp{n}=A\sp n \, \uHat_j\sp{0},\qquad \ub\sp{n}=A\sp n \, \ub\sp{0},\qquad a_u\sp{n}=A\sp n \, a_u\sp{0},
  \label{eq:uamps}
\end{equation}
where $A$ is an amplification factor. It is also convenient to define 
\begin{equation}
  \uHat_j\sp\psup=A\sp{n+1}\bar u_j\sp\psup, \qquad \ub^\ksup=A\sp{n+1}\bar u_b^\ksup, \qquad a_u^\ksup=A\sp{n+1}\bar a_u^\ksup,\qquad \hbox{$k=e$ or $p$,}
  \label{eq:intermediate}
\end{equation}
for the intermediate values in Algorithm~\ref{alg:ampRBad}. 
 We now seek solutions of the equations with $\vert A\vert>1$ so that if no solutions of this form
can be found, then the time-stepping scheme is stable.  
Substitution of~\eqref{eq:uamps} and~\eqref{eq:intermediate} into Step~5 yields the following difference equation 
\begin{equation}
  A\uHat_j\sp{0}=\uHat_j\sp{0}+{\nu\dt\over2}(A+1)D_+D_-\uHat_j\sp{0},\quad j=1,2,\ldots,\qquad \uHat_0\sp{0}=\bar u_b\sp\psup,\quad \lim_{j\rightarrow\infty}\vert \uHat_j\sp{0}\vert<\infty,
  \label{eq:uhat0}
\end{equation}
with solutions of the form $\uHat_j\sp{0}=K\,\xi\sp{j}$, where $\xi$ is a complex constant. Using the assumption that $\vert A\vert>1$ and the boundary conditions gives 
\begin{equation}
  \uHat_j\sp{0}=\bar u_b\sp\psup\xi\sp{j},
  \label{eq:uhatsoln}
\end{equation}
where $\vert\xi\vert<1$ is given by
\begin{equation}
  \xi=q-\sqrt{q\sp2-1},\qquad q=1+{\delta\sp2\over2}\left({A-1\over A+1}\right),
  \label{eq:xiconst}
\end{equation}
with $\delta$ defined in~\eqref{eq:delta} for $\alpha=1/2$.
Similarly, substitution of~\eqref{eq:uamps} and~\eqref{eq:intermediate} into Step~2 yields the difference equation
\begin{equation}
  A\left(1-{\nu\dt\over2}D_+D_-\right)\bar u_j\sp\psup=\left(1+{\nu\dt\over2}D_+D_-\right)\uHat_j\sp{0},\quad j=1,2,\ldots,\qquad \bar u_0\sp\psup=\bar u_b\sp\esup,\quad \lim_{j\rightarrow\infty}\vert \bar u_j\sp\psup\vert<\infty,
  \label{eq:ubar}
\end{equation}
with $\uHat_j\sp{0}$ given in~\eqref{eq:uhatsoln}. The solution of~\eqref{eq:ubar} has the form
\begin{equation}
  \bar u_j\sp\psup=\bigl(\bar u_b\sp\esup-\bar u_b\sp\psup\bigr) \, \eta\sp{j} \, +\uHat_j\sp{0},
  \label{eq:ubarsoln}
\end{equation}
where $\eta$ is defined as before in~\eqref{eq:eta}. 

The solutions for $\uHat_j\sp{0}$ in~\eqref{eq:uhatsoln} and $\bar u_j\sp\psup$ in~\eqref{eq:ubarsoln} are inserted into Steps~3 and~6 of Algorithm~\ref{alg:ampRBad} to eliminate the fluid velocities in favor of the velocities of the body. The remaining steps of the algorithm are then given in terms of velocities and accelerations of the body alone.  For example, the equations in Steps~6 and~7 with the forms in~\eqref{eq:uamps} and~\eqref{eq:intermediate}, and using the solution in~\eqref{eq:uhatsoln}, can be written as the following system of equations
\begin{equation}
  \left[\begin{array}{cc}
    0&\mb+\beta_d\dt\Dun \medskip\\
    A-1 & -{\dt\over2}(A+1)
  \end{array}\right]
  \left[\begin{array}{c}
    \ub\sp{0} \medskip\\
    a_u\sp{0}
  \end{array}\right]=
  \left[\begin{array}{cc}
    \mu L\bigl.(D_{yh}\xi\sp{j})\bigr\vert_{j=0}&\beta_d\dt\Dun \medskip\\
    0 & 0
  \end{array}\right]
  \left[\begin{array}{c}
    \bar u_b\sp\psup \medskip\\
    \bar a_u\sp\psup
  \end{array}\right].
  \label{eq:sysA}
\end{equation}
 A similar analysis of the equations in Steps~3 and~4 gives
\begin{equation}
  \begin{array}{l} 
    \displaystyle{
    \left[\begin{array}{cc}
    \mu L\bigl[\bigl.(D_{yh}\eta\sp{j})\bigr\vert_{j=0}-\bigl.(D_{yh}\xi\sp{j})\bigr\vert_{j=0}\bigr]&\mb+\beta_d\dt\Dun \medskip\\
    A & -{\dt\over2}A
  \end{array}\right]
  \left[\begin{array}{c}
    \bar u_b\sp\psup \medskip\\
    \bar a_u\sp\psup
  \end{array}\right]
  } \bigskip\\
  \displaystyle{
  \qquad\qquad=\left[\begin{array}{cc}
    \mu L\bigl.(D_{yh}\eta\sp{j})\bigr\vert_{j=0}&\beta_d\dt\Dun \medskip\\
    0 & 0
  \end{array}\right]
  \left[\begin{array}{c}
    \bar u_b\sp\esup \medskip\\
    \bar a_u\sp\esup
  \end{array}\right]+
  \left[\begin{array}{cc}
    0 & 0 \medskip\\
    1 & {\dt\over2}
  \end{array}\right]
  \left[\begin{array}{c}
    \bar u_b\sp{0} \medskip\\
    \bar a_u\sp{0}
  \end{array}\right].
  }
  \end{array}
  \label{eq:sysB}
\end{equation}
Finally, the two equations in Step~1 can be written as the system
\begin{equation}
  \left[\begin{array}{c}
    \bar u_b\sp\esup \medskip\\
    \bar a_u\sp\esup
  \end{array}\right]=
  {1\over A\sp2}\left[\begin{array}{cc}
    1 & 2\dt A \medskip\\
    0 & 2A-1
  \end{array}\right]
  \left[\begin{array}{c}
    \bar u_b\sp{0} \medskip\\
    \bar a_u\sp{0}
  \end{array}\right].
  \label{eq:sysC}
\end{equation}
The intermediates states, $\bigl[\bar u_b\sp\ksup\,,\,\bar a_u\sp\ksup\bigr]\sp{T}$, $k=e$ and $p$, can be eliminated from the linear systems in~\eqref{eq:sysA}, \eqref{eq:sysB} and~\eqref{eq:sysC} to give a homogeneous system of the form
\begin{equation}
  M\left[\begin{array}{c}
    \bar u_b\sp{0} \medskip\\
    \bar a_u\sp{0}
  \end{array}\right]=0,
  \label{eq:sysD}
\end{equation}
where the $2\times2$ matrix $M$ involves the amplification factor $A$ and the various parameters of the problem.  
Nontrivial solutions of~\eqref{eq:sysD} exist if ${\rm det}(M)=0$, which leads to the following theorem:
\begin{theorem}
   The \ampRB~algorithm without the velocity-correction, Step~8,
given in Algorithm~\ref{alg:ampRBad} for the model problem~{\adplanar}, is stable in the 
sense of Godunov and Ryabenkii if and only if there are no roots with $|A|>1$ to the following equation:
\begin{equation}
\NAb(A) \eqdef \gamma_b(A-1)\sp3+\gamma_0\bigl[\mbBar(A-1)+\cxi(A+1)\bigr]A\sp2=0,
\label{eq:adconstraint}
\end{equation}
where
\[
\gamma_b \eqdef (2\beta_d\bar{\cal D}\sp{u}-\cxi)(2\beta_d\bar{\cal D}\sp{u}-\ceta),\qquad \gamma_0 \eqdef \mbBar+4\beta_d\bar{\cal D}\sp{u}-\ceta.
\]
Here, $\bar{\cal D}\sp{u}$ and $\mbBar$ are a dimensionless added-damping coefficient and body mass, respectively, given by
\begin{align}
  & \bar{\cal D}\sp{u} \eqdef {\dy\Dun\over\mu L}=1-e\sp{-\delta}, \\
  &  \mbBar \eqdef \left({m_b\over\rho L\dy}\right)\delta\sp2,  \label{eq:mbBarDef}
\end{align}
and $\cxi$ and $\ceta$ are dimensionless Dirichlet-to-Neumann transfer coefficients given by
\begin{align}
\cxi \eqdef -\dy \bigl.(D_{yh}\xi\sp{j})\bigr\vert_{j=0}={3-4\xi+\xi\sp2\over2},   \qquad
\ceta \eqdef -\dy \bigl.(D_{yh}\eta\sp{j})\bigr\vert_{j=0}={3-4\eta+\eta\sp2\over2},  \label{eq:CetaDef}
\end{align}
where $\xi$ and $\eta$, which depend on $A$ and $\delta$, 
are defined in~\eqref{eq:xiconst} and~\eqref{eq:eta}, respectively.
\end{theorem}

A similar stability analysis can be performed for Algorithm~\ref{alg:ampRBad} 
with the optional fluid-velocity correction step (Step~8) included.  This leads to the following theorem:
\begin{theorem}
   The \ampRB~algorithm with the velocity-correction, Step~8,
given in Algorithm~\ref{alg:ampRBad} for the model problem~{\adplanar}, is stable in the 
sense of Godunov and Ryabenkii if and only if there are no roots with $|A|>1$ to the following equation
\begin{equation}
\NAv(A) \eqdef \gamma_v(A-1)\sp3+\gamma_0\bigl[\bar m_b(A-1)+\cxi(A+1)\bigr]A\sp2=0,
\label{eq:advconstraint}
\end{equation}
where $\gamma_b$ in~\eqref{eq:adconstraint} is replaced by
\[
\gamma_v=(2\beta_d\bar{\cal D}\sp{u}-\ceta)\sp2,
\]
in~\eqref{eq:advconstraint} while the other terms are unchanged.  
\end{theorem}
Stability regions for Algorithm~\ref{alg:ampRBad} with, and without, the velocity-correction step,
are defined as follows.
\begin{definition}
  The stability region, $\SRb=\SRb\{\mbBar,\delta,\adc\}$, for Algorithm~\ref{alg:ampRBad} with no velocity-correction is
defined to be the set of values for the non-dimensional parameters $\mbBar\ge 0$, $\delta>0$ and $\adc\ge0$ where there are no
roots $A$ with $|A|>1$ to $\NAb(A)=0$ given in~\eqref{eq:adconstraint}.
\end{definition}
\medskip
\begin{definition}
  The stability region, $\SRv=\SRv\{\mbBar,\delta,\adc\}$, for Algorithm~\ref{alg:ampRBad} with velocity-correction is
defined to be the set of values for the non-dimensional parameters $\mbBar\ge 0$, $\delta>0$ and $\adc\ge0$ where there are no
roots $A$ with $|A|>1$ to $\NAv(A)=0$ given in~\eqref{eq:advconstraint}.
\end{definition}

Portions of the stability regions $\SRb$ and $\SRv$ are presented in 
Figures~\ref{fig:addedDampingNvcModes} and~\ref{fig:addedDampingVcModes}. Generation of these
plots is nontrivial, and so it is appropriate to say a few words about our procedure for doing so.
Consider first the computation of $\SRb$ (the computation of
$\SRv$ being similar). Note that while $\NAb$ in~\eqref{eq:adconstraint} appears to be a cubic polynomial in~$A$, its
coefficients involve square roots of a rational function of~$A$ through the transfer
coefficients~$\cxi$ and $\ceta$ and the definitions of~$\xi$ and $\eta$. 
 However, manipulation of~\eqref{eq:adconstraint}, including
squaring both sides of the equation, reveals a polynomial in $A$ of degree eight. Note that not all
roots of this eighth-degree polynomial are valid since squaring the equation may introduce spurious roots. 
Nonetheless, all eight roots of this polynomial can be found numerically, for example by determining
the eigenvalues of the corresponding companion matrix, and valid roots can be identified as being those that
satisfy the original equation $\NAb=0$. The scheme is determined to be unstable if any of these valid roots satisfy $|A|>1$,
otherwise it is stable.
This procedure can then be carried out with a fine grid of points in the parameter space $(\mbBar,\delta,\adc)$
to reveal the stability region.

Alternately, the boundaries separating stable and unstable regions can be sought directly by using 
the fact that points on this boundary satisfy $|A|=1$ which suggests the transformation 
$A=e^{i\theta}$ with $\theta\in\Real$. 
Furthermore, note that the equation $\NAb=0$ is a quadratic equation in $\adc$.
Given values for $\theta$, $\mbBar$ and $\delta$
one solves this quadratic equation for $\adc$ to give two values $\adc = \adc^\pm( \theta; \mbBar,\delta).$
Values of $\theta$ where $\adc$ is real denote points on the stability boundary. These values can
be determined as the real roots $\theta$ of the scalar equations $\Bc^\pm(\theta)\eqdef\Imag(\adc^\pm(\theta;\,  \mbBar,\delta))=0$.
Note that the case $A \rightarrow 1$ requires some care and reduces to the condition $\gamma_0=0$.

The stability regions $\SRb$ and $\SRv$ depend on the three parameters $(\mbBar,\delta,\adc)$ where for
clarity we recall the expressions for $\mbBar$ and $\delta$, 
\[
  \mbBar = \frac{m_b}{\rho L\dy}\,\delta^2, \quad \delta = \frac{\dy}{\sqrt{\nu\dt/2}} .
\]
For plotting purposes it is convenient to show cuts through the three-dimensional space by 
plotting two-dimensional contours of the regions by holding one of the three parameters fixed
and letting the other two vary.
The lower plots in Figures~\ref{fig:addedDampingNvcModes} and~\ref{fig:addedDampingVcModes}
show such cuts through the stability regions. The lower left plots show the stable and
unstable regions for a typical value of $\delta=0.5$ while the lower right plots show the regions for 
the difficult case $\mbBar=0$ when added damping effects are greatest.
The stable regions are shown in 
white, while {\em instability regions} are shaded in different colors 
which correspond to unstable modes of different character. Characterization of the various 
unstable modes is useful since it sheds light on the origin 
of instabilities that occur in practice (e.g.~$\adc$ is taken to be too small). In order to aid in understanding
these instabilities, the upper panel of plots in 
Figures~\ref{fig:addedDampingNvcModes} and~\ref{fig:addedDampingVcModes}
show the spatial and temporal evolution of the various unstable modes.

\medskip
{\bf \em Instability region I} corresponds to having provided insufficient stabilization of 
added-damping effects (or none at all 
in the case of the traditional scheme with $\adc=0$), and occurs for very light or massless
bodies (a {\em light} body being characterized by the non-dimensional body mass $\mbBar$).
This mode is characterized by $A$ being real
with $A>1$. The behaviour of the 
unstable mode in time is $A^n$, $n=0,1,2,\ldots$ and in space $\xi^j$, $j=0,1,2,\ldots$.
The real parts of these modes are shown at the top of 
Figures~\ref{fig:addedDampingNvcModes} and~\ref{fig:addedDampingVcModes}, 
the $Re(A^n)$ being scaled by its maximum absolute value for graphical purposes.
The unstable modes in region I grow monotonically and {\em exponentially} fast in time, and decay
monotonically in space from a maximal value on the boundary $j=0$.

{\bf \em Instability region II} also corresponds to having provided insufficient stabilization of 
added-damping effects (or none at all in the case of the traditional scheme with $\adc=0$), 
and exists for slightly larger rigid body masses than those characterized by region I. 
This mode is characterized by a complex root $A$, which implies an oscillatory exponential growth
in time. The spatial behaviour is that of a rapidly decaying boundary layer with mild oscillations. 
Note that this root corresponds to a bifurcation of the real root from region I into a complex conjugate pair.

{\bf \em Instability region III} occurs only for the scheme without velocity correction and 
is especially prevalent for small values of $\delta$ (i.e.~{\em large} implicit time-step).
The instability is characterized by a real root with $A<-1$ corresponding to 
an exponentially growing plus-minus type mode in time, and an exponentially decaying mode in space.

{\bf \em Instability region IV} occurs for the scheme both with and without velocity correction and is associated
with having provided excessive stabilization against added-damping effects. This instability is characterized by a 
complex root for $A$, that for large values of $\adc$ 
exhibits rather slow exponential growth in time. 
This mode can be the most difficult to detect as it may take many time-steps to become dominant.

\medskip

{
\newcommand{\labelFont}{\scriptsize}
\newcommand{\figWidth}{4.0cm}
\newcommand{\trimfig}[2]{\trimw{#1}{#2}{.0}{.0}{.0}{.0}}
\newcommand{\figWidtha}{8.5cm}
\begin{figure}[htb]
\begin{center}
\begin{tikzpicture}[scale=1]
  \useasboundingbox (0.0,1) rectangle (16,10);  
  \begin{scope}[xshift=-.25cm]
    \draw(-.05,0.0) node[anchor=south west,xshift=-4pt,yshift=+0pt] {\trimfig{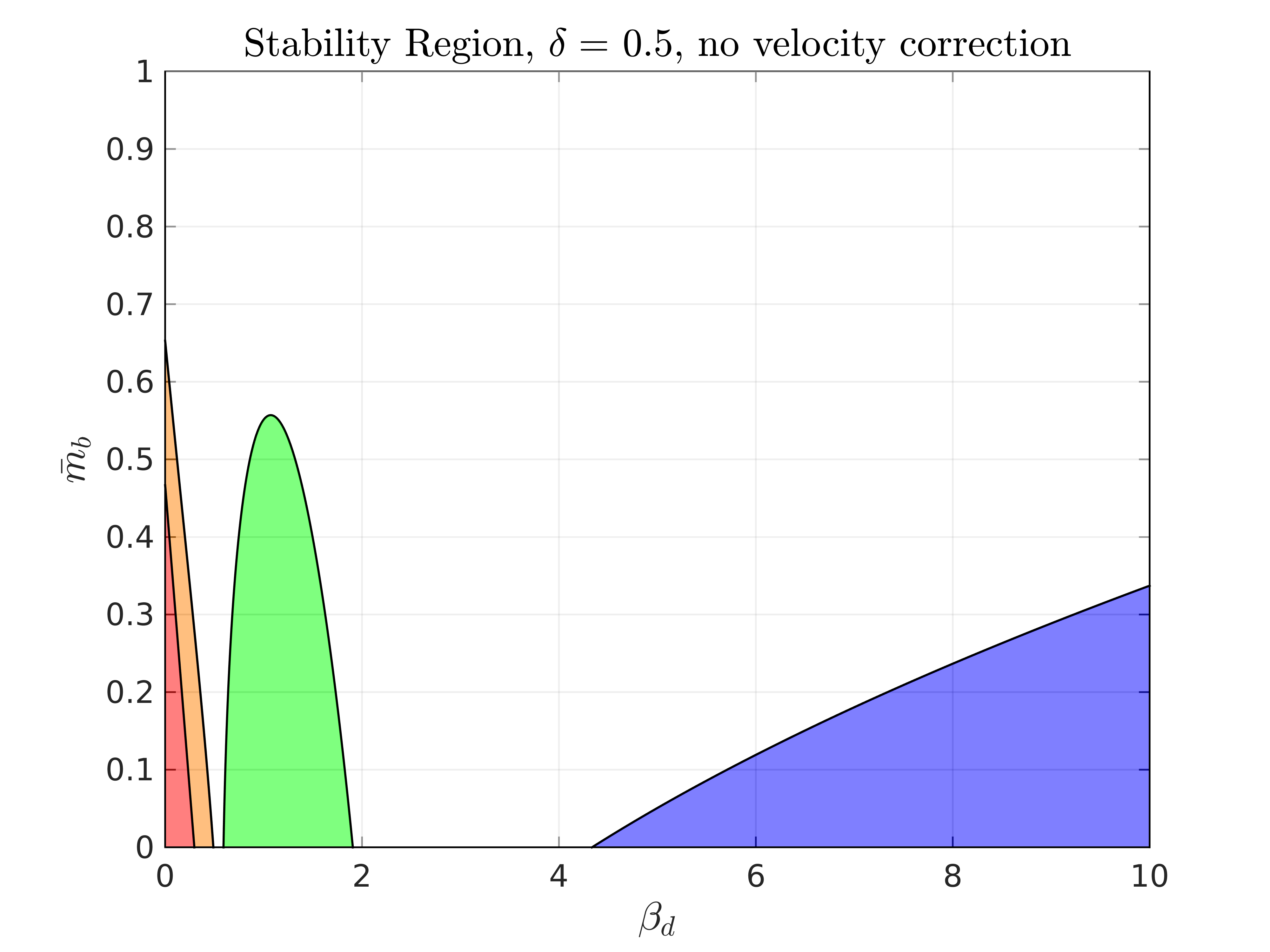}{\figWidtha}};
    \draw(7.75,0.0) node[anchor=south west,xshift=-4pt,yshift=+0pt] {\trimfig{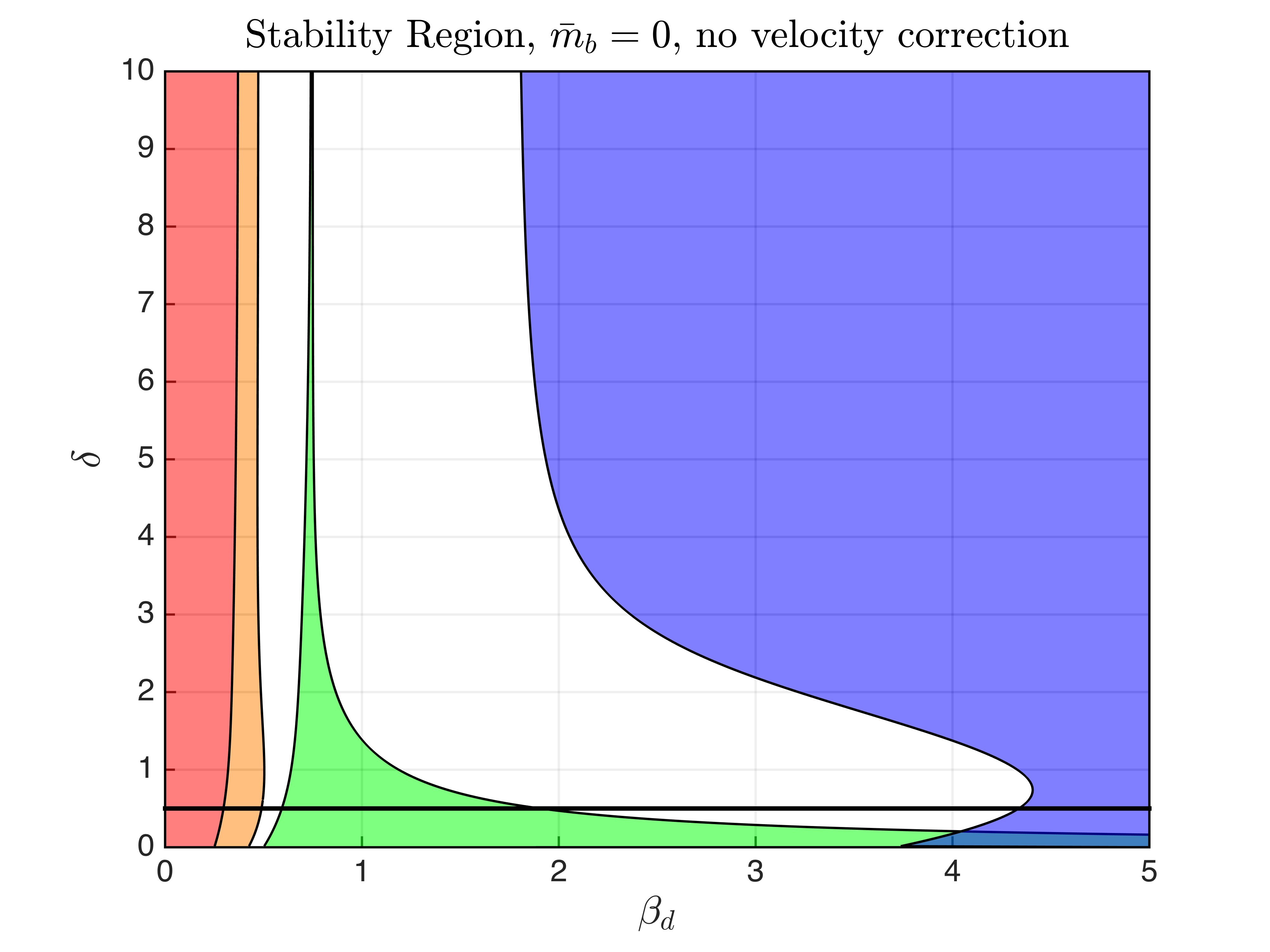}{\figWidtha}};
    \draw(4.5,5.0) node[draw,fill=white] {\labelFont{Stable region}};
    \draw(6.,2.5) node[draw,fill=blue!50] {\labelFont IV};
    \draw(3.2,2.5) node[draw,fill=green!50] {\labelFont III};
    \draw(1.9,4.3) node[draw,fill=orange!50] {\labelFont II};
    \draw(.75,1.5) node[draw,fill=red!50] {\labelFont I};
    \draw(11.9,1.75) node[draw,fill=white] {\labelFont{Stable region}};
  \end{scope}
  \begin{scope}[xshift=-.5cm,yshift=6.5cm]
    \draw(0.0,0.0) node[anchor=south west,xshift=-4pt,yshift=+0pt] {\trimfig{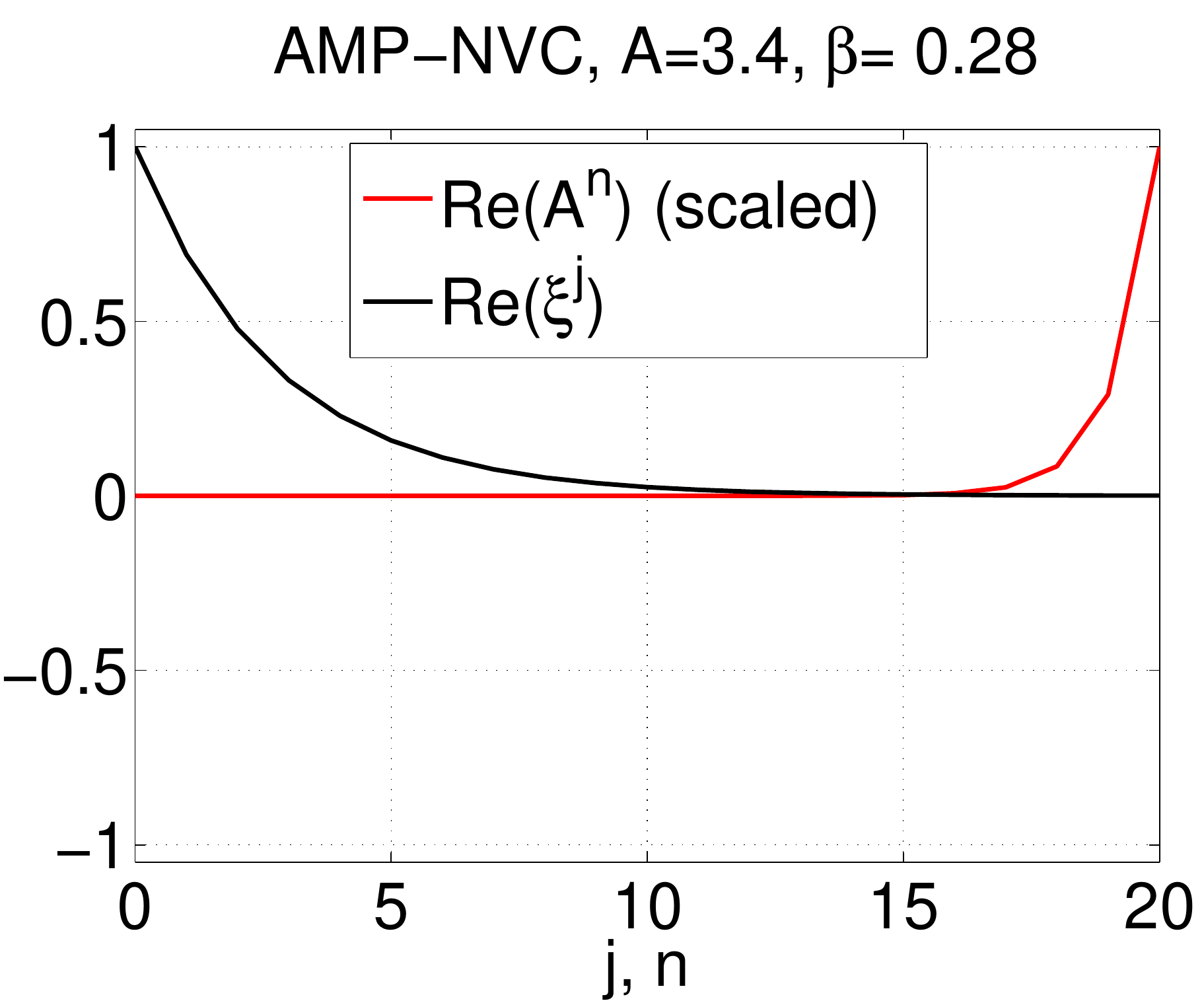}{\figWidth}};
    \draw(2,1.25) node[draw,fill=red!50] {\labelFont I};
    \draw(4.0,0.0) node[anchor=south west,xshift=-4pt,yshift=+0pt] {\trimfig{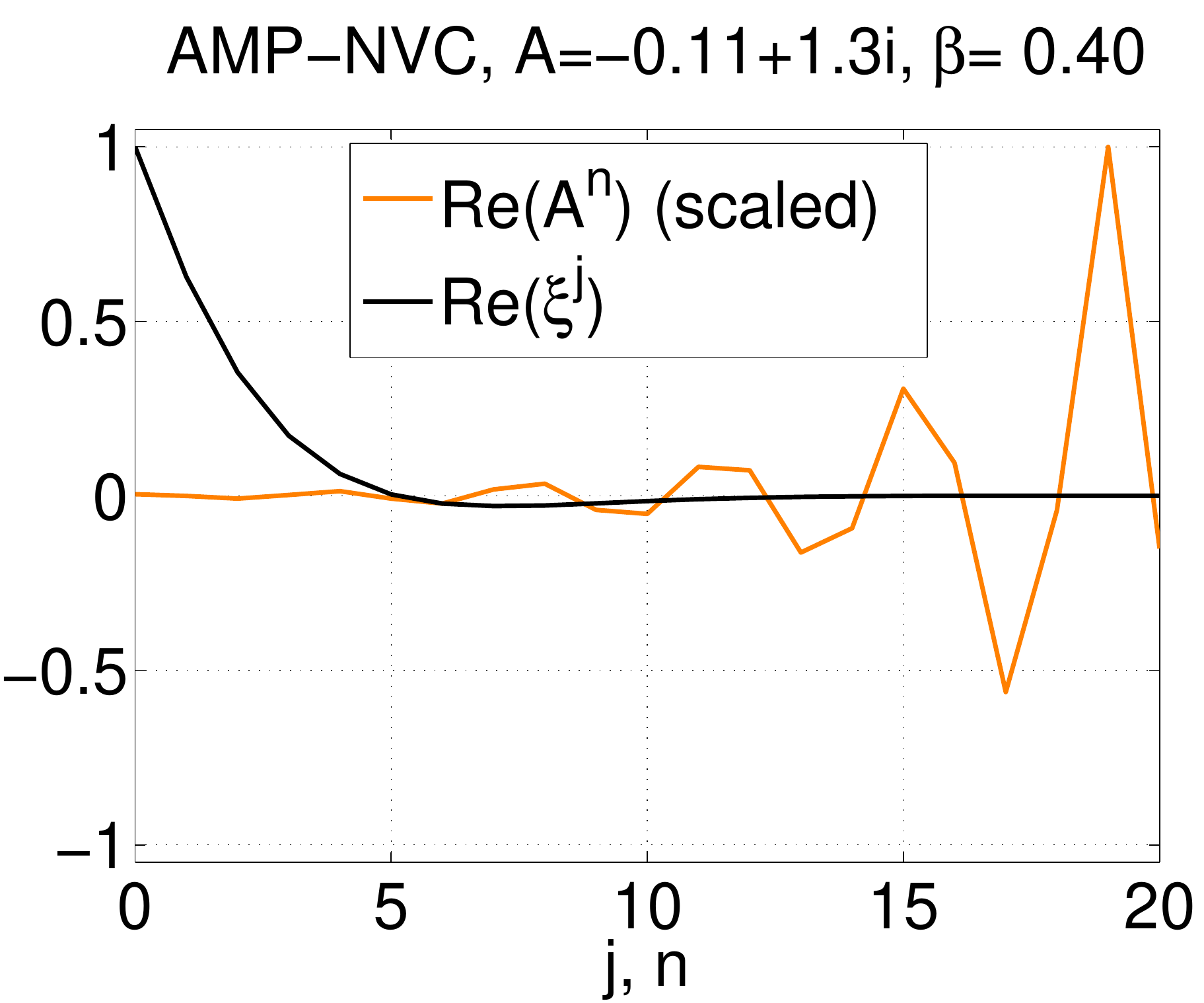}{\figWidth}};
    \draw(6,1.25) node[draw,fill=orange!50] {\labelFont II};
    \draw(8.0,0.0) node[anchor=south west,xshift=-4pt,yshift=+0pt] {\trimfig{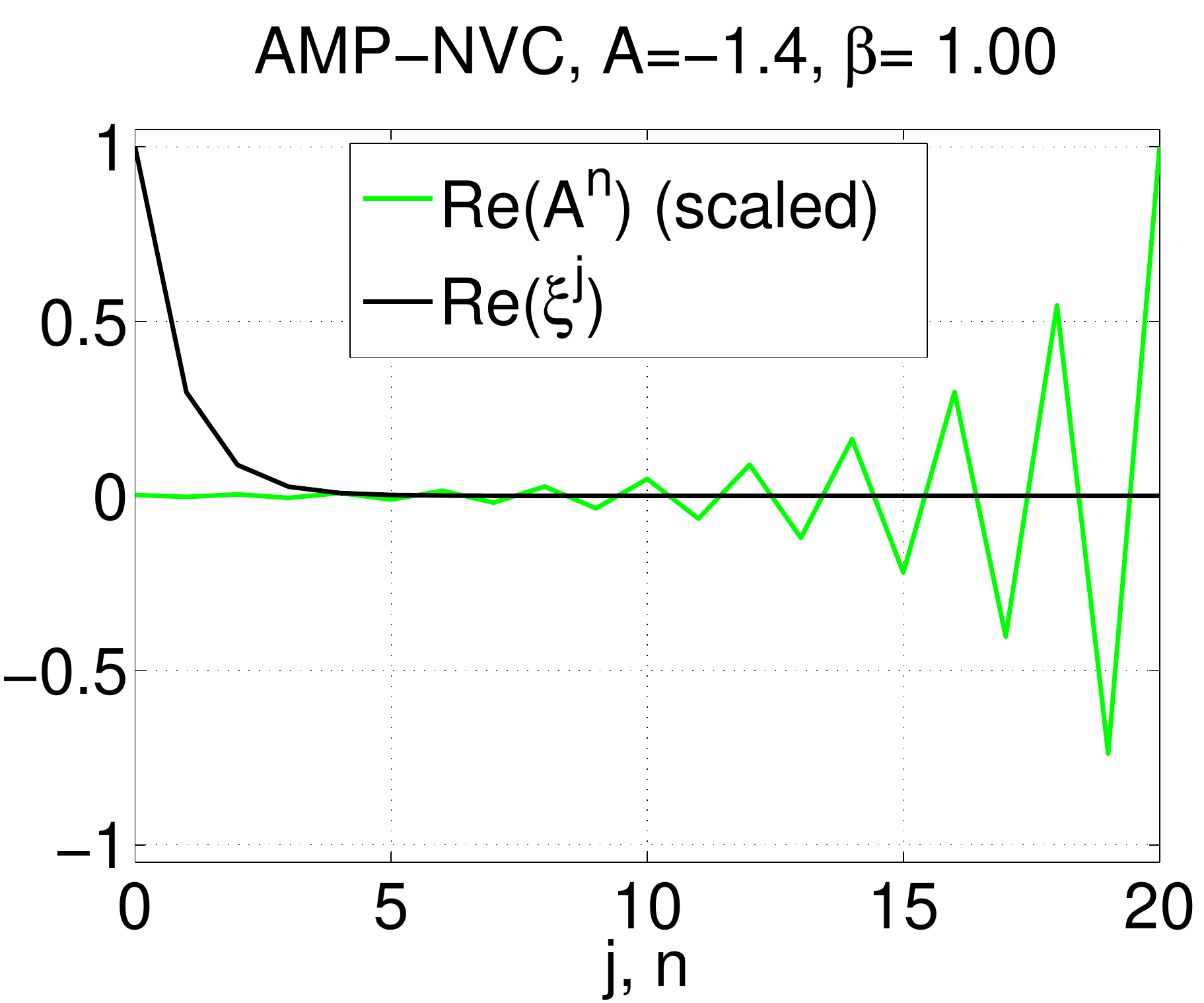}{\figWidth}};
    \draw(10,1.25) node[draw,fill=green!50] {\labelFont III};
    \draw(12.0,0.0) node[anchor=south west,xshift=-4pt,yshift=+0pt] {\trimfig{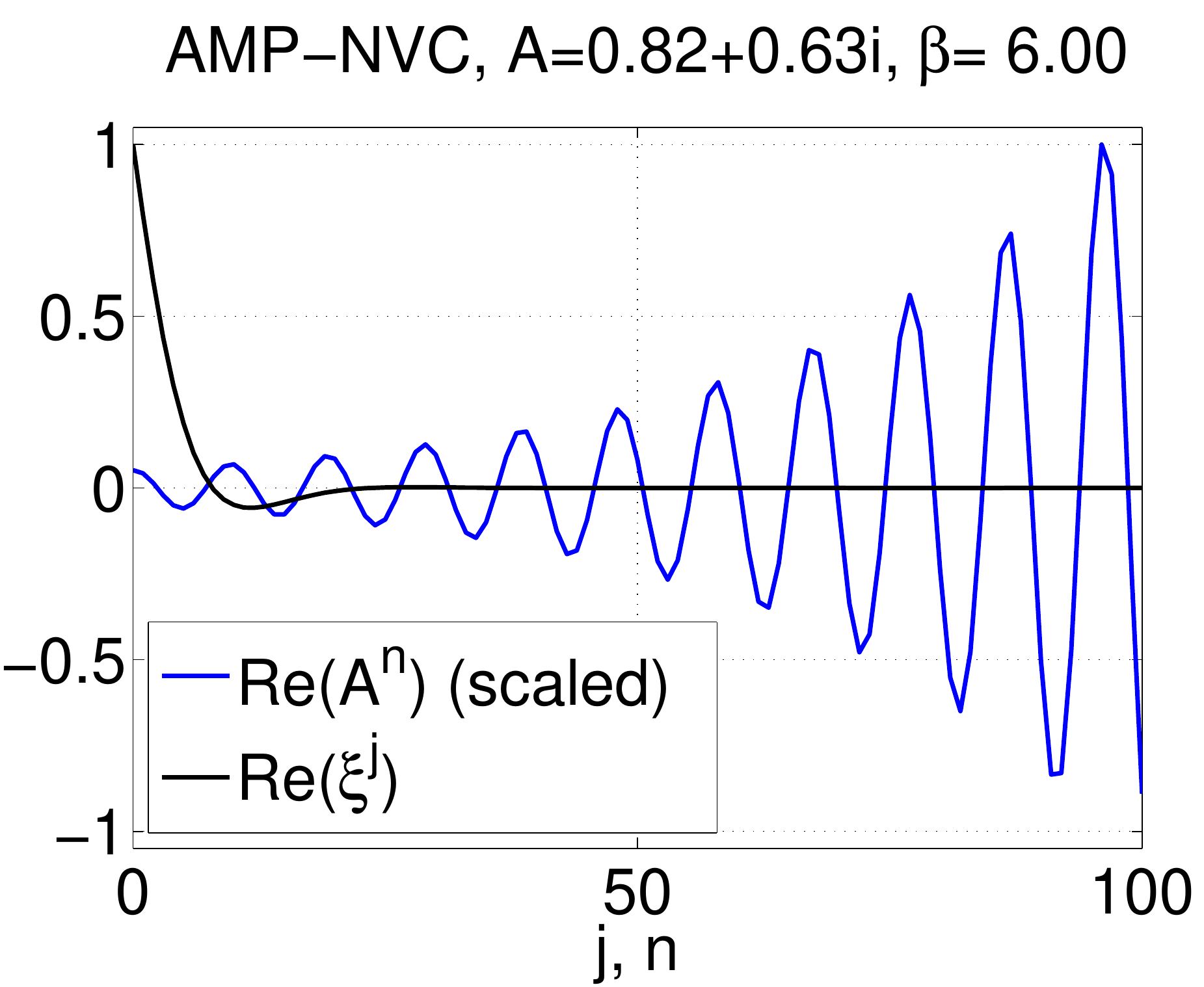}{\figWidth}};
    \draw(14,2.5) node[draw,fill=blue!50] {\labelFont IV};
  \end{scope}
%
\end{tikzpicture}
\end{center}
\caption{Stability regions and unstable modes for the \ampRB-scheme with no velocity correction.
Bottom left: stability region in the $\adc$--$\mbBar$ plane for $\delta=.5$; unstable shaded regions are marked according to the
form of the instability. 
Bottom: right: stability region in the $\adc$--$\delta$ plane for the zero-mass body.
Top: unstable modes in space and time for the four instability types ($\delta=.5$ and $\mbBar=0$).
}
\label{fig:addedDampingNvcModes}
\end{figure}
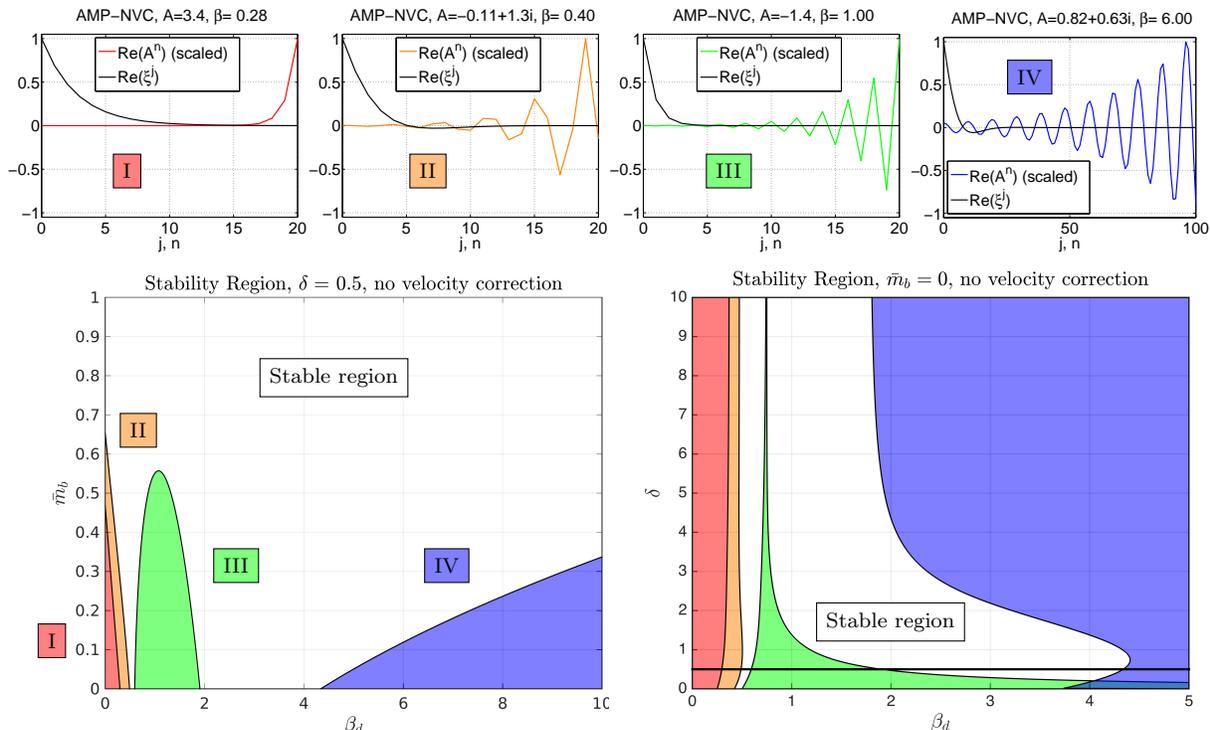
}

Figure~\ref{fig:addedDampingNvcModes} shows the
four instability regions for $\delta=0.5$.  Here we observe that for $\adc$ between~$2$ and~$4$, approximately, the \ampRB-NVC scheme (no velocity correction) is stable for
any value of the mass-ratio $\mbBar\ge 0$. The lower-right plot in Figure~\ref{fig:addedDampingNvcModes}
shows the stable and unstable regions for the limiting case of $\mbBar=0$, where for small $\delta$  
the \ampRB-NVC scheme is unstable except for a very narrow gap near $\beta\approx 0.54$. 
For practical FSI computations, this sensitivity makes the scheme with no velocity correction 
unusable for very light bodies. 

On the other hand, Figure~\ref{fig:addedDampingVcModes} presents the corresponding plots for the \ampRB-VC scheme (with velocity correction).
The additional velocity-correction step is seen to improve the stability of the scheme,
especially for situations when added-damping effects are large. In particular  the green instability region II found in $\SRb$ is
no longer present in $\SRv$. As a result there is range of $\adc$ centered about $\adc=1$ approximately where the scheme 
is stable even as $\mbBar\rightarrow 0$. This makes the scheme more robust when
applied to practical FSI problems. The stability results presented in Figures~\ref{fig:addedDampingNvcModes} 
and~\ref{fig:addedDampingVcModes} are confirmed in the numerical results presented in Section~\ref{sec:numericalVerification}, and in the results of Part~II.

{
\newcommand{\labelFont}{\scriptsize}
\newcommand{\trimfig}[2]{\trimw{#1}{#2}{.0}{.0}{.0}{.0}}
\newcommand{\trimfigz}[2]{\trimw{#1}{#2}{.05}{.05}{.03}{.05}}
\newcommand{\figWidtha}{8.5cm}
\newcommand{\figWidth}{4.0cm}
\newcommand{\figWidthz}{4.cm}
\begin{figure}[htb]
\begin{center}
\begin{tikzpicture}[scale=1]
  \useasboundingbox (0.0,1) rectangle (16,10.25);  
  \begin{scope}[xshift=-.25cm]
  \draw(-.05,0) node[anchor=south west,xshift=-4pt,yshift=+0pt] {\trimfig{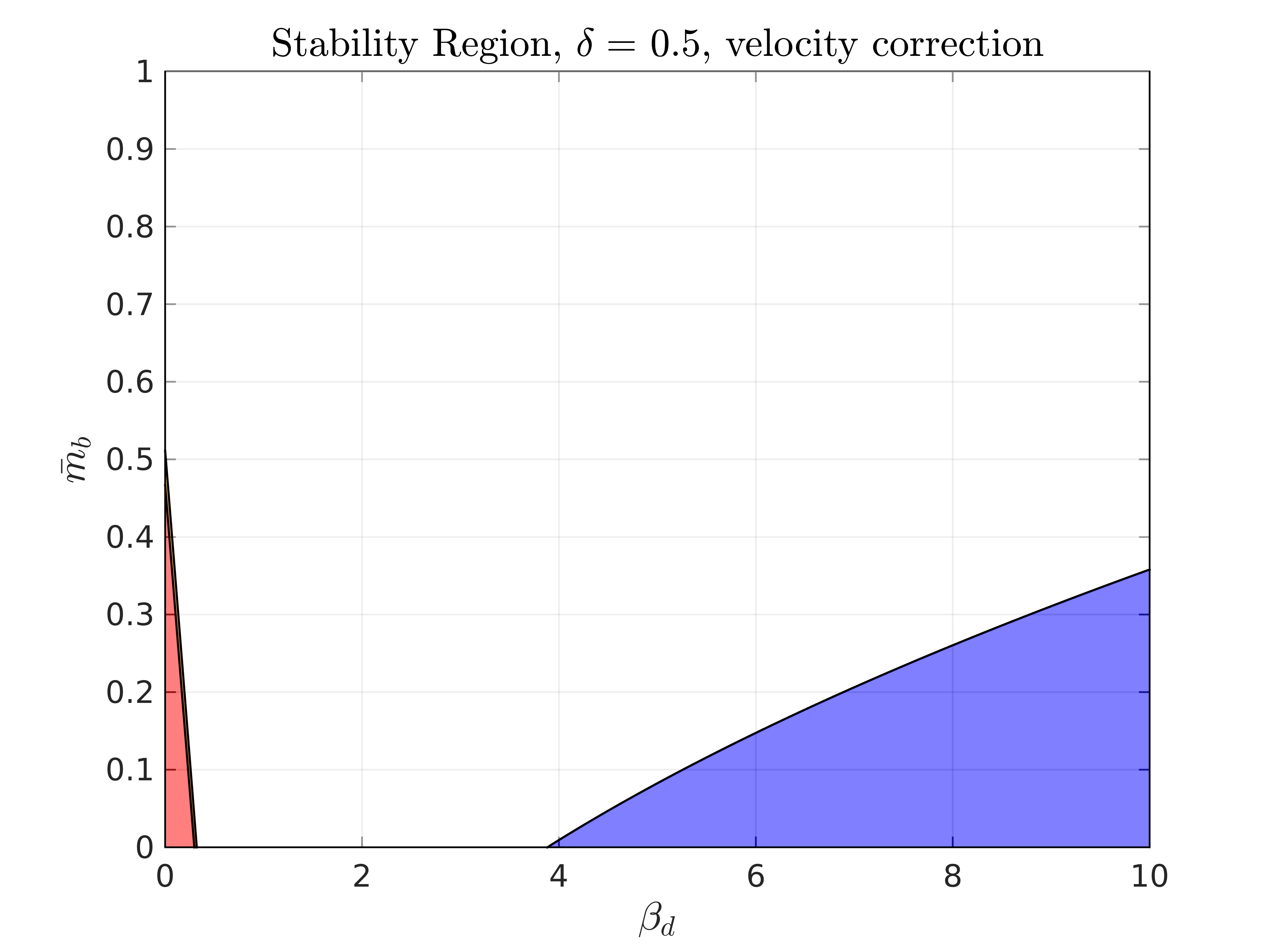}{\figWidtha}};
  \draw(7.75,0) node[anchor=south west,xshift=-4pt,yshift=+0pt] {\trimfig{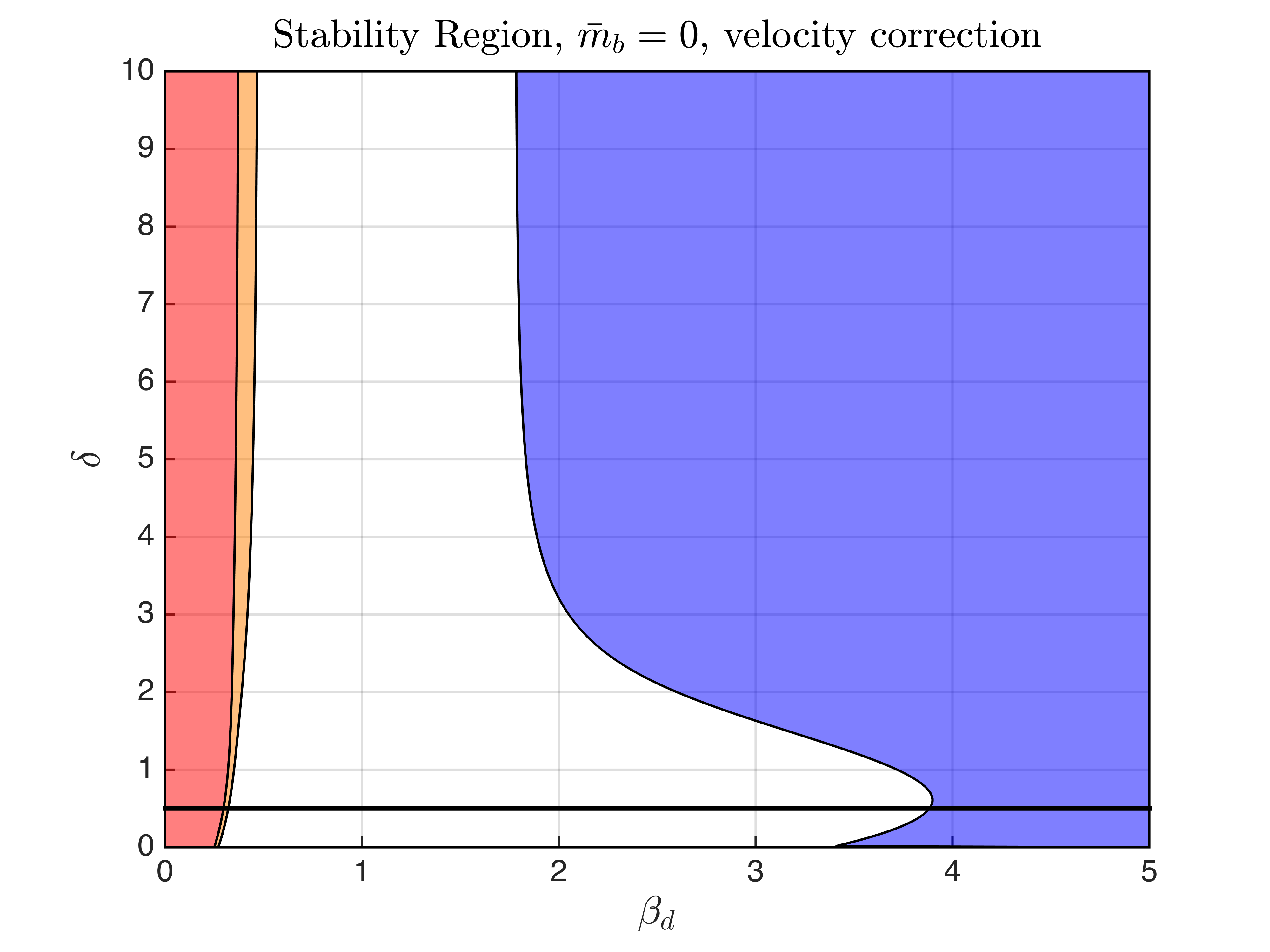}{\figWidtha}};
  \draw(1.5,2.) node[anchor=south west,xshift=-4pt,yshift=+0pt] {\trimfigz{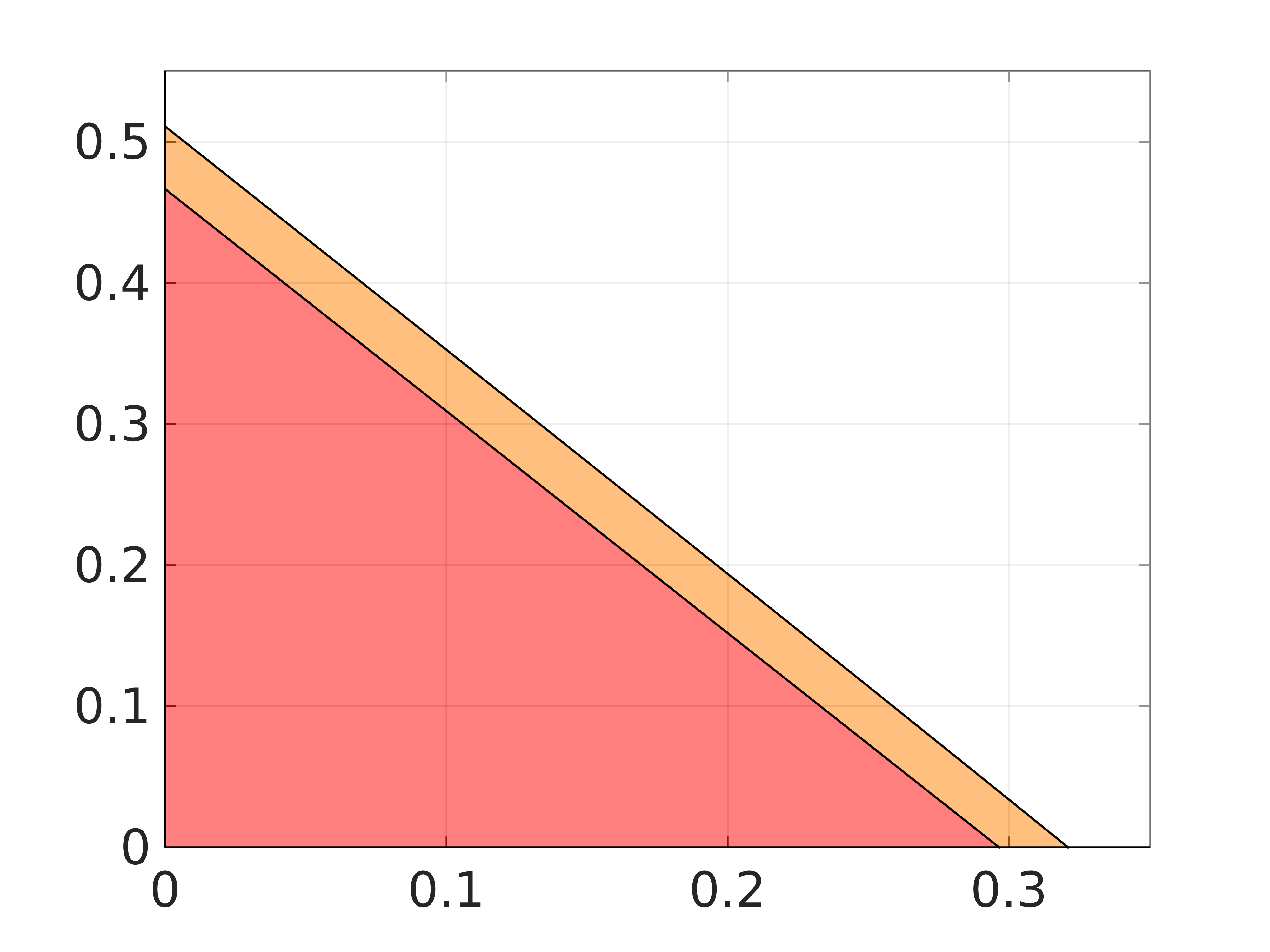}{\figWidthz}};
%
    \draw(5,5.65) node[draw,fill=white] {\labelFont{Stable region}};
    \draw(6.5,2.75) node[draw,fill=blue!50] {\labelFont IV};
    \draw(4.3,4.3) node[draw,fill=orange!50] {\labelFont II};
    \draw(2.15,1.75) node[draw,fill=red!50] {\labelFont I};
    \draw(11.2,1.75) node[draw,fill=white] {\labelFont{Stable region}};
  \end{scope}
  \begin{scope}[xshift=1.5cm,yshift=6.5cm]
    \draw(0.0,0.0) node[anchor=south west,xshift=-4pt,yshift=+0pt] {\trimfig{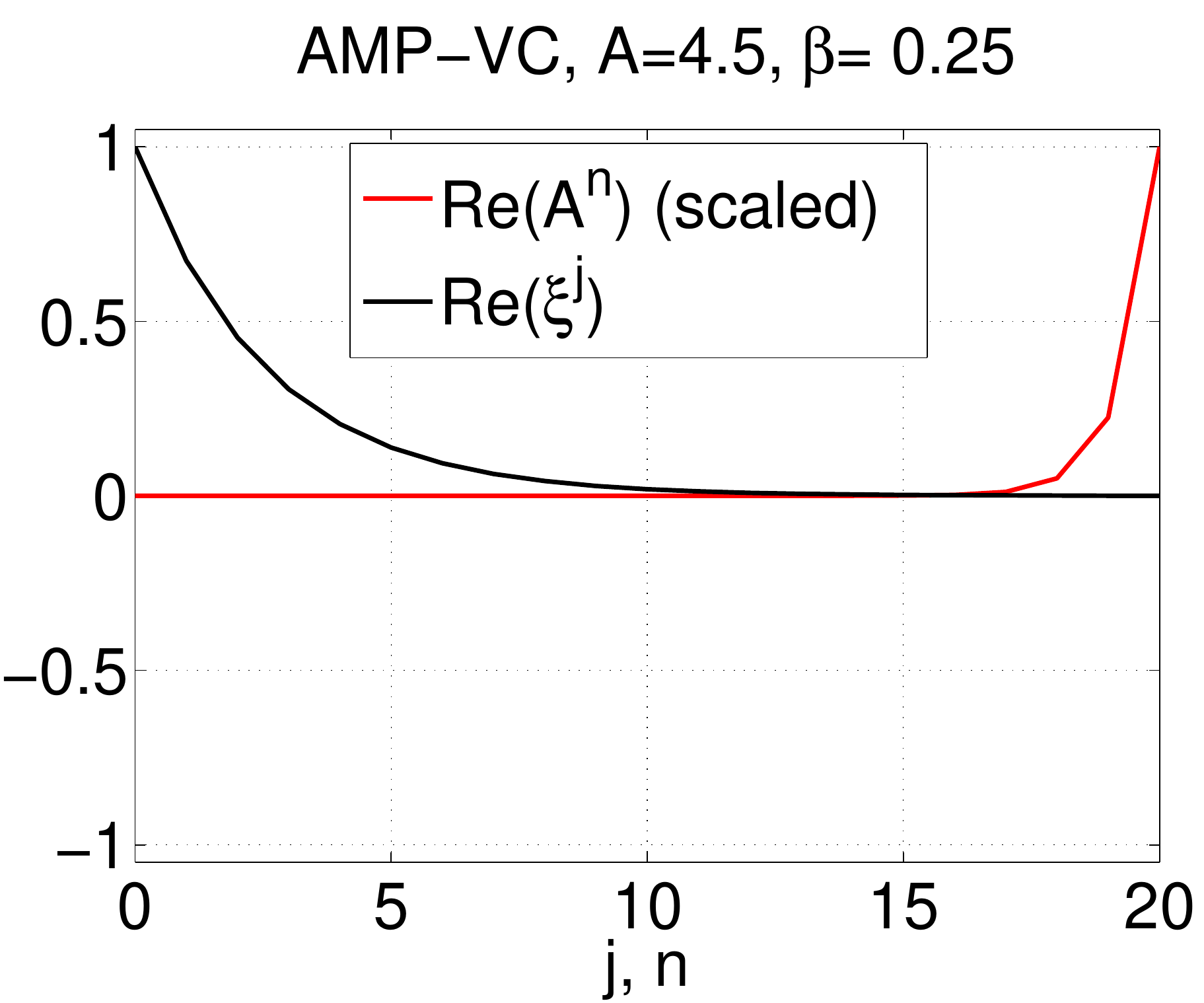}{\figWidth}};
    \draw(2,1.25) node[draw,fill=red!50] {\labelFont I};
    \draw(4.0,0.0) node[anchor=south west,xshift=-4pt,yshift=+0pt] {\trimfig{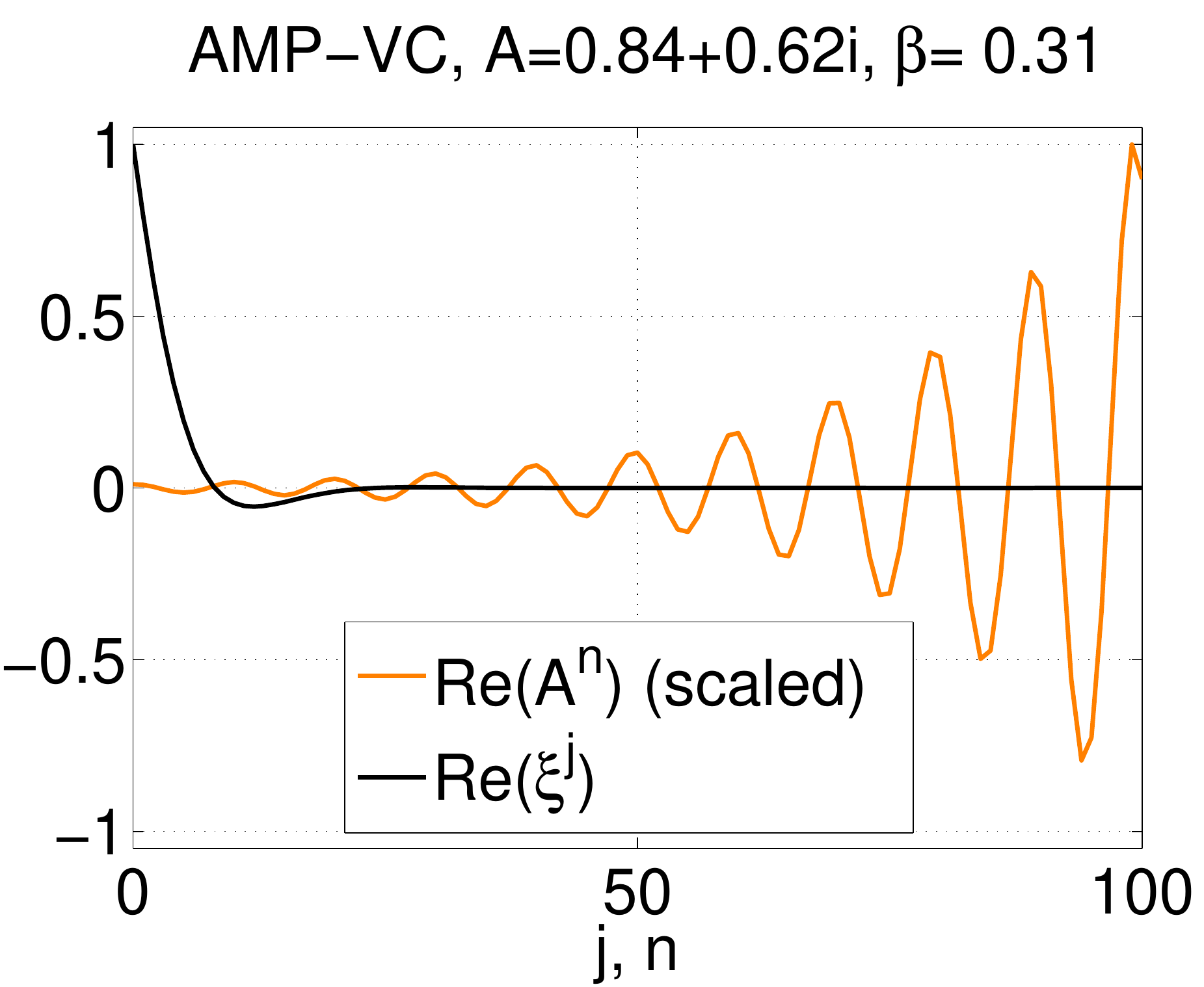}{\figWidth}};
    \draw(6,2.5) node[draw,fill=orange!50] {\labelFont II};
    \draw(8.0,0.0) node[anchor=south west,xshift=-4pt,yshift=+0pt] {\trimfig{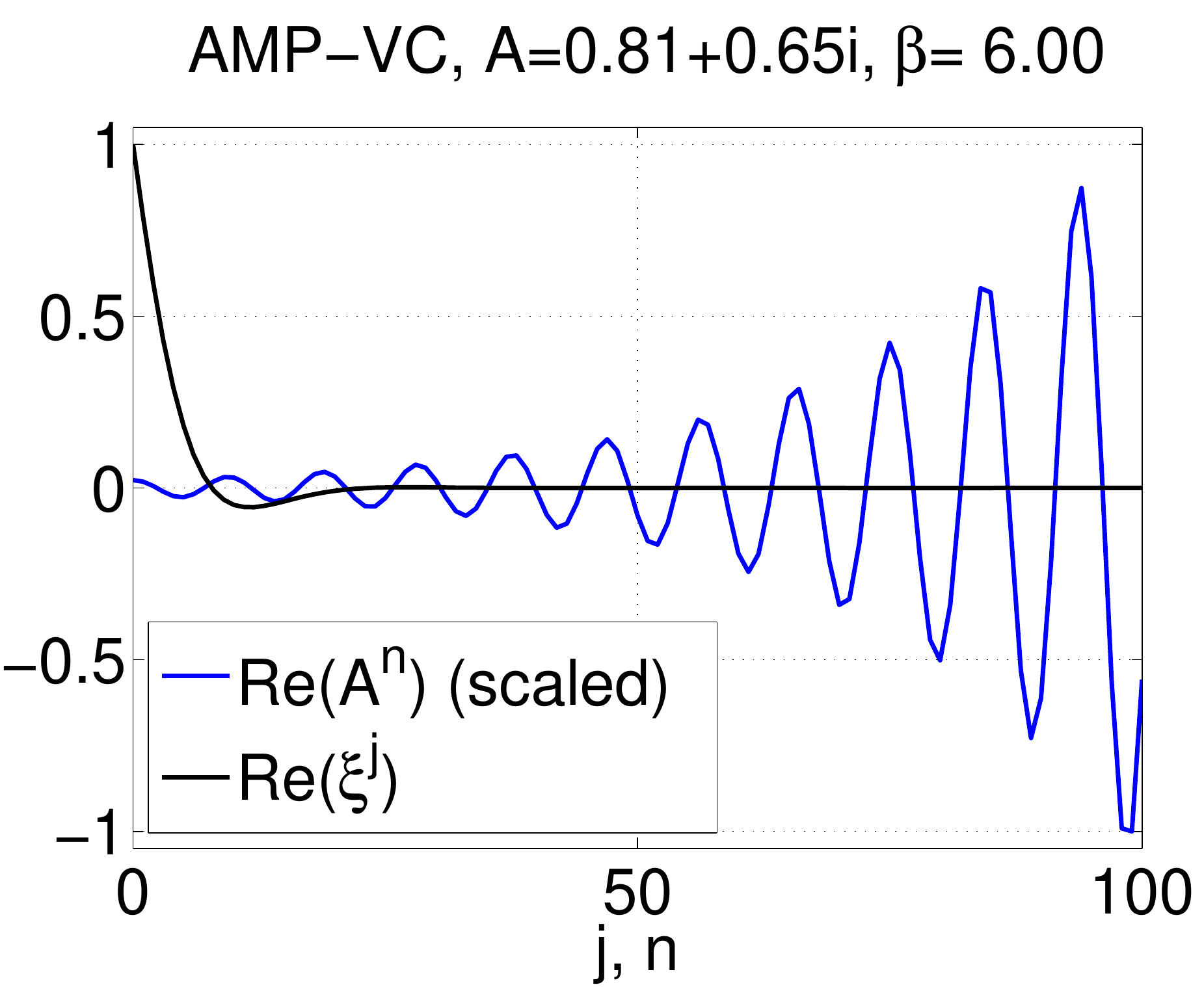}{\figWidth}};
    \draw(10,2.5) node[draw,fill=blue!50] {\labelFont IV};
  \end{scope}
%
\end{tikzpicture}
\end{center}
\caption{
Stability regions and unstable modes for the \ampRB-scheme with velocity correction.
Bottom left: stability region in the $\adc$--$\mbBar$ plane for $\delta=.5$ with zoom of lower left; 
shaded unstable regions are marked according to the
form of the instability. 
Bottom: right: stability region in the $\adc$--$\delta$ plane for the zero-mass body.
Top: unstable modes in space and time for the three instability types ($\delta=.5$ and $\mbBar=0$).
}
\label{fig:addedDampingVcModes}
\end{figure}
}

{
\newcommand{\labelFont}{\scriptsize}
\newcommand{\figWidth}{8.5cm}
\newcommand{\trimfig}[2]{\trimw{#1}{#2}{.0}{.0}{.0}{.0}}
\begin{figure}[htb]
\begin{center}
\resizebox{14cm}{!}{
\begin{tikzpicture}[scale=1]
  \useasboundingbox (0.0,1.) rectangle (16,6.2);  
  \draw(3.2,0) node[anchor=south west,xshift=-4pt,yshift=+0pt] {\trimfig{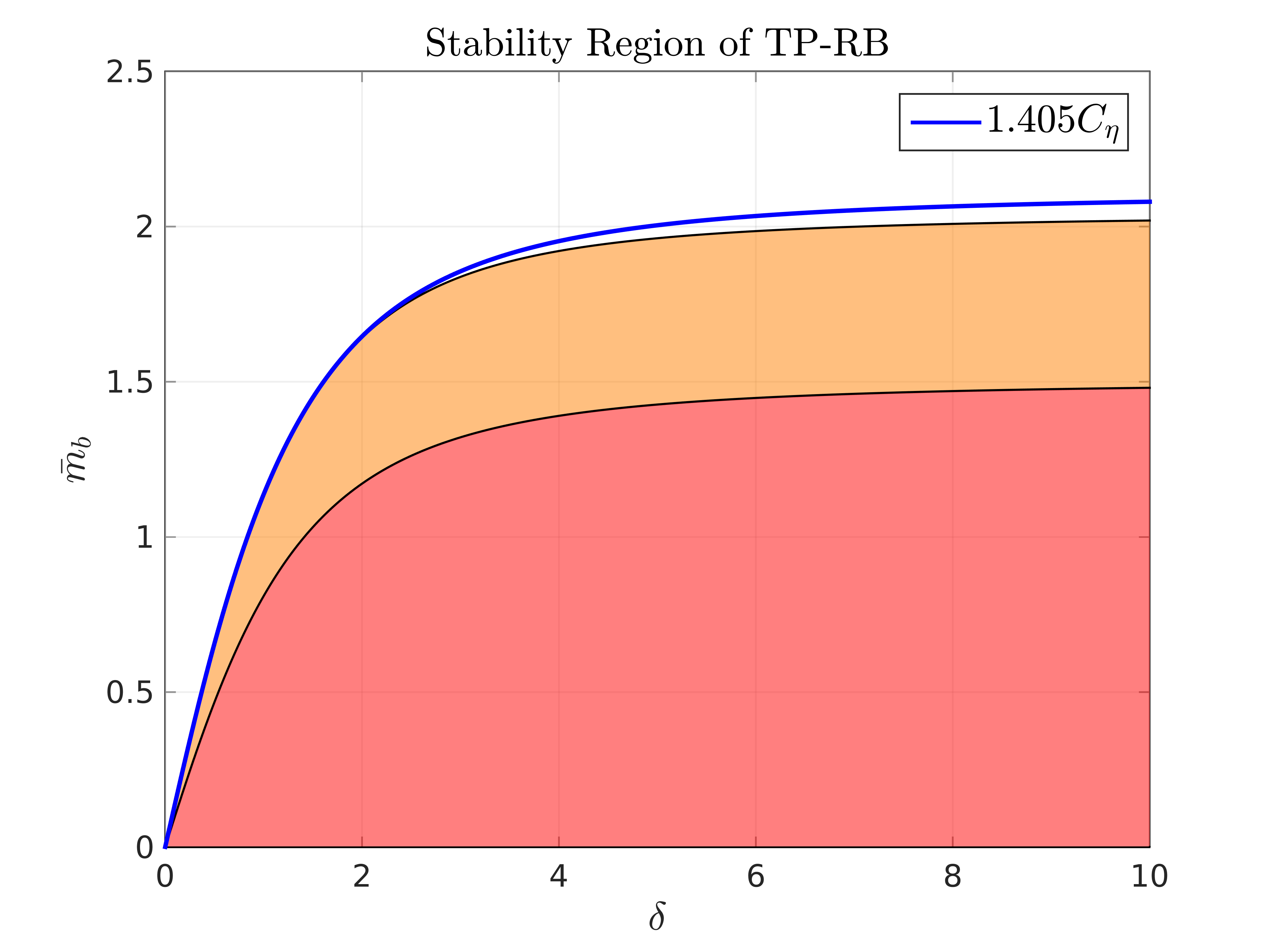}{\figWidth}};
  \draw(6,5.5) node[draw,fill=white] {\labelFont{Stable region}};
%
\end{tikzpicture}
}
\end{center}
\caption{The stability region in the $\delta$--$\mbBar$ plane for the TP-RB scheme for the added-damping model problem.
     The blue curve is an analytic approximation to the boundary of the stable region.}
  \label{fig:addedDampingStabilityTP}
\end{figure}
}

Finally, consider the traditional-partitioned (\tpRB) scheme for the added-damping model problem
corresponding to $\adc=0$ in Algorithm~\ref{alg:ampRBad} (no velocity correction). 
The stability region in this case depends on $\mbBar$ and $\delta$ only, and the numerically
computed stability region is shown in Figure~\ref{fig:addedDampingStabilityTP} as the unshaded (white) region.
The blue curve in the figure is a simple analytical approximation to the stability boundary 
given by the formula 
\newcommand{\const}{K_{\rm tp}}
\begin{align}
	\bar m_b = \const \, C_\eta 	\label{eq:TPRBstabilityRegionApprox}
\end{align}
where $C_\eta$ is defined in~\eqref{eq:CetaDef} and $\const \approx 1.405$ is determined by an approximate fit. 
The analytic form in~\eqref{eq:TPRBstabilityRegionApprox} defines 
a sufficient (and nearly necessary) condition
from which a time-step restriction can be determined.
As can be seen in the subsequent theorem, the {\tpRB} can have a severe
time-step restriction as the mass of the body becomes small.
\begin{theorem}
   The traditional-partitioned scheme for solving the added-damping model problem~{\adplanar}, using Algorithm~\ref{alg:ampRBad}
with $\adc=0$ and no velocity-correction step,
 is stable when $\mbBar>0$ for sufficiently small $\dt$ with $\dt\rightarrow 0$ as $\mb\rightarrow 0$.
The stability region depends on $\mbBar$ and $\delta$ and the numerically
computed stability region is shown in Figure~\ref{fig:addedDampingStabilityTP}. 
Based on the approximate stability boundary~\eqref{eq:TPRBstabilityRegionApprox}, 
a sufficient and nearly necessary condition
for stability is that the time-step must satisfy
\newcommand{\strute}{\rule{0pt}{18pt}}
\begin{equation}
	\dt\,\lesssim\,\frac{2}{\const\, C_\eta}\left(\frac{m_b} {\rho LH}\right)\left(\frac{\dy}{H}\right)\,\frac{H\sp2}{\nu} \quad
  \Longrightarrow \quad
	\dt\,\lesssim
  \begin{cases}
    \displaystyle  \frac{2}{\const^2}\left(\frac{m_b} {\rho L H} \right)^2\frac{H^2}{\nu} & \text{for  $\delta \ll 1$}, \\
     \displaystyle   \frac{4} {3\const} \left(\frac{m_b}{ \rho LH}\right)\left(\frac{\dy}{H}\right)\frac{H^2}{\nu}  \strute & \text{for  $\delta \gg 1$}. 
  \end{cases}  
  \label{eq:TPstabilityBound} 
\end{equation}
\label{theorem:ADTP} 
\end{theorem}
{
\newcommand{\figWidth}{8.6cm}
\newcommand{\figWidthz}{3.25cm}
\newcommand{\trimfig}[2]{\trimw{#1}{#2}{.03}{.05}{.0}{.0}}
\newcommand{\trimfigz}[2]{\trimFigb{#1}{#2}{.15}{.25}{.1}{.15}}
\begin{figure}[htb]
\begin{center}
\begin{tikzpicture}[scale=1]
  \useasboundingbox (0.5,.75) rectangle (16.,6.8);  
  \begin{scope}[xshift=-.8cm,yshift=-.3cm]
  \draw(0,0) node[anchor=south west,xshift=-4pt,yshift=+0pt] {\trimfig{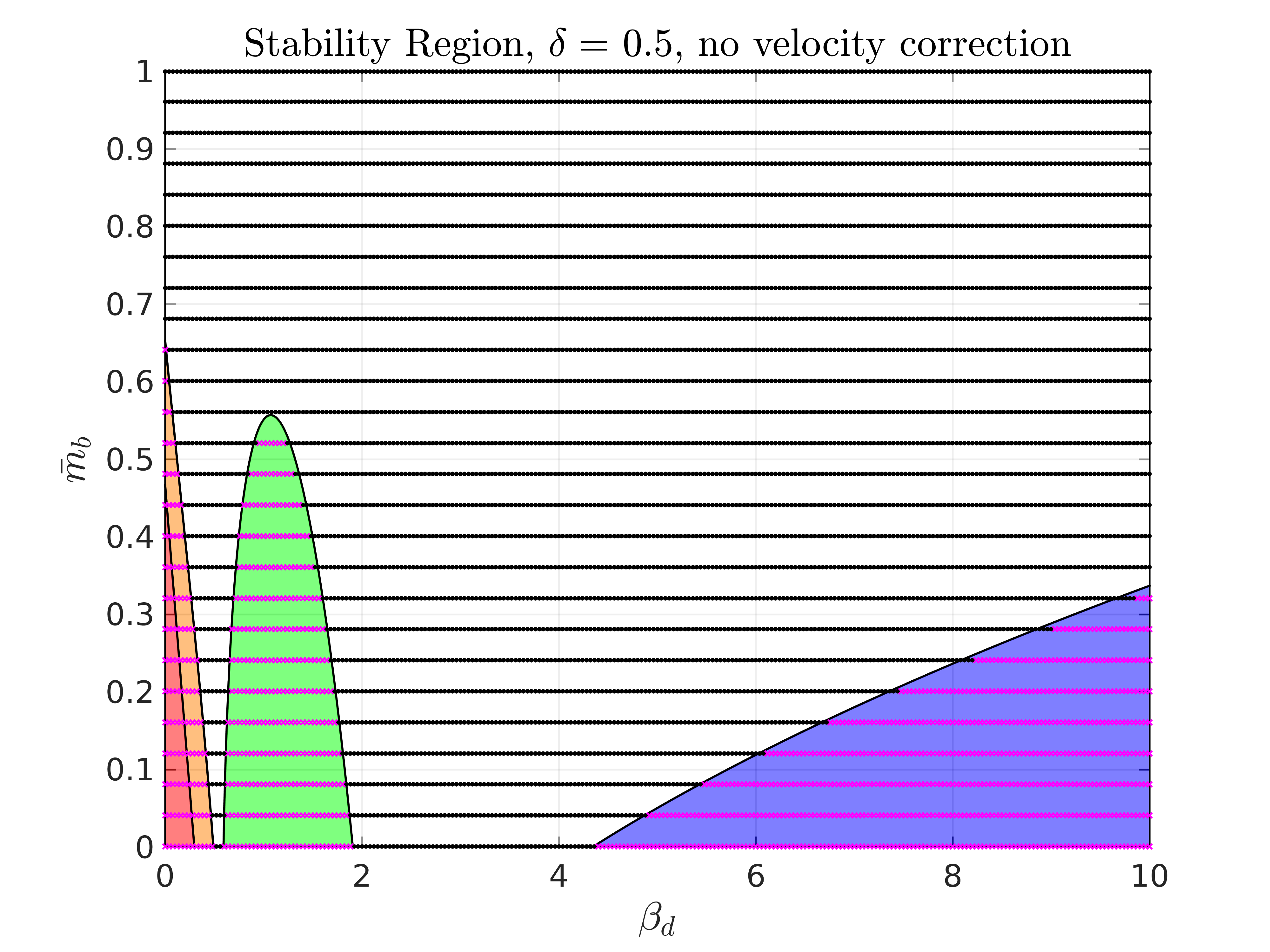}{\figWidth}};
  \draw(1.5,3.8) node[anchor=south west,xshift=-4pt,yshift=+0pt] {\trimfigz{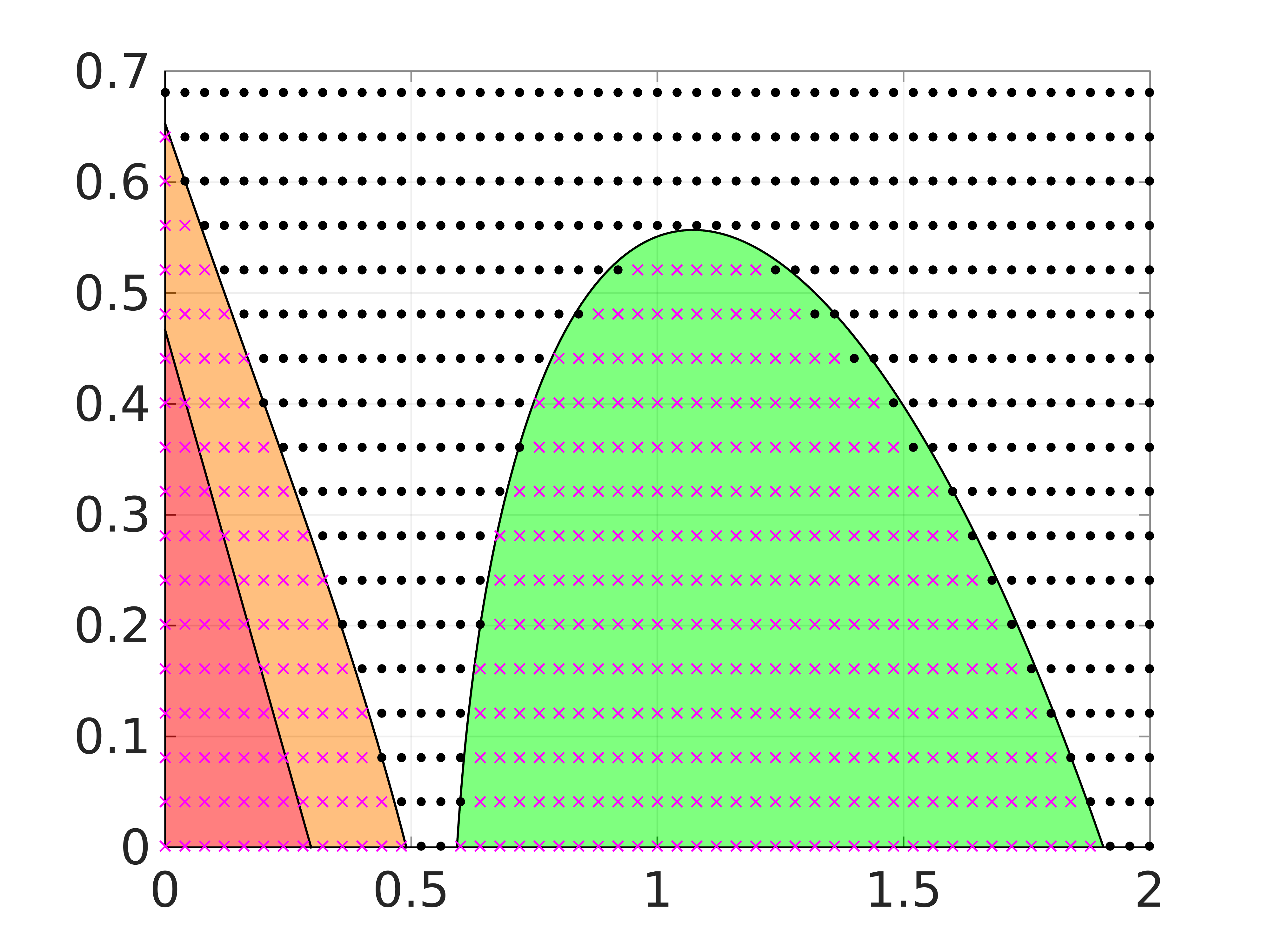}{\figWidthz}};
  \draw(5.1,3.8) node[anchor=south west,xshift=-4pt,yshift=+0pt] {\trimfigz{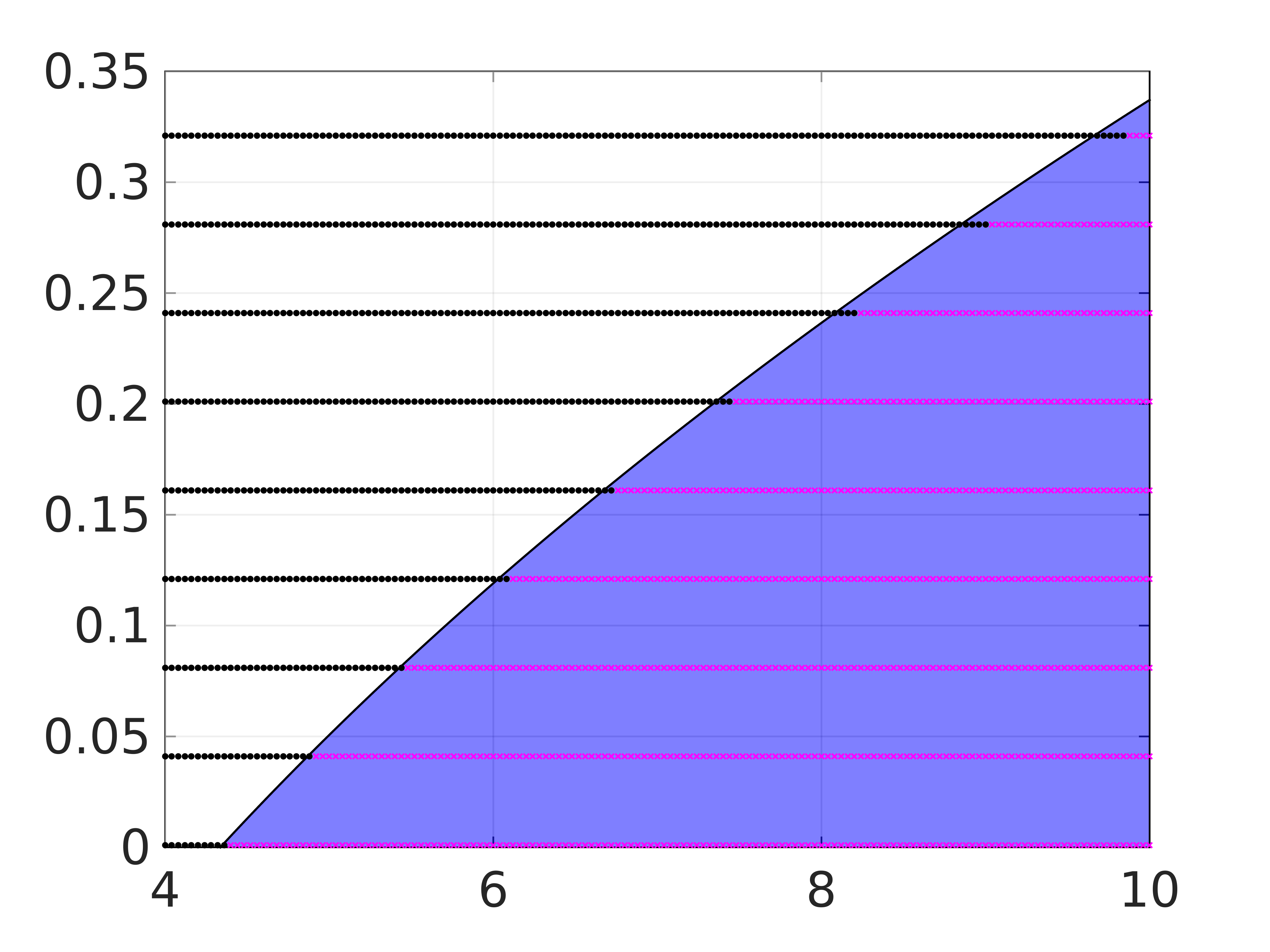}{\figWidthz}};
  \end{scope}

  \begin{scope}[xshift=7.8cm,yshift=-.3cm]
  \draw(0,0) node[anchor=south west,xshift=-4pt,yshift=+0pt] {\trimfig{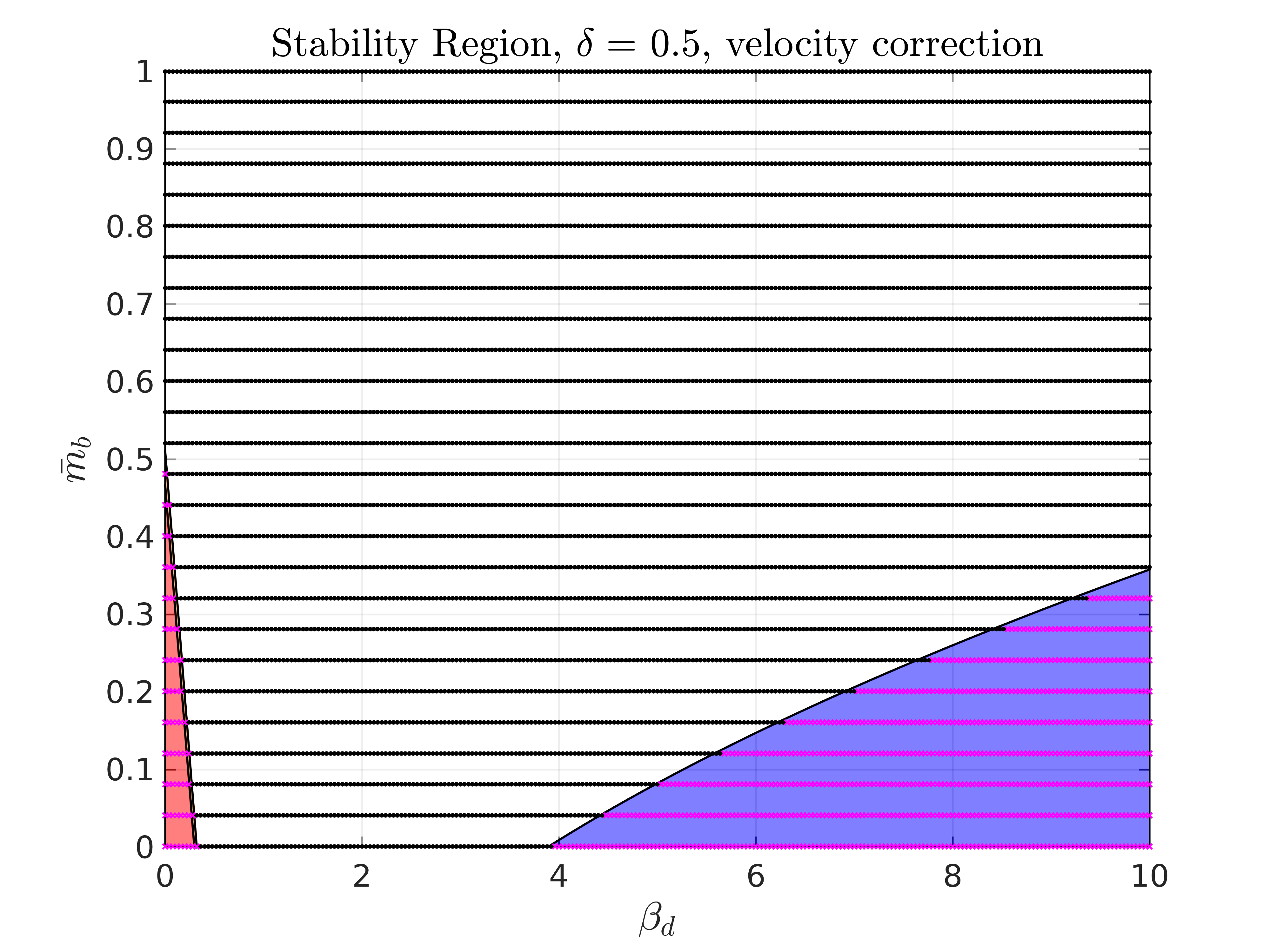}{\figWidth}};
 \draw(1.5,3.8) node[anchor=south west,xshift=-4pt,yshift=+0pt] {\trimfigz{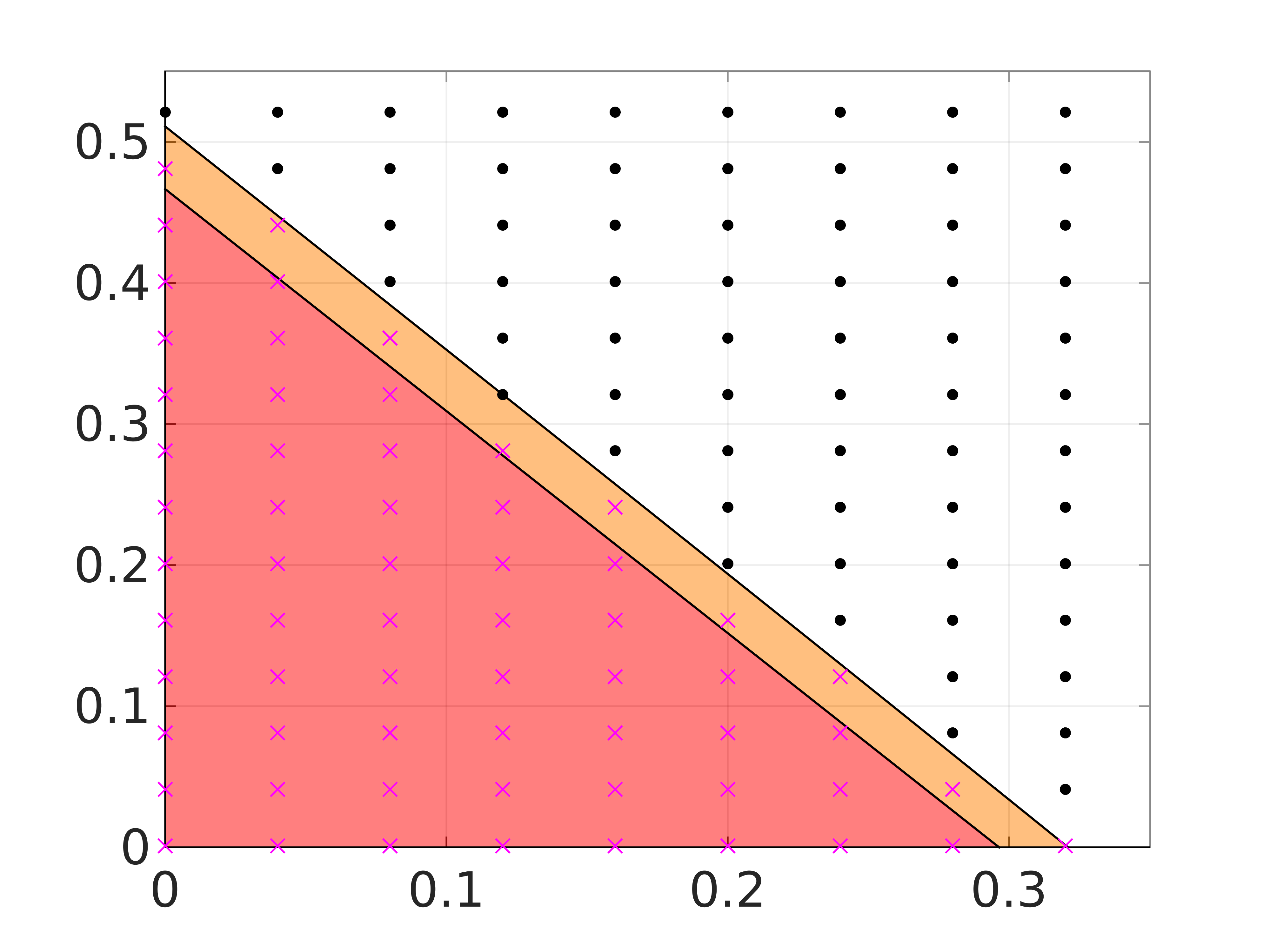}{\figWidthz}};
  \draw(5.1,3.8) node[anchor=south west,xshift=-4pt,yshift=+0pt] {\trimfigz{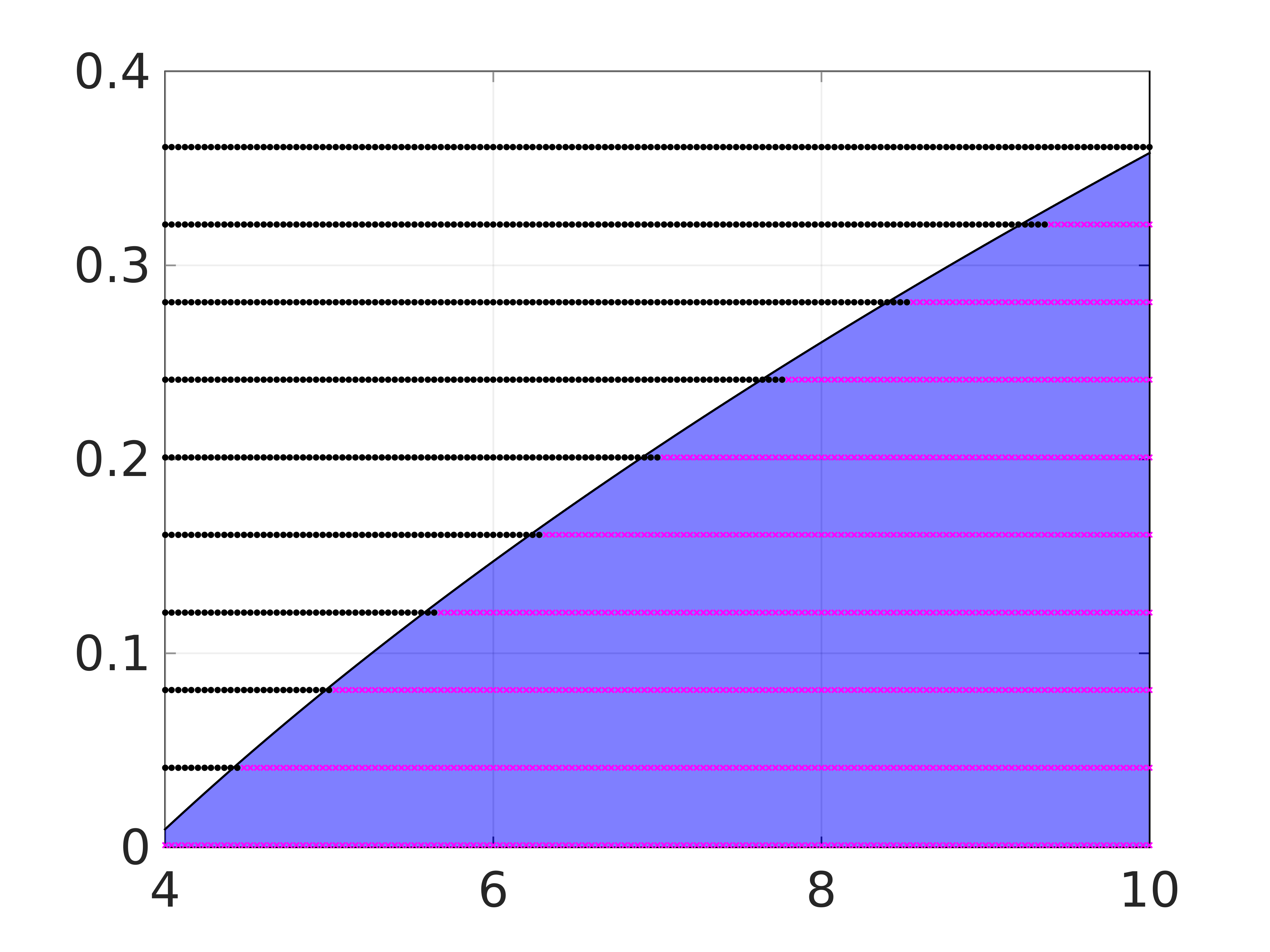}{\figWidthz}};
  \end{scope}
%
\end{tikzpicture}
\end{center}
  \caption{Numerically determined points of stability (black dots for stable and magenta crosses for unstable) compared to the theoretical stability region
     for the added-damping model problem in a rectangular geometry.
  Left: the~\ampRB-NVC scheme. Right:  the \ampRB-VC scheme.
}
  \label{fig:numericalStabilityShearBlock}
\end{figure}
}

\noindent
Equation~\eqref{eq:TPstabilityBound} follows directly from~\eqref{eq:TPRBstabilityRegionApprox} using the definition
of $\mbBar$ and the asymptotic forms for $\ceta$ given by
\[
   \ceta \sim \begin{cases} 
                  \delta      & \text{for $\delta \ll 1$}, \\
                  \frac{3}{2} & \text{for $\delta \gg 1$}.
              \end{cases}
\]
Thus, although the TP-RB scheme can be stabilized with a sufficiently small value of $\dt$ when the mass of the rigid body is finite, the 
maximal stable time-step approaches zero as the square of the mass tends to zero, making 
the TP-RB scheme rapidly impractical for light bodies.

The results of the stability analysis can be confirmed numerically by solving the model problem
and checking whether the computed solutions appear to be stable or unstable.
Figure~\ref{fig:numericalStabilityShearBlock} shows values of the parameters $(\adc,\mbBar)$ with $\delta=0.5$ held fixed where
the~\ampRB~scheme, with and without velocity correction, is determined to be stable or unstable.  Solutions are
computed using  the grid $\Gs^{(4)}$, as defined
in Section~\ref{sec:shearBlock}, with $\dt = 0.05$.  The black
dots in the figure indicate points where the numerical solution appears to be stable up to a final time of $t = 50$, while the
magenta crosses denote points where the solution is unstable.  The numerically computed stability region is seen 
to be in excellent agreement with the theory. There are minor discrepancies near the boundaries of
the stability regions but these can be attributed to the final time not being large enough to show the
slow unstable growth that would be expected near the stability boundaries.

\def\Ibar{\bar{I}}
\def\aabar{\bar{\aa}}
\def\bbbar{\bar{\bb}}
\def\rtfA{\xi}
\newcommand{\amannular}{{MP-AMA}}
\newcommand{\adannular}{{MP-ADA}}

\section{Analysis of FSI model problems in an annular geometry} \label{sec:analysisDisk}

In this section, we extend the analysis of the {\ampRB}~algorithm to two FSI model problems in an
annular geometry, one focusing on added-mass effects and the other on added-damping.  
The analysis for the annular geometry, shown in Figure~\ref{fig:modelDiskInDisk}, confirms
the essential conclusions derived 
from the simplified rectangular-geometry model problems, including the formulation of
the added-damping coefficient $\Duu$ and the determination of the added-damping parameter $\adc$.
For the present problem, the incompressible
fluid occupies the annular domain $\OmegaF$ between an inner radius $\aa$ and an outer radius $\bb$,
and surrounds the rigid body with mass $m_b$ and moment of inertial $I_b$ in the domain $\Omega_b$.
As before, we consider small displacements of the rigid body and linearize the motion of the interface separating
the fluid and the body about the fixed circle with radius $\aa$, denoted by $\Gamma_b$ in the
figure.

{
\def\rad{1.5}
\def\radin{1.625}
\def\radout{3.1}
\def\radline{.8}
\newcommand{\plotDiskModel}{
\fill[fill=red!20,draw=red,line width=2pt] 
      plot[samples=100, domain=0.:360] ( {\rad*cos(\x)} , {\rad*sin(\x)} ) -- cycle ;
}
\newcommand{\plotDiskoutModel}{
\filldraw[fill=blue!20,even odd rule, draw=blue,line width=2pt] 	
	(0,0) circle (\radout)
	(0,0) circle (\radin);
}
\begin{figure}[htb]
\begin{center}
\resizebox{13cm}{!}{
\begin{tikzpicture}[scale=.9]
\useasboundingbox (0,0) rectangle (16.,5.5);
  \begin{scope}[xshift=8cm,yshift=2.5cm]  
	  \plotDiskoutModel
	  \plotDiskModel
	  \draw[thick,black,yshift=-4pt] (0,2.4) node {fluid: $\OmegaF$};
	  \draw[thick,black,xshift=-.7cm,yshift=-2pt] (0,-2.25) node[anchor=south] {interface: $\GammaB$};
	  \draw[thick,black,xshift=4pt] (0,-.4) node {rigid body};
	  \draw[thick,black,xshift=4pt,yshift=-14pt] (0,-.4) node {$\OmegaB$};
	  \draw[thick,black,xshift=2pt] (-.4,.7) node {$\aa$};
	  \draw[thick,black,xshift=3pt,yshift=1pt] (1.6,1.2) node {$\bb$};
	  \draw[very thick,->,black] (0,0) -- (-\radin*1.4/2,\radin*1.4/2); 
	  \draw[very thick,->,black] (0,0) -- (\radout*1.4/2,\radout*1.4/2); 
  \end{scope}
\end{tikzpicture}
}
\end{center}
\caption{The geometry for the annular geometry FSI model problems.}
\label{fig:modelDiskInDisk}
\end{figure}
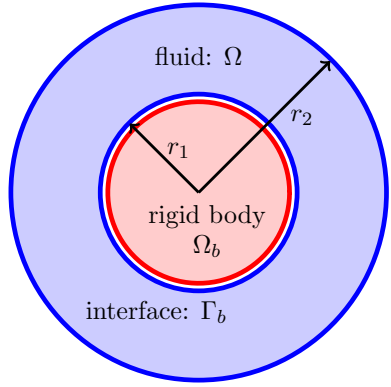
}

It is convenient to describe the velocity of the fluid in terms of its components in the radial and
circumferential directions, denoted by $u(r,\theta,t)$ and $v(r,\theta,t)$, respectively, where
$\xv=(r,\theta)$ are polar coordinates, and the fluid pressure is given by $p(r,\theta,t)$.  It is
sufficient for the purposes of the analysis to consider rotational motion of the body with angular
velocity $\omega_b(t)$ and translational motion of the body only in the horizontal direction with
velocity $u_b(t)$, for example.  Under these assumptions, the governing equations are
\begin{alignat}{2}
 \text{Fluid:}\quad & \rho{\partial u\over\partial t}+{\partial p\over\partial r}=\mu\left[{1\over r}{\partial\over\partial r}\left(r{\partial u\over\partial r}\right)+{1\over r\sp2}{\partial\sp2u\over\partial\theta\sp2}-{u\over r\sp2}-{2\over r\sp2}{\partial v\over\partial\theta}\right] , \quad&& \xv\in\OmegaF, \label{eq:velocityEquation1A}\\
                    & \rho{\partial v\over\partial t}+{1\over r}{\partial p\over\partial\theta}=\mu\left[{1\over r}{\partial\over\partial r}\left(r{\partial v\over\partial r}\right)+{1\over r\sp2}{\partial\sp2v\over\partial\theta\sp2}-{v\over r\sp2}+{2\over r\sp2}{\partial u\over\partial\theta}\right], \quad&& \xv\in\OmegaF, \label{eq:velocityEquation2A}\\
  &  {\partial\over\partial r}(ru)+{\partial v\over\partial\theta} = 0   , \quad&&\xv\in\OmegaF, \label{eq:continuityEquationA}\\
\text{Rigid body:}\quad  &   m_b\,{d u_b\over dt} = \int_{\Gamma_b} \bigl(\ev_x\sp{T}\sigmav\nv\bigr)\,ds + g_u(t) , \quad && \label{eq:bodyEquation1A}\\
  &   I_b\,{d\omega_b\over dt} = \int_{\Gamma_b} \aa\bigl(\tv\sp{T}\sigmav\nv\bigr)\,ds + g_\omega(t) ,   \label{eq:bodyEquation2A} \\
 \text{Interface:}\quad
&   u=u_b\cos\theta,\qquad v=-u_b\sin\theta+\aa\omega_b, \quad&& \xv\in\Gamma_b , \label{eq:matchVA} \\
 \text{Fluid BCs:} \quad 
&  \vv(\bb,\theta,t) =0, \quad&& \theta\in[-\pi,\pi] . \label{eq:modelProblemLastA}
\end{alignat}
Here, $\sigmav$ is the fluid stress tensor, $\nv$ and $\tv$ are the unit normal and tangent vectors on $\Gamma_b$, respectively, $\ev_x$ is the unit normal in the positive $x$-direction, and $g_u(t)$ and $g_\omega(t)$ are an applied force and torque on the body, respectively.  As before, we note that an elliptic equation can be obtained for the fluid pressure using the continuity equation in~\eqref{eq:continuityEquationA} to eliminate the components of velocity in the momentum equations in~\eqref{eq:velocityEquation1A} and~\eqref{eq:velocityEquation2A}.  For this annular geometry, the pressure equation is
\begin{equation}
{\partial\over\partial r}\left(r{\partial p\over\partial r}\right)+{1\over r}{\partial\sp2p\over\partial\theta\sp2}=0,
\end{equation}
which, along with appropriate boundary conditions, determines the pressure 
in the fractional-step method used to solve the two model problems described below.

\subsection{Model problem {\amannular} for added-mass in an annular geometry}

\newcommand{\ama}{{\cal M}_a}

A model problem for added mass can be obtained from the governing equations in~\eqref{eq:velocityEquation1A}--\eqref{eq:modelProblemLastA} by considering a Fourier-mode decomposition of the fluid variables in the $\theta$-direction.  Setting
\begin{equation}
u(r,\theta,t)=\uHat(r,t)\cos\theta,\qquad v(r,\theta,t)=-\vHat(r,t)\sin\theta,\qquad p(r,\theta,t)=\pHat(r,t)\cos\theta,
\label{eq:transformation}
\end{equation}
and dropping the decoupled rotational motion of the body gives
the {\amannular} model problem of the form
\begin{equation}
 \text{{\amannular}:}\;  \left\{ 
   \begin{alignedat}{3}
  & \rho{\partial\uHat\over\partial t}+{\partial\pHat\over\partial r} = 0 , \quad&&  r\in(\aa,\bb),  \\
  & \rho{\partial\vHat\over\partial t}+{\pHat\over r} = 0 , \quad&&  r\in(\aa,\bb),  \\
  & {\partial\over\partial r}(r\uHat)-\vHat = 0   , \quad&& r\in(\aa,\bb), \\
  &   m_b \frac{du_b}{dt} = - \aa\pi\, \pHat(\aa,t)+ g_u(t), \\
  &  \uHat(\aa,t)=u_b, \\
  &  \uHat(\bb,t)=0.
   \end{alignedat}  \right. 
  \label{eq:MP-AMA}
\end{equation}
Here the viscous terms have been dropped 
since it was shown in~\cite{Conca1997} that the added-mass for a cylindrical rigid body immersed 
in an incompressible fluid does not depend on viscosity.
As a consequence, the boundary condition for the circumferential component of the fluid 
velocity at $r=\aa$ and $\bb$ have been dropped as well.
The exact solution of the {\amannular} model problem is given by
\begin{align}
u_b(t)&={1\over m_b+\ama}\int_0\sp{t}g_u(\tau)\,d\tau+u_b(0), \label{eq:ubsolna} \\
\pHat(r,t)&={\rho\aa a_u(t)\over \bb/\aa-\aa/\bb}\left({r\over\bb}+{\bb\over r}\right), \label{eq:pHatsolna} \\
\uHat(r,t)&={\bb\sp2-r\sp2\over\bb\sp2-\aa\sp2}\left({\aa\over r}\right)\sp2u_b(t), \label{eq:uHatsolna} \\
\vHat(r,t)&={\bb\sp2+r\sp2\over\bb\sp2-\aa\sp2}\left({\aa\over r}\right)\sp2u_b(t), \label{eq:vHatsolna}
\end{align}
where
\begin{equation}
\ama=\rho\pi\aa\sp2\left[{1+\bigl(\aa/\bb\bigr)\sp2\over1-\bigl(\aa/\bb\bigr)\sp2}\right],
\label{eq:addedmassannular}
\end{equation}
is the added-mass for this model problem, and the acceleration of the rigid body is given by
\begin{equation}
a_u(t)={d\ub\over dt}={g_u(t)\over m_b+\ama}.
\label{eq:ausolna}
\end{equation}
Note that the added mass $\ama$ approaches the mass of fluid displaced by the rigid body as $\aa/\bb$
tends to zero, and that the added mass becomes large as the fluid gap become small (i.e.~as
$\aa/\bb\rightarrow1$). 
The added-mass approaches infinity as the gap narrows 
since the radial component of the velocity is given as a boundary condition at $r=\bb$ and thus the total volume of fluid is conserved. 
In the prior {\amplanar} rectangular geometry case, the corresponding boundary condition
was on the pressure, which allowed flow through the boundary. 
In fact, if the far boundary in {\amplanar} was taken as a given velocity, then the added-mass would be infinite
and the traditional scheme would be unconditionally unstable regardless of the rigid-body mass. 

The stability of the {\ampRB}~algorithm for the {\amannular} model problem is addressed by
considering the semi-discrete time-stepping scheme given in Algorithm~\ref{alg:ampRBama}.  Here, we
have discretized in time which retains the essential features of the {\ampRB}~algorithm, while the
fluid variables remain continuous in space for convenience in the analysis.  The superscripts
attached to the variables in the algorithm indicate the time level,
e.g.~$\uHat\sp{n}\approx\uHat(r,t\sp{n})$, and whether the variables are extrapolated or predicted
(following the notation used in previous algorithms).

\begin{algorithm}\caption{\ampRB~scheme for {\amannular}  
                         \\ \protect\hphantom{Algorithm 5~~} {\em Added-mass model problem in an annular geometry.}}\small
\[
\begin{array}{l}
\hbox{// \textsl{{\bodyStepIComment}}}\smallskip\\
1.\quad a_u^\esup = 2 a_u^n - a_u^{n-1} , \qquad \ub^\esup = \ub^{n-1} + 2\dt\, a_u^n,\bigskip\\
\hbox{// \textsl{Prediction steps}}\smallskip\\
2.\quad\left\{
\begin{array}{l}
\uHat\sp\psup = \uHat^n - \frac{3\dt}{2\rho}\partial_r\pHat^n +\frac{\dt}{2\rho}\partial_r\pHat^{n-1}, \smallskip\\
\vHat\sp\psup = \vHat^n - \frac{3\dt}{2\rho r}\pHat^n +\frac{\dt}{2\rho r}\pHat^{n-1},
\end{array}\right.   \quad r\in(\aa,\bb), \qquad \uHat^\psup(\aa) = \ub^\esup, \quad \uHat^\psup(\bb) = 0, \smallskip\\
3.\quad \partial_r\bigl(r\partial_r\pHat\sp\psup\bigr)-\pHat\sp\psup/r = 0, \quad r\in(\aa,\bb), \qquad
\left\{\begin{array}{l} \partial_r\pHat^\psup(\aa) + \rho a_u^\psup=0 \smallskip\\ \mb a_u^\psup + \pi\aa\pHat\sp\psup(\aa) = \gunp \end{array}\right., \quad \partial_r\pHat^\psup(\bb)=0, \smallskip\\
4.\quad \ub^\psup = \ub^{n} + \frac{\dt}{2}( a_u^\psup + a_u^n ),\bigskip\\
\hbox{// \textsl{Correction steps}}\smallskip\\
5.\quad\left\{
\begin{array}{l}
\uHat\sp{n+1} = \uHat^n - \frac{\dt}{2\rho}\bigl(\partial_r\pHat^\psup + \partial_r\pHat^{n}\bigr), \smallskip\\
\vHat\sp{n+1} = \vHat^n - \frac{\dt}{2\rho r}\bigl(\pHat^\psup + \pHat^{n}\bigr),
\end{array}\right.   \quad r\in(\aa,\bb), \qquad \uHat^{n+1}(\aa) = \ub^\psup, \quad \uHat^{n+1}(\bb) = 0, \smallskip\\
6.\quad \partial_r\bigl(r\partial_r\pHat\sp{n+1}\bigr)-\pHat\sp{n+1}/r = 0, \quad r\in(\aa,\bb), \quad  \left\{\begin{array}{l} \partial_r\pHat^{n+1}(\aa) + \rho a_u^{n+1}=0 \smallskip\\ \mb a_u^{n+1} + \pi\aa\pHat^{n+1}(\aa) = \gunp \end{array}\right., \quad \partial_r\pHat^{n+1}(\bb)=0, \smallskip\\
7.\quad \ub^{n+1} = \ub^{n} + \frac{\dt}{2}( a_u^{n+1} + a_u^n ).
\end{array}
\]
\label{alg:ampRBama}
\end{algorithm}

The preliminary values for the velocity and acceleration of the body
obtained in Step~1 of Algorithm~\ref{alg:ampRBama} and the predicted values for the components of
the fluid velocity computed in Step~2 are not needed in subsequent steps of the algorithm, and so we
begin with the elliptic problem in Step~3, whose solution is given by
\[
\pHat\sp\psup(r)=\pHat(r,t\sp{n+1}),\qquad a_u\sp\psup=a_u(t\sp{n+1}),
\]
where $\pHat(r,t)$ and $a_u(t)$ are the exact pressure and acceleration of the body given
in~\eqref{eq:pHatsolna} and~\eqref{eq:ausolna}, respectively.  The updates of the velocity of the
body in Step~4 and the components of the fluid velocity in Step~5 use trapezoidal-rule quadratures
of the exact body acceleration and fluid pressure which implies these values are unconditionally stable.  The elliptic
problem for the fluid pressure and the acceleration of the body appears again in Step~6, and the
solutions for $\pHat\sp{n+1}(r)$ and $a_u\sp{n+1}$ agree with the exact solutions at $t=t\sp{n+1}$
as before.  The final step in Algorithm~\ref{alg:ampRBama} is simply a trapezoidal-rule update of
the velocity of the body using the exact body acceleration.  Thus provided $m_b+\ama$ is bounded
away from zero, the discrete solution can be bounded in terms of the data which leads to the following result:
\begin{theorem} The \ampRB~algorithm given in Algorithm~\ref{alg:ampRBama} for the model problem~{\amannular}
is unconditionally stable provided there exists a constant $K>0$ such that 
\begin{align*}
    m_b+\ama \ge K.
\end{align*}
\end{theorem}

The analysis of the traditional-partitioned scheme for the added-mass
model problem~{\amannular} closely follows the previous analysis for the rectangular geometry and we simply state the primary
result without proof. 
\begin{theorem} The traditional-partitioned (\tpRB) algorithm for the added-mass
model problem~{\amannular} is stable if and only if the mass of the body is greater than the added-mass $\ama$,
\begin{align*}
    {\mb\over\ama } > 1. 
\end{align*}
\end{theorem}

\subsection{Model problem {\adannular} for added-damping in an annular geometry}

\newcommand{\linop}{{\cal L}}

A model problem suitable to study added-damping for the annular geometry can be extracted from the
set of equations in~\eqref{eq:velocityEquation1A}--\eqref{eq:modelProblemLastA} by considering the
zeroth Fourier mode of the fluid variables, which is independent of~$\theta$.  In this model, the
circumferential component of the fluid velocity is coupled to the rotation of the rigid body due to
the shear stress acting uniformly along the interface $\Gamma_b$.  The radial component of the 
fluid velocity and the fluid pressure decouple and are therefore dropped from the system.
The result is the 
{\adannular} model problem for added damping given by
\begin{equation}
 \text{{\adannular}:}\;  \left\{ 
   \begin{alignedat}{3}
  & {\partial\vHat\over\partial t}=\nu\linop\vHat, \quad&&  r\in(\aa,\bb),  \\
  &   I_b\,{d\omega_b\over dt} = 2\pi\aa\sp2\mu\left[r{\partial\over\partial r}\left({\vHat\over r}\right)\right]_{r=\aa} + g_\omega(t), \\
  &  \vHat(\aa,t)=\aa \omega_b, \\
  &  \vHat(\bb,t)=0,
   \end{alignedat}  \right. 
  \label{eq:MP-ADA}
\end{equation}
where
\[
\linop\vHat={1\over r}{\partial\over\partial r}\left(r{\partial\vHat\over\partial r}\right)-{\vHat\over r\sp2}.
\]

As in the previous model problem for the added-mass, we consider the semi-discrete time-stepping scheme
given in Algorithm~\ref{alg:ampRBada} for the {\adannular}~problem.  
This algorithm uses the
angular component of the acceleration of the rigid body, which is denoted $b_\omega\sp{n}\approx
b_w(t\sp{n})$.
The added-damping
coefficient, $\Duua$, for the annular geometry takes the form
\newcommand{\deltat}{{\tilde\delta}}
\[
    \Duua=   \mu(2\pi\aa)\aa^2  ~ \frac{1-e^{-\deltat}}{\dr} ,\qquad \deltat \eqdef {\Delta r\over\sqrt{\nu\dt/2}},
\]
which is derived following an analysis of a variational problem similar to that discussed in Section~\ref{sec:addedDampingCoefficient} for the rectangular geometry.  
To assess the stability of the algorithm, it is important to include a choice for the discrete derivative in the AMP interface condition used in Steps~3 and~6. In the analysis to follow, we use the one-sided difference operator, $D_{rh}$, analogous to~\eqref{eq:funkyD}, defined by
\[
D_{rh}f(r)={-3f(r)+4f(r+\Delta r)-f(r+2\Delta r)\over 2\Delta r},
\]
where $f(r)$ is a generic function of $r$.

\begin{algorithm}\caption{\tpRB~scheme ($\beta_d=0$) and \ampRB~scheme ($\beta_d>0$) for {\adannular} 
                         \\ \protect\hphantom{Algorithm 5~~} {\em Added-damping model problem in an annular geometry.}}
\[
\begin{array}{l}
\hbox{// \textsl{{\bodyStepIComment}}}\smallskip\\
1.\quad b_\omega^\esup = 2 b_\omega^n - b_\omega^{n-1} , \qquad \omega_b^\esup = \omega_b^{n-1} + 2\dt\, b_\omega^n,\bigskip\\
\hbox{// \textsl{Prediction steps}}\smallskip\\
2.\quad \vHat\sp\psup = \vHat^n + \frac{\nu\dt}{2}\bigl(\linop\vHat\sp\psup+\linop\vHat\sp{n}\bigr), \quad r\in(\aa,\bb) \qquad \vHat^\psup(\aa) = \aa\omega_b^\esup, \quad \vHat\sp\psup(\bb)=0, \smallskip\\
3.\quad (I_b+\beta_d\dt\,\Duua) b_\omega^\psup=2\pi\aa\sp2\mu \bigl[rD_{rh}(\vHat\sp\psup/r)\bigr]_{r=\aa}+\beta_d\dt\,\Duua b_\omega\sp\esup+g_\omega(t\sp{n+1}), \medskip\\
4.\quad \omega_b^\psup = \omega_b^{n} + \frac{\dt}{2}( b_\omega^\psup + b_\omega^n ),\bigskip\\
\hbox{// \textsl{Correction steps}}\smallskip\\
5.\quad \vHat\sp{n+1} = \vHat^n + \frac{\nu\dt}{2}\bigl(\linop\vHat\sp{n+1}+\linop\vHat\sp{n}\bigr), \quad r\in(\aa,\bb) \qquad \vHat^{n+1}(\aa) = \aa\omega_b^\psup, \quad \vHat\sp{n+1}(\bb)=0, \smallskip\\
6.\quad (I_b+\beta_d\dt\,\Duua) b_\omega^{n+1}=2\pi\aa\sp2\mu \bigl[rD_{rh}(\vHat\sp{n+1}/r)\bigr]_{r=\aa}+\beta_d\dt\,\Duua b_\omega\sp\psup+g_\omega(t\sp{n+1}), \medskip\\
7.\quad \omega_b^{n+1} = \omega_b^{n} + \frac{\dt}{2}( b_\omega^{n+1} + b_\omega^n ), \bigskip\\
\hbox{// \textsl{Fluid-velocity correction step (optional)}}\smallskip\\
8.\quad \vHat\sp{n+1} = \vHat^n + \frac{\nu\dt}{2}\bigl(\linop\vHat\sp{n+1}+\linop\vHat\sp{n}\bigr), \quad r\in(\aa,\bb) \qquad \vHat^{n+1}(\aa) = \aa\omega_b^{n+1}, \quad \vHat\sp{n+1}(\bb)=0,
\end{array}
\]
\label{alg:ampRBada}
\end{algorithm}

The stability analysis of the time-stepping scheme for the {\adannular}~model problem follows a similar approach to that discussed previously in Section~\ref{sec:stabilityMP-AD}. 
Therefore, the discussion here is brief and attention is given primarily to the elements of the analysis that are new.
 We consider the homogeneous problem and set $g_\omega(t)=0$, $\vHat\sp{n}=A\sp{n}\vHat\sp{0}$, $\omega_b\sp{n}=A\sp{n}\omega_b\sp{0}$ and $b_\omega\sp{n}=A\sp{n}b_\omega\sp{0}$, where $A$ is an amplification factor. Furthermore, assign $\vHat\sp\psup=A\sp{n+1}\bar v\sp\psup$, etc., similar to the definitions in~\eqref{eq:uamps} and~\eqref{eq:intermediate}.  The solutions of the boundary-value problems in Steps~2 and~5 take the form
\begin{equation}
\bar v\sp\psup(r)=\aa\bigl(\bar\omega_b\sp\esup-\bar\omega_b\sp\psup\bigr)\phi_1(r)+\bar v\sp{0}(r),\qquad \bar v\sp{0}(r)=\aa\bar\omega_b\sp\psup\phi_2(r),
\label{eq:annulussolutions}
\end{equation}
where
\[
\phi_m(r)={I_1(\zeta_mr)K_1(\zeta_m\bb)-I_1(\zeta_m\bb)K_1(\zeta_mr)\over I_1(\zeta_m\aa)K_1(\zeta_m\bb)-I_1(\zeta_m\bb)K_1(\zeta_m\aa)},\qquad\hbox{$m=1$ or 2,}
\]
and
\[
\zeta_1=\sqrt{{2\over\nu\dt}},\qquad \zeta_2=\zeta_1\sqrt{{A-1\over A+1}}.
\]
Here, $I_1(z)$ and $K_1(z)$ are modified Bessel functions of the first and second kind of order one.
The solutions for $\bar v\sp{0}(r)$ and $\bar v\sp\psup(r)$ in~\eqref{eq:annulussolutions} are
analogous to the discrete solutions in~\eqref{eq:uhatsoln} and~\eqref{eq:ubarsoln} found previously.

The solutions for $\bar v\sp{0}(r)$ and $\bar v\sp\psup(r)$ may now be used in the AMP interface conditions in Steps~3 and~6 to eliminate the circumferential velocities of the fluid so that the remaining steps of the algorithm involve the angular velocities and accelerations of the body alone.  The analysis of the these remaining steps leads to linear systems analogous to those in~\eqref{eq:sysA}, \eqref{eq:sysB} and~\eqref{eq:sysC}, and then to a homogeneous system of the form
\newcommand{\Mt}{{\widetilde M}}
\begin{equation}
\Mt\left[\begin{array}{c}
\bar \omega_b\sp{0} \medskip\\
\bar b_\omega\sp{0}
\end{array}\right]=0.
\label{eq:sysDa}
\end{equation}
Nontrivial solutions of~\eqref{eq:sysDa} are obtained if the determinant of the $2\times2$ matrix $\Mt$ is zero which leads 
to the following theorem.
\begin{theorem}
   The \ampRB~algorithm without the velocity-correction Step 8,
given in Algorithm~\ref{alg:ampRBada} for the model problem~{\adannular}, is stable in the 
sense of Godunov and Ryabenkii if and only if there are no roots with $|A|>1$ to the following equation
\begin{equation}
  \NAbt \eqdef \tilde\gamma_b(A-1)\sp3+\tilde\gamma_0\bigl[\IbBar(A-1)+\ctwo(A+1)\bigr]A\sp2=0,
\label{eq:adconstraintAnnular}
\end{equation}
which is similar to the constraint in~\eqref{eq:adconstraint} derived for the {\adplanar}~model problem.  The coefficients in~\eqref{eq:adconstraintAnnular} are given by
\[
\tilde\gamma_b \eqdef (2\beta_d\bar{\cal D}\sp{\omega}-\ctwo)(2\beta_d\bar{\cal D}\sp{\omega}-\cone),\qquad \tilde\gamma_0 \eqdef \IbBar+4\beta_d\bar{\cal D}\sp{\omega}-\cone,
\]
where $\bar{\cal D}\sp{\omega}$ and $\IbBar$ are a dimensionless added-damping coefficient and moment of inertia of the body, respectively, given by
\begin{align*}
  & \bar{\cal D}\sp{\omega} \eqdef{\Delta r\Duua\over\mu(2\pi\aa)\aa\sp2}=1-e\sp{-\tilde\delta}, \\
  &  \IbBar \eqdef {I_b\over\rho(2\pi\aa)\dr \, \aa\sp2} \, {\tilde\delta\sp2\over2}, 
\end{align*}
and $\cone$ and $\ctwo$ are dimensionless Dirichlet-to-Neumann transfer coefficients given by
\[
\cmmm \eqdef - \Delta r \left[rD_{rh}\left({\phi_m(r)\over r}\right)\right]_{r=\aa}={\aa\over2}\left({3\phi_m(\aa)\over\aa}-{4\phi_m(\aa+\Delta r)\over\aa+\Delta r}+{\phi_m(\aa+2\Delta r)\over\aa+2\Delta r}\right),\quad\text{$m=1$ or $2$}.
\]
\end{theorem}

The constraint equation for $A$ in~\eqref{eq:adconstraintAnnular} applies for Algorithm~\ref{alg:ampRBada} without the optional fluid-velocity correction step.  However, as observed in Section~\ref{sec:stabilityMP-AD} for the {\adplanar}~model problem, including this step improves the stability of the AMP algorithm and then the amplification factor is determined by the following. 
\begin{theorem}
   The \ampRB~algorithm including the velocity-correction Step 8,
given in Algorithm~\ref{alg:ampRBada} for the model problem~{\adannular}, is stable in the 
sense of Godunov and Ryabenkii if and only if there are no roots with $|A|>1$ to the following equation
\begin{equation}
  \NAvt  \eqdef \tilde\gamma_v(A-1)\sp3+\tilde\gamma_0\bigl[\IbBar(A-1)+\ctwo(A+1)\bigr]A\sp2=0,
\label{eq:advconstrainta}
\end{equation}
where $\tilde\gamma_b$ in~\eqref{eq:adconstraintAnnular} is replaced by
\[
\tilde\gamma_v=(2\beta_d\bar{\cal D}\sp{\omega}-\cone)\sp2.
\]
in~\eqref{eq:advconstrainta} while the other terms are unchanged. 
\end{theorem}

{
\newcommand{\labelFont}{\scriptsize}
\newcommand{\figWidth}{8.5cm}
\newcommand{\trimfig}[2]{\trimw{#1}{#2}{.0}{.0}{.0}{.0}}
\begin{figure}[htb]
\begin{center}
\begin{tikzpicture}[scale=1]
  \useasboundingbox (0.0,1) rectangle (16,6.5);  
  \begin{scope}[xshift=-.25cm]
  \draw(-.05,0) node[anchor=south west,xshift=-4pt,yshift=+0pt] {\trimfig{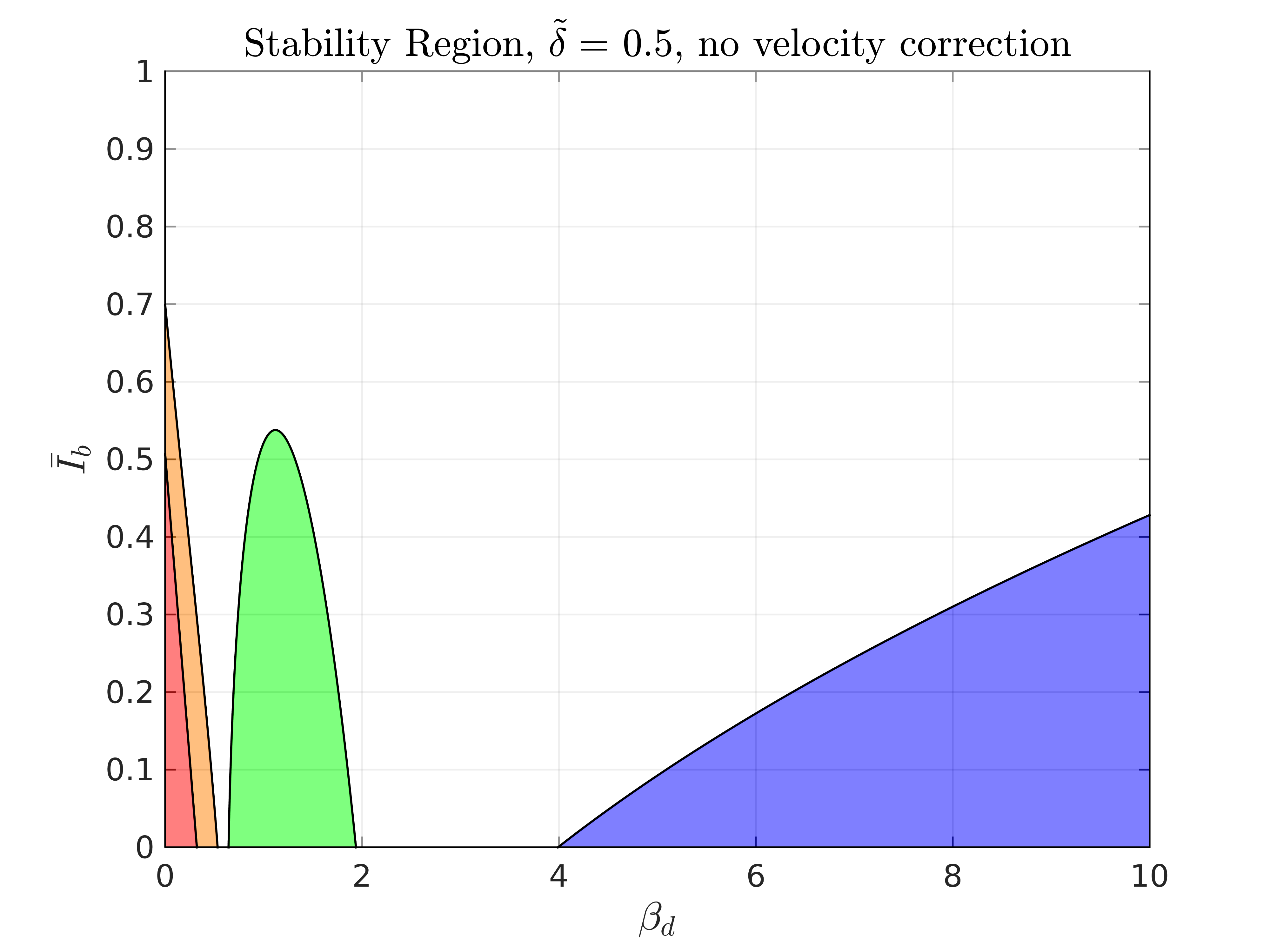}{\figWidth}};
  \draw(7.75,0) node[anchor=south west,xshift=-4pt,yshift=+0pt] {\trimfig{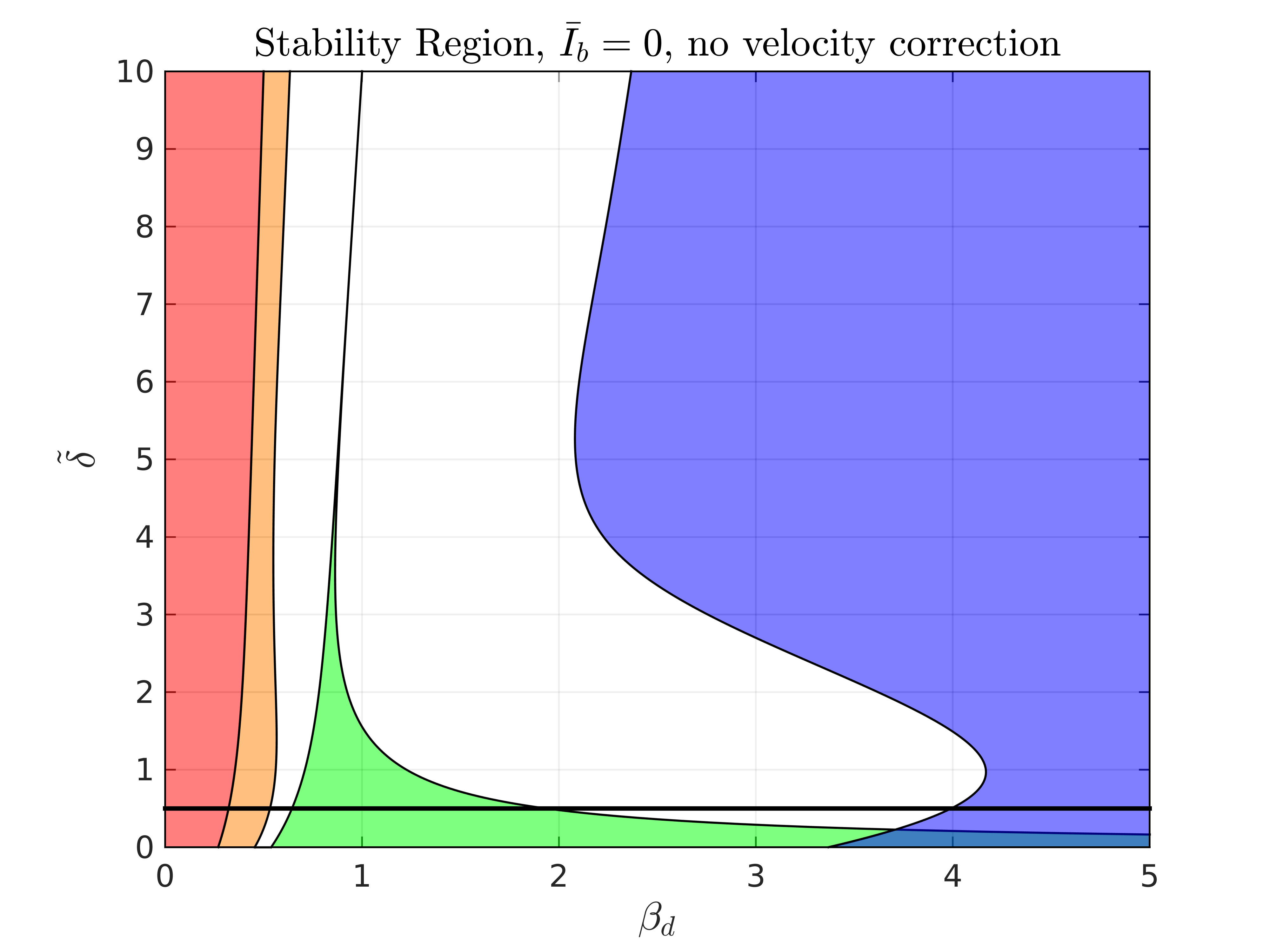}{\figWidth}};
    \draw(4.5,5.0) node[draw,fill=white] {\labelFont{Stable region}};
    \draw(6.,2.5) node[draw,fill=blue!50] {\labelFont IV};
    \draw(3.2,2.5) node[draw,fill=green!50] {\labelFont III};
    \draw(1.9,4.3) node[draw,fill=orange!50] {\labelFont II};
    \draw(.75,1.5) node[draw,fill=red!50] {\labelFont I};
    \draw(11.9,1.75) node[draw,fill=white] {\labelFont{Stable region}}; 
  \end{scope}
%
\end{tikzpicture}
\end{center}
\caption{
Stability regions for the \ampRB~scheme without velocity correction for the added-damping model problem in an annular geometry.
Left: stability region in the $\adc$--$\IbBar$ plane for $\deltat=.5$; 
shaded unstable regions are marked according to the form of the instability. 
Right: stability region in the $\adc$--$\deltat$ plane for the zero-mass body.
Here $\aa = 20 \sqrt{\nu \dt/2} $ and $\bb = 2 \, \aa$.}
\label{fig:addeddampingStablityV1Rotate}
\end{figure}
}

{
\newcommand{\labelFont}{\scriptsize}
\newcommand{\figWidth}{8.5cm}
\newcommand{\figWidthz}{4.0cm}
\newcommand{\trimfig}[2]{\trimw{#1}{#2}{.0}{.0}{.0}{.0}}
\newcommand{\trimfigz}[2]{\trimw{#1}{#2}{.05}{.05}{.03}{.05}}
\begin{figure}[htb]
\begin{center}
\begin{tikzpicture}[scale=1]
  \useasboundingbox (0.0,1) rectangle (16,6.5);  
  \begin{scope}[xshift=-.25cm]
  \draw(-.05,0) node[anchor=south west,xshift=-4pt,yshift=+0pt] {\trimfig{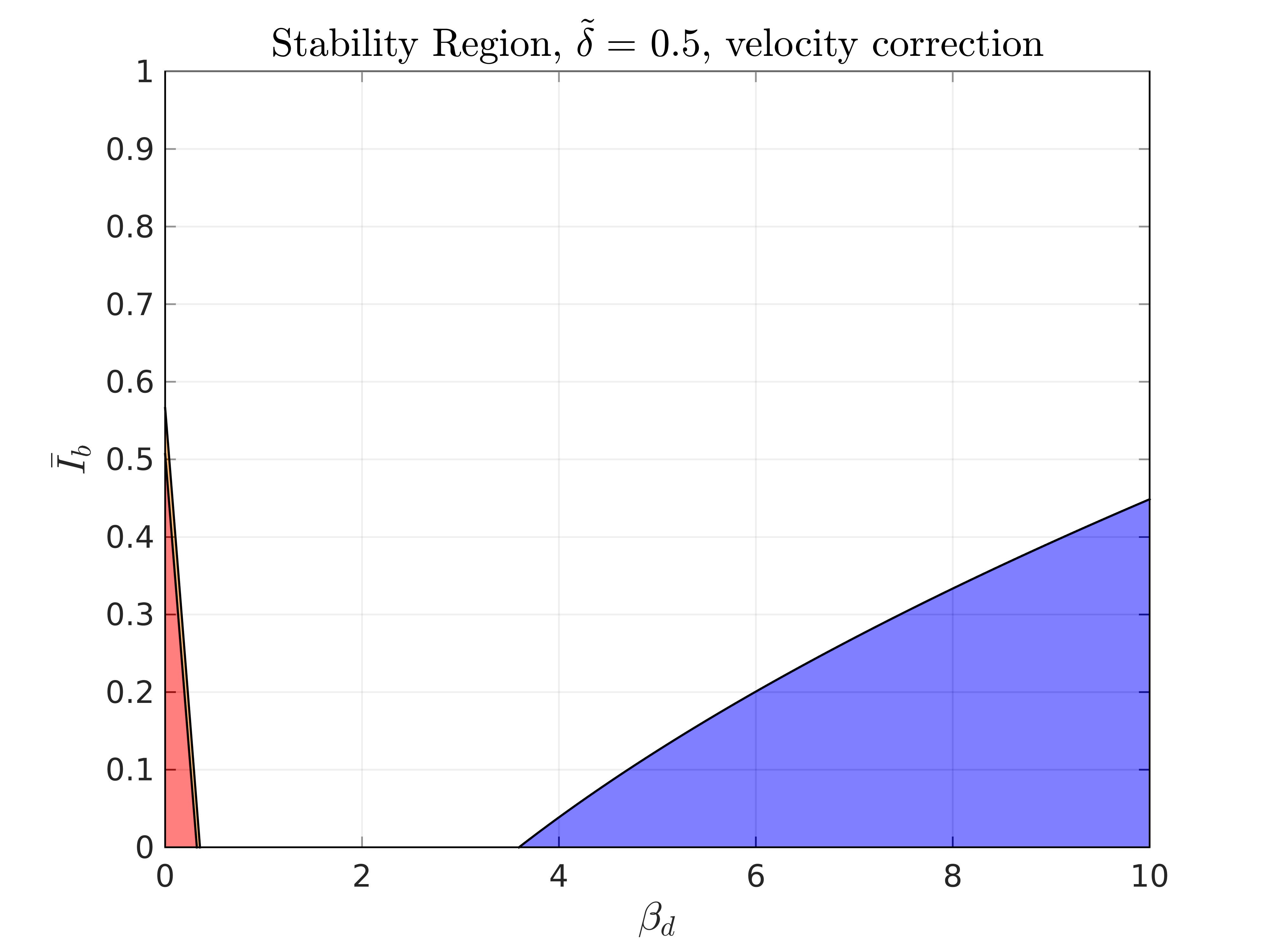}{\figWidth}};
  \draw(1.45,2.) node[anchor=south west,xshift=-4pt,yshift=+0pt] {\trimfigz{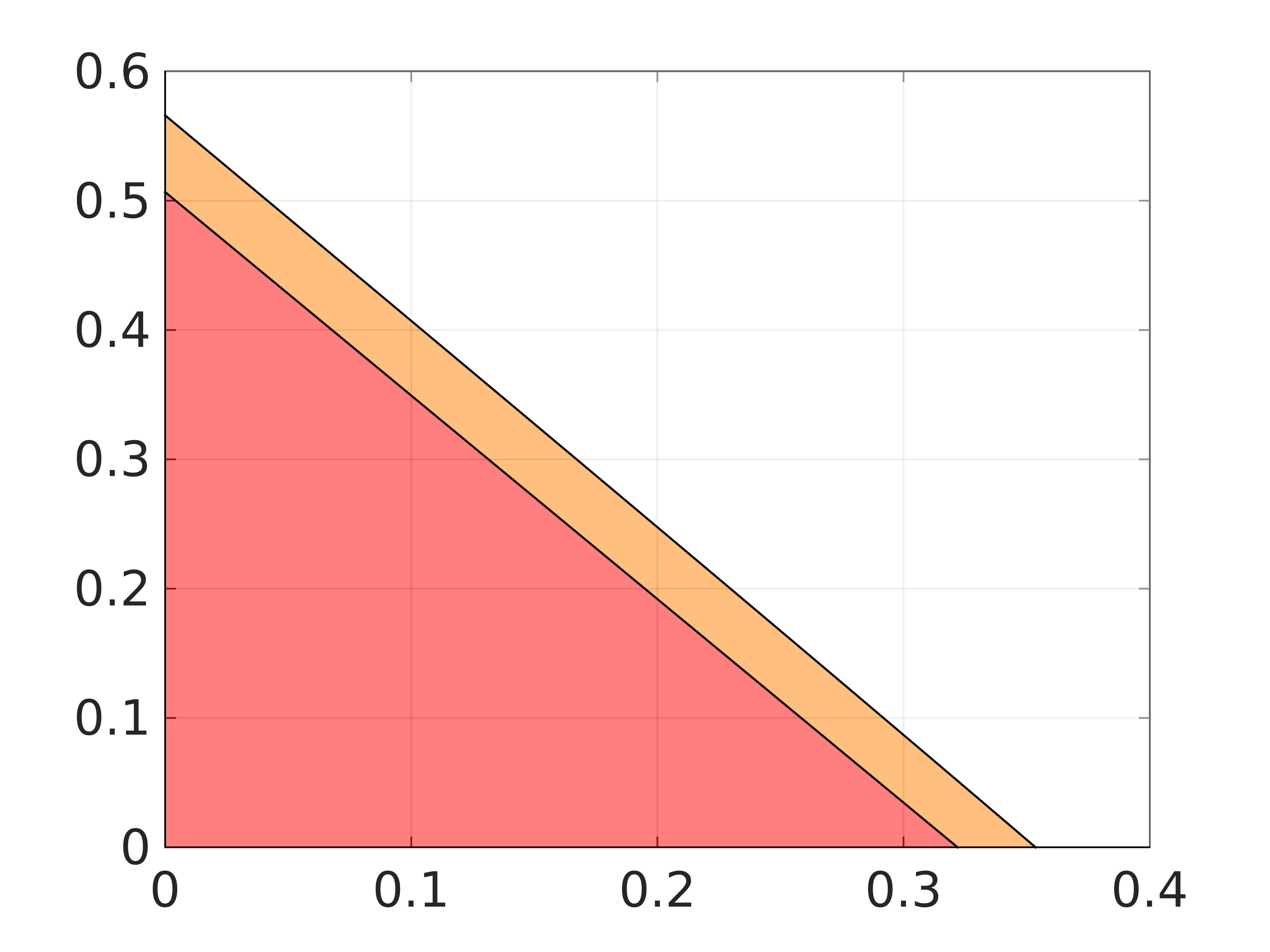}{\figWidthz}};
  \draw(7.75,0) node[anchor=south west,xshift=-4pt,yshift=+0pt] {\trimfig{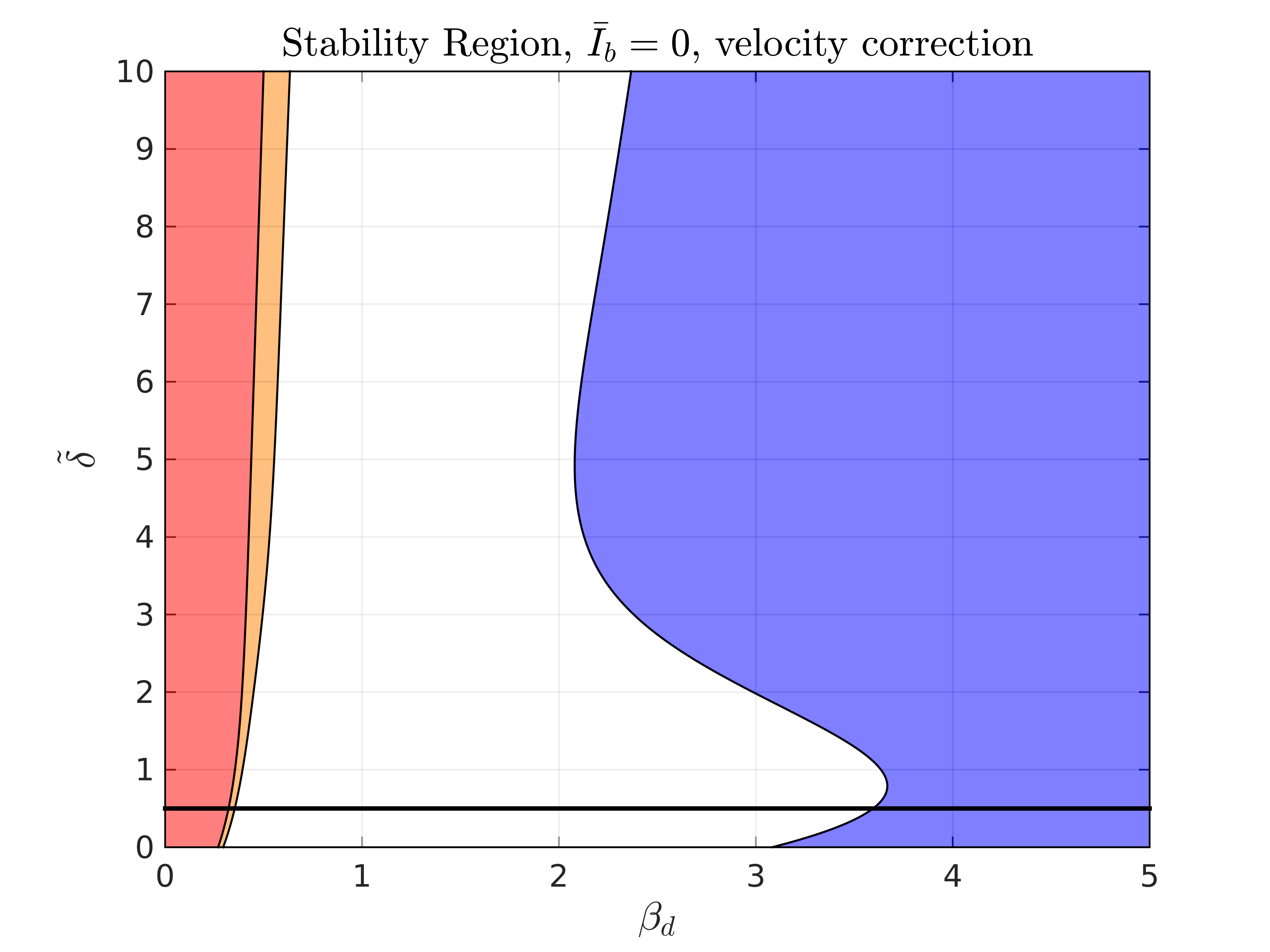}{\figWidth}};
    \draw(5,5.65) node[draw,fill=white] {\labelFont{Stable region}};
    \draw(6.5,2.9) node[draw,fill=blue!50] {\labelFont IV};
    \draw(4.3,4.3) node[draw,fill=orange!50] {\labelFont II};
    \draw(2.15,1.75) node[draw,fill=red!50] {\labelFont I};
    \draw(11.2,1.75) node[draw,fill=white] {\labelFont{Stable region}};
  \end{scope}
%
\end{tikzpicture}
\end{center}
\caption{
	Stability regions for the \ampRB~scheme with velocity correction for the added-damping model problem in an annular geometry.
Left: stability region in the $\adc$--$\IbBar$ plane for $\deltat=.5$; 
shaded unstable regions are marked according to the form of the instability. 
Right: stability region in the $\adc$--$\deltat$ plane for the zero-mass body.
Here $\aa= 20 \sqrt{\nu \dt/2}  $ and $\bb = 2 \, \aa$.}
  \label{fig:addeddampingStablityV2Rotate}
\end{figure}
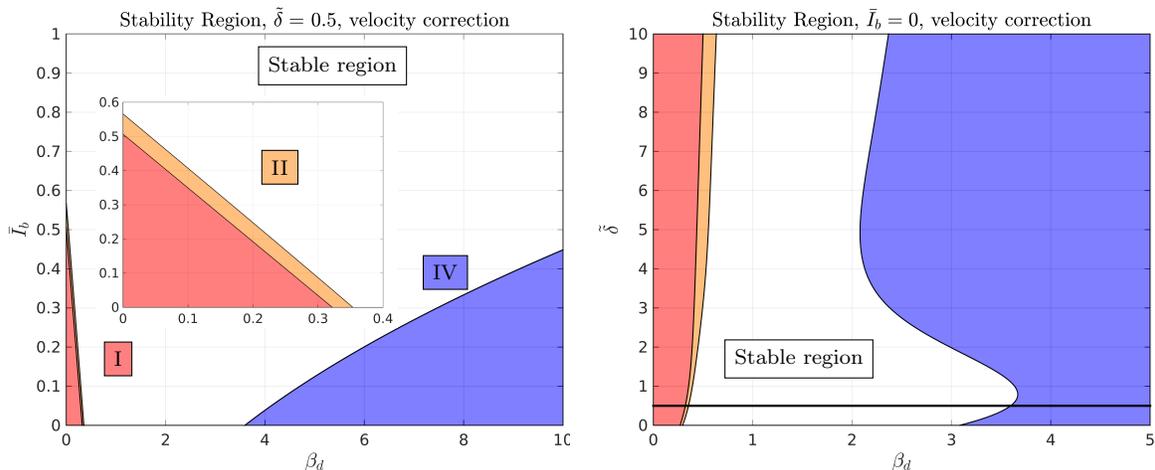
}

The stability regions for the annular case
are not as straightforward to investigate as the prior rectangular case because~\eqref{eq:adconstraintAnnular} 
and~\eqref{eq:advconstrainta}  are transcendental functions, and so the number of roots is not known a-priori. Nonetheless,
the use of the argument principle reveals that there are a finite number of roots with $|A|>1$ for a given choice
of parameters $(\IbBar,\deltat,\adc)$. Using this information we are able to identify locations where the number of
unstable roots changes, and from these points we perform a simple numerical continuation to trace out the 
stability boundaries.

Plots of the stability region for the~\ampRB~scheme without the velocity-correction step
are shown in Figure~\ref{fig:addeddampingStablityV1Rotate}, 
while those for the scheme with velocity-correction are shown in Figure~\ref{fig:addeddampingStablityV2Rotate}.
Comparing to the corresponding plots for the rectangular geometry in 
Figures~\ref{fig:addedDampingNvcModes} and~\ref{fig:addedDampingVcModes}, it is clear that 
the stability regions for the annular geometry are very similar in form to those for the rectangular geometry.
For example, the \ampRB-NVC scheme has four regions of instability, and the 
stability region for the \ampRB-VC scheme shows similar improvements as in the rectangular geometry
with the troublesome instability region III no longer present.
Perhaps the most important conclusion to be drawn from Figure~\ref{fig:addeddampingStablityV2Rotate}
is that there is again a large range for $\adc$, near $\adc=1$, where the \ampRB-VC
is stable for any non-negative value of the scaled moment of inertia $\IbBar$. This allows one to
pick $\adc=1$ and obtain a robustly stable FSI algorithm. Finally we note that the traditional-partitioned 
scheme corresponds to $\adc=0$, and has a similar behaviour to that found 
for the rectangular geometry.

{
\newcommand{\figWidth}{8.6cm}
\newcommand{\figWidthz}{3.25cm}
\newcommand{\trimfig}[2]{\trimw{#1}{#2}{.03}{.05}{.0}{.0}}
\newcommand{\trimfigz}[2]{\trimFigb{#1}{#2}{.15}{.25}{.1}{.15}}
\begin{figure}[htb]
\begin{center}
\begin{tikzpicture}[scale=1]
  \useasboundingbox (0.5,.75) rectangle (16.,6.8);  
  \begin{scope}[xshift=-.8cm,yshift=-.3cm]
    \draw(0,0) node[anchor=south west,xshift=-4pt,yshift=+0pt] {\trimfig{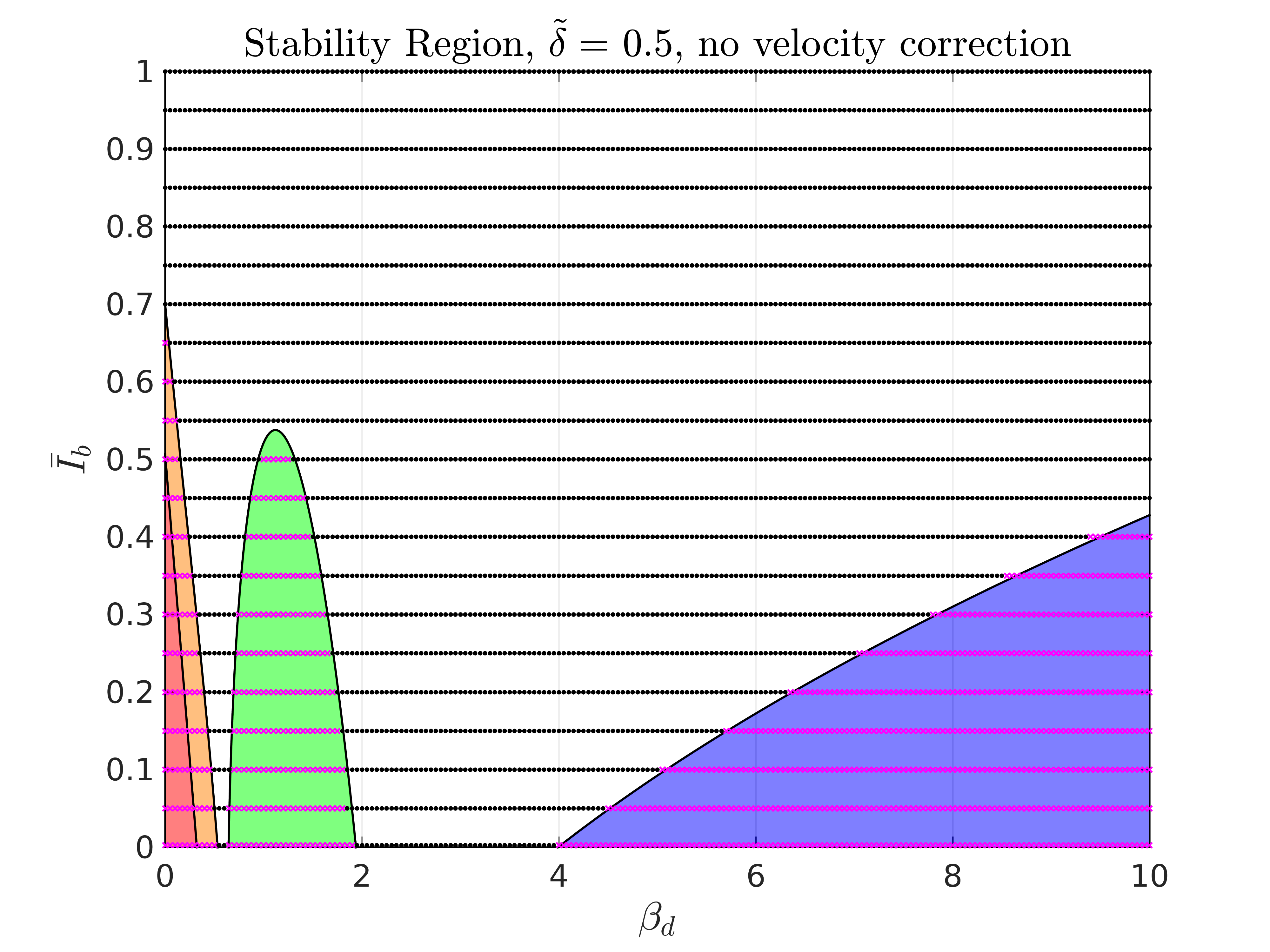}{\figWidth}};
    \draw(1.5,3.8) node[anchor=south west,xshift=-4pt,yshift=+0pt] {\trimfigz{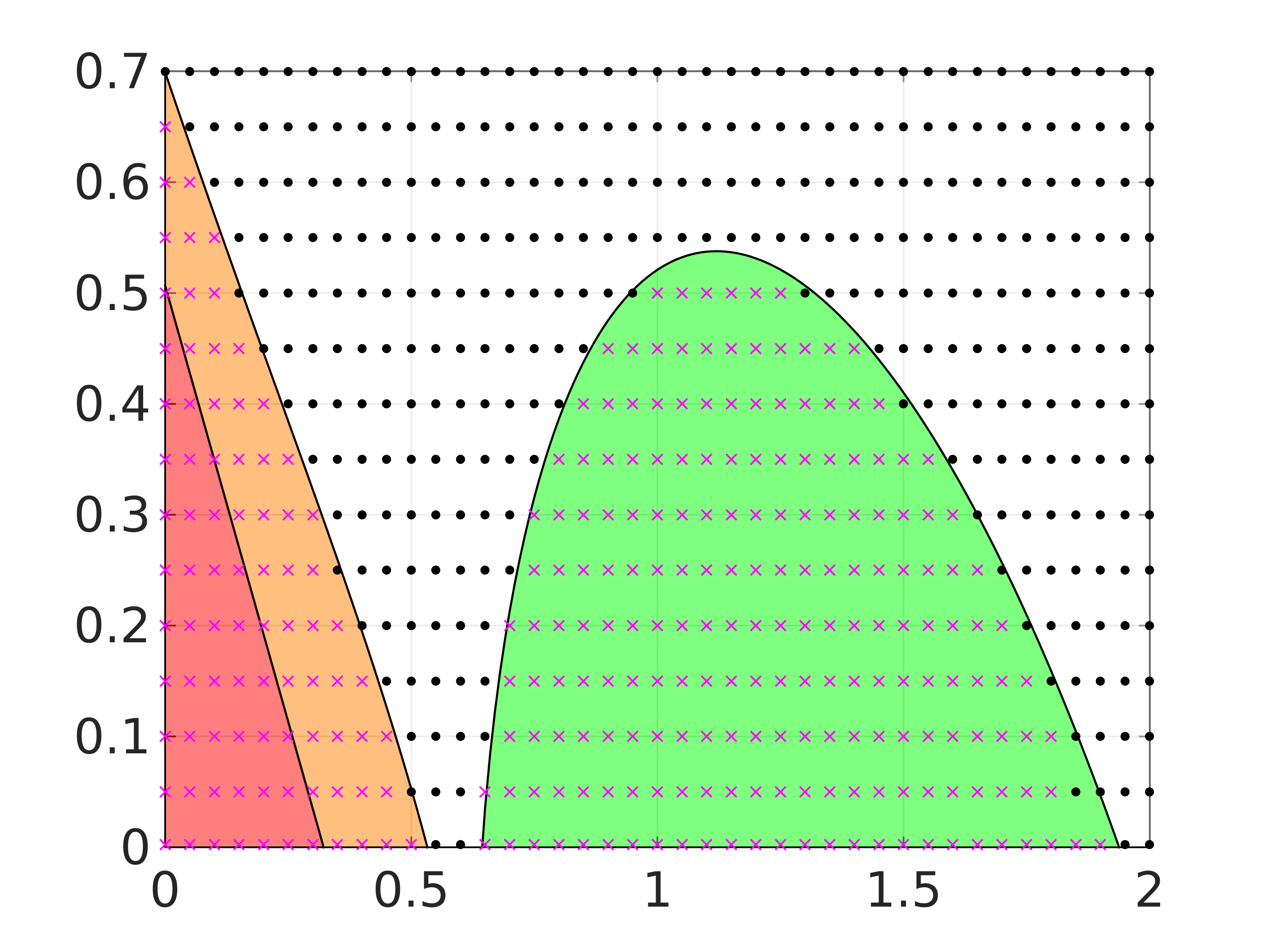}{\figWidthz}};
    \draw(5.1,3.8) node[anchor=south west,xshift=-4pt,yshift=+0pt] {\trimfigz{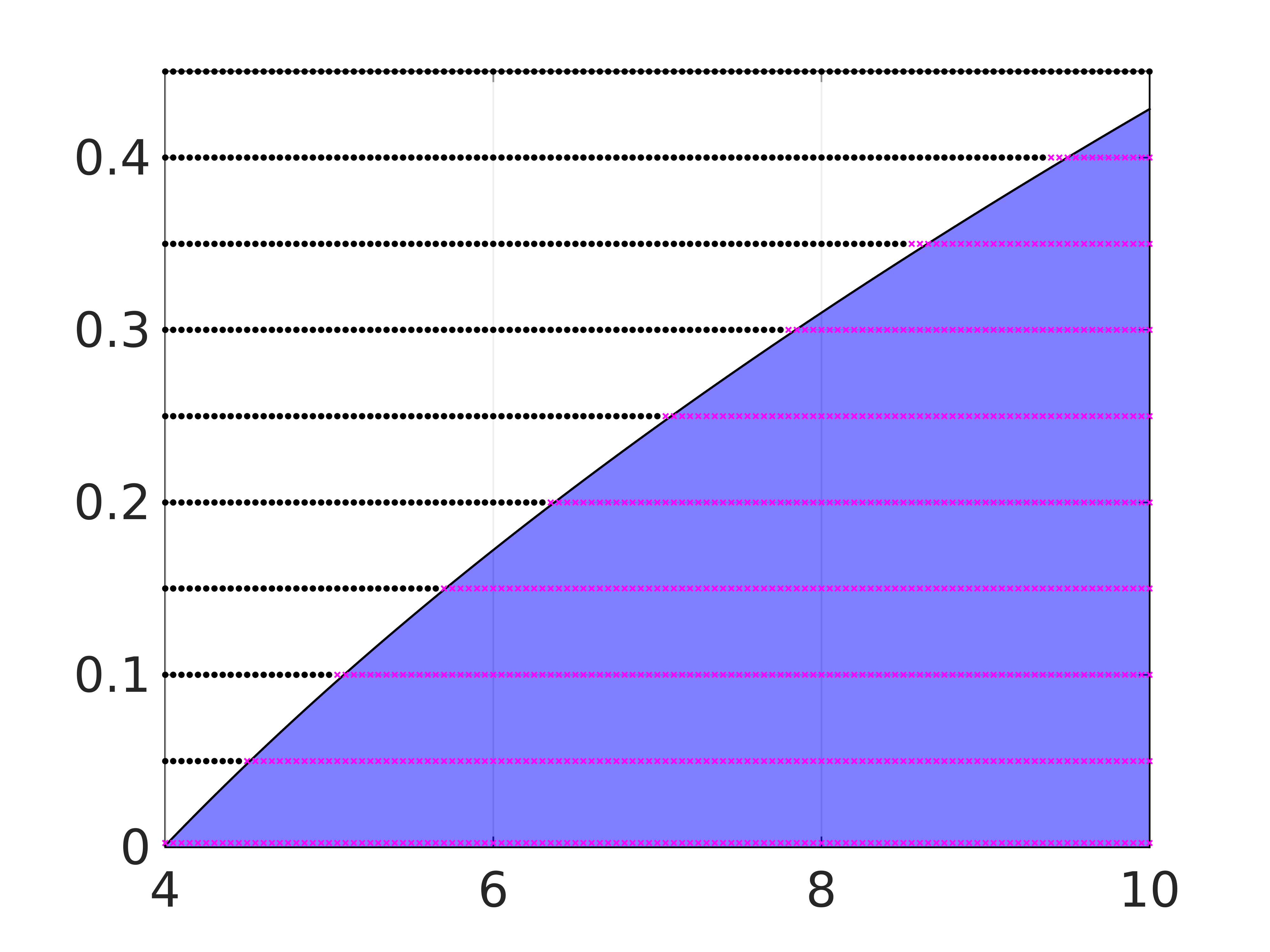}{\figWidthz}};
  \end{scope}
%
  \begin{scope}[xshift=7.8cm,yshift=-.3cm]
    \draw(0,0) node[anchor=south west,xshift=-4pt,yshift=+0pt] {\trimfig{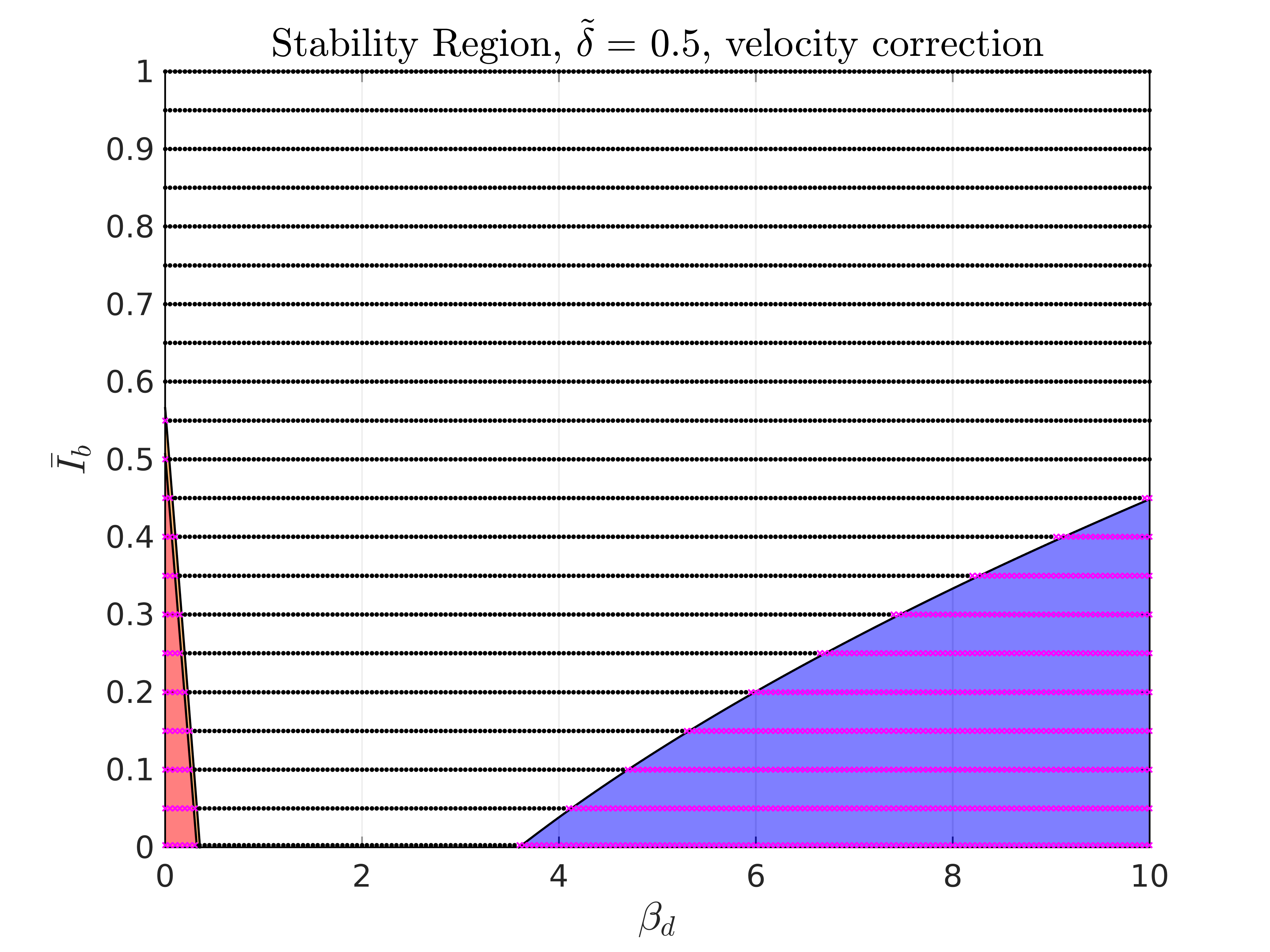}{\figWidth}};
    \draw(1.5,3.8) node[anchor=south west,xshift=-4pt,yshift=+0pt] {\trimfigz{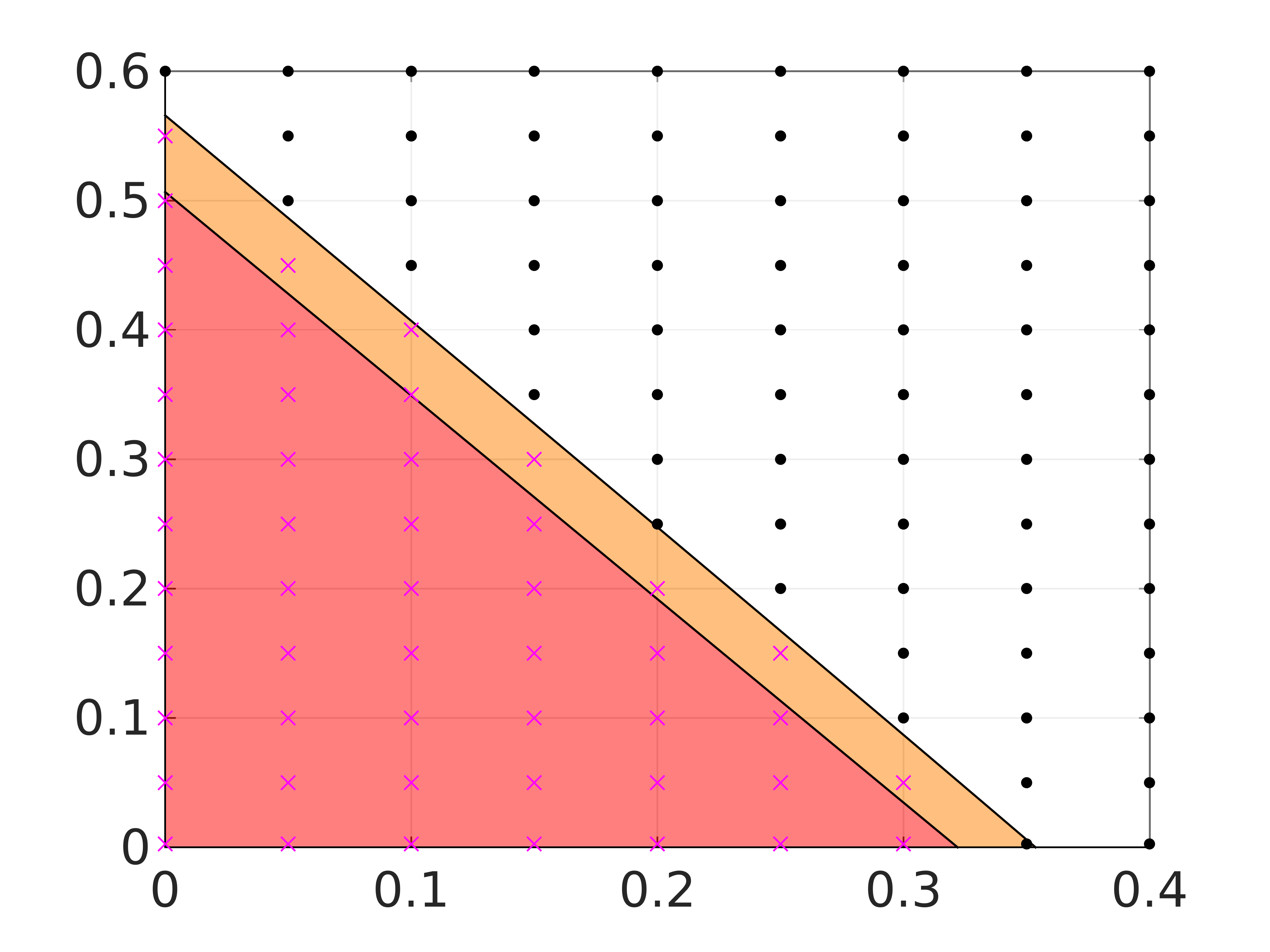}{\figWidthz}};
    \draw(5.1,3.8) node[anchor=south west,xshift=-4pt,yshift=+0pt] {\trimfigz{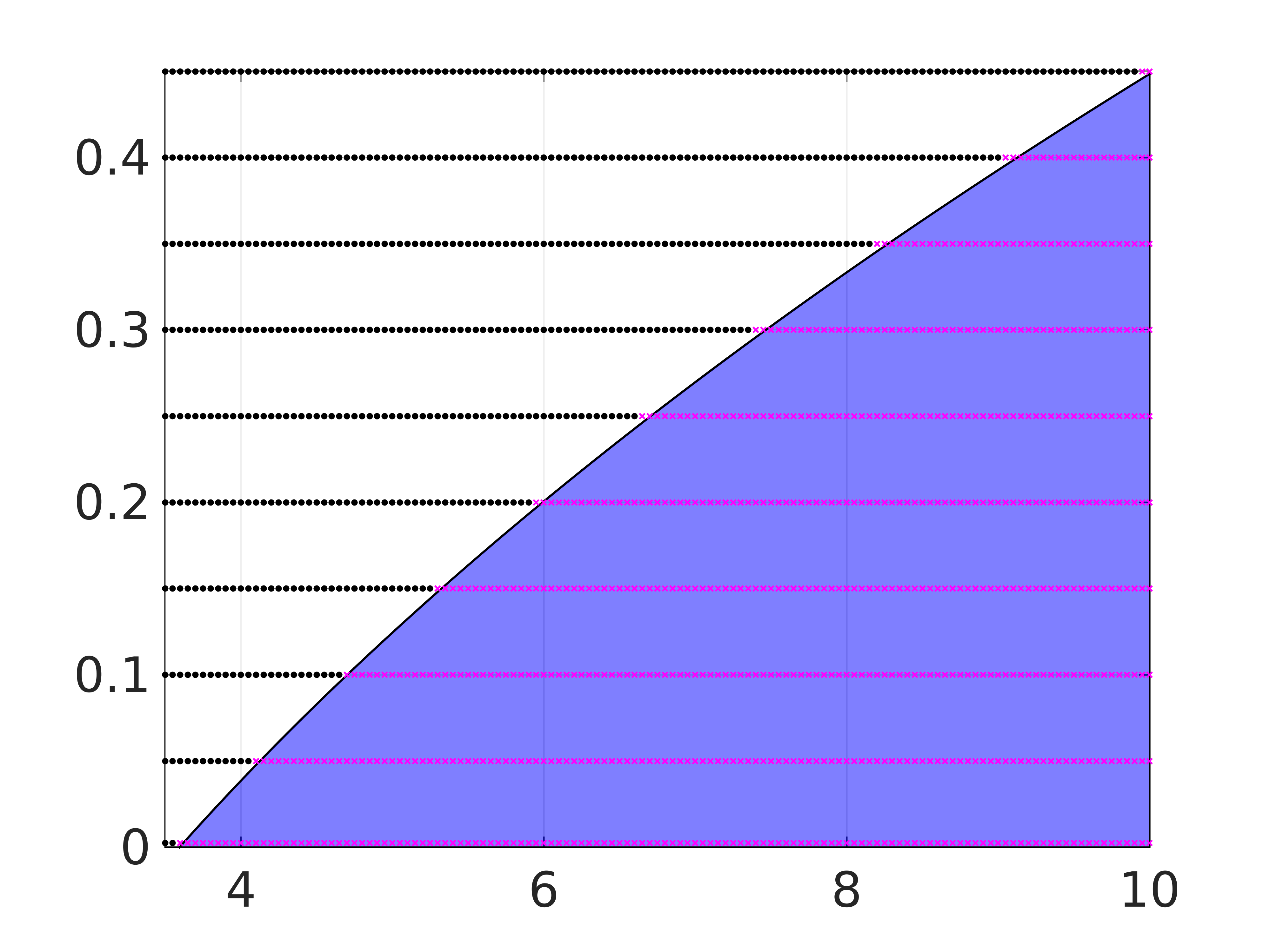}{\figWidthz}};
  \end{scope}
%
\end{tikzpicture}
\end{center}
  \caption{Numerically determined points of stability (black dots for stable and magenta crosses for unstable) compared to the theoretical stability region
     for the added-damping model problem in an annular geometry.
  Left:the~\ampRB-NVC scheme. Right:  the \ampRB-VC scheme.
}
  \label{fig:numericalStabilityRotatingDisk}
\end{figure}
}

Figure~\ref{fig:numericalStabilityRotatingDisk} shows values of the parameters where
the~\ampRB~schemes are numerically determined to be stable or unstable for the annular geometry.  
The numerical stability results are computed in a similar fashion to those for the
rectangular grid case (see Section~\ref{sec:stabilityMP-AD}).
The computational grid is chosen to be the curvilinear grid $\Gs^{(4)}$, as defined in Section~\ref{sec:rotatingDisk},
and the time-step is taken as $\dt = 0.05$. 
In this annular geometry case, the theoretical results used continuous spatial differential operators except
in the interface conditions, while the numerical results use a fully discrete approximation. 
Despite this difference, the points where the numerical scheme is stable or unstable are in excellent agreement with the theory.

\section{Numerical verification} \label{sec:numericalVerification}

In this section, numerical simulations are used to confirm the results of the stability analyses, and
to demonstrate the second-order accuracy of the {\ampRB} scheme by comparison to some exact solutions. 
Calculations with rigid bodies having a
range of masses are presented, but the focus is on light or zero-mass bodies where 
added-damping and added-mass effects are most significant. 
We note that the zero-mass case is representative of cases in which the mass of the body is very
small since the results are nearly identical.

The numerical results are organized into four sections.
Section~\ref{sec:pistonMotion} considers the one-dimensional motion of a rigid rectangular {\em piston}
moving into an adjacent incompressible fluid, This case corresponds to the added-mass  model problem {\amplanar}, 
and the results demonstrate stability and second-order accuracy for a range of rigid-body masses. 
Then in Section~\ref{sec:shearBlock}, results for the one-dimensional added-damping model problem {\adplanar}
of a rectangular body sliding against an incompressible fluid
are presented, and again the results show stability and second-order accuracy for a range of rigid body
masses. 
In addition, by adjusting the added damping parameter, $\adc$, we are able to numerically 
demonstrate the existence of the different instability modes that were predicted by the analysis discussed in
Section~\ref{sec:stabilityMP-AD}.
Subsequently in Sections~\ref{sec:translatingDisk}
and~\ref{sec:rotatingDisk}, the stability and accuracy of the {\ampRB}~scheme in a two-dimensional
annular geometry are demonstrated for two problems involving a translating disk and a rotating disk. In these last two cases the motion of the disk in the numerical solution is free to translate and rotate, and thus both added-mass and added-damping effects can be important and are tested in the scheme.

\newcommand{\massRatio}{M_r}
\subsection{Piston motion of a rigid body and an incompressible fluid} \label{sec:pistonMotion}

Consider the added-mass model problem in a rectangular geometry ({\amplanar})
with fluid domain $\Omegaf=[0,L]\times[0,H]$, where $L = H = 1$, and the fluid density is $\rho = 1$.
The fluid domain is discretized using a uniform grid, denoted by 
$\Gc^{(j)}$, with time-step and grid spacing given by 
$\dt^{(j)}=\dy^{(j)} = 1/(10 j)$. 
The fluid pressure boundary condition at $y=H$ is $p(H,t) = p_H(t)$,
and there is no external body forcing, i.e., $g_v(t) = 0$. The rigid body has 
mass $m_b$, which is varied in the numerical results.
The exact solution to this system is given in $\eqref{eq:vbsoln}$--$\eqref{eq:vHatsoln}$,
and the applied pressure $p_H(t)$ is chosen from~\eqref{eq:vbsoln}  
 so that the rigid body undergoes a specified motion given by 
 $y_b(t) = \frac{1}{4}\sin(2 \pi t)$.
 
Figure~\ref{fig:pistonCompare} shows the position, velocity and acceleration of the rigid body 
computed on the grid $\Gc^{(4)}$ using the~\ampRB~scheme with $\mb = 1$ and $\mb = 0$. 
The exact solution is independent of the rigid-body mass, and so the two numerical
approximations are nearly identical.  Figure~\ref{fig:pistonCompare} also presents the errors
in the approximations of the position, velocity and acceleration of the rigid body for $\Gc^{(4)}$.

{
\newcommand{\figWidth}{7.5cm}
\newcommand{\figWidths}{7.3cm}
\newcommand{\trimfig}[2]{\trimFig{#1}{#2}{.0}{.0}{.0}{.0}}
\begin{figure}[htb]
\begin{center}
\resizebox{14cm}{!}{
\begin{tikzpicture}[scale=1]
  \useasboundingbox (0.0,0) rectangle (16.,6);  
  \draw(0,-.75) node[anchor=south west,xshift=0pt,yshift=+0pt] {\trimfig{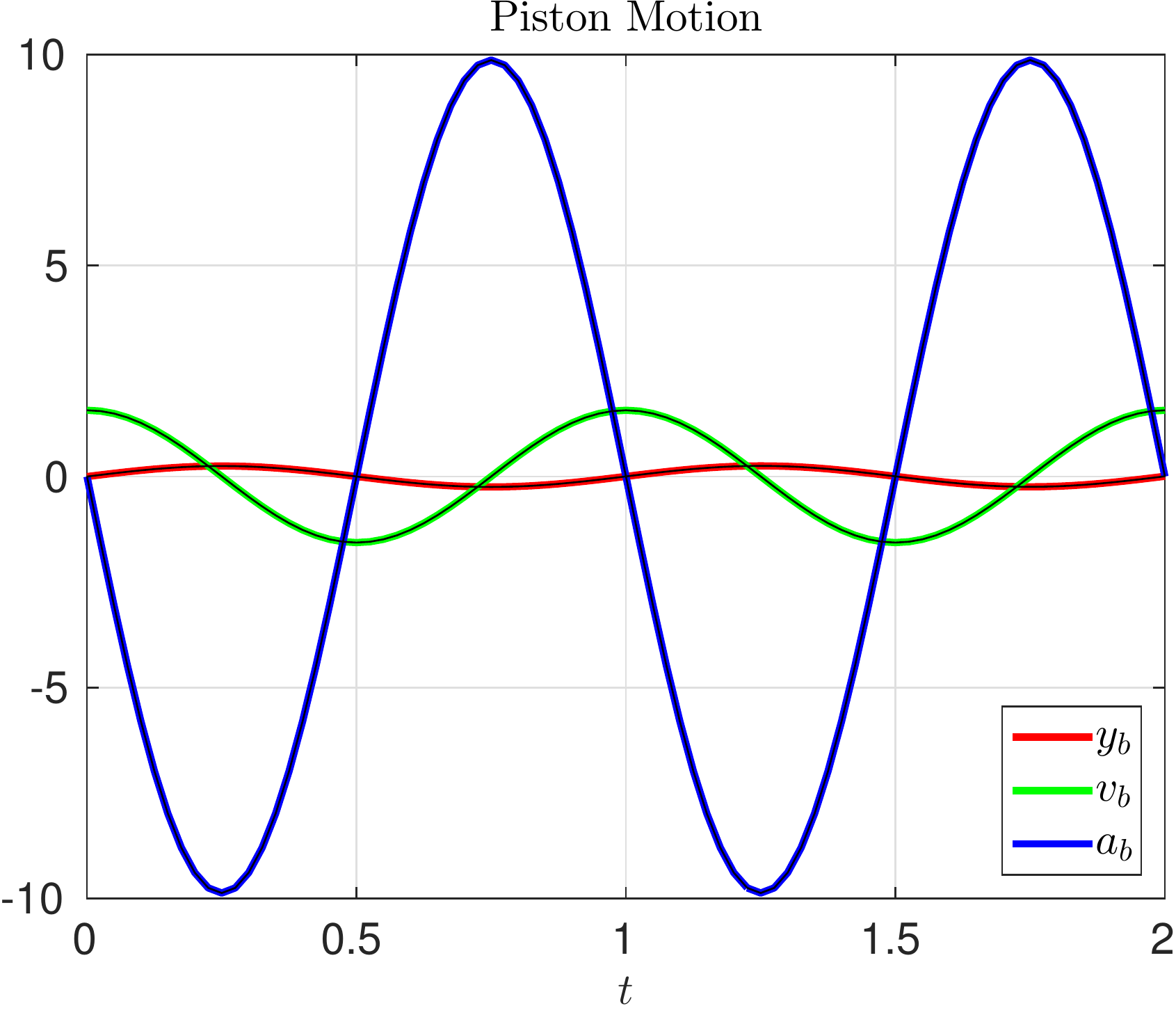}{\figWidth}};
  \draw(8.15,-.75) node[anchor=south west,xshift=0pt,yshift=+0pt] {\trimfig{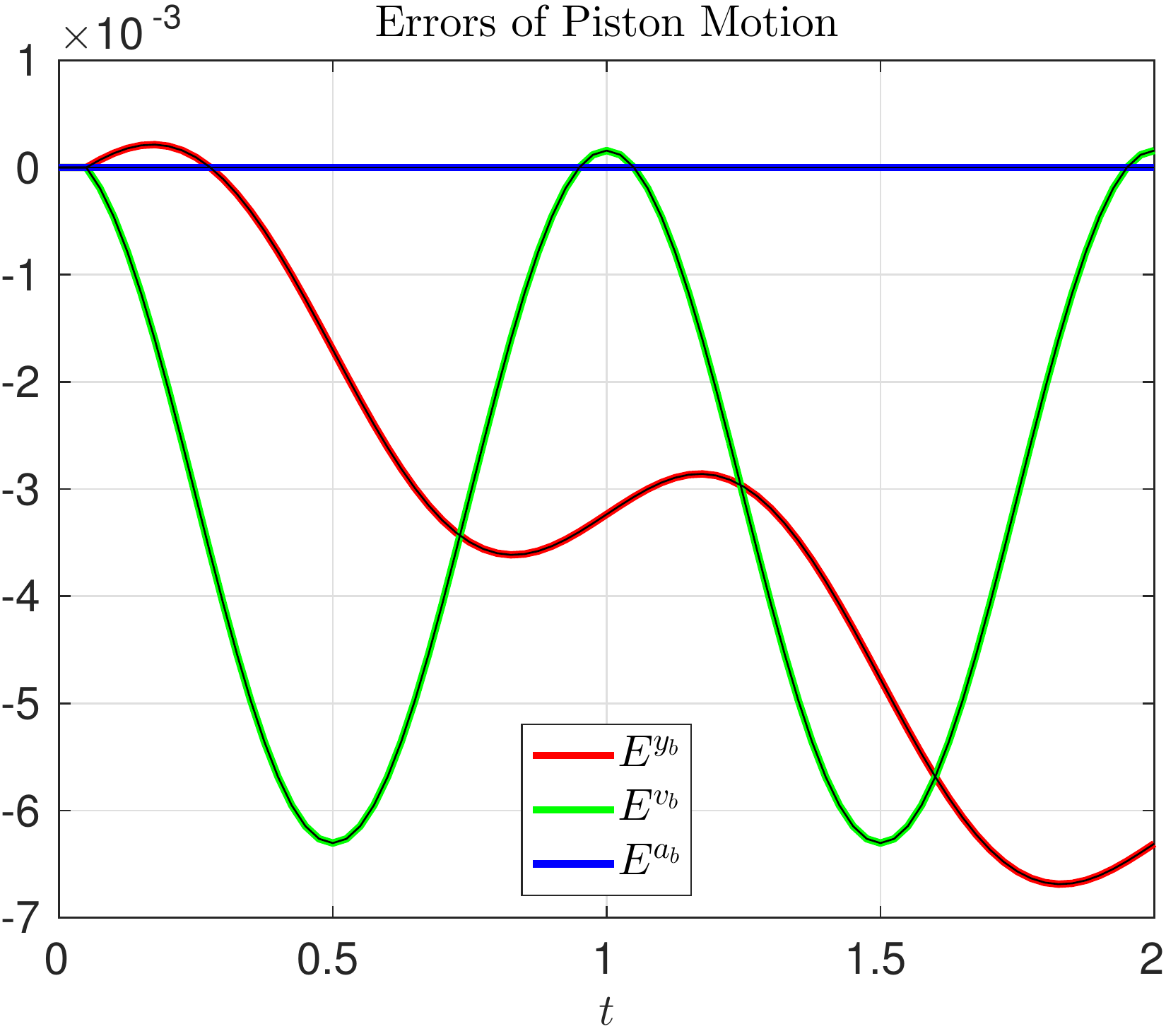}{\figWidths}};
\end{tikzpicture}
}
\end{center}
  \caption{Piston motion.  Stable solutions and errors for the motion of a medium, $\mb=1$, and massless, $\mb=0$ (black lines), rigid body.  The grid used is $\Gc^{(4)}$.}
  \label{fig:pistonCompare}
\end{figure}
}

A grid convergence study is performed for the~\ampRB~scheme, and the results
are presented in Table~\ref{table:piston}. 
In the tables,  $E_j^{p}$, $E_j^{\vv}$, etc., denote the maximum-norm errors in $p$, $\vv$, etc., on grid $\Gc^{(j)}$. 
The convergence rates in the final row of each table are computed using a least squares fit to the log of the errors. 
Three cases are presented in Table~\ref{table:piston}, corresponding to 
a heavy body, $\mb=10$, medium body, $\mb=1$, and massless body, $\mb=0$. In all cases, second-order accuracy
is demonstrated for the position and velocity of the body and the velocity of the fluid, while the fluid pressure and the acceleration of the body are exact
to machine precision, as predicted by the theory. Note also that the observed errors in the position
and velocities are nearly identical for all rigid-body masses because the exact solution is the same for all cases.
Finally note that the second-order accurate {\tpRB} scheme is found to be unstable when the added-mass
ratio $\massRatio \lesssim 1$, which 
matches the theory for the {\tpRB} scheme as presented in Section~\ref{sec:stability}.
 
{ 
\newcommand{\errp}{$E_j^{p}$}
\newcommand{\erruv}{$E_j^{\vv}$}
\newcommand{\errxA}{$E_j^{\xvb}$}
\newcommand{\errvA}{$E_j^{\vvb}$}
\newcommand{\erraA}{$E_j^{\avb}$}
%
\newcommand{\pistonTableI}{%
\begin{tabular}{|c|c|c|c|c|c|c|c|c|c|c|} \hline 
  \multicolumn{11}{|c|}{Piston motion, $\mb=10$} \\ \hline 
\strutt$\dy^{(j)}$&     \errp     &  r   &     \erruv     &  r   &     \errxA     &  r   &     \errvA     &  r   &     \erraA     &  r    \\[3pt] \hline 
1/20 & \num{3.1}{-13} &  & \num{6.5}{-3} &      & \num{1.4}{-2} &      & \num{6.5}{-3} &      & \num{3.0}{-14} &       \\ \hline 
1/40 & \num{8.8}{-13} &  & \num{2.1}{-3} &  3.1 & \num{3.6}{-3} &  3.8 & \num{2.1}{-3} &  3.1 & \num{8.7}{-14} &   \\ \hline 
1/80 & \num{2.6}{-13} &  & \num{5.5}{-4} &  3.8 & \num{9.0}{-4} &  4.0 & \num{5.5}{-4} &  3.8 & \num{2.5}{-14} &   \\ \hline 
1/160& \num{2.5}{-12} &  & \num{1.4}{-4} &  3.9 & \num{2.2}{-4} &  4.0 & \num{1.4}{-4} &  3.9 & \num{2.5}{-13} &   \\ \hline 
rate &        &       &         1.85 &       &         1.98 &       &         1.85 &       &         &      \\ \hline
\end{tabular}
}
\newcommand{\pistonTableII}{%
\begin{tabular}{|c|c|c|c|c|c|c|c|c|c|c|} \hline 
  \multicolumn{11}{|c|}{Piston motion, $\mb=1$} \\ \hline 
\strutt$\dy^{(j)}$&     \errp     &  r   &     \erruv     &  r   &     \errxA     &  r   &     \errvA     &  r   &     \erraA     &  r    \\[3pt] \hline 
1/20 & \num{7.1}{-15} &  & \num{6.5}{-3} &      & \num{1.4}{-2} &      & \num{6.5}{-3} &      & \num{1.8}{-15} &       \\ \hline 
1/40 & \num{5.2}{-14} &  & \num{2.1}{-3} &  3.1 & \num{3.6}{-3} &  3.8 & \num{2.1}{-3} &  3.1 & \num{5.2}{-14} &  \\ \hline 
1/80 & \num{2.8}{-13} &  & \num{5.5}{-4} &  3.8 & \num{9.0}{-4} &  4.0 & \num{5.5}{-4} &  3.8 & \num{2.1}{-13} &  \\ \hline 
1/160& \num{1.4}{-13} &  & \num{1.4}{-4} &  3.9 & \num{2.2}{-4} &  4.0 & \num{1.4}{-4} &  3.9 & \num{1.2}{-13} &  \\ \hline 
rate &         &       &         1.85 &       &         1.98 &       &         1.85 &       &         &      \\ \hline
\end{tabular}
}
\newcommand{\pistonTableIV}{%
\begin{tabular}{|c|c|c|c|c|c|c|c|c|c|c|} \hline 
  \multicolumn{11}{|c|}{Piston motion, $\mb=0$} \\ \hline 
\strutt$\dy^{(j)}$&     \errp     &  r   &     \erruv     &  r   &     \errxA     &  r   &     \errvA     &  r   &     \erraA     &  r    \\[3pt] \hline 
1/20 & \num{1.8}{-15} &   & \num{6.5}{-3} &      & \num{1.4}{-2} &      & \num{6.5}{-3} &      & \num{0.0}{0} &   \\ \hline 
1/40 & \num{3.1}{-14} &   & \num{2.1}{-3} &  3.1 & \num{3.6}{-3} &  3.8 & \num{2.1}{-3} &  3.1 & \num{6.0}{-14} &   \\ \hline 
1/80 & \num{2.6}{-14} &   & \num{5.5}{-4} &  3.8 & \num{9.0}{-4} &  4.0 & \num{5.5}{-4} &  3.8 & \num{9.4}{-14} &   \\ \hline 
1/160& \num{4.8}{-14} &   & \num{1.4}{-4} &  3.9 & \num{2.2}{-4} &  4.0 & \num{1.4}{-4} &  3.9 & \num{1.5}{-13} &   \\ \hline 
rate &				 &       &         1.85 &      &         1.98  &       &         1.85 &       &          &    \\ \hline
\end{tabular}
}
{
\begin{table}[hbt]\tableFont 
\begin{center}
  \pistonTableI \\
\bigskip
  \pistonTableII \\
\bigskip	
  \pistonTableIV
\caption{Piston motion. Maximum errors and estimated convergence rates at $t=0.8$ computed using the~\ampRB~scheme
   for a heavy, $\mb=10$, medium, $\mb=1$, and massless, $\mb=0$, rigid body. 
    The column labeled ``r'' provides the ratio of the errors at the current grid spacing to that on the next coarser grid.
 }
\label{table:piston}
\end{center}
\end{table}
}
} 

\subsection{A rigid body sliding against an incompressible fluid} \label{sec:shearBlock}

Turning now to the situation where added-damping effects are dominant,
consider the  added-damping model problem ({\adplanar}) for the same square domain used
in the previous case. For this problem the fluid boundary condition at $y=H$ 
is taken as a no-slip condition $u(H,t)=0$, and there is no external forcing, $g_u(t)=0$, as before.
In this case, an exact solution is given by
\begin{align*}
 &   u = -\frac{ \sin( \lambda(y-H))}{\sin(\lambda H)} \, e^{-\nu \lambda^2 t},
\end{align*}
where $\lambda$ is an eigenvalue determined by a solution of the transcendental equation
\begin{align*}
 &   \lambda H \, {\tan(\lambda H)} = \frac{\rho L H} {\mb} = \frac{1}\massRatio. 
\end{align*}
Note that there are infinitely many solutions to this equation, and the numerical tests use the smallest
positive eigenvalue, as given in Table~\ref{table:lambdaH} for a few values of $\massRatio$.

\begin{table}[H]\tableFont 
  \begin{centre}
    \begin{tabular}{|c|c|c|c|c|}\hline
&			 $\massRatio=10$	  & $\massRatio=1$	   & $\massRatio=10^{-3}$ & $\massRatio=0$ \\ \hline 
$\lambda H$  & 0.311052848200298  &  0.860333589019380 & 1.56922710098197     & 1.57079632679489  \\ \hline 
    \end{tabular}
  \caption{Scaled eigenvalues, $\lambda H$, in the exact solution of the sliding block problem for different values of $\massRatio$.}
  \label{table:lambdaH}
  \end{centre}
\end{table}

{
\newcommand{\figWidth}{7.5cm}
\newcommand{\figWidths}{7.25cm}
\newcommand{\trimfig}[2]{\trimFig{#1}{#2}{.0}{.0}{.0}{.0}}
\begin{figure}[htb]
\begin{center}
\resizebox{14cm}{!}{
\begin{tikzpicture}[scale=1]
  \useasboundingbox (0, 0) rectangle (16,5.5);  
  \draw(0,-0.75) node[anchor=south west,xshift=0pt,yshift=0pt] {\trimfig{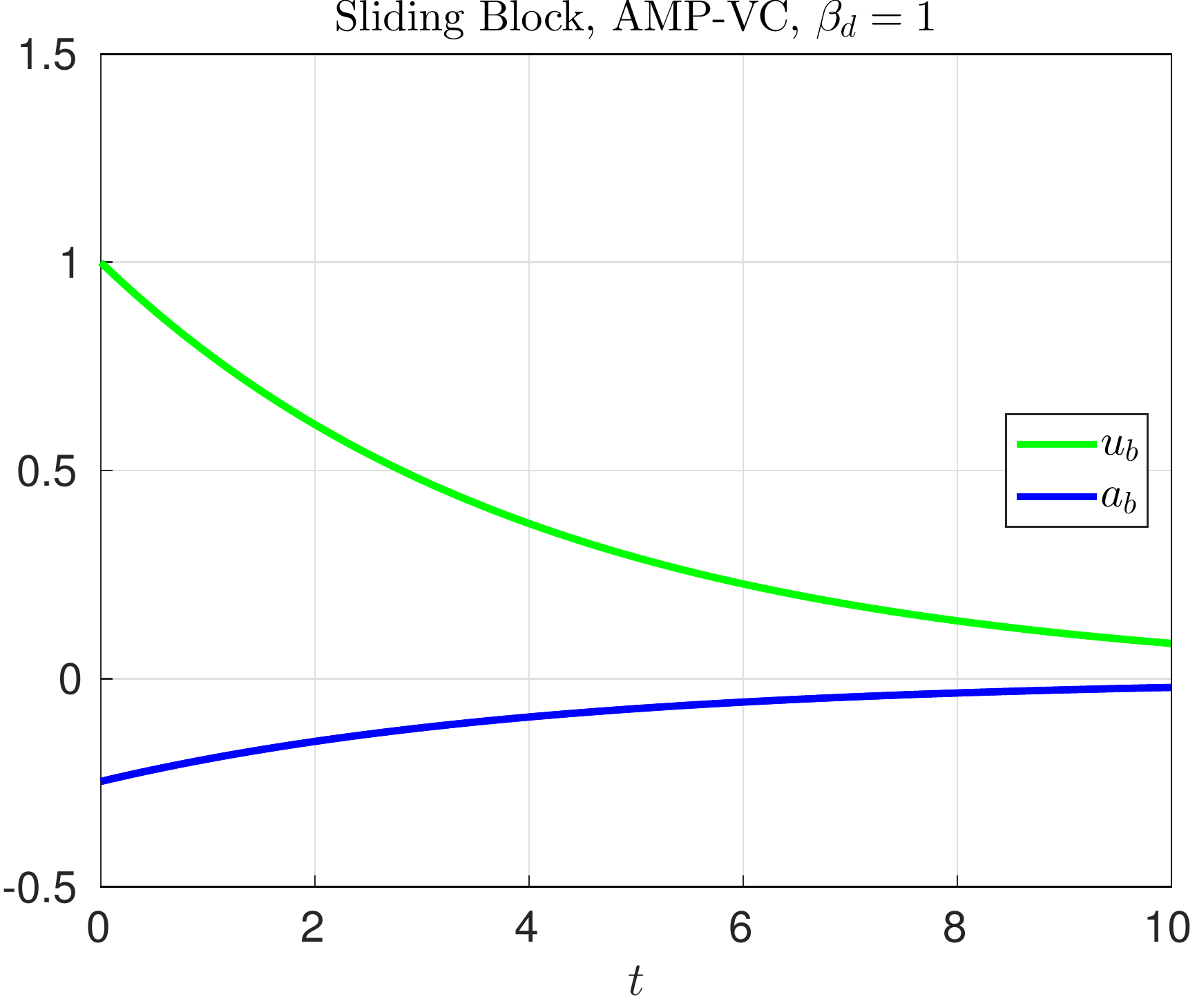}{\figWidth}};
  \draw(8.2,-0.75) node[anchor=south west,xshift=0pt,yshift=0pt] {\trimfig{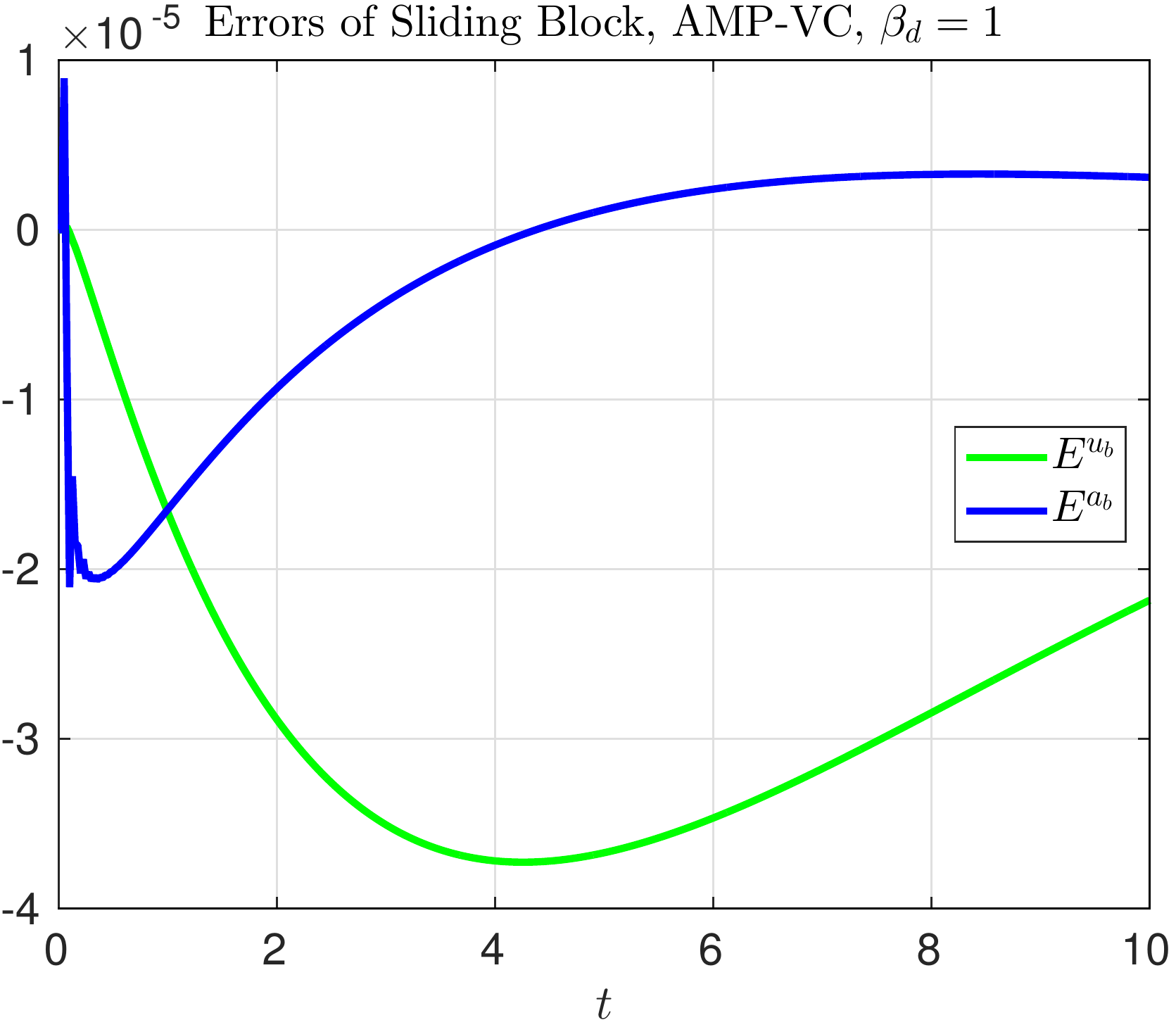}{\figWidths}};
\end{tikzpicture}
}
\end{center}
  \caption{Sliding block.  Stable solutions and errors for the motion of a massless, $\mb=0$, body.
   The~\ampRB-VC scheme on the grid $\Gs^{(4)}$ is used with $\adc=1$.
    }
\label{fig:shearBlockStable}
\end{figure}
}

{ 
\newcommand{\errp}{$E_j^{p}$}
\newcommand{\erruv}{$E_j^{\vv}$}
\newcommand{\errxA}{$E_j^{\xvb}$}
\newcommand{\errvA}{$E_j^{\vvb}$}
\newcommand{\erraA}{$E_j^{\avb}$}
%
\newcommand{\shearBlockTableI}{%
\begin{tabular}{|c|c|c|c|c|c|c|} \hline 
  \multicolumn{7}{|c|}{Sliding block, $\mb=10$} \\ \hline 
\strutt$\dy^{(j)}$&     \erruv     &  r   &     \errvA     &  r   &     \erraA     &  r    \\[3pt] \hline 
1/20  & \num{7.0}{-7} &      & \num{7.0}{-7} &      & \num{7.4}{-7} &       \\ \hline 
1/40  & \num{1.9}{-7} &  3.8 & \num{1.9}{-7} &  3.8 & \num{1.9}{-7} &  4.0  \\ \hline 
1/80  & \num{4.6}{-8} &  4.0 & \num{4.6}{-8} &  4.0 & \num{4.6}{-8} &  4.0  \\ \hline 
1/160 & \num{1.2}{-8} &  4.0 & \num{1.2}{-8} &  4.0 & \num{1.2}{-8} &  4.0  \\ \hline 
1/320 & \num{2.9}{-9} &  4.0 & \num{2.9}{-9} &  4.0 & \num{2.9}{-9} &  4.0  \\ \hline 
1/640 & \num{7.4}{-10} &  4.0 & \num{7.4}{-10} &  4.0 & \num{7.2}{-10} &  4.0  \\ \hline 
rate & 1.98 &          & 1.98 &          & 2.00 &     \\ \hline
\end{tabular}
}
\newcommand{\shearBlockTableII}{%
\begin{tabular}{|c|c|c|c|c|c|c|} \hline 
  \multicolumn{7}{|c|}{Sliding block, $\mb=1$} \\ \hline 
\strutt$\dy^{(j)}$&     \erruv     &  r   &     \errvA     &  r   &     \erraA     &  r    \\[3pt] \hline 
1/20  & \num{3.3}{-5} &      & \num{3.3}{-5} &      & \num{3.0}{-5} &       \\ \hline 
1/40  & \num{8.5}{-6} &  3.9 & \num{8.5}{-6} &  3.9 & \num{7.2}{-6} &  4.1  \\ \hline 
1/80  & \num{2.1}{-6} &  4.1 & \num{2.1}{-6} &  4.1 & \num{1.8}{-6} &  4.0  \\ \hline 
1/160 & \num{5.2}{-7} &  4.0 & \num{5.2}{-7} &  4.0 & \num{4.4}{-7} &  4.0  \\ \hline 
1/320 & \num{1.3}{-7} &  4.0 & \num{1.3}{-7} &  4.0 & \num{1.1}{-7} &  4.0  \\ \hline 
1/640 & \num{3.3}{-8} &  4.0 & \num{3.3}{-8} &  4.0 & \num{2.8}{-8} &  4.0  \\ \hline 
rate & 1.99 &          & 1.99 &          & 2.01 &         \\ \hline
\end{tabular}
}
\newcommand{\shearBlockTableIV}{%
\begin{tabular}{|c|c|c|c|c|c|c|} \hline 
  \multicolumn{7}{|c|}{Sliding block, $\mb=0$} \\ \hline 
\strutt$\dy^{(j)}$&     \erruv     &  r   &     \errvA     &  r   &     \erraA     &  r    \\[3pt] \hline 
1/20  & \num{6.1}{-5} &      & \num{3.7}{-5} &      & \num{5.6}{-5} &       \\ \hline 
1/40  & \num{1.9}{-5} &  3.2 & \num{1.7}{-5} &  2.1 & \num{1.6}{-5} &  3.4  \\ \hline 
1/80  & \num{5.2}{-6} &  3.6 & \num{5.1}{-6} &  3.3 & \num{4.4}{-6} &  3.7  \\ \hline 
1/160 & \num{1.4}{-6} &  3.7 & \num{1.4}{-6} &  3.7 & \num{1.1}{-6} &  3.9  \\ \hline 
1/320 & \num{3.6}{-7} &  3.9 & \num{3.6}{-7} &  3.8 & \num{2.9}{-7} &  3.9  \\ \hline 
1/640 & \num{9.3}{-8} &  3.9 & \num{9.3}{-8} &  3.9 & \num{7.3}{-8} &  4.0  \\ \hline 
rate & 1.88 &          & 1.76 &          & 1.92 &							  \\ \hline
\end{tabular}
}
{
\begin{table}[htb]\tableFont 
\begin{center}
  \shearBlockTableI \quad 
  \shearBlockTableII \\
\bigskip
  \shearBlockTableIV
\caption{Sliding block. Maximum errors and estimated convergence rates at $t=1$ computed using the~\ampRB-VC~scheme
   for a heavy, $\mb=10$, medium, $\mb=1$, and massless, $\mb=0$, rigid body. }
\label{table:shearBlock}
\end{center}
\end{table}
}
} 

Figure~\ref{fig:shearBlockStable} shows the velocity and acceleration (and errors in the same) of 
a massless rigid body, $\mb=0$,  with  $\nu=0.1$ 
on grid $\Gs^{(4)}$ as computed by the~\ampRB-VC with $\adc=1$. 
Results from grid convergence studies using the \ampRB-VC scheme with $\adc = 1$ are presented in Table~\ref{table:shearBlock} and demonstrate second-order accuracy for massless, medium, and heavy rigid bodies. As in the prior convergence study the time step is taken as $\dt^{(j)} = \dy^{(j)} = 1/(10j)$. Also note that the second-order accurate~\tpRB~scheme is found to be unstable when $\mb \lesssim 0.653$ at $\delta = 0.5$, which agrees with the theoretical prediction of $\mb \le 0.656$ from~\eqref{eq:TPRBstabilityRegionApprox}.

{
\newcommand{\labelFont}{\scriptsize}
\newcommand{\trimfig}[2]{\trimw{#1}{#2}{.0}{.0}{.0}{.0}}
\newcommand{\trimfigz}[2]{\trimw{#1}{#2}{.05}{.05}{.03}{.05}}
\newcommand{\figWidth}{4.0cm}
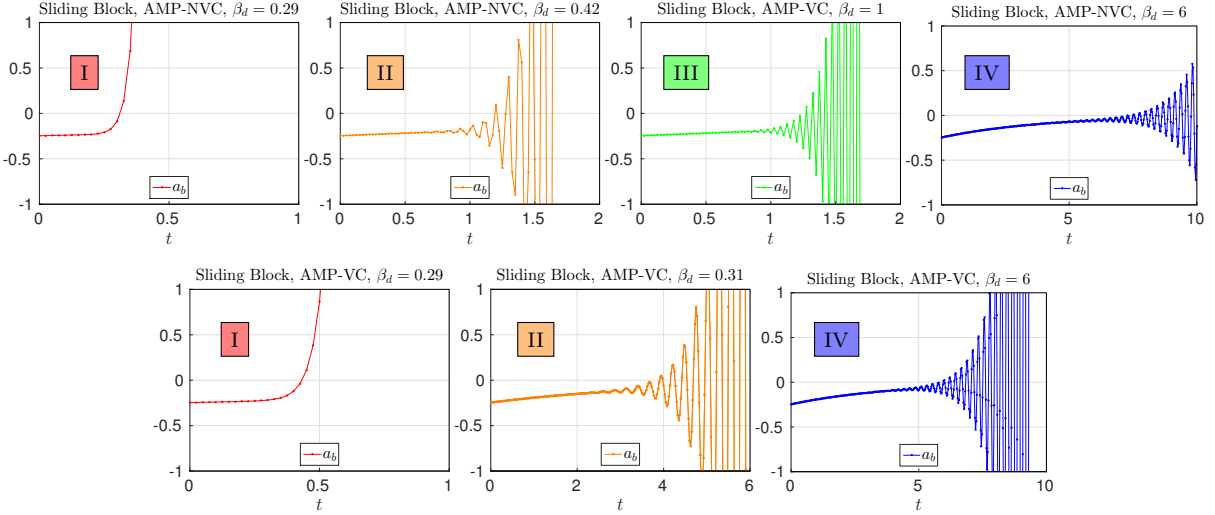
\begin{figure}[htb]
\begin{center}
\begin{tikzpicture}[scale=1]
  \useasboundingbox (0.0,1) rectangle (16,7.2);  
\begin{scope}[xshift=-.5cm,yshift=3.8cm]
    \draw(0.0,0.0) node[anchor=south west,xshift=-4pt,yshift=+0pt] {\trimfig{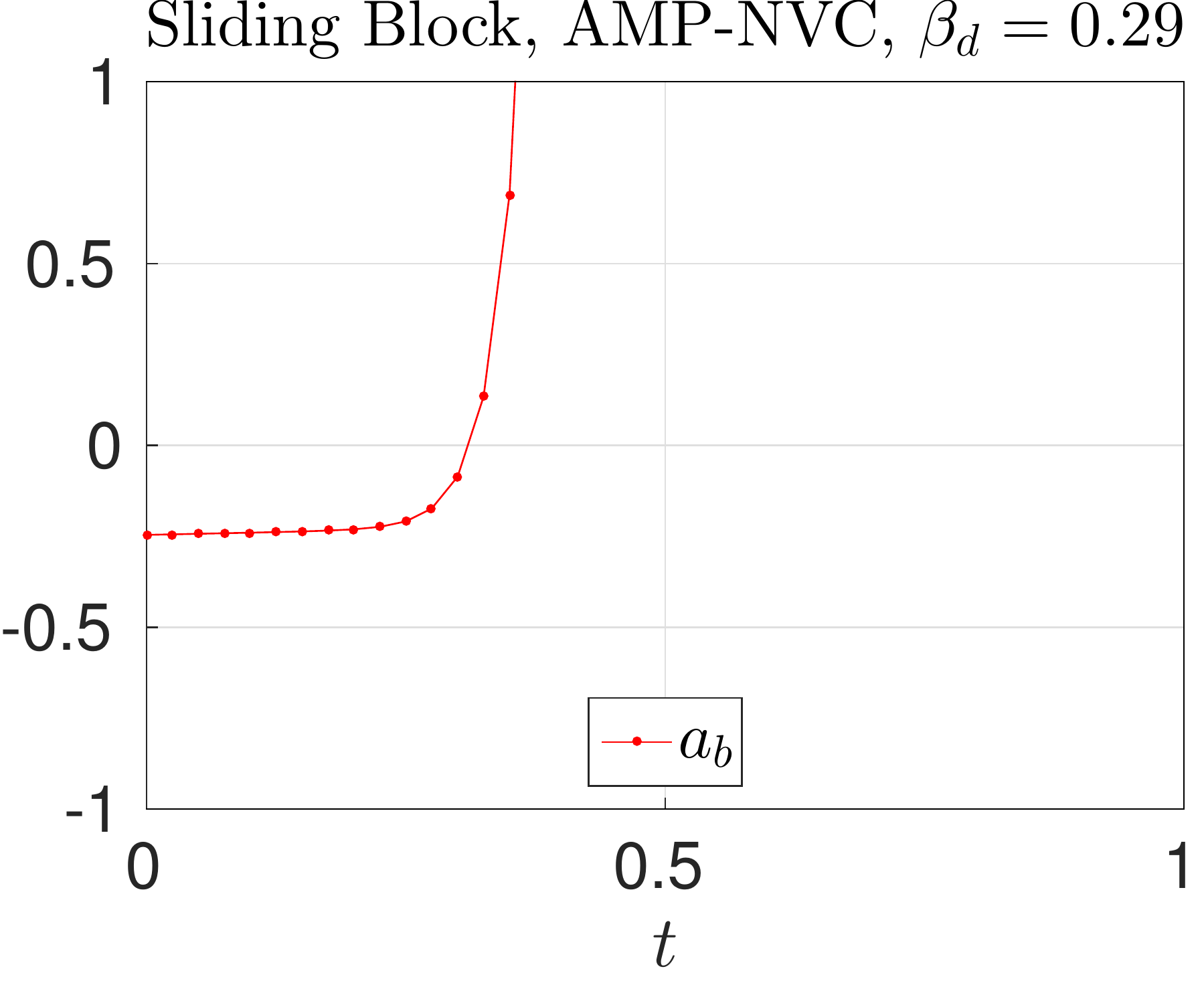}{\figWidth}};
    \draw(1.5,2.5) node[draw,fill=red!50] {\labelFont I};
    \draw(4.0,0.0) node[anchor=south west,xshift=-4pt,yshift=+0pt] {\trimfig{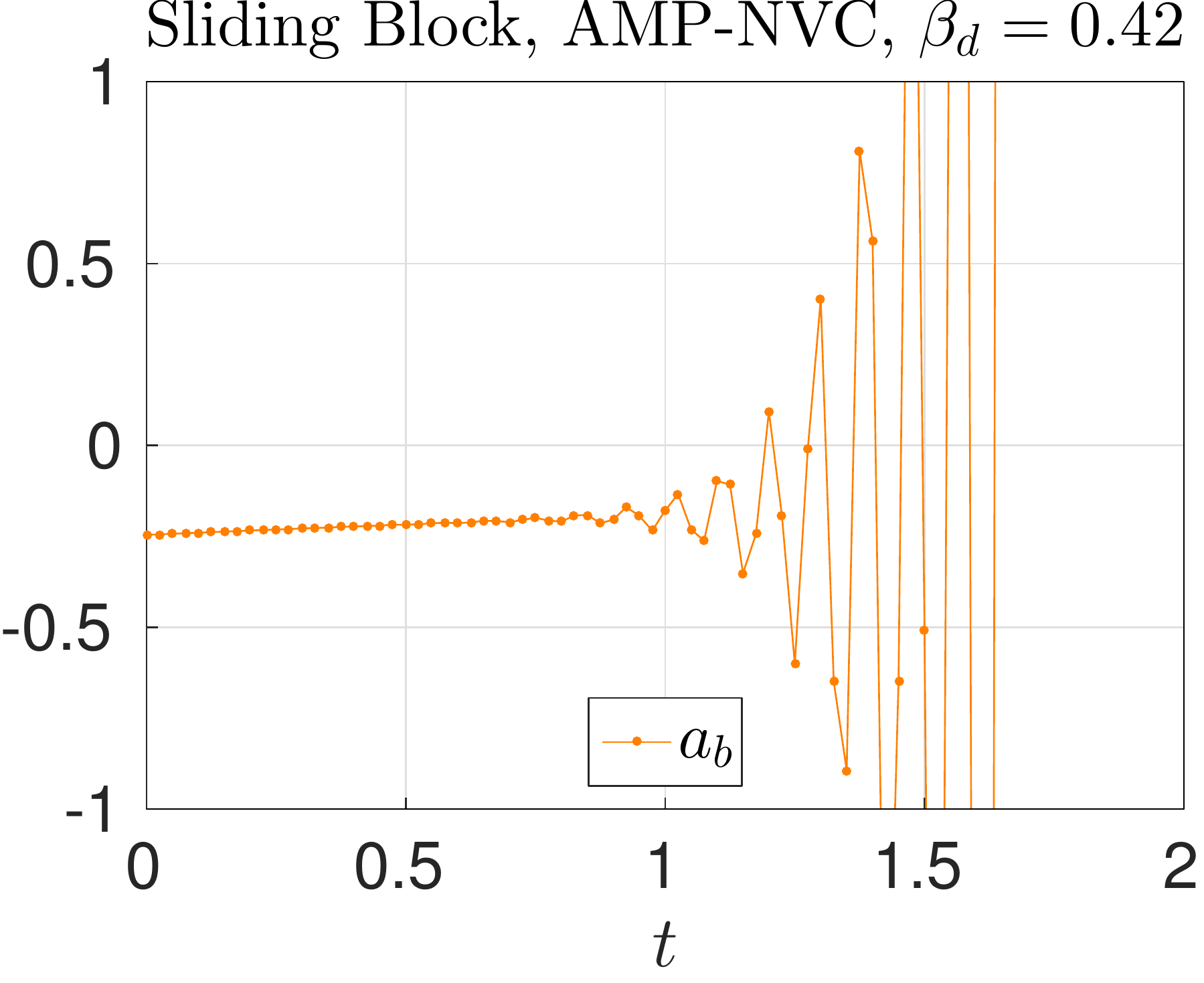}{\figWidth}};
    \draw(5.5,2.5) node[draw,fill=orange!50] {\labelFont II};
	\draw(8.0,0.0) node[anchor=south west,xshift=-4pt,yshift=+0pt] {\trimfig{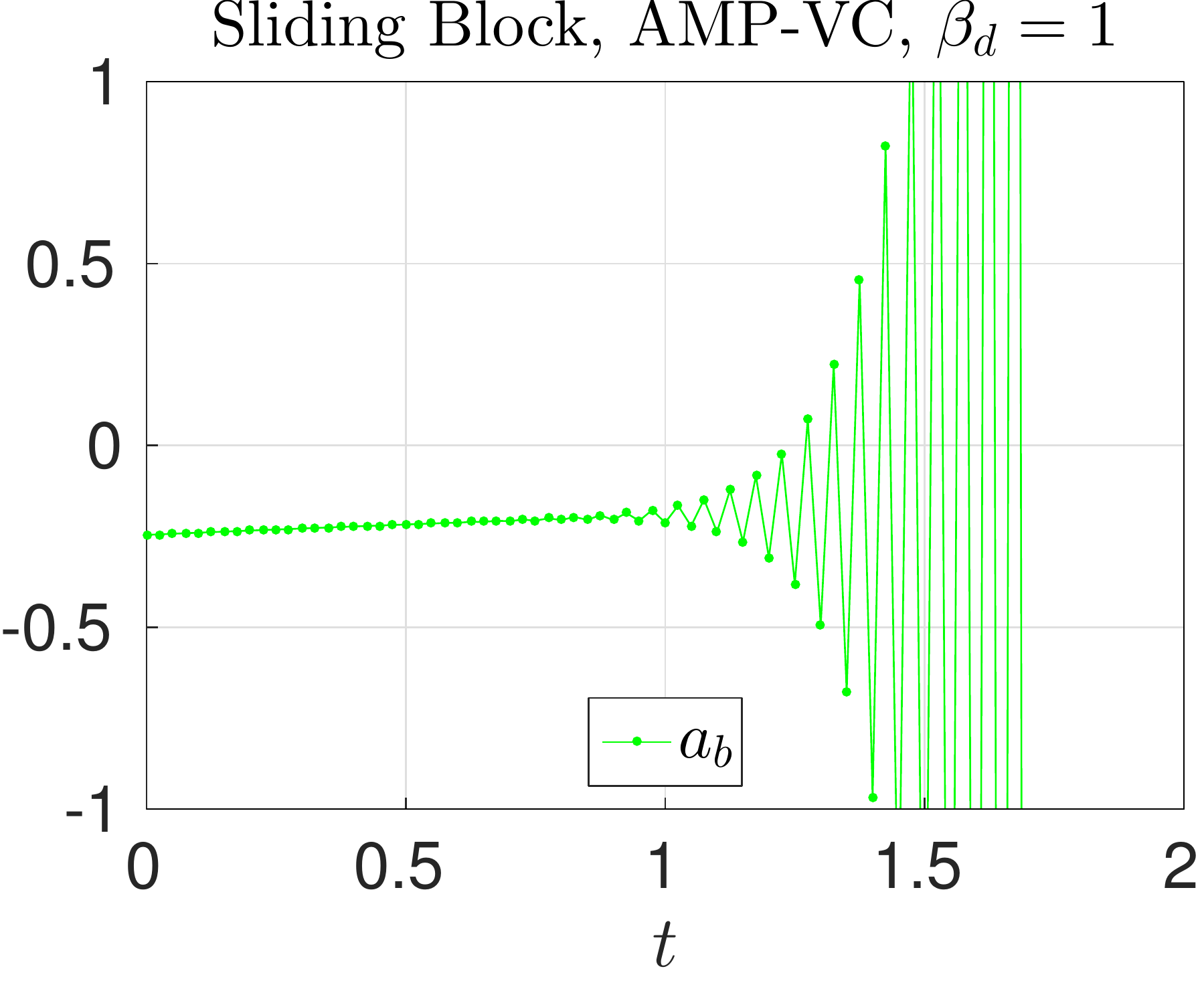}{\figWidth}};
    \draw(9.5,2.5) node[draw,fill=green!50] {\labelFont III};
    \draw(12.0,0.0) node[anchor=south west,xshift=-4pt,yshift=+0pt] {\trimfig{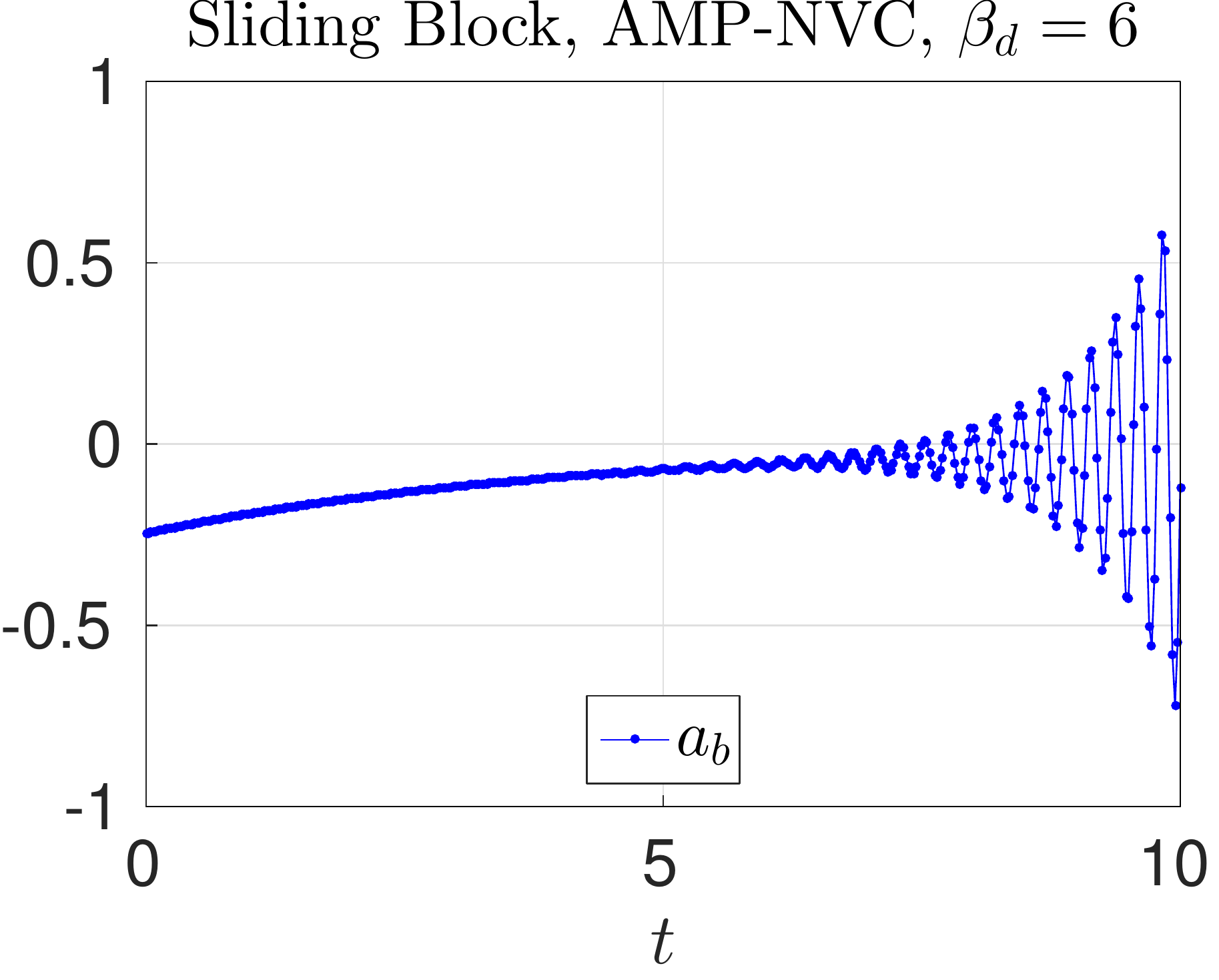}{\figWidth}};
    \draw(13.5,2.5) node[draw,fill=blue!50] {\labelFont IV};
 \end{scope}
 \begin{scope}[xshift=1.5cm,yshift=0.25cm]
    \draw(0.0,0.0) node[anchor=south west,xshift=-4pt,yshift=+0pt] {\trimfig{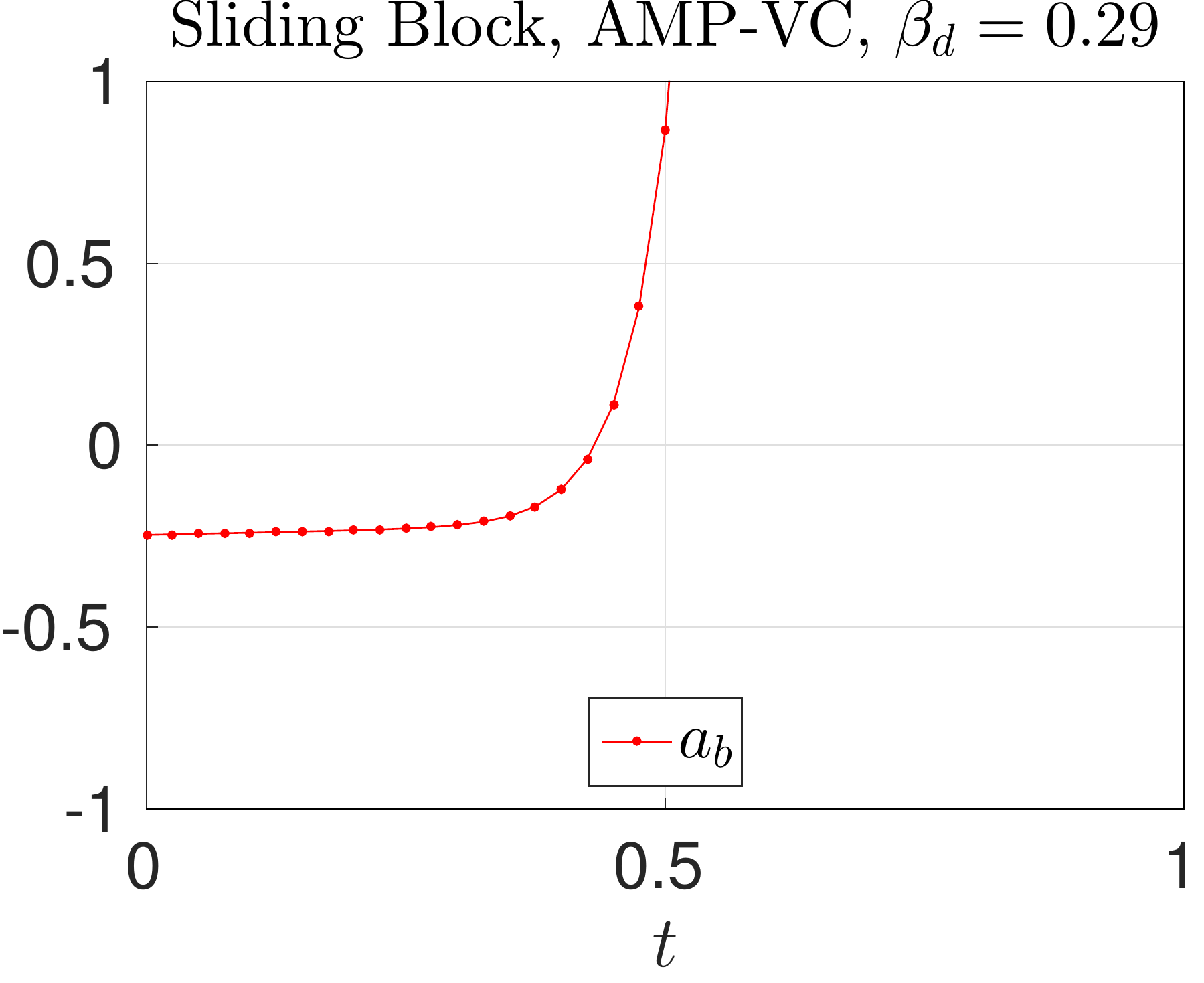}{\figWidth}};
    \draw(1.5,2.5) node[draw,fill=red!50] {\labelFont I};
    \draw(4.0,0.0) node[anchor=south west,xshift=-4pt,yshift=+0pt] {\trimfig{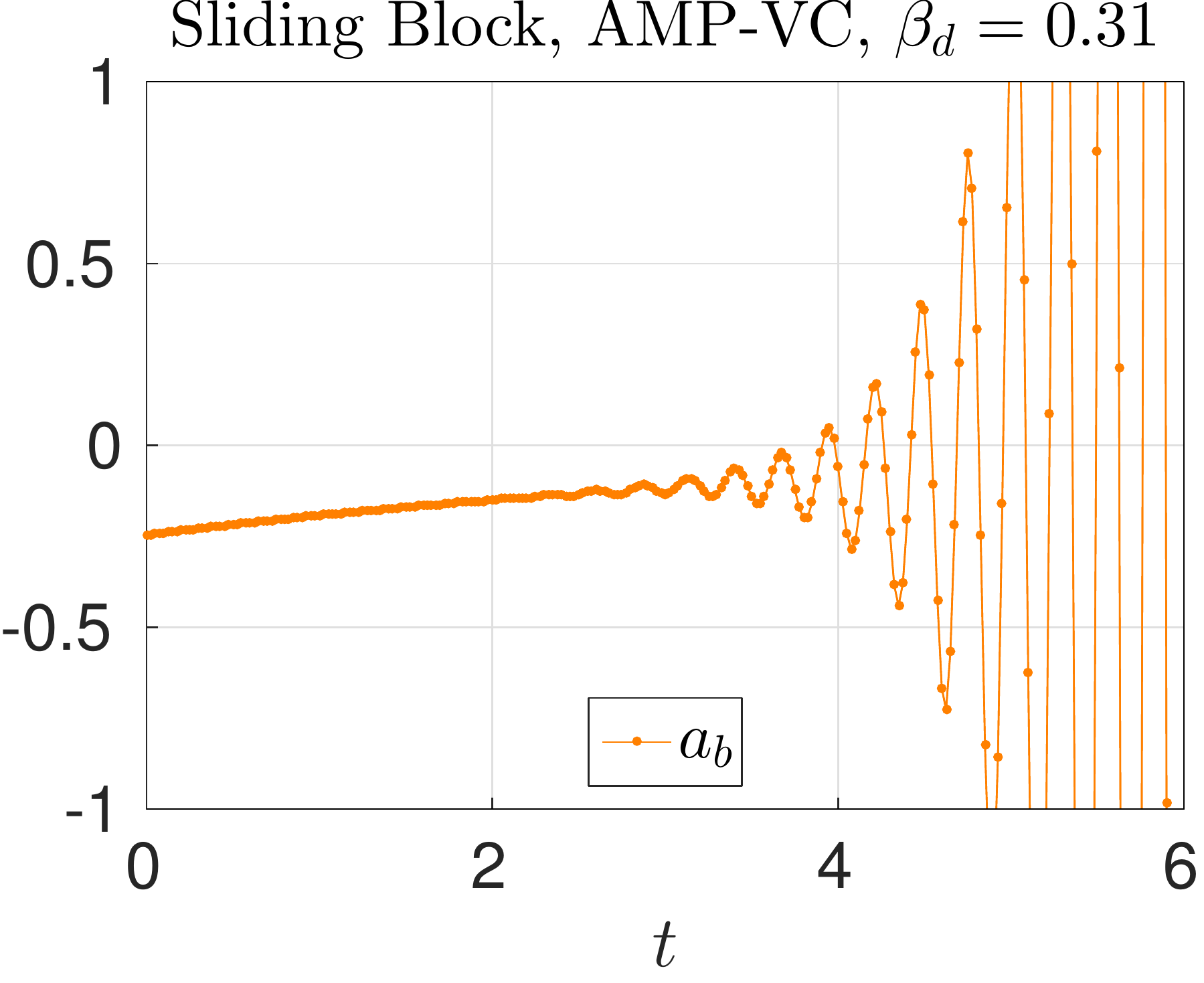}{\figWidth}};
    \draw(5.5,2.5) node[draw,fill=orange!50] {\labelFont II};
    \draw(8.0,0.0) node[anchor=south west,xshift=-4pt,yshift=+0pt] {\trimfig{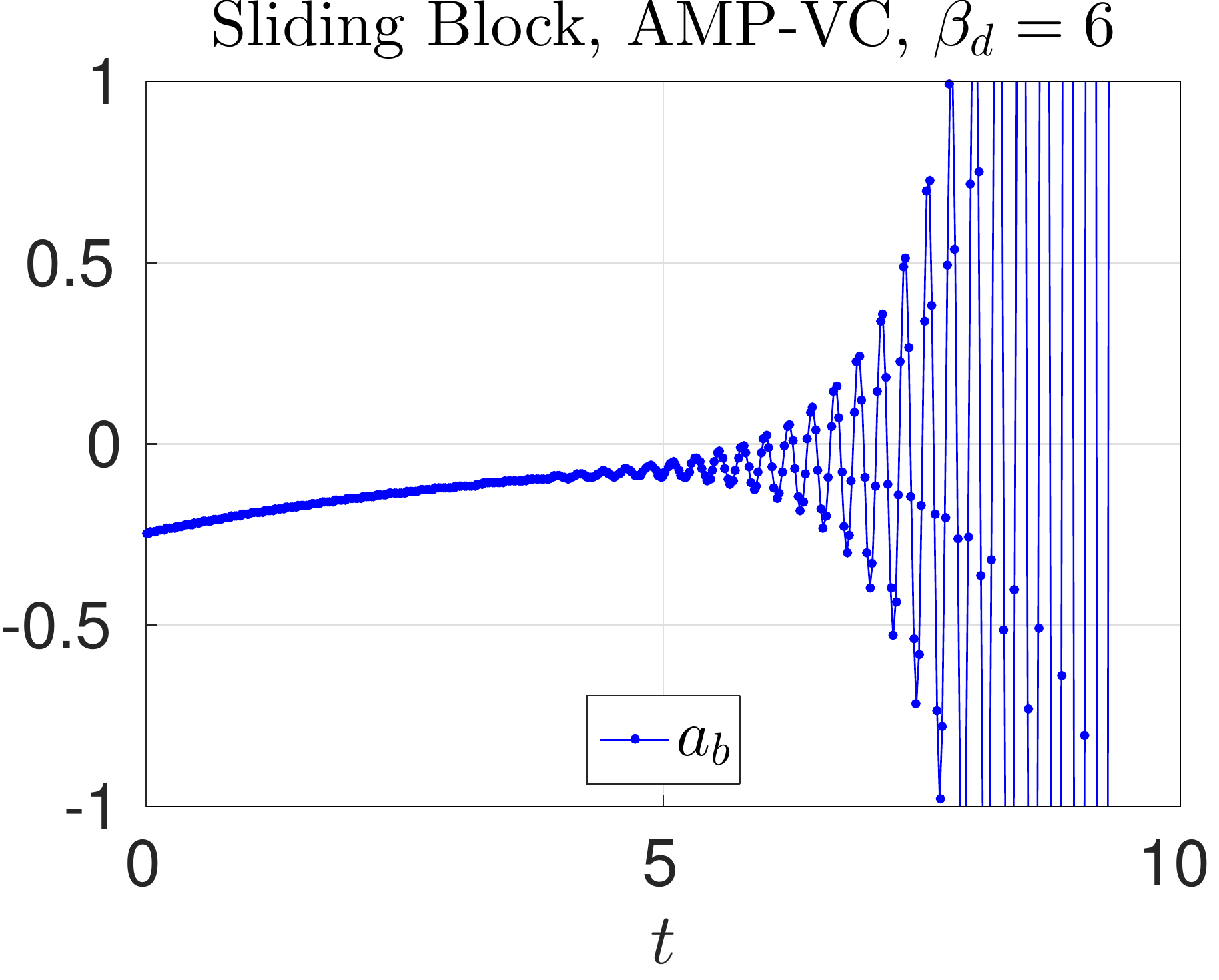}{\figWidth}};
    \draw(9.5,2.5) node[draw,fill=blue!50] {\labelFont IV};
  \end{scope}
 %
\end{tikzpicture}
\end{center}
  \caption{
	  Sliding block. Acceleration of a light rigid body with $\mb=0.001$.
   Top: typical instabilities in the {\ampRB}-NVC~scheme:
    (I)~for $\adc = 0.29$   there is a rapid, monotonically growing instability;
	(II)~for $\adc = 0.42$  there is a rapidly growing, high-frequency instability;
    (III)~for $\adc = 1$ 	there is a rapidly growing, step-to-step oscillatory instability; and
    (IV)~for $\adc = 6$ 	there is a slowly growing, low-frequency instability. 
   Bottom: typical instabilities in the {\ampRB}-VC~scheme:
   (I)~for $\adc = 0.29$ 	there is a rapid, monotonically growing instability;
   (II)~for $\adc = 0.31$ 	there is a rapidly growing, high-frequency instability; and
   (IV)~for $\adc = 6$ 	there is a slowly growing, low-frequency instability. 
    }
  \label{fig:shearBlockCurves}
\end{figure}
}

Next we probe the classifications of unstable modes as discussed in
Section~\ref{sec:adSubmodel}. Figure~\ref{fig:shearBlockCurves} shows computed results for
the~\ampRB~scheme where the added-damping parameter $\adc$ is intentionally chosen to elicit the
unstable behaviour of a type typified by one of the four instability regimes.
The plots in the top row of the figure show results for the~\ampRB-NVC scheme, while
the plots in the bottom row show results for the~\ampRB-VC scheme.  By comparing the behaviour of
the computed instabilities with those predicted by the stability analysis (see Figures~\ref{fig:addedDampingNvcModes} and~\ref{fig:addedDampingVcModes}), we find the numerical results are in excellent agreement
with the results of the theory for all modes of instability.

\newcommand{\amplitude}{{\alpha_b}}

\subsection{Translating disk in an incompressible fluid} \label{sec:translatingDisk}

To discuss added-mass effects and the {\ampRB} scheme in a two-dimensional setting, consider the
translation of a rigid disk of radius $\aa=1$ surrounded by an annular fluid region bounded by a radius
$\bb=2$ (see Figure~\ref{fig:modelDiskInDisk}).  The interface matching conditions are applied on the surface of the disk at $r=r_1$ and a no-slip boundary
condition is taken on the outer boundary $r=\bb$.  An exact solution to this problem can be derived assuming an
infinitesimal displacement of the body.  The analysis, as described in~\ref{sec:exactTranslatingDisk}, seeks solutions
with the horizontal velocity of the rigid body taken to be
\[
u_b(t) = \amplitude e^{-\lambda^2 \nu t},
\]
where $\amplitude$ is the (small) amplitude of the motion, $\lambda$ is an eigenvalue and $\nu=\mu/\rho$ is the kinematic viscosity.  The corresponding fluid velocity and pressure in this solution have the form
\begin{align*}
   \vv(r,\theta,t) &=   (\hat Q(r)/r)\begin{bmatrix} \cos\sp2\theta \\ \sin\theta \cos\theta\end{bmatrix}e\sp{-\lambda\sp2\nu t} 
                      + \hat Q\sp\prime(r)\begin{bmatrix} \sin\sp2\theta \\ -\sin\theta\cos\theta \end{bmatrix}e\sp{-\lambda\sp2\nu t} , \\
   p(r,\theta,t) &= \mu\left(A_pr+B_p/r\right)e\sp{-\lambda\sp2\nu t}\cos\theta+ p_0 ,
\end{align*}
where $(r,\theta)$ are polar coordinates, $\hat Q(r)$ is given in~\eqref{eq:hatQ}, $A_p$ and $B_b$ are constant specified in the exact solution, and $p_0$ is an arbitrary constant.  The eigenvalues $\lambda$ are determined from~\eqref{eq:transDiskEigenvalues} and the smallest (positive) values are given in Table~\ref{table:cylEig} for various choices of the density $\rho_b$ of the solid disk (see~\ref{sec:exactTranslatingDisk}).

In the numerical tests to follow, the amplitude of the horizontal motion is taken as $\amplitude = 10^{-7}$ (to approximately satisfy the assumption of infinitesimal displacement), and the arbitrary constant $p_0$ is set to zero.
The inner and outer radii of the annular fluid domain are taken to be $\aa=1$ and $\bb=2$, respectively, and the
density and viscosity of the fluid are $\rho=1$ and $\mu=0.1$.
The values of $\lambda$ used for the exact solution, and the corresponding initial conditions for the simulations, are
taken from Table~\ref{table:cylEig}.
The annular grids used in the computations are denoted by $\Gs^{(j)}$ and have grid spacing $\dr^{(j)} = 1/(10 j)$ and $\Delta \theta^{(j)} = 2 \pi/ (60 j)$.  The time step is taken as $\dt^{(j)} = 1/(10 j)$.

{
\def\rad{1.25}
\newcommand{\plotDisk}{
\fill[fill=red!20,draw=red,line width=2pt] 
      plot[samples=100, domain=0.:360] ( {\rad*cos(\x)} , {\rad*sin(\x)} ) -- cycle ;
}
\newcommand{\figWidth}{5.5cm}
\newcommand{\trimfig}[2]{\trimFig{#1}{#2}{.0}{.0}{.0}{.0}}
\newcommand{\figWidtha}{5.5cm}
\newcommand{\trimfiga}[2]{\trimFig{#1}{#2}{.045}{.31}{.18}{.2}}
\newcommand{\trimfigb}[2]{\trimFig{#1}{#2}{0.47}{0.4}{0.12}{0.12}}
\begin{figure}[htb]
\begin{center}
\resizebox{14cm}{!}{
\begin{tikzpicture}[scale=1]
  \useasboundingbox (0.0,.25) rectangle (16,6.);  
  \draw(-0.2,0) node[anchor=south west,xshift=-4pt,yshift=+10pt]  {\trimfig{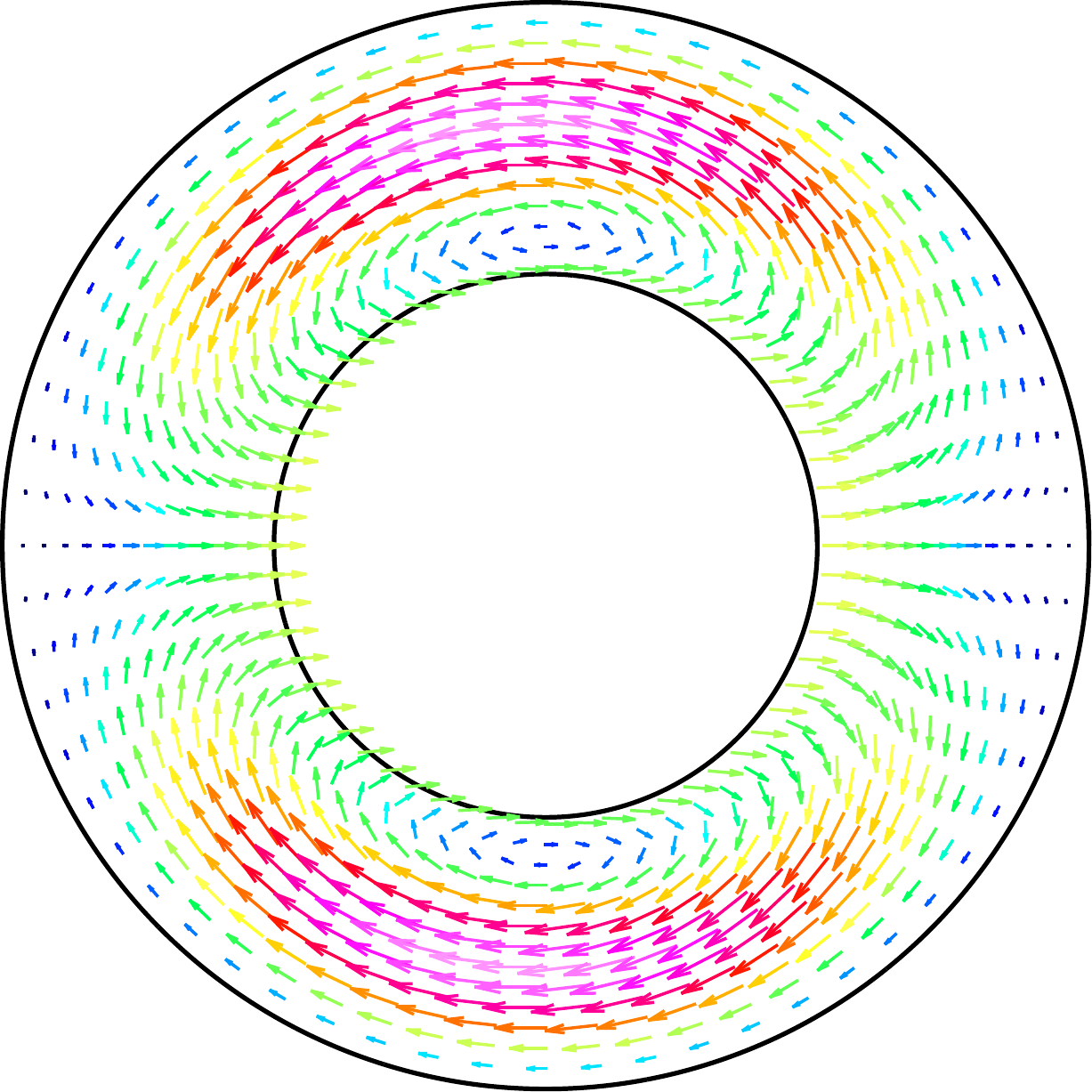}{\figWidth}};
  \draw(5.3, 0) node[anchor=south west,xshift=-4pt,yshift=-20pt] {\trimfigb{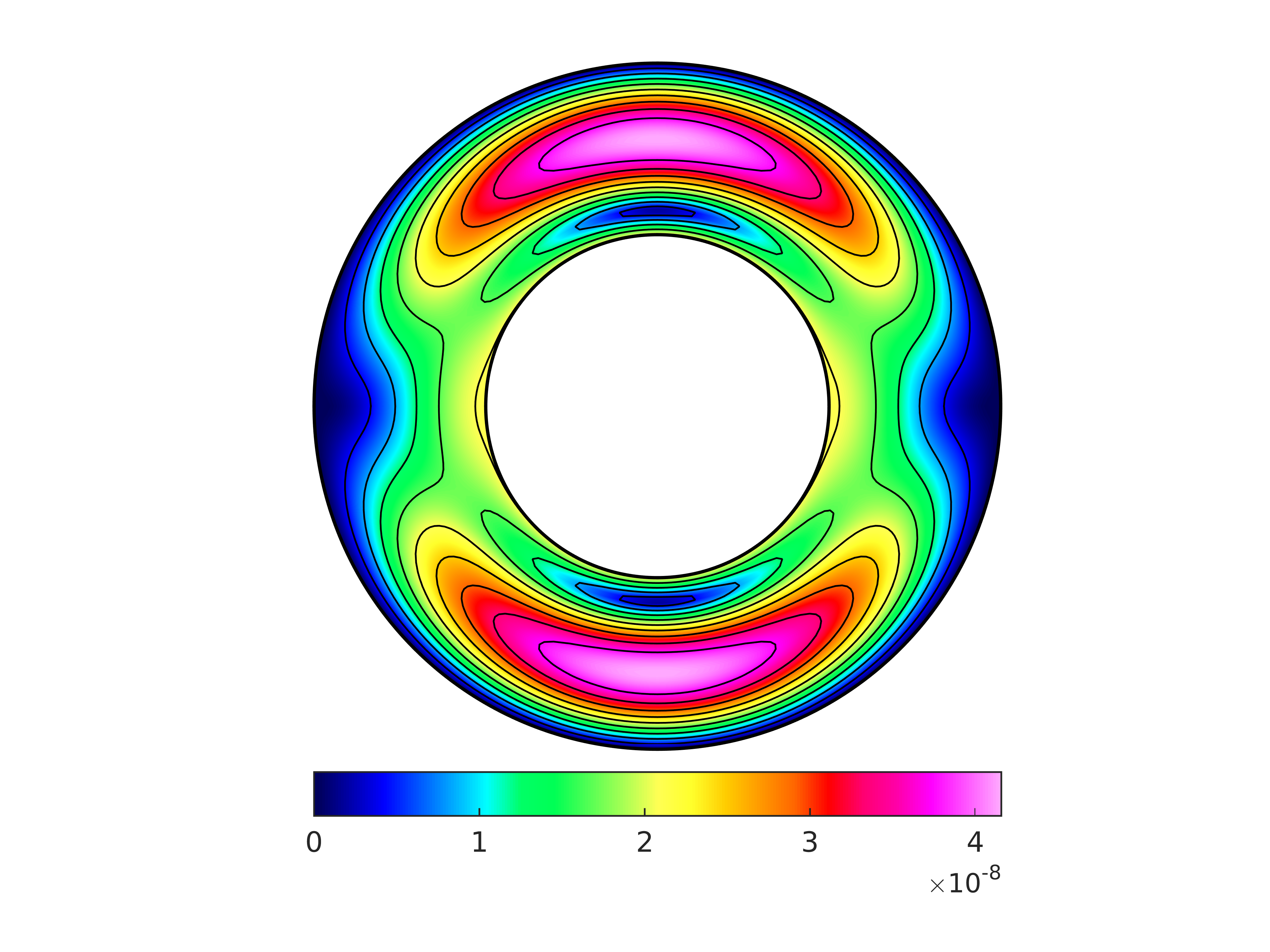}{\figWidtha}};
  \draw(10.8,0) node[anchor=south west,xshift=-4pt,yshift=-20pt] {\trimfigb{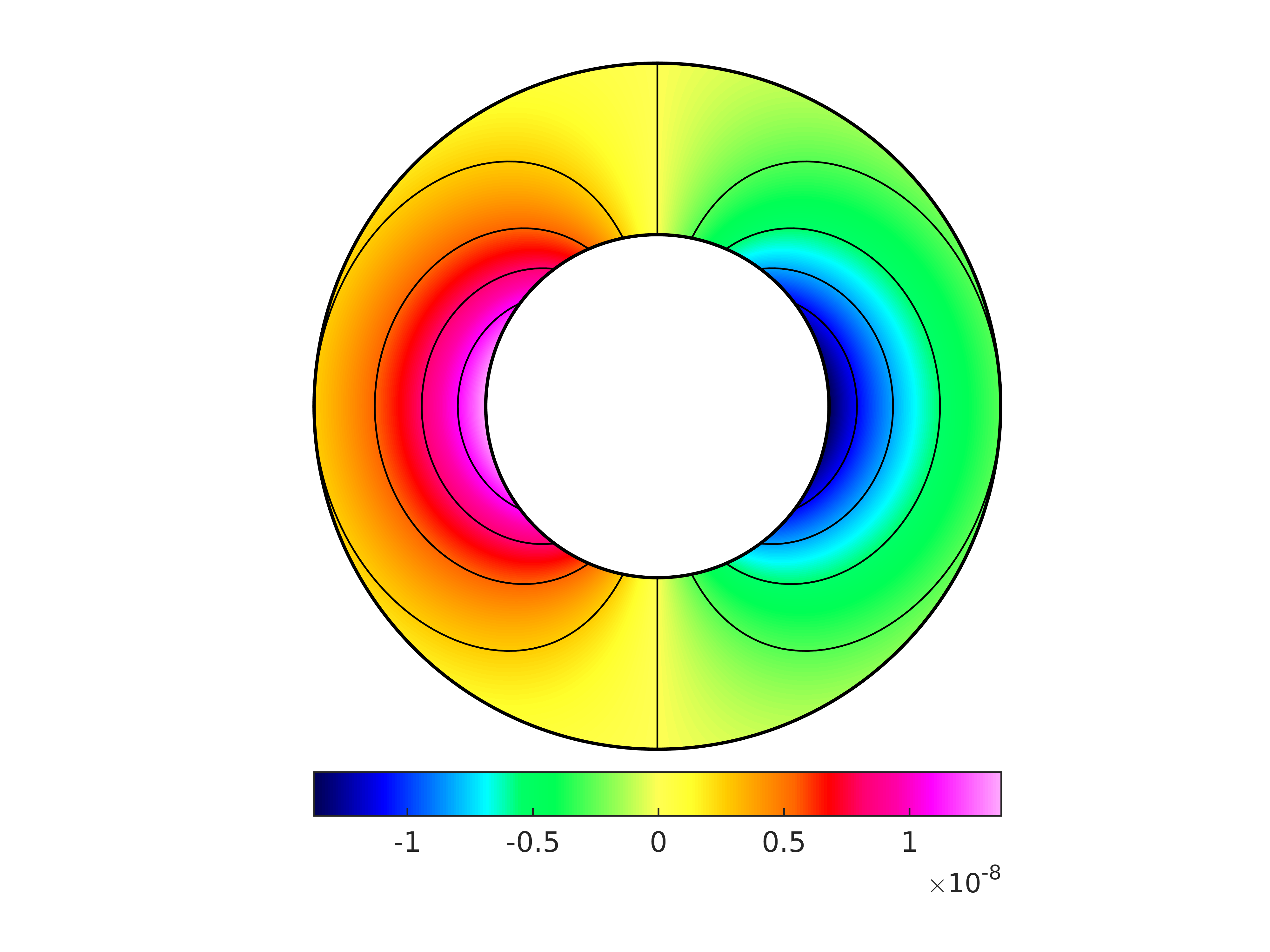}{\figWidtha}};
  \begin{scope}[xshift=2.58cm,yshift=3.28cm]
    \plotDisk
    \draw[very thick,->,red] (-.3,0) -- (.3,0); 
  \end{scope}
\end{tikzpicture}
} 
\end{center}
  \caption{Translating disk.  Computed velocity vectors (left), speed (middle) and pressure (right) using the grid $\Gc^{(4)}$
  for $\rhos=0$. The color bar below the middle plot also applied to the left plot.
   }
  \label{fig:translatingDiskInDisk}
\end{figure}
}

{
\newcommand{\figWidth}{7.5cm}
\newcommand{\figWidths}{7.25cm}
\newcommand{\trimfig}[2]{\trimFig{#1}{#2}{.0}{.0}{.0}{.0}}
\begin{figure}[htb]
\begin{center}
\resizebox{14cm}{!}{
\begin{tikzpicture}[scale=1]
  \useasboundingbox (0, 0) rectangle (16.,6);  
  \draw(0,-.75) node[anchor=south west,xshift=0pt,yshift=+0pt] {\trimfig{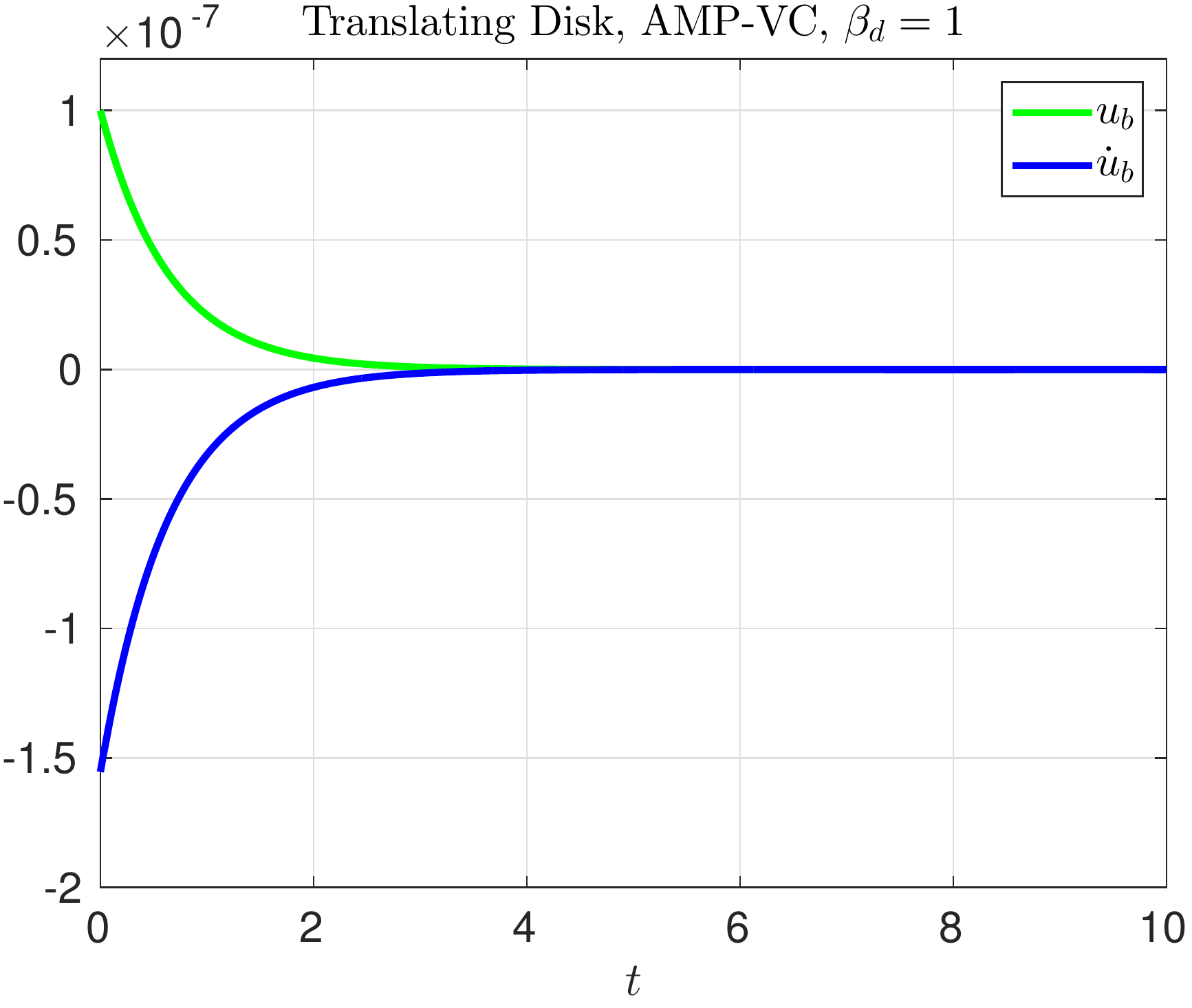}{\figWidth}};
  \draw(8,-.75) node[anchor=south west,xshift=0pt,yshift=+0pt] {\trimfig{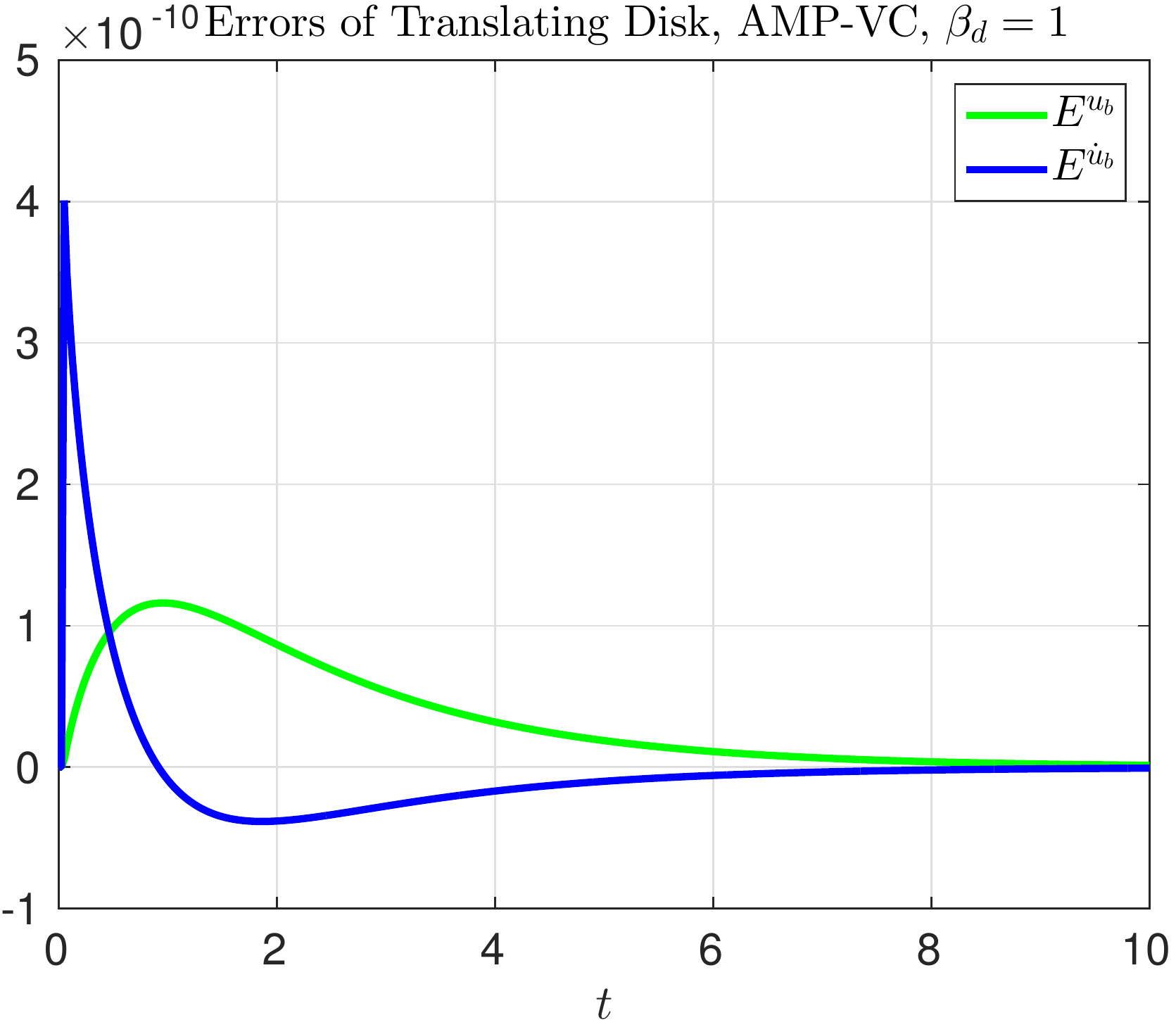}{\figWidths}};
\end{tikzpicture}
} 
\end{center}
  \caption{
  Translating disk. Velocity and acceleration of the body (left) and errors (right) for a massless, $\rhob=0$, rigid disk using the grid $\Gc^{(4)}$. 
The~\ampRB-VC scheme is used with $\adc=1$.
  }
  \label{fig:translatingDiskCurves}
\end{figure}
}

{ 
\newcommand{\errp}{$E_j^{p}$}
\newcommand{\erruv}{$E_j^{\vv}$}
\newcommand{\errvb}{$E_j^{\vv_b}$}
\newcommand{\errvbt}{$E_j^{\dot\vv_b}$}

\newcommand{\translatingDiskTableI}{%
\begin{tabular}{|c|c|c|c|c|c|c|c|c|} \hline 
  \multicolumn{9}{|c|}{Translating disk, $\rhob=10$} \\ \hline 
\strutt  $\dr^{(j)}$  &     \errp     &  r   &     \erruv     &  r   &     \errvb     &  r   &     \errvbt     &  r    \\[3pt] \hline 
1/10 & \num{4.1}{-9} &      & \num{2.0}{-9} &      & \num{6.1}{-10} &      & \num{6.0}{-10} &       \\ \hline 
1/20 & \num{1.2}{-9} &  3.5 & \num{5.8}{-10} &  3.5 & \num{1.7}{-10} &  3.5 & \num{1.7}{-10} &  3.5  \\ \hline 
1/40 & \num{3.0}{-10} &  3.9 & \num{1.6}{-10} &  3.7 & \num{4.9}{-11} &  3.6 & \num{4.5}{-11} &  3.8  \\ \hline 
1/80 & \num{7.6}{-11} &  4.0 & \num{4.1}{-11} &  3.8 & \num{1.3}{-11} &  3.8 & \num{1.2}{-11} &  3.9  \\ \hline 
1/160& \num{1.9}{-11} &  4.0 & \num{1.0}{-11} &  3.9 & \num{3.3}{-12} &  3.9 & \num{2.9}{-12} &  3.9  \\ \hline 
rate &      1.94     &      &      1.90     &      &      1.88     &      &      1.92     &      \\ \hline
\end{tabular}
}
\newcommand{\translatingDiskTableII}{%
\begin{tabular}{|c|c|c|c|c|c|c|c|c|} \hline 
  \multicolumn{9}{|c|}{Translating disk, $\rhob=1$} \\ \hline 
\strutt  $\dr^{(j)}$  &     \errp     &  r   &     \erruv     &  r   &     \errvb     &  r   &     \errvbt     &  r    \\[3pt] \hline 
1/10 & \num{2.1}{-10} &      & \num{9.7}{-10} &      & \num{9.7}{-10} &      & \num{2.7}{-10} &       \\ \hline 
1/20 & \num{6.3}{-11} &  3.3 & \num{2.6}{-10} &  3.8 & \num{2.6}{-10} &  3.8 & \num{7.4}{-11} &  3.7  \\ \hline 
1/40 & \num{1.8}{-11} &  3.5 & \num{6.7}{-11} &  3.9 & \num{6.7}{-11} &  3.9 & \num{1.6}{-11} &  4.6  \\ \hline 
1/80 & \num{4.8}{-12} &  3.7 & \num{1.7}{-11} &  3.9 & \num{1.7}{-11} &  3.9 & \num{3.5}{-12} &  4.6  \\ \hline 
1/160& \num{1.2}{-12} &  3.9 & \num{4.3}{-12} &  4.0 & \num{4.3}{-12} &  4.0 & \num{8.6}{-13} &  4.1  \\ \hline 
rate &      1.86     &      &      1.96     &      &      1.96     &      &      2.10     &      \\ \hline 
\end{tabular}
}
\newcommand{\translatingDiskTableIV}{%
\begin{tabular}{|c|c|c|c|c|c|c|c|c|} \hline 
  \multicolumn{9}{|c|}{Translating disk, $\rhob=0$} \\ \hline 
\strutt  $\dr^{(j)}$  &     \errp     &  r   &     \erruv     &  r   &     \errvb     &  r   &     \errvbt     &  r    \\[3pt] \hline 
1/10 & \num{2.9}{-10} &      & \num{1.7}{-9} &      & \num{1.7}{-9} &      & \num{8.0}{-11} &       \\ \hline 
1/20 & \num{8.8}{-11} &  3.3 & \num{4.5}{-10} &  3.7 & \num{4.5}{-10} &  3.7 & \num{1.7}{-11} &  4.7  \\ \hline 
1/40 & \num{2.3}{-11} &  3.8 & \num{1.2}{-10} &  3.9 & \num{1.2}{-10} &  3.9 & \num{9.1}{-12} &  1.9  \\ \hline 
1/80 & \num{5.8}{-12} &  3.9 & \num{2.9}{-11} &  4.0 & \num{2.9}{-11} &  4.0 & \num{3.2}{-12} &  2.9  \\ \hline  
1/160& \num{1.5}{-12} &  3.9 & \num{7.4}{-12} &  4.0 & \num{7.4}{-12} &  4.0 & \num{8.3}{-13} &  3.8  \\ \hline 
rate &      1.92     &      &      1.96     &      &      1.96     &      &      1.56     &      \\ \hline
\end{tabular}
}
{
\begin{table}[hbt]\tableFont 
\begin{center}
  \translatingDiskTableI \\
\medskip
  \translatingDiskTableII \\
\medskip
  \translatingDiskTableIV
\caption{Translating disk. Maximum errors and estimated convergence rates at $t=1$ computed using the~\ampRB-VC~scheme
   for a heavy, $\rhob=10$, medium, $\rhob=1$, and massless, $\rhob=0$, rigid disk. }
\label{table:translatingDisk}
\end{center}
\end{table}
}
} 

Figures~\ref{fig:translatingDiskInDisk} and~\ref{fig:translatingDiskCurves} illustrate the behaviour of
the solution for a
massless rigid body ($\rhos=0$) as determined by the {\ampRB-VC}~scheme using $\Gs^{(4)}$. Note that although the exact solution to this test case has no rotation, the numerical
implementation allows the full set of motions including both $x$ and $y$ translations as well as
rotation. As such it is important to include stabilization against added-damping effects and this is done
by including the added-damping terms with $\adc=1$ . 
Results from a convergence study
are provided in Table~\ref{table:translatingDisk}, which demonstrate second-order accuracy for all
choices of body densities, $\rho_b$.

\subsection{Rotating disk in a disk of incompressible fluid} \label{sec:rotatingDisk}

In this final example we consider a rotating solid disk surrounded by an annular fluid region.
Rotational added-damping effects play an important role in this problem.
The geometry of the problem and the grid for the fluid domain are identical to those given in Section~\ref{sec:translatingDisk}. 
The interface matching conditions are applied on the surface of the disk, $\aa=1$, and no-slip boundary conditions
are imposed on the fixed outer boundary, $\bb=2$.  The fluid density and viscosity are taken to be $\rho=1$ and $\mu=0.1$.

An exact solution can be found for a disk rotating about its centre assuming there is no displacement of the disk in the $x$ or $y$ directions.  If there is no external torque on the disk, then a straightforward application of separation of variables leads to a solution of the form
\[
\vv(r,\theta,t) = v_\theta(r,t) \begin{bmatrix} -\sin\theta \\ \cos\theta \end{bmatrix},\qquad
p(r,t) = \rho\int_{\aa}^{r}  \frac{[v_\theta(s,t)]^2}{s} \, ds + p_0,
\]
where $p_0$ is a constant and
\[
v_\theta(r,t) =  \amplitude \left[
           \frac{J_1(\lambda r) Y_1(\lambda \bb) - J_1(\lambda \bb) Y_1(\lambda r)}
                {J_1(\lambda \aa) Y_1(\lambda \bb) - J_1(\lambda \bb) Y_1(\lambda \aa)}\right]
           \,e^{-\lambda^2 \nu t},\]
is the circumferential component of the velocity.  Here, $J_1$ and $Y_1$ are Bessel functions of order one and $\amplitude$ is the velocity at $r=\aa$ and $t=0$.  The angular velocity of the rigid body in the exact solution is given by
\[
\omega_b(t)={v_\theta(\aa,t)\over \aa}={\amplitude e^{-\lambda^2 \nu t}\over\aa}.
\]
We note that the solution decays at a rate determined by the eigenvalues~$\lambda$ for which there are an infinite number.  These eigenvalues satisfy a constraint determined from the equation for the angular acceleration of the body.  This constraint has the form
\begin{align*}
  \left[ \frac{\partial v}{\partial r}  - \left[ \frac{1}{\aa} - \frac{\lambda^2 I_b}{2\pi\rho \aa^3}\right] \, v \right]_{r=\aa} = 0 ,
\end{align*}
where $I_b$ is the moment of inertia of the solid disk.  Assuming a disk with uniform density, $\rho_b$, we have $I_b = \rhos \pi \aa^4 /2$, and the smallest (positive) values of $\lambda$ are given in Table~\ref{table:cylRotEig} for various values of $\rhob$.  We note that unlike the solution described in~\ref{sec:exactTranslatingDisk} for the infinitesimal motion of a translating disk, the solution given here for the rotating disk is an exact solution of the full Navier-Stokes equations for any value of the amplitude~$\amplitude$.

\begin{table}[H]\tableFont
\begin{center}
\begin{tabular}{|c|c|c|c|c|}\hline
             & $\rhob=10$ & $\rhob=1$ & $\rhob=10^{-3}$ & $\rhob=0$\\ \hline 
$\lambda$    & 0.969674911038943 & 1.97045369767466 & 2.27309495331973 & 2.27340628201175  \\ \hline 
\end{tabular}
\caption{Eigenvalues, $\lambda$, in the exact solution of the rotating disk for different values of $\rhob$.}
\label{table:cylRotEig}
\end{center}
\end{table}

{
\def\rad{1.22}
\newcommand{\plotDisk}{
\fill[fill=red!20,draw=red,line width=2pt] 
      plot[samples=100, domain=0.:360] ( {\rad*cos(\x)} , {\rad*sin(\x)} ) -- cycle ;
}
\newcommand{\figWidth}{5.6cm}
\newcommand{\trimfig}[2]{\trimFig{#1}{#2}{.42}{.31}{.12}{.12}}
\newcommand{\figWidtha}{5.5cm}
\newcommand{\trimfiga}[2]{\trimFig{#1}{#2}{.045}{.31}{.18}{.2}}
\newcommand{\trimfigb}[2]{\trimFig{#1}{#2}{0.47}{0.4}{0.12}{0.12}}
\begin{figure}[htb]
\begin{center}
\resizebox{11cm}{!}{
\begin{tikzpicture}[scale=1]
  \useasboundingbox (5.0,.25) rectangle (16.9,6);  
  \draw(5.3,0) node[anchor=south west,xshift=-4pt,yshift=-20pt] {\trimfigb{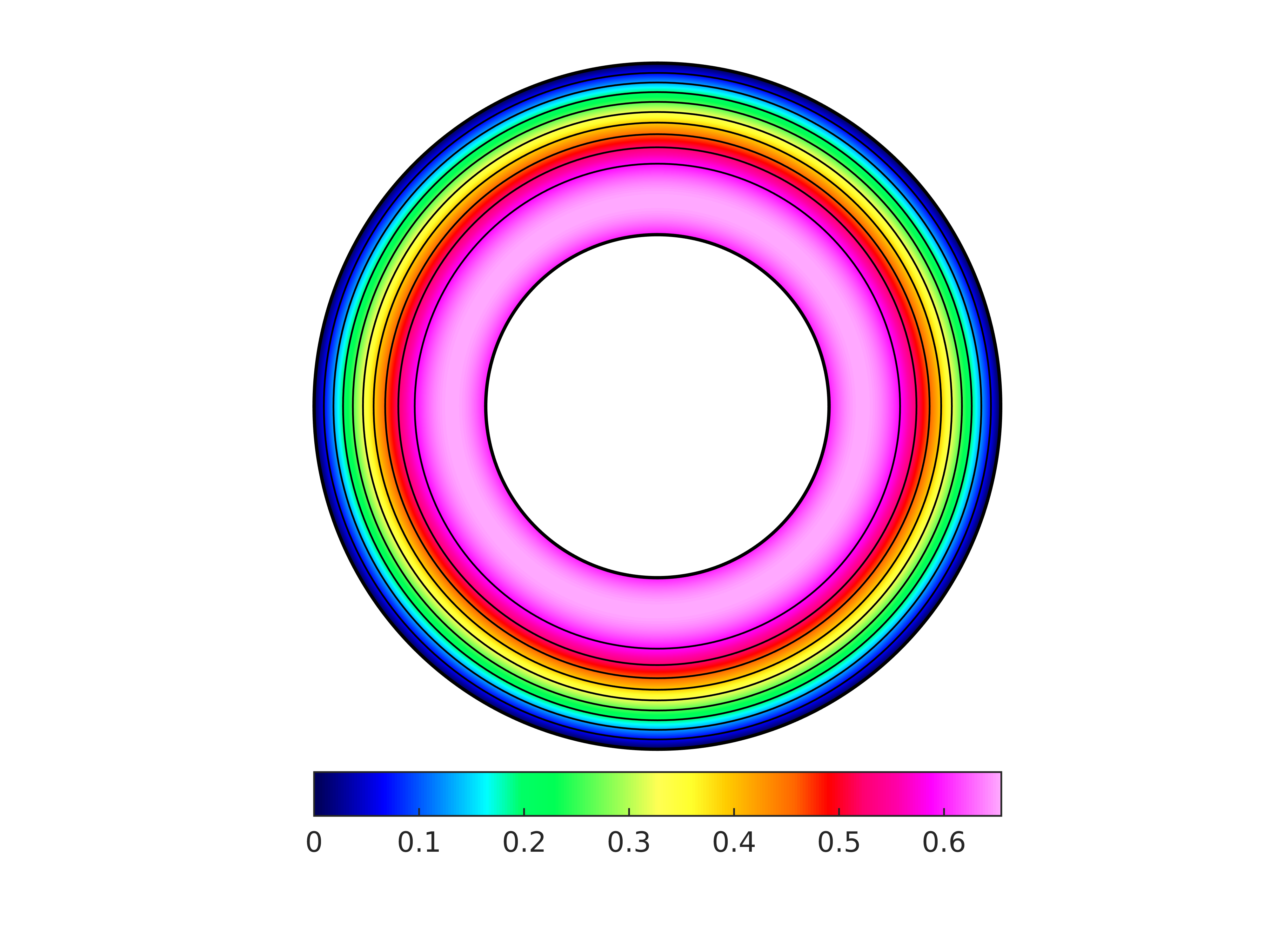}{\figWidtha}};
  \draw(11.3,0) node[anchor=south west,xshift=-4pt,yshift=-20pt] {\trimfigb{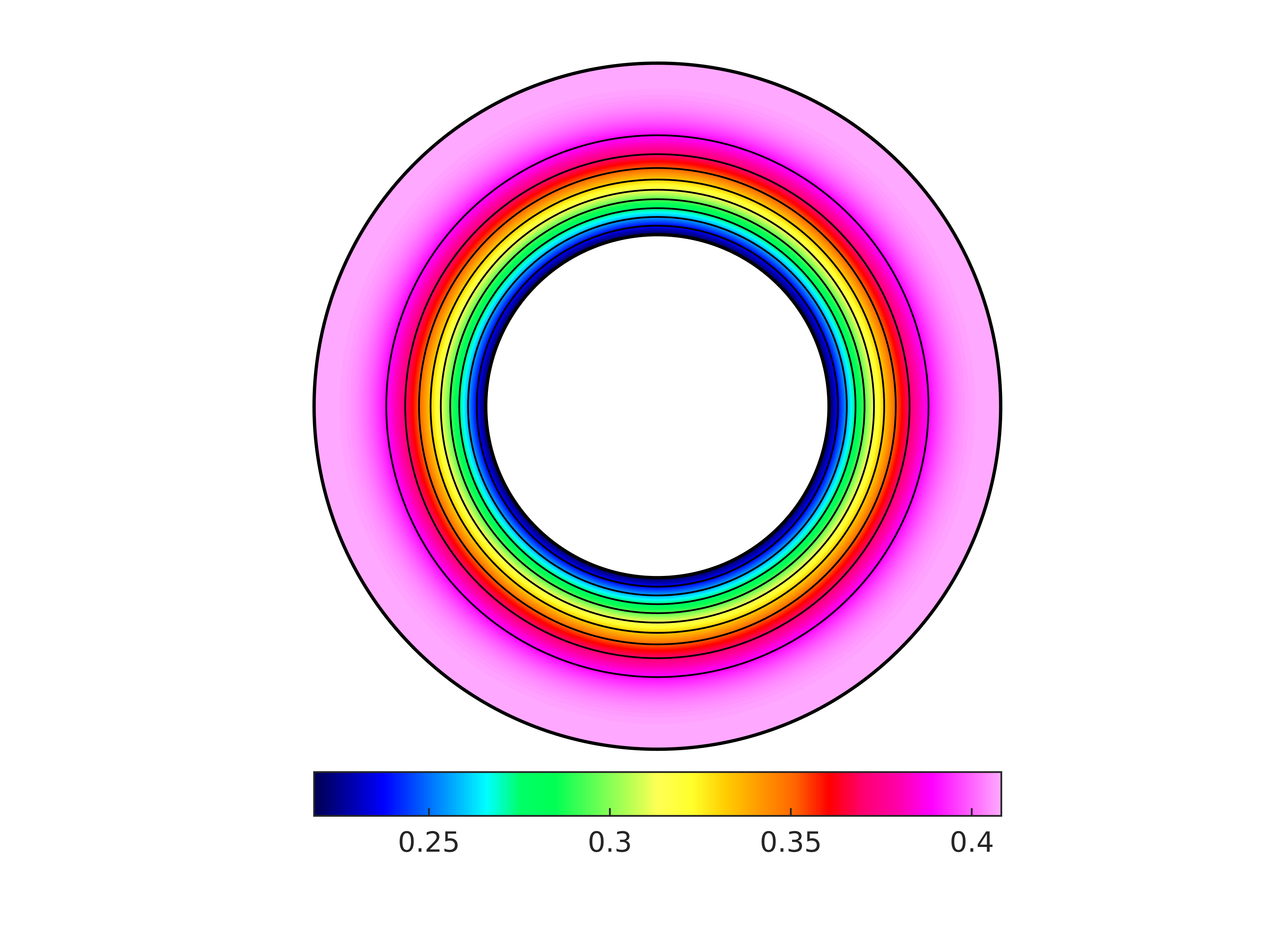}{\figWidtha}};
  \begin{scope}[xshift=8.07cm,yshift=3.3cm]
    \plotDisk
    \draw[very thick,->,red] (\rad*.5-.1,0) arc (0:320:.5cm); %
  \end{scope}
\end{tikzpicture}
} 
\end{center}
  \caption{Rotating disk. Computed speed (left) and pressure (right) for a massless disk, $\rhob=0$, using the grid $\Gc^{(4)}$.  The AMP-RB-VC scheme is used with $\adc=1$.}
  \label{fig:diskInDiskGrid}
\end{figure}
}

{
\newcommand{\figWidth}{7.5cm}
\newcommand{\figWidths}{7.25cm}
\newcommand{\trimfig}[2]{\trimFig{#1}{#2}{.0}{.0}{.0}{.0}}
\begin{figure}[htb]
\begin{center}
\resizebox{14cm}{!}{
\begin{tikzpicture}[scale=1]
  \useasboundingbox (0, 0) rectangle (16,6);  
  \draw(0,-0.75) node[anchor=south west,xshift=0pt,yshift=0pt] {\trimfig{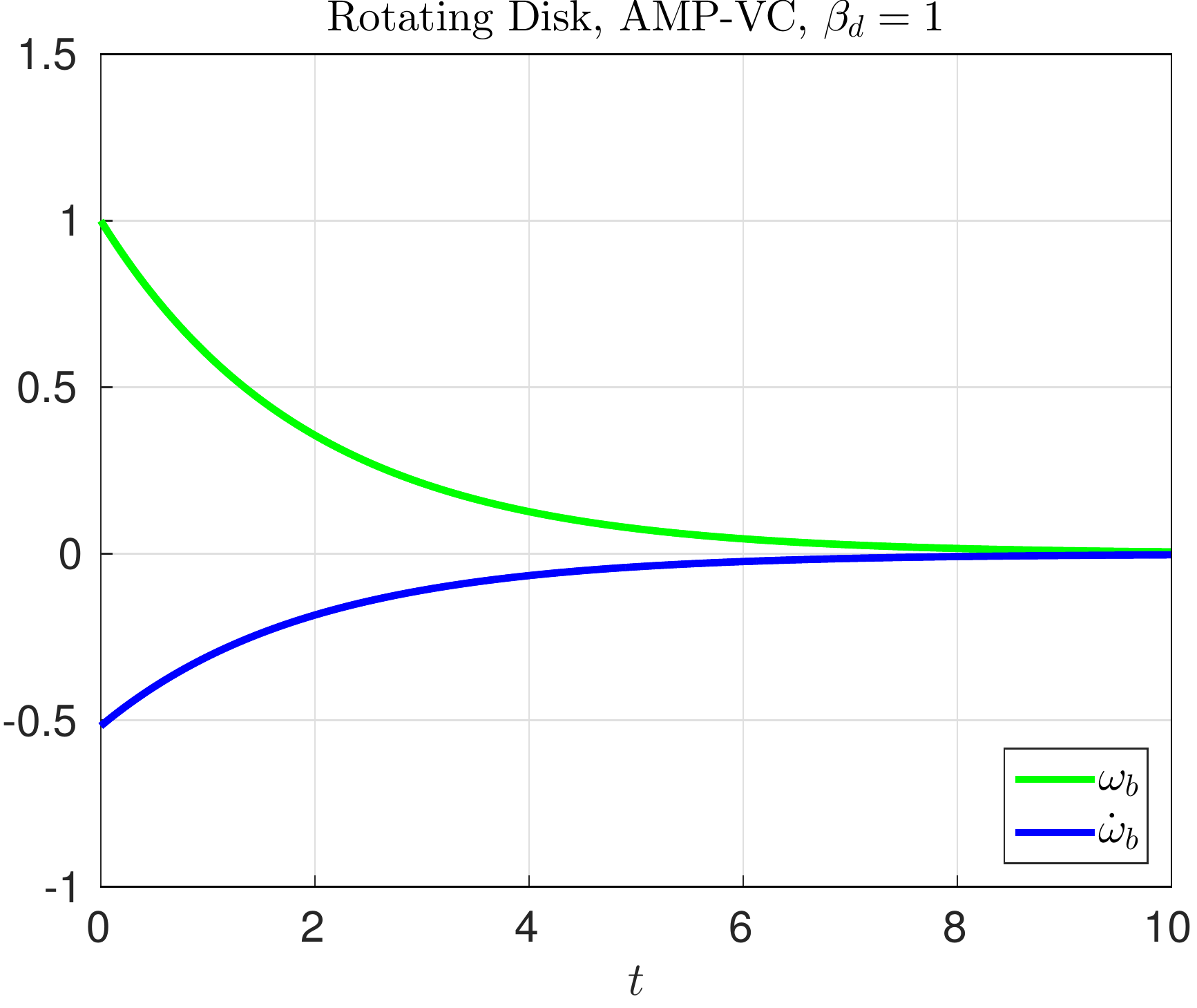}{\figWidth}};
  \draw(8,-0.75) node[anchor=south west,xshift=0pt,yshift=0pt] {\trimfig{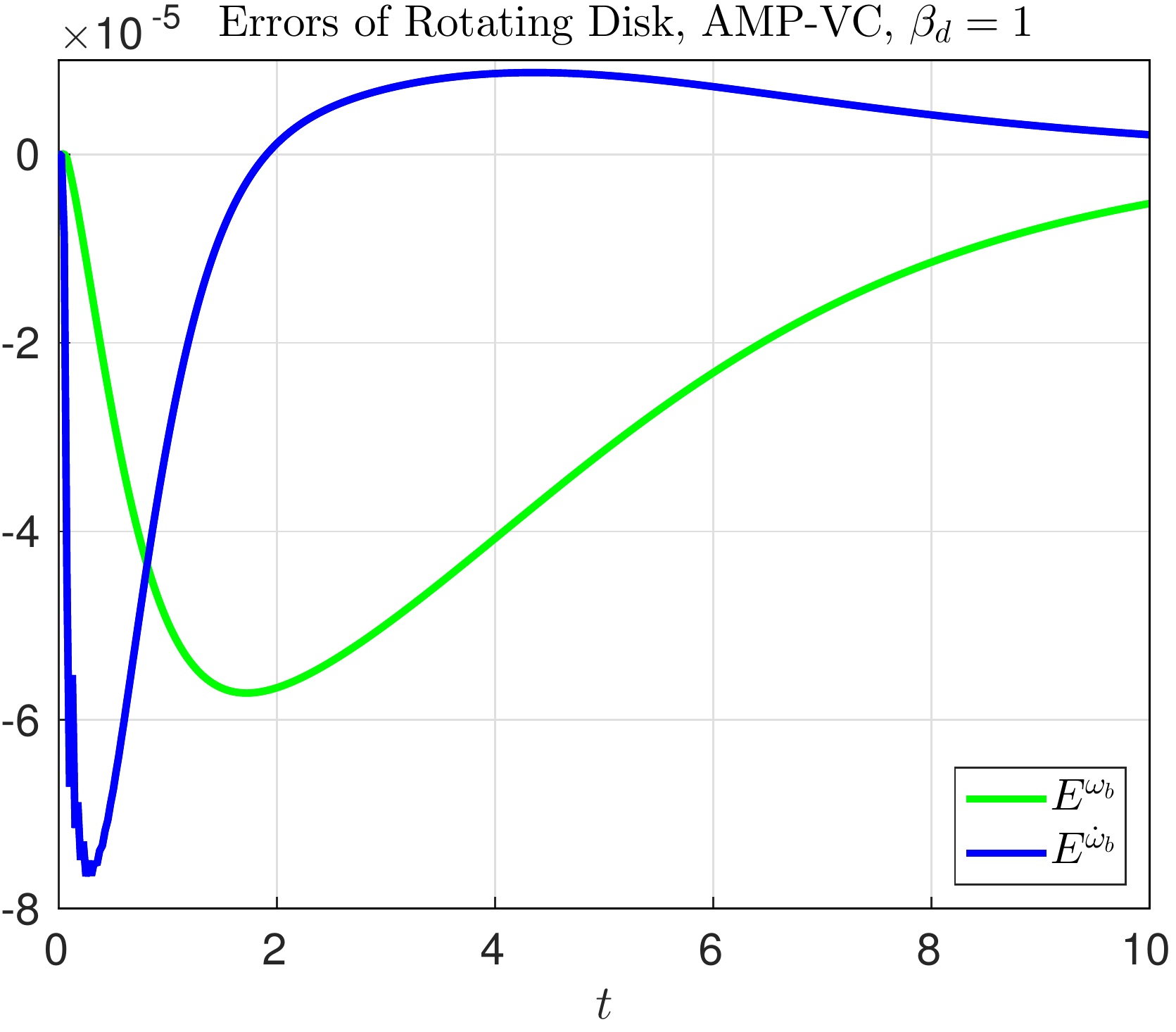}{\figWidths}};
\end{tikzpicture}
}
\end{center}
  \caption{Rotating disk.  Angular velocity and acceleration (left) and errors (right) of a massless disk, $\rhob=0$, using the grid $\Gc^{(4)}$.
   The~\ampRB-VC scheme is used with $\adc=1$.
    }
\label{fig:rotatingDiskStable}
\end{figure}
}

{
\newcommand{\errp}{$E_j^{p}$}
\newcommand{\erruv}{$E_j^{\vv}$}
\newcommand{\errwC}{$E_j^{\omega_3}$}
\newcommand{\errwtC}{$E_j^{\dot\omega_3}$}

\newcommand{\rotatingDiskTableI}{%
\begin{tabular}{|c|c|c|c|c|c|c|c|c|} \hline 
  \multicolumn{9}{|c|}{Rotating disk, $\rhob=10$} \\ \hline 
\strutt  $\dr^{(j)}$  &     \errp     &  r   &     \erruv     &  r   &     \errwC     &  r   &     \errwtC     &  r    \\[3pt] \hline 
1/10 & \num{3.9}{-3} &      & \num{1.5}{-3} &      & \num{4.3}{-4} &      & \num{4.4}{-4} &       \\ \hline 
1/20 & \num{1.1}{-3} &  3.6 & \num{4.3}{-4} &  3.4 & \num{1.2}{-4} &  3.5 & \num{1.3}{-4} &  3.4  \\ \hline 
1/40 & \num{2.8}{-4} &  3.9 & \num{1.2}{-4} &  3.6 & \num{3.5}{-5} &  3.5 & \num{3.5}{-5} &  3.7  \\ \hline 
1/80 & \num{7.1}{-5} &  4.0 & \num{3.2}{-5} &  3.8 & \num{9.3}{-6} &  3.7 & \num{9.1}{-6} &  3.8  \\ \hline 
1/160& \num{1.8}{-5} &  4.0 & \num{8.1}{-6} &  3.9 & \num{2.4}{-6} &  3.9 & \num{2.3}{-6} &  3.9  \\ \hline 
rate &      1.95     &      &      1.88     &      &      1.87     &      &      1.90     &      \\ \hline 
\end{tabular}
}
\newcommand{\rotatingDiskTableII}{%
\begin{tabular}{|c|c|c|c|c|c|c|c|c|} \hline 
  \multicolumn{9}{|c|}{Rotating disk, $\rhob=1$} \\ \hline 
\strutt  $\dr^{(j)}$  &     \errp     &  r   &     \erruv     &  r   &     \errwC     &  r   &     \errwtC     &  r    \\[3pt] \hline 
1/10 & \num{3.7}{-4} &      & \num{1.1}{-3} &      & \num{2.7}{-4} &      & \num{5.8}{-5} &       \\ \hline 
1/20 & \num{1.0}{-4} &  3.7 & \num{3.1}{-4} &  3.6 & \num{4.1}{-5} &  6.6 & \num{2.2}{-5} &  2.6  \\ \hline 
1/40 & \num{2.6}{-5} &  3.8 & \num{8.1}{-5} &  3.9 & \num{6.7}{-6} &  6.0 & \num{6.7}{-6} &  3.3  \\ \hline 
1/80 & \num{6.6}{-6} &  3.9 & \num{2.1}{-5} &  3.9 & \num{1.2}{-6} &  5.5 & \num{1.8}{-6} &  3.7  \\ \hline 
1/160& \num{1.7}{-6} &  3.9 & \num{5.2}{-6} &  4.0 & \num{2.5}{-7} &  4.9 & \num{4.6}{-7} &  3.8  \\ \hline 
rate &      1.95     &      &      1.95     &      &      2.51     &      &      1.76     &      \\ \hline 
\end{tabular}
}
\newcommand{\rotatingDiskTableIV}{%
\begin{tabular}{|c|c|c|c|c|c|c|c|c|} \hline 
  \multicolumn{9}{|c|}{Rotating disk, $\rhob=0$} \\ \hline 
\strutt  $\dr^{(j)}$  &     \errp     &  r   &     \erruv     &  r   &     \errwC     &  r   &     \errwtC     &  r    \\[3pt] \hline 
1/10 & \num{1.2}{-3} &      & \num{1.2}{-3} &      & \num{4.2}{-4} &      & \num{5.3}{-4} &       \\ \hline 
1/20 & \num{3.6}{-4} &  3.2 & \num{3.4}{-4} &  3.6 & \num{1.6}{-4} &  2.5 & \num{1.3}{-4} &  4.0  \\ \hline 
1/40 & \num{9.5}{-5} &  3.8 & \num{9.0}{-5} &  3.8 & \num{5.0}{-5} &  3.3 & \num{3.0}{-5} &  4.4  \\ \hline 
1/80 & \num{2.4}{-5} &  4.0 & \num{2.3}{-5} &  3.9 & \num{1.4}{-5} &  3.7 & \num{7.0}{-6} &  4.3  \\ \hline 
1/160& \num{6.1}{-6} &  3.9 & \num{5.8}{-6} &  4.0 & \num{3.5}{-6} &  3.9 & \num{1.7}{-6} &  4.0  \\ \hline 
rate &      1.91     &      &      1.94     &      &      1.74     &      &      2.07     &      \\ \hline 
\end{tabular}
}
{
\begin{table}[htb]\tableFont 
\begin{center}
  \rotatingDiskTableI \\
\medskip
  \rotatingDiskTableII \\
\medskip
  \rotatingDiskTableIV
\caption{Rotating disk. Maximum errors and estimated convergence rates at $t=1$ computed using the~\ampRB~scheme
   for a heavy, $\rhob=10$, medium, $\rhob=1$ and massless, $\rhob=0$, rigid disk. }
\label{table:rotatingDisk}
\end{center}
\end{table}
}
} 

{
\newcommand{\labelFont}{\scriptsize}
\newcommand{\trimfig}[2]{\trimw{#1}{#2}{.0}{.0}{.0}{.0}}
\newcommand{\trimfigz}[2]{\trimw{#1}{#2}{.05}{.05}{.03}{.05}}
\newcommand{\figWidth}{4.0cm}
\newcommand{\figWidths}{3.9cm}
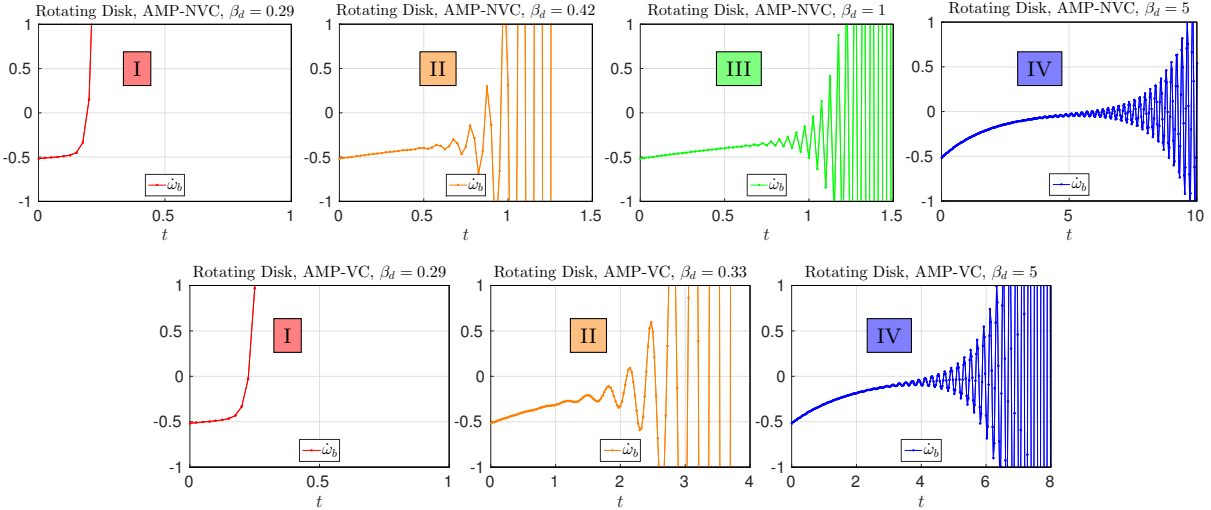
\begin{figure}[htb]
\begin{center}
\begin{tikzpicture}[scale=1]
  \useasboundingbox (0.0,1) rectangle (16,7.2);  
 \begin{scope}[xshift=1.5cm,yshift=.25cm]
    \draw(0.0,0.0) node[anchor=south west,xshift=-4pt,yshift=+0pt] {\trimfig{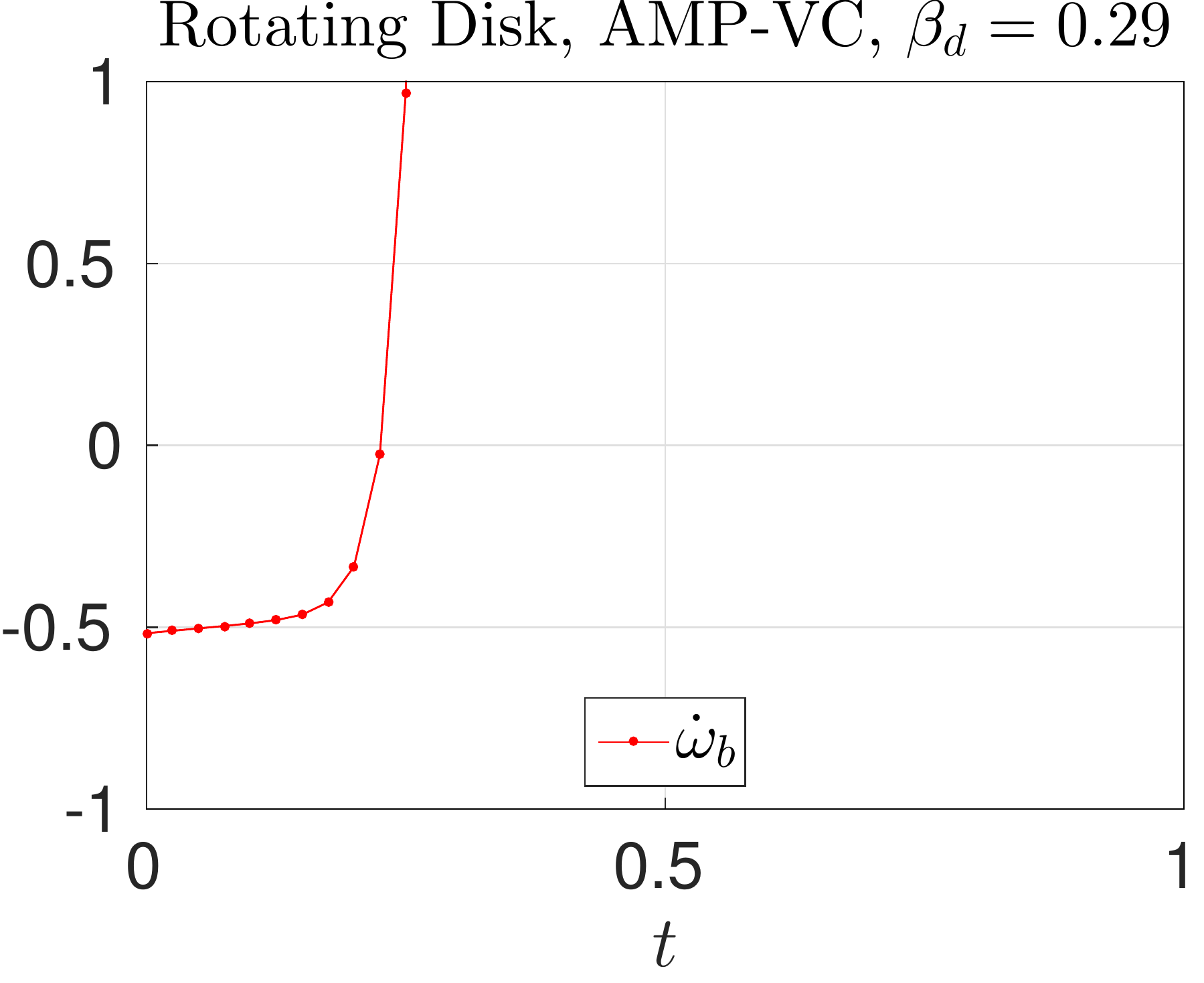}{\figWidth}};
    \draw(2.2,2.5) node[draw,fill=red!50] {\labelFont I};
    \draw(4.0,0.0) node[anchor=south west,xshift=-4pt,yshift=+0pt] {\trimfig{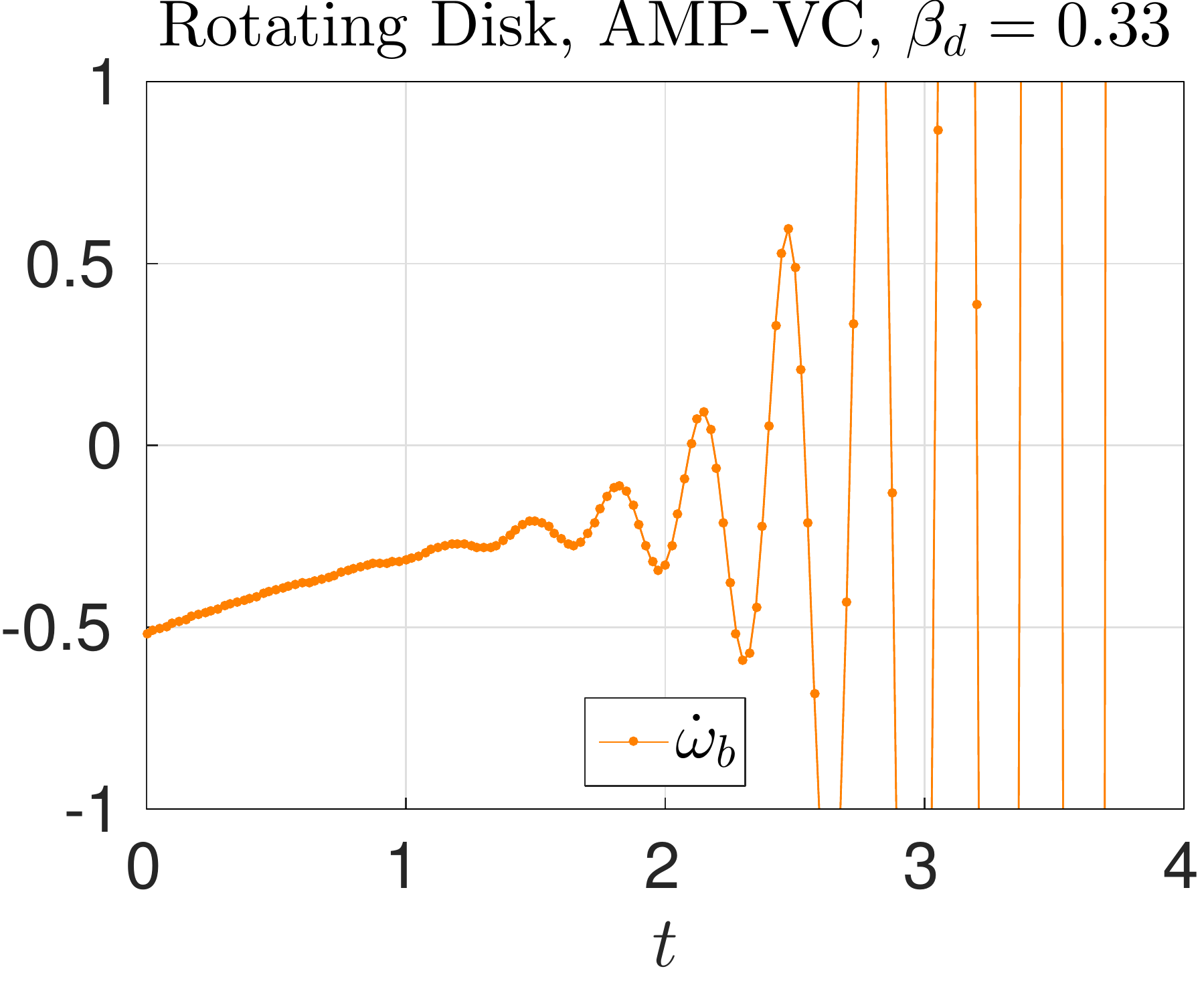}{\figWidth}};
    \draw(6.2,2.5) node[draw,fill=orange!50] {\labelFont II};
    \draw(8.0,0.0) node[anchor=south west,xshift=-4pt,yshift=+0pt] {\trimfig{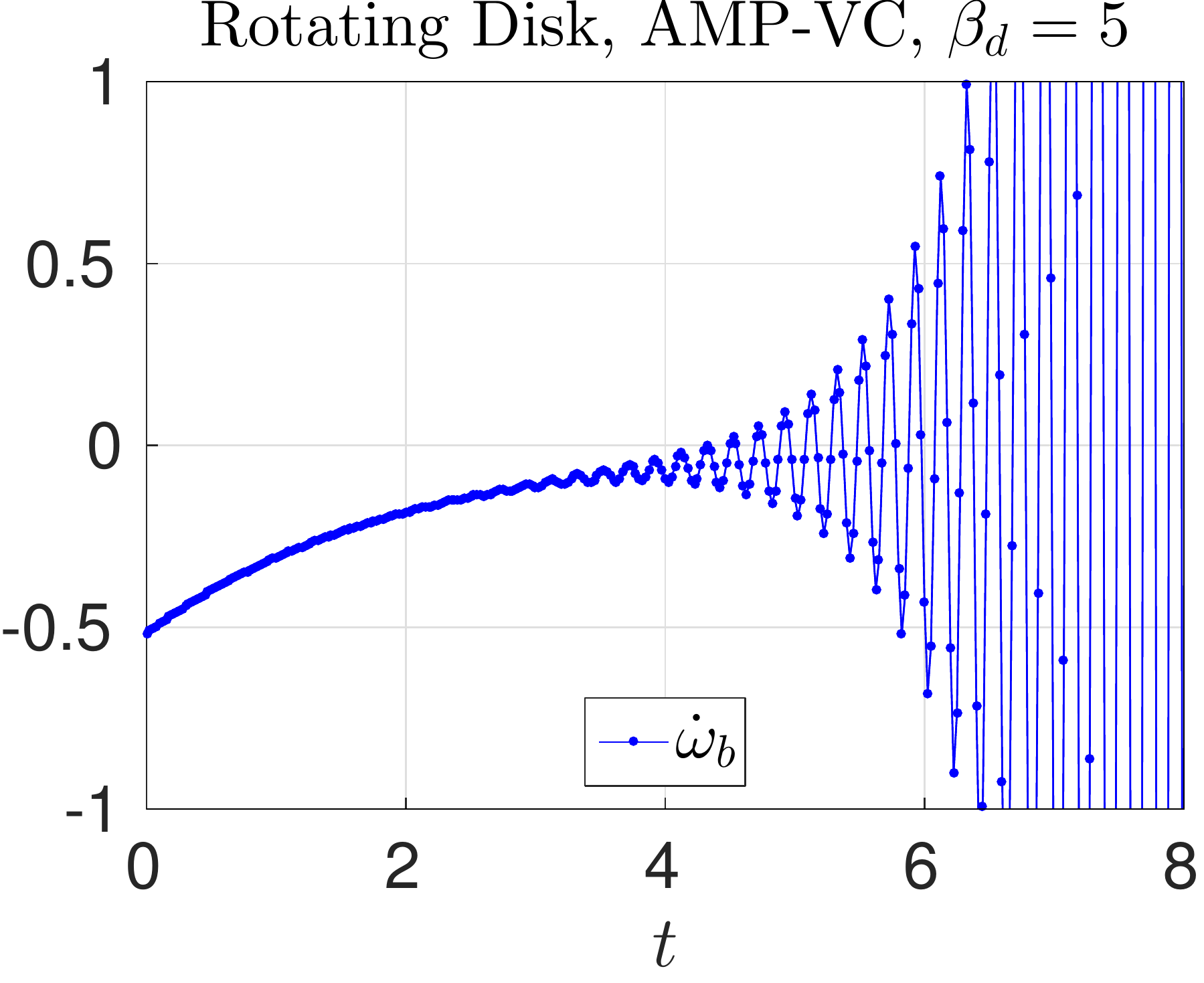}{\figWidth}};
    \draw(10.2,2.5) node[draw,fill=blue!50] {\labelFont IV};
  \end{scope}

  \begin{scope}[xshift=-.5cm,yshift=3.8cm]
    \draw(0.0,0.0) node[anchor=south west,xshift=-4pt,yshift=+0pt] {\trimfig{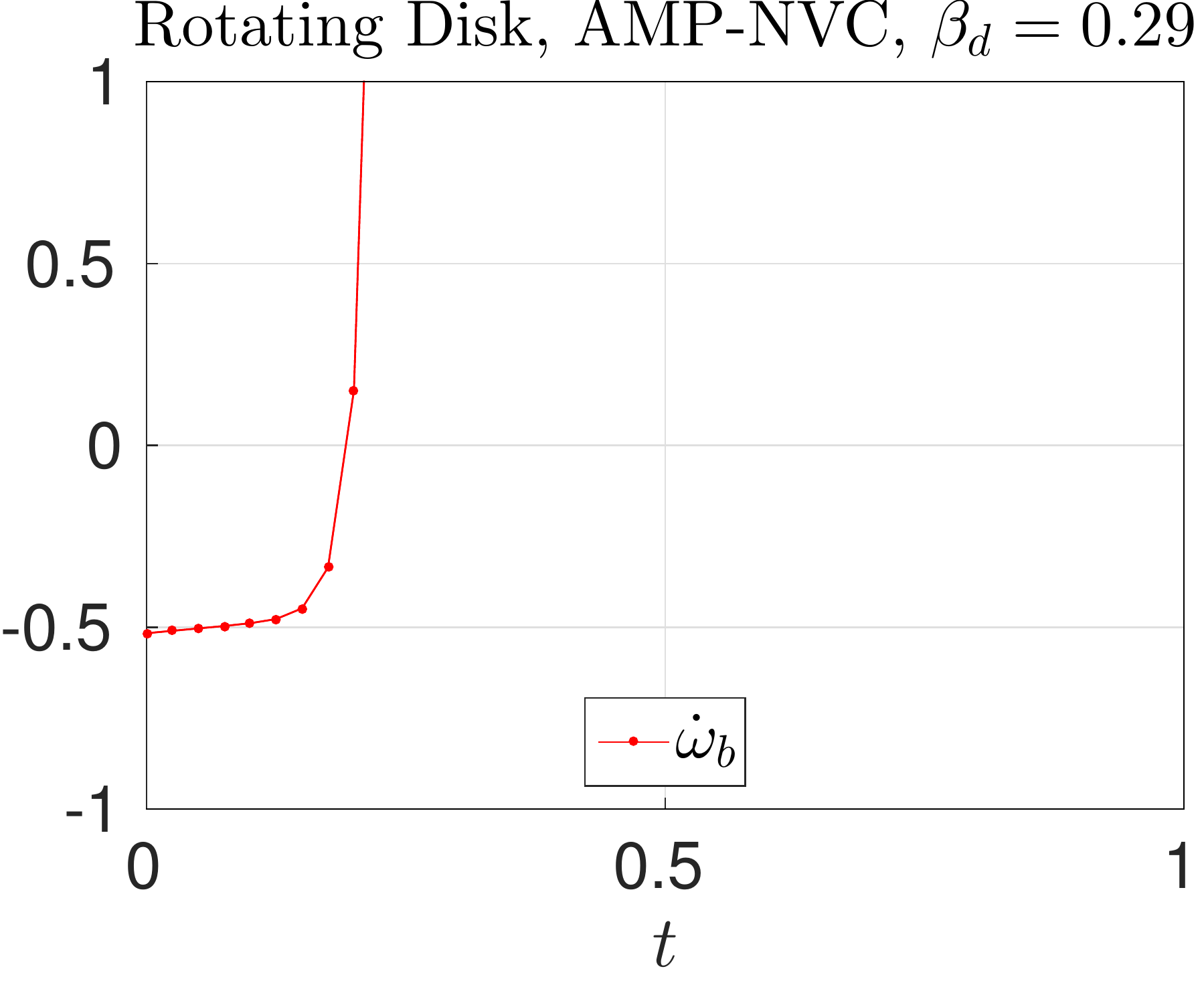}{\figWidths}};
    \draw(2.2,2.5) node[draw,fill=red!50] {\labelFont I};
    \draw(4.0,0.0) node[anchor=south west,xshift=-4pt,yshift=+0pt] {\trimfig{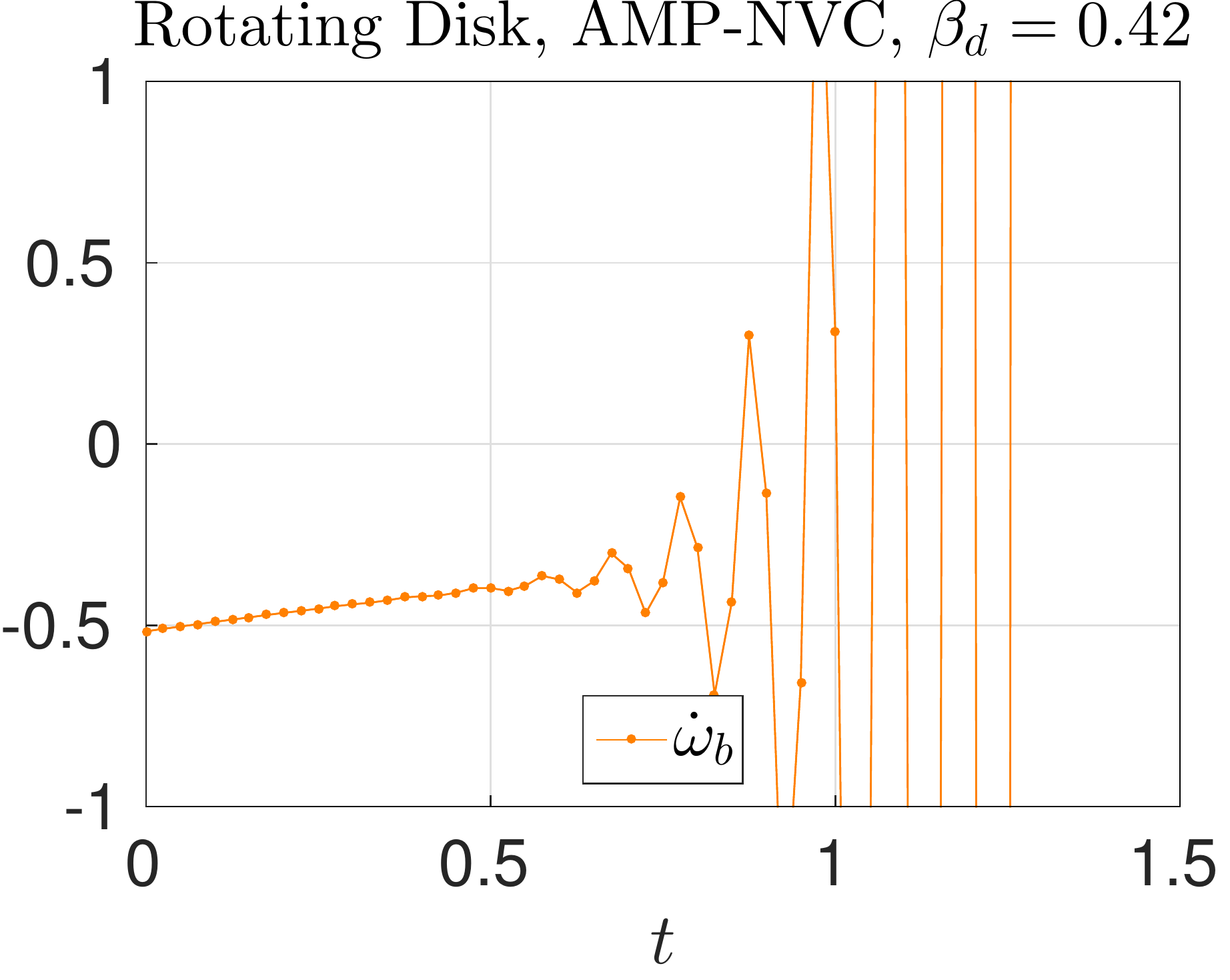}{\figWidth}};
    \draw(6.2,2.5) node[draw,fill=orange!50] {\labelFont II};
	\draw(8.0,0.0) node[anchor=south west,xshift=-4pt,yshift=+0pt] {\trimfig{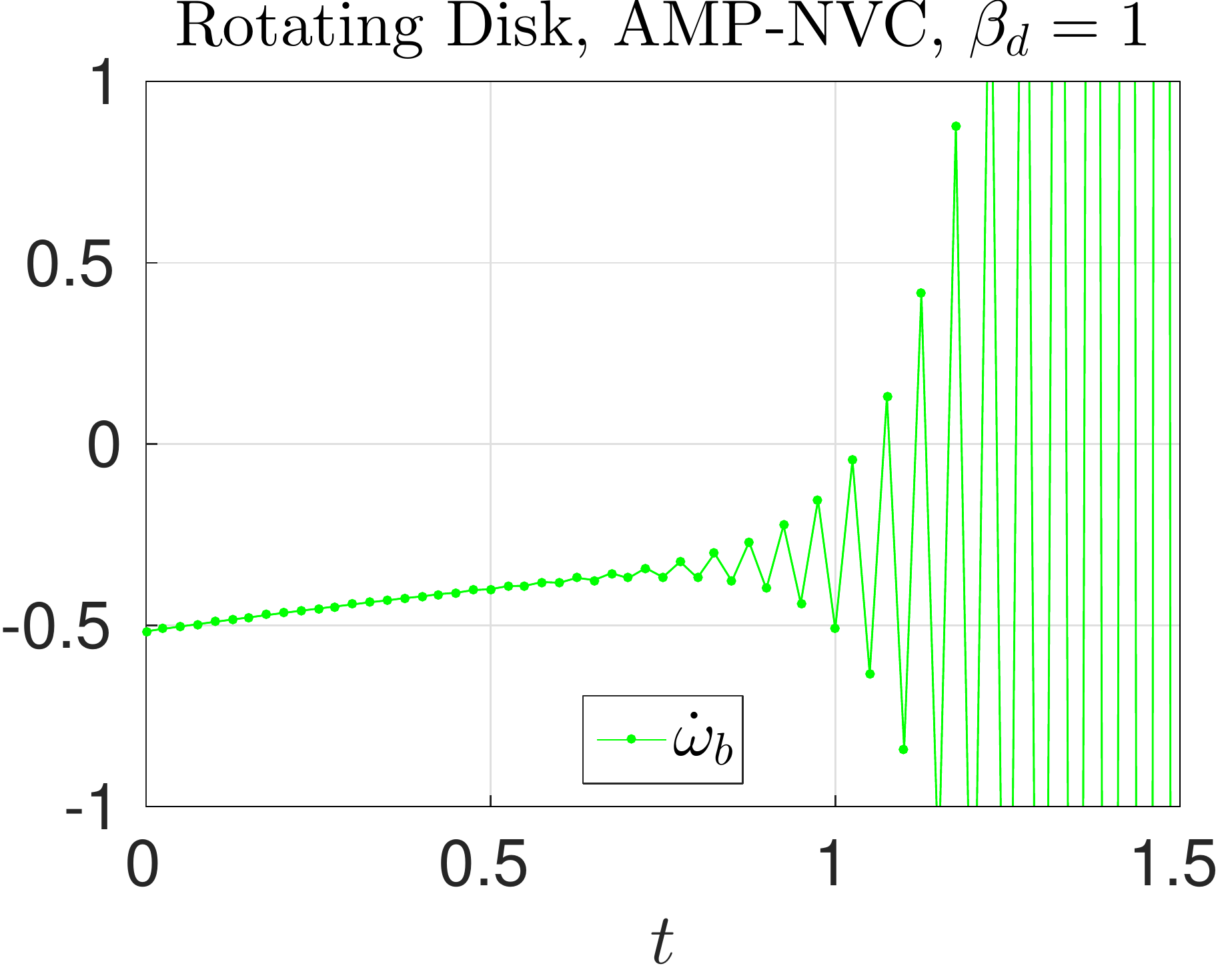}{\figWidth}};
    \draw(10.2,2.5) node[draw,fill=green!50] {\labelFont III};
    \draw(12.0,0.0) node[anchor=south west,xshift=-4pt,yshift=+0pt] {\trimfig{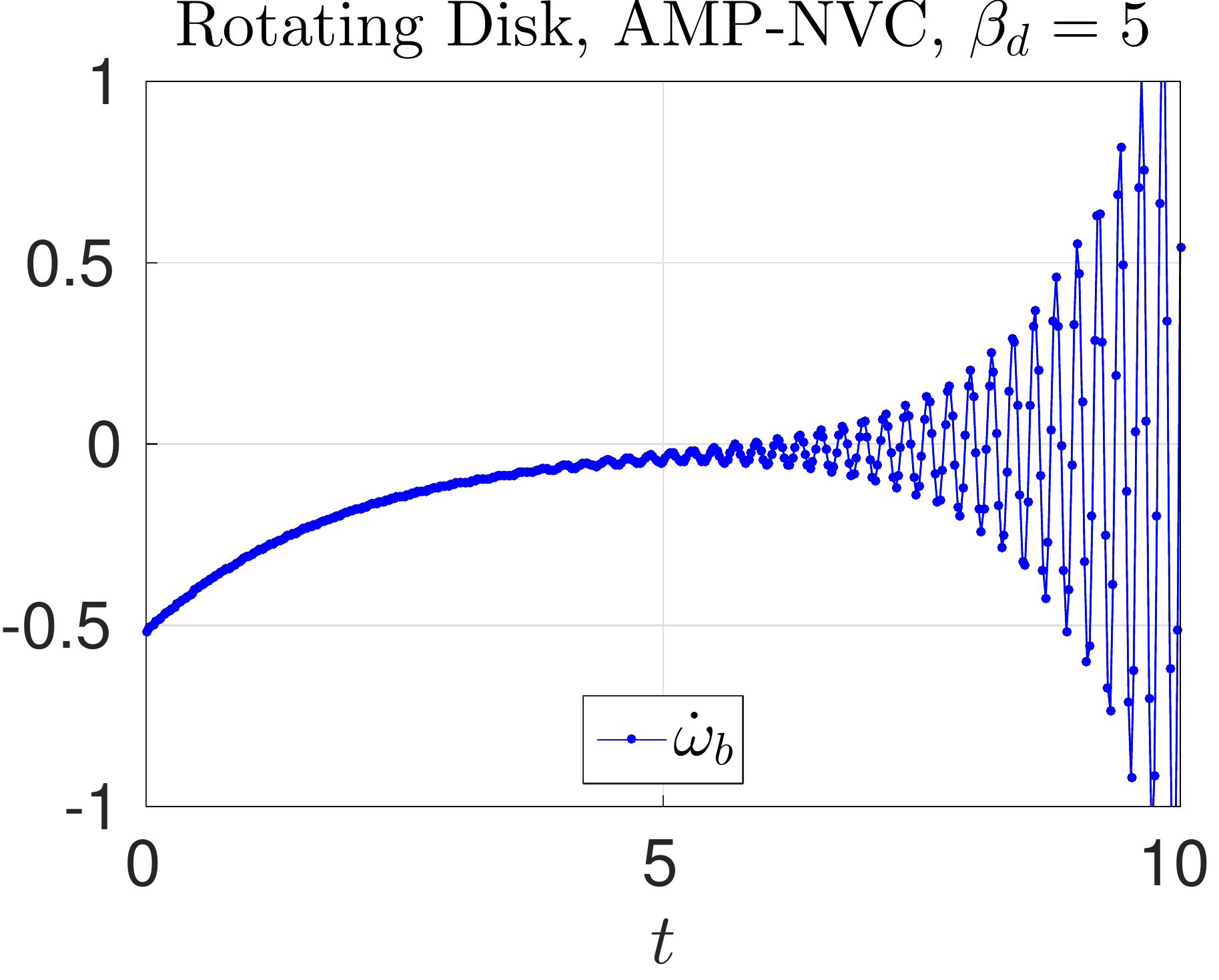}{\figWidth}};
    \draw(14.2,2.5) node[draw,fill=blue!50] {\labelFont IV};
  \end{scope}
%
\end{tikzpicture}
\end{center}
  \caption{Rotating disk.  Angular acceleration of a rotating light rigid body with $\rhob=0.001$ using the grid $\Gc^{(4)}$.
   Top: typical instabilities in the {\ampRB}-NVC~scheme: 
   (I)~for $\adc = 0.29$    there is a rapid,  monotonically growing instability;
   (II)~for $\adc = 0.42$   there is a rapidly growing, high-frequency instability;
   (III)~for $\adc = 1$     there is a rapidly growing, step-to-step instability; and
   (IV)~for $\adc = 5$      there is a slowly growing, low-frequency instability. 
   Bottom:  typical instabilities in the {\ampRB}-VC~scheme:
   (I)~for $\adc =0.29$ 	there is a rapid, monotonically growing instability;
   (II)~for $\adc = 0.33$ 	there is a rapidly growing, high-frequency instability; and
   (IV)~for $\adc = 5$   	there is a slowly growing, low-frequency instability. 
    }
  \label{fig:shearDiskCurves}
\end{figure}
}

Figures~\ref{fig:diskInDiskGrid} shows shaded contours of the speed and pressure of the fluid at
$t=1$, while Figure~\ref{fig:rotatingDiskStable} shows the behaviour of the angular
velocity and acceleration of the disk along with the errors in these quantities.  The solution illustrated in these two figures
is computed for a massless disk ($\rhob=0$) using the {\ampRB}-VC~scheme on the grid $\Gc^{(4)}$ with $\adc = 1$.
In addition, Table~\ref{table:rotatingDisk} presents the errors at $t=1$ and estimated 
convergence rates for simulations with $\rhob=0$, $1$, and
$10$. These convergence rates demonstrate the second-order accuracy of the \ampRB~approach even
for the difficult case of a massless body in which added-damping effects are strong.
While the exact solution for the rigid body
assumes only rotation, the numerical implementation does not make this assumption
so that both rotation and translation of the body are allowed.
Therefore, it is important to include stabilization against
added-mass effects, as discussed previously, although the computed
translational velocities are very close to zero (in good agreement with the exact solution).

Finally, we probe the classifications of unstable modes, similar to what was done in the the prior
case of the sliding block in Section~\ref{sec:shearBlock}. Figure~\ref{fig:shearDiskCurves}
presents time traces of the angular acceleration for a very light rigid body with $\rhos=0.001$
computed using the grid $\Gs^{(4)}$. Different values of the added-damping parameter~$\adc$ are used to elicit
unstable behaviours of the various types. Even though the \ampRB-VC scheme is stable for a range of values of
$\adc$ about $\adc=1$ as predicted by the theory and demonstrated in previous figures, instabilities occur
in the \ampRB-VC and \ampRB-NVC~schemes for values of $\adc$ outside the stability region.  We observe that the behaviours of
the instabilities shown in Figure~\ref{fig:shearDiskCurves} are in good agreement with those predicted
by the analysis of Section~\ref{sec:analysisDisk}.

\section{Conclusions} \label{sec:conclusions}

In this first part of a two-part series,
a partitioned scheme for simulating fluid structure interactions between
incompressible flows and rigid bodies is developed and analyzed.
This new~\ampRB~scheme is stable, without sub-time-step
iterations, for bodies of any mass, including bodies of zero mass.
A fractional-step time-stepping scheme
is used to evolve the incompressible fluid where the velocity and pressure are solved in separate stages.
Coupling of the rigid-body equations occurs in the pressure equation where the degrees of freedom
corresponding to rigid-body accelerations are solved in conjunction with the fluid pressure, and 
both added-mass and added-damping effects are effectively addressed by the {\ampRB} scheme.
Added-mass effects are treated with a generalized Robin interface 
condition that couples the fluid pressure and body accelerations.
Added-damping effects are handled by introducing an added-damping 
coefficient into the equations of motion for the rigid body.
The form of this added-damping coefficient (which in general geometries becomes
an added-damping tensor) is derived from a linearization of
the integrals that define the force and torques on the body.
An approximate form for the added-damping coefficient, suitable for
use with general bodies, is determined from a local analysis and variational problem.

The stability properties of the proposed interface conditions are studied for 
a rectangular-geometry and an annular-geometry model-problem
using mode analysis. 
The analyses of the model problems for added-mass effects show that the~{\ampRB} scheme 
is stable for bodies of any mass including bodies of zero mass. 
On the other hand, the analysis of the added-damping model problem reveals that for 
very light bodies a final velocity-correction stage should be used to ensure stability. The 
mode analysis is used to determine the values for the stabilizing
added-damping parameter $\adc$ (which scales the added-damping coefficient)
to ensure stability of the scheme. It is shown that there is range of values for $\adc$ 
where the scheme remains stable independent of the problem parameters $\nu$, $\dt$, grid spacing, body mass, and body moment-of-inertia. As a result the parameter $\adc$ can chosen once and for all and does not need to be adjusted from problem to problem.
With this additional correction phase, the {\ampRB} scheme is shown to be stable against added-damping 
effects for any body mass, including the zero mass case.

Numerical results are provided for model problems in rectangular and annular geometries 
that confirm the results of the theoretical analyses. 
These results verify stability and second-order accuracy of the {\ampRB} scheme, 
including for bodies of zero mass.
Part II of this work considers the
extension of the {\ampRB} scheme to general geometries.

\appendix
\section{Exact solution for a translating disk in an annular Stokes flow} \label{sec:exactTranslatingDisk}

An exact solution to an FSI problem involving the horizontal motion of a solid disk in an incompressible Stokes flow is described.  The solution is used in Section~\ref{sec:translatingDisk} to access the accuracy of the {\ampRB}~scheme for an FSI problem in an annular geometry.  The configuration of the problem is shown in Figure~\ref{fig:modelDiskInDisk} and the equations governing the motion of the body and the velocity and pressure of the fluid are given in~\eqref{eq:velocityEquation1A}--\eqref{eq:modelProblemLastA}.  As before, we assume that the amplitude of the motion is small and linearize the problem about the fixed interface at $r=\aa$.  Using the transformation of the fluid variables in~\eqref{eq:transformation}, the governing equations become
\begin{alignat}{2}
 \text{Fluid:}\quad & \rho{\partial\hat u\over\partial t}+{\partial\hat p\over\partial r}=\mu\left[{1\over r}{\partial\over\partial r}\left(r{\partial\hat u\over\partial r}\right)-{2\hat u\over r\sp2}+{2\hat v\over r\sp2}\right] , \quad&& \aa<r<\bb, \label{eq:hatVelocityEquation1A}\\
                    & \rho{\partial\hat v\over\partial t}+{\hat p\over r}=\mu\left[{1\over r}{\partial\over\partial r}\left(r{\partial\hat v\over\partial r}\right)-{2\hat v\over r\sp2}+{2\hat u\over r\sp2}\right], \quad&& \aa<r<\bb, \label{eq:hatVelocityEquation2A}\\
  &  {\partial\over\partial r}(r\hat u)-\hat v = 0   , \quad&&\aa<r<\bb, \label{eq:hatContinuityEquationA}\\
\text{Rigid body:}\quad  &   m_b\,{d u_b\over dt} = \aa\pi\left.\left[-\hat p +\mu{\partial\hat v\over\partial r}\right]\right|_{r=\aa}, \quad && \label{eq:hatBodyEquation1A}\\
 \text{Interface:}\quad
&   \hat u(\aa,t)=\hat v(\aa,t)=u_b(t), \label{eq:hatMatchVA} \\
 \text{Fluid BCs:} \quad 
&  \hat u(\bb,t)=\hat v(\bb,t)=0 . \label{eq:hatModelProblemLastA}
\end{alignat}
Here, we have assumed that the body does not rotate so that $\omega_b(t)=0$ and that there is no external force on the body, i.e.~$g_u(t)=0$.  The equations are homogeneous and the motion is initiated by the (non-zero) initial conditions which we specify later.

Separable solutions of the linear equations for the fluid in~\eqref{eq:hatVelocityEquation1A}--\eqref{eq:hatContinuityEquationA} have the form
\begin{equation}
\hat p(r,t)=\mu\left(A_pr+{B_p\over r}\right)e\sp{-\lambda\sp2\nu t},
\label{eq:hatPressureSoln}
\end{equation}
and
\begin{equation}
\hat u(r,t)=\left({\hat Q(r)\over r}\right)e\sp{-\lambda\sp2\nu t},\qquad \hat v(r,t)=\hat Q\sp\prime(r)e\sp{-\lambda\sp2\nu t},
\label{eq:hatVelocitySoln}
\end{equation}
where
\begin{equation}
\hat Q(r)=A_1J_1(\lambda r)+B_1Y_1(\lambda r)+{1\over\lambda\sp2}\left(A_pr-{B_p\over r}\right).
\label{eq:hatQ}
\end{equation}
Here, $J_1$ and $Y_1$ are Bessel functions of order one, $\nu=\mu/\rho$ is the kinematic viscosity, and $A_p$, $B_p$, $A_1$ and $B_1$ are constants.  For later purposes, define the following constants:
\begin{alignat*}{2}
b_{11}& ={\pi \bb \over 2 \aa}\Bigl(2 Y_1(\lambda \bb)-\lambda \bb Y_0(\lambda \bb)\Bigr),
\qquad &&
b_{12}={\pi \over2}\lambda \aa Y_0(\lambda \bb), \\
b_{21} & ={\pi \bb\over 2\aa}\Bigl(2J_1(\lambda \bb)-\lambda \bb J_0(\lambda \bb)\Bigr), &&
b_{22}={\pi \over2}\lambda \aa J_0(\lambda \bb),
\end{alignat*}
and
\begin{alignat*}{2}
c_{11} & ={b_{11}J_1(\lambda \aa)-b_{21}Y_1(\lambda \aa)+1\over\lambda\sp2 \aa\sp2}, &&
c_{12}={b_{12}J_1(\lambda \aa)-b_{22}Y_1(\lambda \aa)-1\over\lambda\sp2 \aa\sp2}, \\
c_{21}& ={b_{11}\lambda \aa J_1\sp\prime(\lambda \aa)-b_{21}\lambda \aa Y_1\sp\prime(\lambda \aa)+1\over\lambda\sp2 \aa \sp 2} ,\qquad &&
c_{22}={b_{12}\lambda \aa J_1\sp\prime(\lambda \aa)-b_{22}\lambda \aa Y_1\sp\prime(\lambda \aa)+1\over\lambda\sp2 \aa\sp2}.
\end{alignat*}
Using~\eqref{eq:hatVelocitySoln} in the boundary conditions~\eqref{eq:hatMatchVA} and~\eqref{eq:hatModelProblemLastA}, together with identities for Bessel functions, gives
\[
A_1 = \frac{1}{\aa \lambda^2} \left[ \aa^2 b_{11} A_p + b_{12} B_p \right], \qquad
B_1 = -\frac{1}{\aa \lambda^2} \left[ \aa^2 b_{21} A_p + b_{22} B_p \right],
\]
and
\[
A_p= {\amplitude\over\aa\sp2}\left[{c_{22}-c_{12}\over c_{11}c_{22}-c_{12}c_{21}}\right],\qquad
B_p= \amplitude\left[{c_{11}-c_{21}\over c_{11}c_{22}-c_{12}c_{21}}\right],
\]
where $\amplitude$ is the (small) amplitude of the velocity of the body given by
\begin{equation}
u_b(t)=\amplitude e\sp{-\lambda\sp2\nu t}.
\label{eq:bodySoln}
\end{equation}
Finally, using~\eqref{eq:hatPressureSoln}, \eqref{eq:hatVelocitySoln} and~\eqref{eq:bodySoln} in the equation governing the motion of the rigid body~\eqref{eq:hatBodyEquation1A} gives the transcendental equation
\begin{align}
&\frac{\mb}{\rho \pi } \lambda\sp2  + \bigl[-1+b_{11}J_1\sp{\prime\prime}(\lambda \aa)-b_{21}Y_1\sp{\prime\prime}(\lambda \aa)\bigr]\left[{c_{22}-c_{12}\over c_{11}c_{22}-c_{12}c_{21}}\right] \nonumber \\
&\qquad + \left[-1+b_{12}J_1\sp{\prime\prime}(\lambda \aa)-b_{22}Y_1\sp{\prime\prime}(\lambda \aa)-{2\over\lambda\sp2\aa\sp2}\right]\left[{c_{11}-c_{21}\over c_{11}c_{22}-c_{12}c_{21}}\right] = 0,
\label{eq:transDiskEigenvalues}
\end{align}
which determines the eigenvalues $\lambda$ of the separable solutions.

There are an infinite number of solutions of~\eqref{eq:transDiskEigenvalues} for a given set of parameters of the problem.  We are interested in the smallest (positive) value of $\lambda$ which corresponds to the solution with the slowest rate of decay.  For example, if $\mb=\rho_b\pi\aa\sp2$, $\rho=1$, $\aa=1$ and $\bb=2$, then the smallest values of $\lambda$ can be computed by solving~\eqref{eq:transDiskEigenvalues} numerically for different values of $\rho_b$.  Some of these solutions are collected in Table~\ref{table:cylEig}.  Finally, the initial conditions for the problem are taken from~\eqref{eq:hatVelocitySoln} and~\eqref{eq:bodySoln} for a computed value of $\lambda$ upon setting $t=0$.

\begin{table}[H]\tableFont
\begin{center}
\begin{tabular}{|c|c|c|c|c|}\hline
             & $\rhos=10$ & $\rhos=1$ & $\rhos=10^{-3}$ & $\rhos=0$  \\ \hline 
$\lambda$    & 1.84604165628237 & 3.40692142656753 & 3.94155013131292 & 3.94219755667612 \\ \hline 
\end{tabular}
\caption{Eigenvalues, $\lambda$, in the exact solution of the translating disk for different values of $\rhob$.}
\label{table:cylEig}
\end{center}
\end{table}

\bibliographystyle{elsart-num}
\bibliography{journal-ISI,jwb,henshaw,henshawPapers,fsi}

\begin{thebibliography}{10}
\expandafter\ifx\csname url\endcsname\relax
  \def\url#1{\texttt{#1}}\fi
\expandafter\ifx\csname urlprefix\endcsname\relax\def\urlprefix{URL }\fi

\bibitem{takashi1992}
N.~Takashi, T.~J. Hughes, An arbitrary {Lagrangian-Eulerian} finite element
  method for interaction of fluid and a rigid body, Comput. Method. Appl. Mech.
  Engrg. 95~(1) (1992) 115--138.

\bibitem{HuPatankarZhu2001}
H.~H. Hu, N.~A. Patankar, M.~Y. Zhu, Direct numerical simulations of
  fluid-solid systems using the arbitrary langrangian-eulerian technique, J.
  Comput. Phys. 169~(2) (2001) 427--462.

\bibitem{vierendeels2005analysis}
J.~Vierendeels, K.~Dumont, E.~Dick, P.~Verdonck, Analysis and stabilization of
  fluid-structure interaction algorithm for rigid-body motion, AIAA J. 43~(12)
  (2005) 2549--2557.

\bibitem{coquerelle2008vortex}
M.~Coquerelle, G.-H. Cottet, A vortex level set method for the two-way coupling
  of an incompressible fluid with colliding rigid bodies, J. Comput. Phys.
  227~(21) (2008) 9121--9137.

\bibitem{gibou2012efficient}
F.~Gibou, C.~Min, Efficient symmetric positive definite second-order accurate
  monolithic solver for fluid/solid interactions, J. Comput. Phys. 231~(8)
  (2012) 3246--3263.

\bibitem{glowinski1999distributed}
R.~Glowinski, T.-W. Pan, T.~I. Hesla, D.~D. Joseph, A distributed {Lagrange}
  multiplier/fictitious domain method for particulate flows, Int. J. Multiphase
  Flow 25~(5) (1999) 755--794.

\bibitem{glowinski2000distributed}
R.~Glowinski, T.-W. Pan, T.~I. Hesla, D.~D. Joseph, J.~P\'eriaux, A distributed
  {Lagrange} multiplier/fictitious domain method for the simulation of flow
  around moving rigid bodies: application to particulate flow, Comput. Method.
  Appl. Mech. Engrg. 184~(2) (2000) 241--267.

\bibitem{GlowinskiPanHelsaJosephPeriaux2001}
R.~Glowinski, T.~W. Pan, T.~I. Hesla, D.~D. Joseph, J.~P\'eriaux, A fictitious
  domain approach to the direct numerical simulation of incompressible viscous
  flow past moving rigid bodies: Applications to particulate flow, J. Comput.
  Phys. 169 (2001) 363--426.

\bibitem{costarelli2016embedded}
S.~D. Costarelli, L.~Garelli, M.~A. Cruchaga, M.~A. Storti, R.~Ausensi, S.~R.
  Idelsohn, An embedded strategy for the analysis of fluid structure
  interaction problems, Comput. Method. Appl. Mech. Engrg. 300 (2016) 106--128.

\bibitem{kajishima2002interaction}
T.~Kajishima, S.~Takiguchi, Interaction between particle clusters and
  particle-induced turbulence, Int. J. Heat Fluid Flow 23~(5) (2002) 639--646.

\bibitem{uhlmann2005immersed}
M.~Uhlmann, An immersed boundary method with direct forcing for the simulation
  of particulate flows, J. Comput. Phys. 209~(2) (2005) 448--476.

\bibitem{kim2006immersed}
D.~Kim, H.~Choi, Immersed boundary method for flow around an arbitrarily moving
  body, J. Comput. Phys. 212~(2) (2006) 662--680.

\bibitem{lee2008immersed}
T.-R. Lee, Y.-S. Chang, J.-B. Choi, D.~W. Kim, W.~K. Liu, Y.-J. Kim, Immersed
  finite element method for rigid body motions in the incompressible
  {Navier}--{Stokes} flow, Comput. Method. Appl. Mech. Engrg. 197~(25) (2008)
  2305--2316.

\bibitem{borazjani2008curvilinear}
I.~Borazjani, L.~Ge, F.~Sotiropoulos, Curvilinear immersed boundary method for
  simulating fluid structure interaction with complex {3D} rigid bodies, J.
  Comput. Phys. 227~(16) (2008) 7587--7620.

\bibitem{breugem2012second}
W.-P. Breugem, A second-order accurate immersed boundary method for fully
  resolved simulations of particle-laden flows, J. Comput. Phys. 231~(13)
  (2012) 4469--4498.

\bibitem{kempe2012improved}
T.~Kempe, J.~Fr{\"o}hlich, An improved immersed boundary method with direct
  forcing for the simulation of particle laden flows, J. Comput. Phys. 231~(9)
  (2012) 3663--3684.

\bibitem{yang2012simple}
J.~Yang, F.~Stern, A simple and efficient direct forcing immersed boundary
  framework for fluid--structure interactions, J. Comput. Phys. 231~(15) (2012)
  5029--5061.

\bibitem{bhalla2013unified}
A.~P.~S. Bhalla, R.~Bale, B.~E. Griffith, N.~A. Patankar, A unified
  mathematical framework and an adaptive numerical method for fluid--structure
  interaction with rigid, deforming, and elastic bodies, J. Comput. Phys. 250
  (2013) 446--476.

\bibitem{yang2015non}
J.~Yang, F.~Stern, A non-iterative direct forcing immersed boundary method for
  strongly-coupled fluid--solid interactions, J. Comput. Phys. 295 (2015)
  779--804.

\bibitem{wang2015strongly}
C.~Wang, J.~D. Eldredge, Strongly coupled dynamics of fluids and rigid-body
  systems with the immersed boundary projection method, J. Comput. Phys. 295
  (2015) 87--113.

\bibitem{kim2016penalty}
Y.~Kim, C.~S. Peskin, A penalty immersed boundary method for a rigid body in
  fluid, Phys. Fluids 28~(3) (2016) 033603.

\bibitem{lacis2016stable}
U.~L{\=a}cis, K.~Taira, S.~Bagheri, A stable fluid--structure-interaction
  solver for low-density rigid bodies using the immersed boundary projection
  method, J. Comput. Phys. 305 (2016) 300--318.

\bibitem{Conca1997}
C.~Conca, A.~Osses, J.~Planchard, Added mass and damping in fluid--structure
  interaction, Comput. Method. Appl. Mech. Engrg. 146~(3-4) (1997) 387--405.

\bibitem{Robinson-MosherSchroederFedkiw2011}
A.~Robinson-Mosher, C.~Schroeder, R.~Fedkiw, A symmetric positive definite
  formulation for monolithic fluid structure interaction, J. Comput. Phys.
  230~(4) (2011) 1547--1566.

\bibitem{BadiaQuainiQuarteroni2008}
S.~Badia, A.~Quaini, A.~Quarteroni, Splitting methods based on algebraic
  factorization for fluid--structure interaction, SIAM J. Sci. Comput. 30~(4)
  (2008) 1778--1805.

\bibitem{ICNS}
W.~D. Henshaw, A fourth-order accurate method for the incompressible
  {N}avier-{S}tokes equations on overlapping grids, J. Comput. Phys. 113~(1)
  (1994) 13--25\citeCount{154}.

\bibitem{splitStep2003}
W.~D. Henshaw, N.~A. Petersson, A split-step scheme for the incompressible
  {Navier-Stokes} equations, in: M.~M. Hafez (Ed.), Numerical Simulation of
  Incompressible Flows, World Scientific, 2003, pp. 108--125\citeCount{30}.

\bibitem{rbins2016r}
J.~W. Banks, W.~D. Henshaw, D.~W. Schwendeman, Q.~Tang, A stable partitioned
  {FSI} algorithm for rigid bodies and incompressible flow. {Part II}: General
  formulation, preprint arXiv:??, arXiv, submitted for publication (2016).

\bibitem{mog2006}
W.~D. Henshaw, D.~W. Schwendeman, Moving overlapping grids with adaptive mesh
  refinement for high-speed reactive and non-reactive flow, J. Comput. Phys.
  216~(2) (2006) 744--779\citeCount{54}.

\bibitem{lrb2013}
J.~W. Banks, W.~D. Henshaw, B.~Sj{\"o}green, A stable {FSI} algorithm for light
  rigid bodies in compressible flow, J. Comput. Phys. 245 (2013)
  399--430\citeCount{1}.

\bibitem{fsi2012}
J.~W. Banks, W.~D. Henshaw, D.~W. Schwendeman, Deforming composite grids for
  solving fluid structure problems, J. Comput. Phys. 231~(9) (2012)
  3518--3547\citeCount{5}.

\bibitem{flunsi2016}
J.~W. Banks, W.~D. Henshaw, A.~Kapila, D.~W. Schwendeman, An added-mass
  partitioned algorithm for fluid-structure interactions of compressible fluids
  and nonlinear solids, J. Comput. Phys. 305 (2016) 1037--1064\citeCount{0}.

\bibitem{fib2014}
J.~W. Banks, W.~D. Henshaw, D.~W. Schwendeman, An analysis of a new stable
  partitioned algorithm for {FSI} problems. {Part I}: Incompressible flow and
  elastic solids, J. Comput. Phys. 269 (2014) 108--137\citeCount{2}.

\bibitem{fis2014}
J.~W. Banks, W.~D. Henshaw, D.~W. Schwendeman, An analysis of a new stable
  partitioned algorithm for {FSI} problems. {Part II}: Incompressible flow and
  structural shells, J. Comput. Phys. 268 (2014) 399--416\citeCount{2}.

\bibitem{beamins2016}
L.~Li, W.~D. Henshaw, J.~W. Banks, D.~W. Schwendeman, G.~A. Main, A stable
  partitioned {FSI} algorithm for incompressible flow and deforming beams, J.
  Comput. Phys. 312 (2016) 272--306.

\bibitem{Petersson00}
N.~A. Petersson, Stability of pressure boundary conditions for {S}tokes and
  {N}avier-{S}tokes equations, J. Comput. Phys. 172~(1) (2001) 40--70.

\bibitem{GKSII}
B.~Gustafsson, H.-O. Kreiss, A.~Sundstr\"om, Stability theory of difference
  approximations for mixed initial boundary value problems. {II}, Mathematics
  of Computation 26~(119) (1972) 649--686.

\end{thebibliography}

\end{document}